\documentclass[12pt,leqno]{article}

\usepackage{amsmath}
\usepackage{amscd}
\usepackage{amsmath}
\usepackage{amssymb}
\usepackage{latexsym}
\makeatletter

\@addtoreset{equation}{section}
\makeatother

\newtheorem{thm}{Theorem}[subsection]
\newtheorem{pr}[thm]{Proposition}
\newtheorem{df}[thm]{Definition}
\newtheorem{lm}[thm]{Lemma}
\newtheorem{cor}[thm]{Corollary}
\newtheorem{cn}[thm]{Conjecture}

\newtheorem{ex}[thm]{Example}

\newcommand{\qed}{\hfill{\rule{7pt}{7pt}}\smallskip}
\input xy \xyoption{all} \CompileMatrices

\begin{document}

\title{Ramification theory
for varieties over a local field}
\author{{\sc Kazuya Kato and Takeshi Saito}}
\maketitle

\begin{abstract}
We define generalizations of
classical invariants of 
wild ramification
for coverings on a variety
of arbitrary dimension over a local field.
For an $\ell$-adic sheaf,
we define its Swan class
as a 0-cycle class supported on
the wild ramification locus.
We prove a formula
of Riemann-Roch type
for the Swan conductor
of cohomology 
together with 
its relative version,
assuming
that the local field is
of mixed characteristic.

We also prove the integrality
of the Swan class
for curves over
a local field
as a generalization
of the Hasse-Arf theorem.
We derive a proof 
of a conjecture of Serre
on the Artin character
for a group action
with an isolated fixed point
on a regular local ring,
assuming the dimension is 2.
\end{abstract}

\section*{Introduction}

\subsection*{0.1. The goal of this paper}

Let $K$ be a complete
discrete valuation field
and 
${\cal O}_K$ be the valuation ring.
We assume that the residue field $F$ 
is perfect of 
characteristic $p>0$. 
In this article, we generalize the classical ramification theory 
of extensions of $K$ briefly recalled in 0.2 to the ramification theory for varieties over $K$ as is described in 0.3--0.5
below. We also prove a conductor formula of Riemann-Roch type stated in 0.6.

We fix a prime number
$\ell$ different
from $p$.
Let $U$ be a separated
smooth scheme
of finite type over
$K$ and ${\cal F}$
be a smooth $\ell$-adic
sheaf on $U$.
The alternating
sum ${\rm Sw}_K
H^*(U_{\bar K},
{\cal F})$
of the
Swan conductor
is defined as
an invariant measuring
the wild ramification
of the $\ell$-adic
representation
$H^*(U_{\bar K},
{\cal F})$
of the absolute
Galois group $G_K$ of $K$.
We define
an element
${\rm Sw}_U{\cal F}$
(see 0.3--0.5
in the introduction and
Definition \ref{dfSwU} in the text)
called the Swan class
as a certain 0-cycle
class supported
on the closed
fiber
of a compactification of $U$ over ${\cal O}_K$
and prove a
conductor formula
$${\rm Sw}_KH^*_c
(U_{\bar K},{\cal F})
=
{\rm rank}\ {\cal F}
\cdot
{\rm Sw}_KH^*_c
(U_{\bar K},{\mathbb Q}_\ell)
+
\deg {\rm Sw}_U{\cal F},
\leqno{\rm (0.1)}
$$
assuming that
$K$ is of characteristic 0 in
Corollary \ref{corUsm}.
We also prove
a relative version
(see (0.6) below)
of the conductor formula
in Theorem \ref{thmcf}.

The formula (0.1)
is an arithmetic analogue of
the higher dimensional
generalization
of the Grothendieck-Ogg-Shafarevich
established in \cite{KSA}.
The term
${\rm Sw}_KH^*_c
(U_{\bar K},{\mathbb Q}_\ell)$
has been computed
in the case where $U$
is further assumed
proper over $K$
by the conductor formula
of Bloch, proved 
under some mild assumption
in \cite{KSI}.
In this paper, we will focus
on a mixed characteristic case.
Another approach in a geometric
equal characteristic case
is studied in \cite{Tsu}.

\subsection*{0.2 Invariants in classical ramification theory}

We first recall the classical ramification theory. For a finite separable extension 
$L$ of $K$, we have the following invariants of ramification in (i)--(iii) below, which 
are integers $\geq 0$.
In (ii) and (iii),  we assume that $L/K$ is a Galois extension with Galois 
group $G$.

\begin{itemize}
\item[(i)] 
The different $D_{L/K}$ and the logarithmic different $D_{L/K}^{\log} 
=D_{L/K}-e_{L/K}+1$, where $e_{L/K}$ is the ramification index of the extension $L/K$.
\item[(ii)] 
The Lefschetz number $i(\sigma)$ and the logarithmic Lefschetz number 
$j(\sigma)$ for $\sigma\in G\setminus \{1\}$ defined as 
\begin{eqnarray*}
i(\sigma)&=\min\{\text{ord}_L(\sigma(a)-a)\;|\;a\in {\cal O}_L\}, \\
j(\sigma)&=\min\{\text{ord}_L(\sigma(a)/a-1)\;|\;a\in L^\times\}.
\end{eqnarray*}
\item[(iii)] 
The Artin conductor ${\rm Art}(\rho)$ and the Swan conductor ${\rm Sw}(\rho)$ for a finite 
dimensional representation $\rho$ of $G$ over a field of characteristic $0$. They are 
defined by
\begin{eqnarray*}
{\rm{Art}}(\rho) 
&= \displaystyle{\frac{1}{e_{L/K}}  \sum_{\sigma\in G-\{1\}} 
i(\sigma)(\dim(\rho)-{\rm Tr}(\rho(\sigma))),}
\\
{\rm{Sw}}(\rho) &
= 
\displaystyle{\dfrac{1}{e_{L/K}}  
\sum_{\sigma\in G-\{1\}} j(\sigma)(\dim(\rho)-{\rm Tr}(\rho(\sigma))).}
\end{eqnarray*} 
\end{itemize}
The Hasse-Arf theorem
asserts the highly non-trivial fact
that these 
conductors are in fact integers.

These invariants
are linked by several important formulas (see \cite{CL}
for example). 
For example, in the case $L/K$ is Galois with Galois group $G$, we have $$D_{L/K}= 
\sum_{\sigma\in G\setminus \{1\}} i(\sigma), \quad D_{L/K}^{\log}= \sum_{\sigma\in 
G\setminus \{1\}} j(\sigma).$$

The invariants $D_{L/K}^{\log}$, $j(\sigma)$ and ${\rm{Sw}}(\rho)$ are the parts of 
$D_{L/K}$, $i(\sigma)$ and ${\rm{Art}}(\rho)$, respectively,  which handles the wild 
ramification.
We will focus on the wild ramification and introduce generalizations of $D_{L/K}^{\log}$, 
$j(\sigma)$, and $\rm{Sw}(\rho)$.

\subsection*{0.3 Generalization}

In our generalization of ramification theory in \cite{KSA} (resp.\  in this paper), 
in place of $L/K$ in 0.2,
we consider 
a finite \'etale morphism $$f\colon V\to U$$ of non-singular algebraic varieties over a 
perfect field $k$ of characteristic $p$  (resp.\ over $K$) and study the ramification 
of $f$ along the boundary of compactifications of $V$ and $U$ over $k$ (resp.\ ${\cal O}_K$).  
We call the case over $k$ the geometric case (geo) and the case over $K$ the arithmetic 
case (ari). For simplicity, in this introduction, we assume ${\rm{char}}\; K=0$ in 
the case (ari). 
Although the main theme of this paper is the arithmetic case, we describe also 
the geometric case in 0.3 and
0.4 to compare.

In the case (geo) (resp.\ (ari)), for a proper scheme $Y$ over $k$ (resp.\ ${\cal O}_K$) which 
contains $V$ as a dense open subscheme, we define in Definition \ref{dftmT}  the wild ramification locus 
$\Sigma_{V/U}Y$ of $f\colon V\to U$ on $Y$ as a closed subset of $Y$ . The wild ramification 
locus satisfies the following properties:
1. $V\cap \Sigma_{V/U}Y=\emptyset$.
2. If $Y'$ is a proper scheme over $k$ (resp.\ ${\cal O}_K$) containing $V$ as a dense open subscheme and if $Y'\to Y$ is a morphism inducing the identity on $V$, 
then $\Sigma_{V/U}Y$ coincides with the image of $\Sigma_{V/U}Y'$.
3. In the case (ari) (recall that we assume ${\rm{char}}\;K=0$), $\Sigma_{V/U}Y$ is 
contained in the special fiber $Y\otimes_{{\cal O}_K} F$ of $Y$.

\medskip

For a  proper scheme $X$ over $k$ (resp.\ ${\cal O}_K$) which contains $U$ as a dense open 
subscheme, we define the wild ramification locus $\Sigma_{V/U}X$ of $f$ on $X$, which 
is a closed subset of $X$,
to be the image  
$\bar f(\Sigma_{V/U}Y)$ for
a morphism $\bar f\colon Y\to X$ 
of compactifications
extending $f\colon V\to U$.
This also satisfies 
analogous properties 
corresponding to the above 1, 2, 3.

For a commutative ring $R$, let
$$F_0G(\partial_{V/U} V)_R:=\varprojlim_Y (F_0G(\Sigma_{V/U}Y)\otimes_{\mathbb Z} R), 
\quad F_0G(\partial_{V/U} U)_R:=\varprojlim_X (F_0G(\Sigma_{V/U}X)\otimes_{\mathbb 
Z} R)$$ where $Y$ runs through proper integral schemes over $k$ (resp.\
${\cal O}_K$) containing $V$ as a dense open subscheme and $X$ runs through proper integral 
schemes over $k$ (resp.\ ${\cal O}_K$) containing $U$ as a dense open subscheme.
Here
$G(-)$ denotes the Grothendieck group of coherent sheaves, and $F_0G(-)$ denotes the 
part generated by the classes of coherent sheaves with finite supports.

Let $Z(V/U)$ 
denote the free abelian group
on the set of connected components of 
the complement
$V\times_U V\setminus \Delta_V$ 
of the diagonal 
$\Delta_V\subset V\times_U V$.
Note that since $f:V\to U$ is \'etale, 
$\Delta_V$ is open and closed in $V\times_U V$. 
The definition of
generalizations
of invariants of
wild ramification
is based on 
a homomorphism 
$$Z(V/U)
\to F_0G(\partial_{V/U}V)_{\mathbb Q},
\leqno{(0.2)}$$ 
whose definition
will be sketched
in 0.4 below. 
The homomorphism (0.2) is
called the localized intersection product
with logarithmic diagonal and denoted by $(-, \Delta_V)^{\log}$ (resp.\ $((-, 
\Delta_V))^{\log})$
in the case (geo) (reps.\ (ari)).
Though $V\times_U V\setminus \Delta_V$ does not intersect with the diagonal, 
the localized intersection with the log diagonal appears on the boundary of $V$ in 
a compactification $Y$.

\begin{itemize}
\item[(i)] We define $$D_{V/U}^{\log}\in F_0G(\partial_{V/U}V)_{\mathbb Q}$$ by 
$$D_{V/U}^{\log} = ([V\times_U V\setminus \Delta_V], \Delta_V)^{\log} \quad \text{in 
the case (geo)},$$ $$D_{V/U}^{\log} = (([V\times_U V\setminus \Delta_V], 
\Delta_V))^{\log} \quad \text{in the case (ari)}.$$
\item[(ii)]
In the case $V\to U$ is a Galois covering with Galois group $G$, then for $\sigma\in 
G\setminus \{1\}$, we define $$j(\sigma)\in F_0G(\partial_{V/U} V)_{\mathbb Q}$$ by
$$j(\sigma) = ([\Gamma_{\sigma}], \Delta_V)^{\log} \quad \text{in the case 
(geo)},$$ $$j(\sigma) = (([\Gamma_{\sigma}], \Delta_V))^{\log} \quad \text{in the case 
(ari)},$$  where $\Gamma_{\sigma}$ is the graph of $\sigma$.
\item[(iii)] 
For a finite dimensional representation $\rho$ of $G$ over a field of 
characteristic $0$, we define the Swan class
$${\rm{Sw}}(\rho)=  \frac{1}{\sharp(G)}  
\sum_{\sigma\in G-\{1\}} f_*(j(\sigma))(\dim(\rho)-{\rm Tr}(\rho(\sigma)))\in
F_0G(\partial_{V/U}U)_{{\mathbb Q}(\zeta_{p^\infty})}$$ where $f_*$ is the push forward
$F_0G(\partial_{V/U}V)_{\mathbb Q}\to F_0G(\partial_{V/U}U)_{\mathbb Q}$ and 
${\mathbb Q}(\zeta_{p^\infty})=\cup_n {\mathbb Q}(\zeta_{p^n})$ with $\zeta_{p^n}$ a 
primitive $p^n$-th root of unity.
\end{itemize}

In (ii), we have $D_{V/U}^{\log}=\sum_{\sigma\in G\setminus \{1\}} 
j(\sigma)$ simply because $V\times_U V \setminus \Delta_V=\coprod_{\sigma\in 
G\setminus\{1\}} \Gamma_{\sigma}$.
We expect that we can remove  $\otimes\mathbb Q$ and $\otimes{\mathbb 
Q}(\zeta_{p^\infty})$ in the definitions of the above invariants in (i)--(iii).

To formulate a conductor formula given below, we define ${{\rm Sw}}(\rho)$ also 
for a finite dimensional representation $\rho$ of $G$ over a field of characteristic 
$\ell$ by
$${\rm Sw}(\rho)=  \frac{1}{\sharp(G)}  \sum_{\sigma\in G-\{1\}} 
f_*(j(\sigma))(\dim(\rho)-{\rm Tr}^{\rm Br}(\rho(\sigma)))\in
F_0G(\partial_{V/U}U)_{{\mathbb Q}(\zeta_{p^\infty})}$$ using the Brauer trace.
The definition makes sense
because we have
$j(\sigma)=0$ unless
the order of $\sigma$ is not a power of $p$.

\medskip
The relation with classical ramification theory is as follows.

In the case (ari), assume $U={\rm Spec}\ K, V={\rm Spec}\ L,Y={\rm Spec}\ {\cal O}_L$. Then 
if $L/K$ is wildly ramified (resp.\ at most tamely ramified),  $\Sigma_{V/U}Y$ consists of the 
closed point of $Y$ (resp.\ the empty set). If $L/K$ is wildly ramified, we have
$F_0G(\Sigma_{V/U}Y)=\mathbb Z$, and $D_{V/U}^{\log}$ and 
$j(\sigma)$ defined above
recover the classical $D_{L/K}^{\log}$ and $j(\sigma)$, 
respectively.
 In the case (geo), assume that $Y$, $V$, $U$ are smooth curves over $k$, and let 
$K_0$ (resp.\ $L_0$) be the function field of $U$ (resp.\ $V$). Then $\Sigma_{V/U}Y$ consists of the places of $L_0$ 
where the extension $L_0/K_0$ is
wildly ramified and $F_0G(\Sigma_{V/U}Y)$ is the direct sum of $\mathbb 
Z$ indexed by these places. For
  $v\in \Sigma_{V/U}Y$, if $u$ denotes the place of $K_0$ lying under $v$ and if $K$ (resp.\ 
$L$) denotes the  completion of $K_0$ (resp.\
$L_0$) at $u$ (resp.\ $v$), the $v$-components of $D_{V/U}^{\log}$ and 
$j(\sigma)$ defined above coincide with $D_{L/K}^{\log}$ and $j(\sigma)$ in the classical 
ramification theory for $L/K$, respectively.

The revolutionary idea that
the invariant of ramification
should be defined as
a 0-cycle class on the ramification locus is due to S.\ Bloch \cite{bloch}.

\subsection*{0.4 The definition of the localized intersection product with logarithmic diagonal}

Let $V\to U$ be
a finite \'etale morphism
of smooth integral schemes
over $k$ (resp.\ over $K$)
in the case (geo) (resp.\ (ari)).
We put $n=\dim U$
(reps.\ $n=\dim U_K+1$) in the case (geo) (resp.\ (ari)).
The embedding theorem of Nagata
and the theory of alteration of de Jong
gives us a cartesian diagram 
of integral schemes over $k$ 
(resp.\ ${\cal O}_K$) in the case (geo) (resp.\ (ari))
$$\begin{CD}
V@<g<<W\\
@V{\cap}VV@VV{\cap}V\\
Y@<{\bar g}<<Z
\end{CD}
\leqno{(0.3)}$$
where $Y$ and $Z$ are proper over $k$ (resp.\ ${\cal O}_K$)
and satisfy the following properties:
The vertical arrows are open immersions 
with dense images,  the arrow
$\bar g\colon Z\to Y$
is surjective and 
generically finite
and $Z$ 
is regular and contains 
$W$ as the complement
of a divisor with simple 
normal crossings.

In the case (geo) (resp.\ (ari)),
we define
the logarithmic 
self-product $(Z\times_k Z)^\sim$ (resp.\
$(Z\times_{{\cal O}_K} Z)^\sim$)
 as a modification of the usual product $Z\times_k Z$ (resp.\
$Z\times_{{\cal O}_K} Z$).
Let $P$ denote
$(Z\times_k Z)^\sim$ (resp.\
$(Z\times_{{\cal O}_K} Z)^\sim$).
The diagonal map
$Z\to Z\times_k Z$ (resp.\
$Z\to Z\times_{{\cal O}_K} Z$)
is canonically lifted
to a closed immersion 
$Z\to P$ called the log diagonal map.
The scheme $P$ contains 
$W\times_k W$ (resp.\
$W\times_{{\cal O}_K} W$)
as an open subscheme.
Let $A$ be the closure of $W\times_U W\setminus W\times_V W$ in $P$
and $\Sigma$ be
the intersection of $A$ with the logarithmic diagonal
$Z$ in $P$.

We define the intersection product
with the logarithmic diagonal $Z$ in $P$
as a homomorphism
$$G(A) \to G(\Sigma)
\leqno{(0.4)}$$
as follows.
See Proposition \ref{prmap} for detail. 
Regard ${\cal O}_Z$ as an ${\cal O}_P$-module via the log diagonal. 
In the case (geo), 
the map (0.4) is defined as the usual intersection product with the class $[{\cal O}_Z]$
for a smooth scheme $P$.
Namely, it maps the class of
a coherent ${\mathcal O}_P$-module
${\cal F}$ supported on $A$ to the alternating sum: $$[\\{\cal F}]\mapsto 
\sum_{i=0}^{\dim P} 
(-1)^i[{\cal T}or_i^{{\cal O}_P}({\cal F}, {\cal O}_Z)].$$ In the case (ari), it is defined 
as $$[\\{\cal F}]\mapsto [{\cal T}or_{2i}^{{\cal O}_P}({\cal F}, {\cal O}_Z)]-[{\cal 
T}or_{2i-1}^{{\cal O}_P}({\cal F}, {\cal O}_Z)]$$ 
for sufficiently large integer $i$.
In the case (ari), 
it is proved in \cite{KSI}
that the class $[{\cal T}or_j^{{\cal O}_P}({\cal F}, {\cal O}_Z)]\in G(\Sigma)$ depends only on 
the parity of $j$
for sufficiently large $j$. 

The maps (0.4) for various diagrams (0.3) induce (0.2) as follows.
Let $F_\bullet$ denote the
topological filtration on 
the Grothendieck group $G(-)$.
We regard the free abelian group
$Z(V/U)$ as the graded quotient 
${\rm Gr}^F_n
G(V\times_U V\setminus \Delta_V)$
by the canonical surjection defined 
by taking the length at the generic point of each connected component of $V\times_U 
V\setminus \Delta_V$.
We prove in Proposition \ref{prprod}
that the homomorphisms
${\rm Gr}^F_\bullet
G(A)\to {\rm Gr}^F_{\bullet-n}
G(\Sigma)$
induced by (0.4) 
factor through the canonical 
surjection ${\rm Gr}^F_\bullet
G(A)\to {\rm Gr}^F_\bullet
G(W\times_U W\setminus W\times_V W)$
defined by the restriction
if the following condition is satisfied:
\begin{itemize}
\item[($X$)]
There exists a cartesian diagram
$$\begin{CD}
U@<{f\circ g}<<W\\
@V{\cap}VV @VV{\cap}V\\
X@<<<Z
\end{CD}$$
over $k$ (resp.\ over ${\cal O}_K$)
where 
$X$ is a proper scheme
over $k$ (resp.\ over ${\cal O}_K$)
containing
$U$ as the complement of a Cartier divisor.
\end{itemize}
Consequently, we obtain 
$$
Z(V/U)=
{\rm Gr}^F_n
G(V\times_U V \setminus \Delta_V) \to 
{\rm Gr}^F_n
G(W\times_U W\setminus W\times_V 
W) \to F_0G(\Sigma)$$
where the first arrow is the pull-back
by $g\times g$
and the second arrow is
induced by (0.4).

If further the condition
\begin{itemize}
\item[($Y$)]
 $\bar g(\Sigma)
\subset
\Sigma_{V/U}Y$
\end{itemize}
is satisfied,
the composition $$Z(V/U)
\to
F_0G(\Sigma)\overset{\bar g_*}\to
F_0G(\Sigma_{V/U}Y)\otimes_{\mathbb Z} {\mathbb Q}$$
with
the push-forward map $\bar g_*$
divided by the generic 
degree $[Z:Y]$ of $Z$ over $Y$
is defined. 
We prove in Theorem \ref{thmmap} that 
such $Z$ satisfying the conditions ($X$) and ($Y$)
does exist and that
the composition
$Z(V/U)
\to
F_0G(\Sigma_{V/U}Y)\otimes_{\mathbb Z} {\mathbb Q}$ is independent of $Z$
and forms an inverse system
to define the required map
(0.2).

\subsection*{0.5 The Swan class of a constructible sheaf}

In order to formulate
a conductor formula of Riemann-Roch type
in 0.6 below,
we extend the definition of
the Swan class to constructible sheaves.
In the rest of introduction, we
consider the arithmetic case (we assume   ${\rm{char}}\; K=0$).

For a separated scheme $U$ of finite type over $K$ and for a commutative ring $R$, 
let $$F_0G(\partial_F U)_R:= \varprojlim_X  (F_0G(X\otimes_{{\cal O}_K} F)\otimes_{\mathbb 
Z} R)$$ where $X$ runs through proper schemes over ${\cal O}_K$ which contain $U$ as a 
dense open subscheme.

Let $f\colon V\to U$
be a finite \'etale morphism
and let $\bar Z(V/U)$
denote the free abelian group
on the set of connected components of
$V\times_U V$.
It is the direct sum of
$Z(V/U)$ with the free
abelian group of rank 1
generated by the class of $\Delta_V$.
The composition $((-, \Delta_V))^{\log}\colon Z(V/U)\to 
F_0G(\partial_{V/U} V)_{\mathbb Q}\to F_0G(\partial_F V)_{\mathbb Q}$ is naturally 
extended to $$((-, \Delta_V))^{\log}
\colon \bar Z(V/U)\to  F_0G(\partial_F 
V)_{\mathbb Q}.
\leqno{(0.5)}$$
To define the map (0.5),
we proceed similarly as in 0.4.
Namely, we consider a diagram (0.3) and,
letting $A'\subset P$ denote the closure of $W\times_U W$, we define 
a homomorphism $$G(A')\to F_0G(Z_F)\; ;[{\cal F}]\mapsto [{\cal 
T}or_{2i}^{{\cal O}_P}({\cal F}, {\cal O}_Z)]-[{\cal T}or_{2i-1}^{{\cal O}_P}({\cal 
F}, {\cal O}_Z)] \quad (i\gg 0)$$ 
similarly as (0.4).
Then this induces (0.5) in the same way as
(0.4) induces (0.2).
In particular, $j(\sigma)\in F_0G(\partial_F V)_{\mathbb Q}$ is defined even for 
$\sigma=1$ as $((\Delta_V, \Delta_V))^{\log}$.

We extend 
the definition of the Swan class 
of smooth sheaves sketched in 0.2 and 0.3 to constructible sheaves.

\medskip
\noindent{\bf Proposition 1}
{\rm (Proposition \ref{prSwan}, Corollary \ref{corSwb})}
{\em Assume ${\rm{char}}\; K=0$. 
Then, there is a unique way to define 
$${{\rm Sw}}_U {\cal F}, \;  \overline{{\rm Sw}}_U {\cal F}\in F_0G(\partial_F 
U)_{{\mathbb Q}(\zeta_{p^\infty})}$$ for any separated scheme $U$ of finite type over 
$K$ and for any constructible ${\bar {\mathbb F}}_{\ell}$-sheaf $\cal F$ on $U$, 
satisfying the following conditions {\rm (1)--(3)}.

\medskip
{\rm
(1)} Assume $U$ is a non-singular variety and  ${\cal F}$ is locally constant. Let $f\colon V\to U$ be 
a finite \'etale Galois covering of $U$ with Galois group $G$ on which the pull-back 
$f^*{\cal F}$ is a constant sheaf. Then ${\rm{Sw}}_U {\cal F}$ is the image of 
$\rm{Sw}(\rho)$ in {\rm 0.3}
for the ${\bar { \mathbb 
F}}_{\ell}$-representation $\rho$ of $G$  corresponding to $\cal F$.
We also have
$$\overline{{\rm{Sw}}}_U {\cal F}=  {\rm{Sw}}_U {\cal F}-{\rm{rank}\  }{\cal F}\cdot 
((\Delta_U, \Delta_U))^{\log}  = -\frac{1}{\sharp(G)} \sum_{\sigma\in G} 
f_*j(\sigma)\cdot {\rm Tr}^{\rm Br}(\rho(\sigma)).$$ \medskip

{\rm
(2)} For an exact sequence $0\to {\cal F}'\to {\cal F}\to {\cal F}''\to 0$ of 
constructible ${\bar {\mathbb F}}_{\ell}$-sheaves on $U$, we have $${\rm Sw}_U{\cal 
F}= {\rm Sw}_U{\cal F'} +{\rm Sw}_U{\cal F''}, \quad {\overline{\rm Sw}_U}{\cal 
F}= {\overline{\rm Sw}_U}{\cal F'}
+{\overline{\rm Sw}_U}{\cal F''}.$$

\medskip

{\rm 
(3)} If $i:U'\to U$ is an immersion of schemes
over $K$, we have $${\overline 
{\rm Sw}_U}(i_{!}{\cal F}) =i_{!} {\overline {\rm Sw}_{U'}}{\cal F}$$ for a 
constructible ${\bar {\mathbb F}}_{\ell}$-sheaf $\cal F$ on $U'$. Here $i_{!}$ on the 
right hand side is the canonical homomorphism $F_0G(\partial_F {U'})_{{\mathbb 
Q}(\zeta_{p^\infty})}\to F_0G(\partial_F U)_{{\mathbb Q}(\zeta_{p^\infty})}$.
}
\medskip

Here in (1), $((\Delta_U, \Delta_U))^{\log}\in F_0G(\partial_F U)_{\mathbb Q}$ is 
defined by (0.5) for $U=V$.
The key ingredient of the proof of Proposition 1 is an excision formula Theorem
\ref{thmexc}.

For $U$ as in Proposition 1 and for a constructible ${\bar {\mathbb 
Q}}_{\ell}$-sheaf ${\cal F}$ on $U$, 
we define ${{\rm Sw}}_U {\cal F},  \overline{{\rm Sw}}_U {\cal F}\in 
F_0G(\partial_F U)_{{\mathbb Q}(\zeta_{p^\infty})}$ as those of the ${\bar {\mathbb 
F}}_{\ell}$-sheaf which is obtained from $\cal F$ by taking modulo $\ell$. In the case 
$U$ is regular and $\cal F$ is smooth and trivialized by a finite \'etale Galois covering 
$V\to U$ with Galois group $G$, ${{\rm Sw}}_U {\cal F}$ is the image of 
${{\rm Sw}}(\rho)$ in 0.3 where $\rho$ is the representation of $G$ over ${\bar 
{\mathbb Q}}_{\ell}$ corresponding to $\cal F$.

\subsection*{0.6 The conductor formula}

We prove the following conductor formula of Riemann-Roch type.

\medskip
\noindent
{\bf Theorem 2}
{\rm  (Theorem \ref{thmcf})}
{\em
Assume ${\rm{char}}\;K=0$.
Let
$f\colon U\to V$ be a morphism of
separated
schemes of finite
type over $K$ and let ${\cal F}$ be a
constructible ${\bar {\mathbb F}}_{\ell}$-sheaf (resp.\ ${\bar {\mathbb
Q}}_{\ell}$-sheaf)
 on $U$. Then we have
$$
\overline{\rm Sw}_VRf_!{\cal F}
=
f_!
\overline{\rm Sw}_U{\cal F}
\leqno{\rm (0.6)}
$$
where $f_!$ on the right hand side is the canonical homomorphism $F_0G(\partial_F 
U)_{{\mathbb Q}(\zeta_{p^\infty})}\to F_0G(\partial_F V)_{{\mathbb 
Q}(\zeta_{p^\infty})}$.
}
\medskip

In the case where
${\cal F}$ is smooth
and $V={\rm Spec}\ K$,
the equality (0.6)
specializes to
the conductor formula (0.1)
for the alternating
sum of the Swan conductor
(Corollary \ref{corUsm}).
It also gives
$$
{\rm Sw}_KH^*_c
(U_{\bar K},{\mathbb Q}_\ell)
=
-\deg((\Delta_U,
\Delta_U))^{\log},
\leqno{\rm (0.7)}$$
which
is a generalization
of the conductor formula
of Bloch \cite{bloch}
proved under some mild assumption
in \cite{KSI}.
A special case of
$\dim U_K=1$ and $V={\rm Spec}\ K$
has been studied
in \cite{abb1}.
A crucial ingredient
in the proof of the equality
(0.6) is a logarithmic
variant 
Theorem \ref{thmlLTF} of
the Lefschetz trace
formula for
open varieties.

\subsection*{0.7 Integrality}

As a generalization
of the classical theorem
of Hasse-Arf,
we expect
that the Swan class
${\rm Sw}_U{\cal F}$
should have no
denominator,
Conjecture \ref{cnHA}.
By a standard argument
using Brauer induction,
it is reduced to
the rank one case.
Theorem \ref{thmrk1}
comparing
the Swan class
${\rm Sw}_U{\cal F}$
for a smooth sheaf ${\cal F}$
of rank 1
with a cycle class
$c_{\cal F}$,
defined earlier by one
of the authors,
implies the following integrality.

\medskip
\noindent 
{\bf Theorem 3}
{\rm (Corollary \ref{corint}.1)}
{\em Assume 
${\rm{char}}\;K=0$. Let $U$ be a scheme of finite type over $K$ of dimension $\leq 
1$ and let $\cal F$ be a constructible ${\bar {\mathbb F}}_{\ell}$-sheaf (resp.\ ${\bar 
{\mathbb
Q}}_{\ell}$-sheaf)
 on $U$. Then
${\rm{Sw}}_U \cal F$ belongs to the image of $$F_0G(\partial_F U)_{\mathbb Z} \to 
F_0G(\partial_F U)_{{\mathbb Q}(\zeta_{p^\infty})}.$$ }
\smallskip

From this integrality,
we derive 
the two dimensional case of a conjecture of Serre \cite{Sear}, as is announced 
in \cite{ICM}. In \cite{Sear}, Serre conjectures (see Conjecture \ref{cnSe}) that the theory of Artin 
characters in the ramification theory of a discrete valuation ring can be generalized 
to any regular local ring $A$ with a finite group of automorphisms under a condition 
of isolated fixed point.
An equal characteristic case has been proved
earlier in \cite{Artin}
and some special case has been proved
in \cite{abb2}.

\medskip
\noindent 
{\bf Theorem 4}
{\rm (Corollary \ref{corint}.2.)}
{\em The conjecture of Serre {\rm \cite{Sear}} is true in the case $\dim (A)=2$.
}

\subsection*{0.8 Organization of this paper}

We sketch the content
of each section.
The first three sections
are preliminaries.
In Section \ref{sLef},
after preparing
general terminologies
on semi-stable schemes,
log products, etc.,
we prove a logarithmic
Lefschetz trace
formula,
which is a crucial step
in the proof
of the formula (0.6).
The trace formula
is a sort of mixture
of those proved
in \cite{KSI} and
in \cite{KSA}.
In Section \ref{stm},
we study the tame
ramification
of an \'etale morphism
along the boundary,
using log products.
The purpose of studying
tame ramification
first is
to define the
wild ramification
locus and to focus on it.
We give criterions
for tameness
in terms of valuation rings,
using the quasi-compactness
of the limit of
proper modifications.
In Section \ref{slp},
first we compute
certain tor-sheaves,
which is a crucial
step in the proof
of the excision formula.
We also give some
complement on
the localized Chern classes
and the excess intersection
formula studied in \cite{KSI}
as a preliminary
for the computation
of the logarithmic different.

In Sections \ref{siw},
\ref{siw2}
and \ref{sfml},
we define the invariants
of wild ramification
and establish
their properties.
First,
in Section \ref{siw},
we study the local
structures
of log products
of schemes
over $S={\rm Spec}\
{\cal O}_K$.
In Section \ref{siw2},
we define the
invariants and
study its basic properties.
Section \ref{sfml}
is technically the heart
of the article.
We prove the excision
formula for the invariants.
We also give a formula
in some semi-stable case,
which is a crucial
step in the proof of
the formula (0.6).

In Section \ref{sSw},
we define the Swan class
and prove the formula (0.6).
In Section \ref{srk1},
we compute the Swan
class in the case of
rank 1
and deduce
the integrality of
the Swan class
and complete
the proof of the
conjecture of Serre
in the case of
dimension 2.

The logical structure
of the proof of the
formula (0.6)
is summarized as follows.
We deduce a formula
Proposition \ref{prcysd}
in some semi-stable case
from the log Lefschetz
trace formula
Theorem \ref{thmlLTF}.
We prove a formula
Propositions \ref{prcfs},
\ref{prcfs0}
for stable curves
using
Proposition \ref{prcysd}
and a compatibility
with cospecialization
map Proposition \ref{prcc}.
We complete the proof of the formula (0.6)
in Theorem \ref{thmcf}
by deducing it from
a special case
Corollary \ref{corcfs},
by devissage.

\tableofcontents

\newpage 
\section{Log Lefschetz trace
formula}\label{sLef}

We prove a logarithmic
Lefschetz trace
formula 
Theorem \ref{thmlLTF}
for schemes
over a discrete valuation ring
and give a complement
in Section \ref{sscc}.
They play
a crucial role in the proof
of the conductor formula
in Section \ref{sscfs}.
As preliminaries,
we fix terminologies on
semi-stable schemes,
log blow-ups,
log products
and on log stalks
in Sections \ref{ssst},
\ref{sssb},
\ref{sslb},
\ref{sslst}
respectively.

\subsection{Semi-stable schemes
and stable curves}\label{ssst}
We fix some terminology
on semi-stable schemes.

\begin{df}\label{dfsst}
Let $f\colon
X\to S$ be a morphism of
schemes
and $r\ge 0$
be an integer.

{\rm 1.}
We say that $X$ is 
{\rm weakly strictly
semi-stable}
of relative dimension $r$
over
$S$ if the following condition
is satisfied:
\begin{itemize}
\item[
{\rm (\ref{dfsst}.1)}]
For every point $x\in X$,
there exist
an open neighborhood
$x\in U\subset X$,
an affine open neighborhood
$s=f(x)\in {\rm Spec}\ R
\subset S$,
an integer $1\le q\le r+1$,
an element $a\in R$
and an \'etale morphism 
$$\begin{CD}
U@>>> {\rm Spec}\ 
R[T_1,\ldots,T_{r+1}]/
(T_1\cdots T_q-a)
\end{CD}$$
over $S$.
\end{itemize}
If $S={\rm Spec}\ R$
for a discrete valuation
ring $R$, 
we say
a weakly strictly
semi-stable scheme over $S$
is {\rm strictly semi-stable}
if $a\in R$
in {\rm (\ref{dfsst}.1)}
is a uniformizer.

{\rm 2.}
We say that $X$ is 
{\rm weakly semi-stable}
of relative dimension $r$
over
$S$ if,
\'etale locally on $X$ and on $S$,
it is 
weakly strictly 
semi-stable
of relative dimension $r$
over $S$.
Namely,
if the following
condition is satisfied:
\begin{itemize}
\item[
{\rm (\ref{dfsst}.2)}]
For every geometric point
$\bar x\to X$,
there exist
\'etale neighborhoods
$\bar x\to U
\to X$
and 
$\bar s=f(\bar x)\to V
\to S$
and a morphism $U\to V$
compatible with $X\to S$
and with $\bar x\to \bar s$
such that
$U$ is
weakly strictly semi-stable
of relative dimension $r$
over $V$.
\end{itemize}
If $S={\rm Spec}\ R$
for a discrete valuation
ring $R$,
a scheme $X$ over $S$
is said to be
{\rm semi-stable} 
if,
\'etale locally on $X$,
it is 
strictly semi-stable
over $S$.
\end{df}

If $X$
is weakly semi-stable over $S$,
the scheme $X$
is flat over $S$
and is smooth over $S$
on a dense open subscheme
of each fiber.

We show that,
locally on $X$,
the subscheme of $S$
defined
by $a$ 
is well-defined
and
the subschemes of $X$
defined by
$T_1,\ldots,T_q$
are well-defined
up to permutation.

\begin{lm}\label{lmsst}
Let $S={\rm Spec}\ R$
be an affine scheme
and 
$f\colon X\to S$ be a scheme
over $S$.
Assume that
$X$ is \'etale over
$R[T_1,\ldots,T_{r+1}]/
(T_1\cdots T_q-a)$
for an element $a\in R$
and $q\ge 1$.
Let $x$
be a point of $X$
where the morphism
$f\colon X\to S$
is not smooth
and $s=f(x)$.

{\rm 1.}
The annihilator of
the ${\cal O}_{S,s}$-module
$\Omega^{r+1}_{X/S,x}$
is generated by $a$.

{\rm 2.}
Assume that
there exist
$q$ irreducible
components
of the fiber
$X_s$ containing $x$.
Then, the
intersection
${\rm Spec}\ 
{\cal O}_{X,x}
\cap (X\times_S
{\rm Spec}\ R/(a))^{\rm sm}$
with the smooth locus
has $q$ connected
components.
Their schematic closures
in ${\rm Spec}\ 
{\cal O}_{X,x}$
are defined by
$T_1,\ldots,T_q$.
\end{lm}

{\it Proof.}
1.
An explicit computation
in the case where
$X={\rm Spec}\
R[T_1,\ldots,T_{r+1}]/
(T_1\cdots T_q-a)$
shows that
the ${\cal O}_{S,s}$-module
$\Omega^{r+1}_{X/S,x}$
is generated by
one element 
and the annihilator
is generated by
$T_1\cdots 
T_{i-1}T_{i+1}\cdots T_q$
for $i=1,\ldots,q$.
The assertion follows
from this easily.

2.
The irreducible
components
of the fiber
$X_s$
containing $x$
are defined by
$T_1,\ldots,T_q$.
The connected
components of
the intersection
${\rm Spec}\ 
{\cal O}_{X,x}
\cap (X\times_S
{\rm Spec}\ R/(a))^{\rm sm}$
are also defined by
$T_1,\ldots,T_q$.
Thus the assertion follows.
\qed

In Definition \ref{dfsst}.2,
we may take $V=S$
in the condition (\ref{dfsst}.2)
by Lemma \ref{lmsst}.1.

\begin{cor}\label{corwss}
Let
$X$ be a weakly semi-stable scheme
over a scheme $S$.
Then, 
$X$ is weakly strictly semi-stable
over $S$ 
if and only if,
for every point
$s$ of $S$,
each irreducible 
component
of the fiber
$X_s=X\times_Ss$
is smooth
over $s$.
\end{cor}

{\it Proof.}
If $X$ is weakly strictly semi-stable,
each irreducible 
component
of the fiber
$X_s=X\times_Ss$
is clearly smooth over $s$ for
every point $s$ of $S$.
Let $x\in X$ be a point above $s\in S$.
If $X\to S$ is smooth
at $x$,
it is weakly strictly
semi-stable at $x$.
Assume $f\colon
X\to S$ is not smooth at $x$.
Let $\bar x$
be a geometric point
above $x$.
Then the irreducible components
of the strict henselization
${\rm Spec}\
{\cal O}_{X_s,\bar x}$
of the fiber
are defined by
$T_1,\ldots,T_q$
in the notation of Lemma \ref{lmsst}.1.
The pull-back of
an irreducible component
of the fiber $X_s$
is the union of some of them.
Hence, 
each irreducible 
component
of the fiber
$X_s=X\times_Ss$
is smooth at $x$
if and only if
the ideals
$(T_1),\ldots,(T_q)$
are defined in 
${\mathcal O}_{X_s,x}$.
\qed

We may modify a
weakly semi-stable curve
to a
weakly strictly
semi-stable curve,
under an assumption.
This construction will be used
in the proof of Lemma \ref{lmcysd}
in the case (\ref{lmcysd}.1a).

\begin{lm}\label{lmD1}
Let $X$ be a 
weakly semi-stable curve
over a normal scheme $S$
and let $E\subset X$
denote the closed subset
consisting of the points
where $X$ is not smooth
over $S$.
Assume that
$X$ is smooth on
a dense open subscheme
of $S$
and that
the following
condition is satisfied:
\begin{itemize}
\item[{\rm (\ref{lmD1}.1)}]
For every point $x\in E$
and $s=f(x)$,
the element $a\in {\cal O}_{S,s}$
in Lemma {\rm \ref{lmsst}}
is a square up to a unit.
\end{itemize}
Then, 
there exists a quasi-coherent ideal
${\cal I}\subset {\cal O}_X$
such that 
${\cal I}={\cal O}_X$ outside $E$
and that 
the blow-up $X'$ of $X$ at ${\cal I}$ 
is weakly strictly semi-stable
over $S$.
\end{lm}

{\it Proof.}
Let $x$ be a point of $E$.
Then, \'etale locally on a neighborhood
of $x$,
$X$ is \'etale over 
the scheme defined by $T_1T_2-a$
and the ideal $(a)$ is
well-defined
by Lemma {\rm \ref{lmsst}}.1.
We put $a=b^2$.
By the assumption that
$S$ is normal,
the ideal $(b)$
is also well-defined.
Since the ideal $(T_1,T_2)$
is the annihilator
of $\Omega^1_{X/S}$ at $x$,
the ideal ${\cal I}
\subset {\cal O}_X$ \'etale
locally defined by 
$(T_1,T_2,b)$ is well-defined
on $X$.
Then, the blow-up
$X'\to X$
by the ideal ${\cal I}$
satisfies the condition
by Corollary \ref{corwss}.
\qed

\begin{df}\label{dfncd}
Let $f\colon
X\to S$ be a weakly
semi-stable scheme over $S$.

{\rm 1.}
Let 
$D=D_1+\cdots+D_n$
be the sum
of Cartier divisors of $X$.
Then,
we say that $D$ 
has {\rm simple 
normal crossings 
relatively to} $S$
if the following condition
is satisfied:
\begin{itemize}
\item[
{\rm (\ref{dfncd}.1)}]
For every point $x\in X$,
there exist
an open neighborhood
$x\in U\subset X$,
a weakly
semi-stable scheme $Y$
over $S$
and a smooth morphism
$U\to {\mathbf A}^m_Y$
to the affine space 
with coordinate $T_1,\ldots,T_m$
such that,
for each
$i=1,\ldots,n$,
the restriction 
$D_i\times_XU$ is
either empty or
defined by
$T_{j_i}$
for some
$1\le j_i\le m$.
Further, for $1\le i<i'\le n$
such that
$D_i\times_XU$ 
and $D_{i'}\times_XU$ 
are non-empty,
we have $j_i\neq j_{i'}$.
\end{itemize}

{\rm 2.}
Let $D$ be
a Cartier divisor of $X$.
Then,
we say that $D$ 
has {\rm normal crossings 
relatively to} $S$
if, \'etale locally on $X$,
it has simple 
normal crossings
relatively to $S$.
\end{df}

If a Cartier divisor of $X$
has {\rm normal crossings 
relatively to} $S$,
it is flat over $S$.
If $D=D_1+\cdots+D_n$
is a divisor 
with simple normal
crossings
relatively to $S$,
for a subset
$I\subset \{1,\ldots,n\}$,
the intersection
$D_I=
\bigcap_{i\in I}D_i$
is weakly strictly
semi-stable over $S$.
If $X$ is smooth over $S$,
the terminology
on simple normal crossing
divisors
is the same as the usual one
defined in \cite[2.1]{SGA1}.

We recall the following
fact on the tameness
of the direct image
for a proper semi-stable scheme.

\begin{lm}\label{lmsstm}
Let $S$ be a regular noetherian
scheme and $D\subset S$
be a divisor with
normal crossings.
Let $f\colon X\to S$
be a proper weakly semi-stable
scheme
such that
the base change
$X\times_SW\to W=
S\setminus D$
is smooth
and $E\subset X$ be a 
divisor with normal crossings
relatively to $S$.
We put $U=X\setminus E$
and $f_U\colon U\to S$
be the restriction of
$f$.

Then, for an integer $n\ge 1$
invertible on $S$,
the higher direct image
$R^qf_{U!}{\mathbb Z}/
n{\mathbb Z}$
is locally constant on 
$W=S\setminus D$
and is tamely ramified
along
$D$ for every $q\ge 0$.
\end{lm}

{\it Proof.}
By the assumption
that $S$ is regular
and $D$ has normal crossings,
it is reduced to
the case where
$S={\rm Spec}\ {\cal O}_K$
for a discrete valuation
ring and $D$ consists
of the closed point $s$,
by Abhyankhar's lemma
\cite[Proposition 5.2]{SGA1}.
Let $j\colon U\to X$
denote the open immersion.
Then, it suffices
to show that
the action of
the inertia group
$I={\rm Gal}(K^{\rm sep}/
K^{\rm ur})$ on
the sheaf $R^q\psi
j_!{\mathbb Z}/n{\mathbb Z}$
of nearby cycles
is tamely ramified.
If $E=\emptyset$,
it is proved in \cite{RZ}.

We show the general case.
Since the assertion is \'etale local on $X$,
we may assume
$X$ is weakly strictly semi-stable over $S$. 
Let $j\colon U\to X$ denote
the open immersion
and for a finite set $I$
of indices 
of irreducible components
$E_i$ of $E$,
let $i_I\colon E_I\to X$
be the closed immersion
of the intersections.
Then, $E_I$ are semi-stable
over $S$ and we
have an exact sequence
$0\to j_!{\mathbb Z}/n{\mathbb Z}
\to
{\mathbb Z}/n{\mathbb Z}
\to
\bigoplus_{\sharp I=1}i_{I*}{\mathbb Z}/n{\mathbb Z}
\to
\bigoplus_{\sharp I=2}i_{I*}
\pi_*{\mathbb Z}/n{\mathbb Z}
\to \cdots
\to$.
Using this,
the assertion 
is reduced to
the case
$E=\emptyset$.
\qed

We recall the
definition
of a stable curve
\cite{DM}.
Let
$f\colon X\to S$
be a proper
weakly semi-stable scheme
of relative dimension 1
over a scheme $S$
and 
$(s_i)_{
i=1,\ldots,d}$
be a finite family
of sections
$s_i\colon S\to X$.
Let $\omega_{X/S}=
R^{-1}f^!{\cal O}_X$
be the relative
dualizing sheaf.
Then, we say a pair
$(f\colon X\to S,
(s_i)_{i=1,\ldots,d})$
is a $d$ pointed
stable curve
if the following
condition
is satisfied.
\begin{itemize}
\item
The divisor
$D=\sum_{i=1}^d
s_i(S)$
has 
simple normal crossings
relatively to $S$,
the canonical map
${\cal O}_S\to
f_*{\cal O}_X$
is an isomorphism
and
the invertible
${\cal O}_X$-module
$\omega_{X/S}(D)$
is relatively ample.
\end{itemize}
If $(X,(s_i))$ 
is a pointed stable curve
over $S$,
the sections
$s_i(S)$ do not
meet each other
and are contained
in the locus where
$f$ is smooth.
Further the 
${\cal O}_S$-module
$f_*\omega_{X/S}$
is locally free.
The rank of
$f_*\omega_{X/S}$
is called the genus of $X$.
If $(X,(s_i))$
is a $d$ pointed
stable curve of genus $g$,
we have
$2g-2+d>0$.

We recall some facts
on the moduli of pointed
stable curves,
used in the proof of
the conductor formula
for a relative curve
in Proposition \ref{prcfs}
and Corollary \ref{corcfs}.
Let $\bar S=
\bar {\cal M}_{g,d}$
be the moduli stack
of $d$ pointed stable curves
of genus $g$.
It is a proper smooth
Deligne-Mumford stack
over ${\mathbb Z}$
\cite{KM}
and the coarse moduli
scheme 
$\bar M_{g,d}$
is a projective scheme \cite{GIT}.
Let $f\colon X\to \bar S$
be the universal family
and $s_1,\ldots,s_d\colon
\bar S\to X$
be the universal sections.
Let
$S={\cal M}_{g,d}
\subset \bar S$
be the open substack
where $X$ is smooth.
It is the complement
of a divisor with normal crossings
\cite{KM}.

Let $n\ge1$ be
an integer.
The $n$-torsion part
${\rm Jac}_{X[\frac1n]
/S[\frac1n]}[n]
=R^1f_{S[\frac 1n]*}
\mu_n$
of the Jacobian
is a locally constant sheaf
of ${\mathbb Z}/n{\mathbb Z}$-modules
of rank $2g$ on $S[\frac 1n]$.
Let 
${\cal M}_{g,d,n}$
over
${\cal M}_{g,d}[\frac 1n]=
S[\frac 1n]$
be the moduli of
an isomorphism
$({\mathbb Z}/n{\mathbb Z})^{2g}
\to R^1f_{S[\frac 1n]*}
\mu_n$.
If $n\ge 3$,
then ${\cal M}_{g,d,n}$
is represented
by a scheme $S_n=M_{g,d,n}$
smooth over ${\mathbb Z}[\frac1n]$.
Further,
the normalization
$\bar S_n=
\bar{\cal M}_{g,d,n}$
of
$\bar{\cal M}_{g,d}[\frac 1n]$
in $S_n=M_{g,d,n}$
is a projective
scheme over
${\mathbb Z}[\frac 1n]$
\cite{Gab}.
See also \cite[2.24]{dJI}.

\subsection{Semi-stable schemes
and log blow-up}\label{sssb}

We briefly
recall the log blow-up
and apply it
to give some constructions
related to semi-stable schemes.
For terminologies on log blow-up,
we refer to
\cite[Section 4.2]{KSI}.
Let $P$ be a finitely 
generated commutative integral
saturated torsion free monoid,
called a torsion free fs-monoid for short.
In other words,
the associated group
$P^{\rm gp}$ is
a finitely generated
free abelian group
and there exists
a finitely many
elements
$f_1,\ldots,f_m$
of the dual group 
$P^{{\rm gp}* }=
{\rm Hom}(P^{\rm gp},{\mathbb Z})$
such that
$P$
is identified with the submonoid
$\{x\in P^{\rm gp}\mid
f_i(x)\ge 0$ for $i=1,\ldots,m\}
\subset P^{\rm gp}$,
see \cite[Proposition 1.1]{O}.
We identify
the dual monoid
$P^*=
{\rm Hom}_{\rm monoid}(P,{\mathbb N})$
with the submonoid
$\{f\in P^{{\rm gp}* }\mid
f(x)\ge 0$ for $x\in P\}$.
If $P^\times=
\{x\in P\mid x^{-1}\in P\}$
is trivial,
the abelian group
$P^{{\rm gp}* }$
is generated
by the submonoid $P^*$.
Further
in this case,
$P^*$
is the intersection
of $P^{{\rm gp}* }$
in $P^{{\rm gp}* }
\otimes_{\mathbb Z}
{\mathbb Q}$
with 
$\{a_1f_1+
\cdots+
a_mf_m\mid
a_i\in 
{\mathbb Q},
a_i\ge 0\}$
(loc.\ cit.).

Let $X$ be a log scheme
and 
$f\colon P\to
\Gamma(X,{\cal O}_X)$
be a chart.
It defines a
strict morphism
$f^*\colon
X\to {\rm Spec}\ {\mathbb Z}[P]$
of log schemes.
Recall that a morphism
$X\to Y$
is strict if
the log structure
of $X$ is the pull-back
of that of $Y$.
Let $V_\sigma\subset
P^{{\rm gp}* }
\otimes_{\mathbb Z}
{\mathbb Q}$
be a ${\mathbb Q}$-linear
subspace.
Then the intersection
$N_\sigma
=P^*\cap V_\sigma$
is a finitely generated
saturated submonoid.
Let $P_\sigma
\supset P$
be the 
finitely generated
saturated monoid
defined by
$\{x\in P^{\rm gp}\mid
f(x)\ge 0$ for $f\in N_\sigma\}$.
Then, we define
a scheme $X_\sigma$
by
$X_\sigma
=X\times_
{{\rm Spec}\ 
{\mathbb Z}[P]}
{\rm Spec}\ 
{\mathbb Z}[P_\sigma]$.
Let $g\colon
P\to \Gamma(X,{\cal O}_X)$
be another chart
such that
there exists
a morphism
$u\colon
P\to \Gamma(X,{\cal O}_X^\times)$
satisfying
$g=f\cdot u$.
Then, the schemes
$X_\sigma$ over $X$
defined by $f$
and by $g$
are canonically
isomorphic to each other.

Let $\Sigma$ be
a subdivision
of the dual monoid
$P^*$.
Recall that 
a subdivision
$\Sigma$ consists
of finite family
of submonoids
$N_\sigma=
P^*\cap V_\sigma$ of
the dual monoid $P^*$
indexed by $\sigma\in \Sigma$.
Recall also that
$\Sigma$ is regular
means that the monoid
$N_\sigma\subset P^*$ for every
$\sigma\in \Sigma$
is isomorphic to ${\mathbb N}^r$
for some $r\ge 0$
and hence
$P_\sigma$
is isomorphic
to ${\mathbb N}^r
\times
{\mathbb Z}^{n-r}$
where $n$ is the rank of
$P^{\rm gp}$.
By patching
the schemes
$X_\sigma$ over $X$,
we obtain a scheme
$X_\Sigma$ over $X$.
Recall that
if $\Sigma$
is a proper subdivision,
the scheme $X_\Sigma$
is proper over $X$.
In this case,
we call $X_\Sigma$
a log blow-up of $X$.

Let $S$ be a regular noetherian
scheme and $D\subset S$
be a divisor with
normal crossings.
Let $j_W\colon
W=S\setminus D
\to S$ denote
the open immersion
and we regard $S$
as a log scheme
defined by 
the log structure
${\cal M}_S=
{\cal O}_S\cap j_{W*}
{\cal O}_W^\times$.
We consider a weakly semi-stable
scheme $f\colon X\to S$
and a divisor $E\subset X$
with normal crossings
relatively to $S$
such that
the base change
$X\times_SW\to W$
is smooth.
Let $j_U\colon
U=X\setminus 
(f^{-1}(D)\cup E)
\to X$ denote
the open immersion
and we regard $X$
as a log scheme
defined by 
the log structure
${\cal M}_X=
{\cal O}_X\cap j_{U*}
{\cal O}_U^\times$.
Then the map
$f\colon X\to S$
is log smooth.

We construct
proper modifications
of weakly semi-stable
schemes using log blow-ups.
This will be used
at the end of the proof
of Theorem \ref{thmlLTF}.

\begin{lm}\label{lmsOK}
Let
${\cal O}_K$ be
a discrete valuation ring 
and $X$ be
a weakly semi-stable scheme
over $S={\rm Spec}\ {\cal O}_K$
with smooth generic fiber $X_K$.
Then, there exists
a proper modification
$X'\to X$
such that
$X'_K\to X_K$
is an isomorphism and that
$X'$ is semi-stable
over $S$.
\end{lm}

{\it Proof.}
Let $\pi$ be
a prime element of ${\cal O}_K$.
First, we consider
the case where
there exists an 
\'etale morphism
$X\to {\rm Spec}\ 
{\cal O}_K
[T_1,\ldots,T_{r+1}]/
(T_1\cdots T_q-
\pi^e)$
for an integer $e\ge 1$.
Let $P_{q,e}$
be the monoid
${\mathbb N}^q+
\langle (\frac 1e,
\ldots,\frac1e)
\rangle\subset
{\mathbb Q}^q$.
The uniformizer
$\pi$ and
the pull-backs
of $T_1,
\ldots, T_q$
define
a morphism
$X\to
{\rm Spec}\ {\mathbb Z}[P_{q,e}]
=
{\rm Spec}\ {\mathbb Z}
[T_1,\ldots,T_q,S]/
(T_1\cdots T_q-S^e)$.

We identify the dual monoid
$N_{q,e}=P^*_{q,e}$
with 
$\{(a_1,\ldots,a_q)
\in {\mathbb N}^q\mid
a_1+\cdots+a_q
\equiv 0\bmod e\}$.
Let $B_{q,e}
\subset N_{q,e}$
be the finite set
$\{
(a_1,\ldots,a_q)
\in {\mathbb N}^q\mid
a_1+\cdots+a_q=e\}$
and define $\Sigma_{q,e}$
by
$\{\sigma\subset
B_{q,e}\mid
(a_1,\ldots,a_q),
(b_1,\ldots,b_q)
\in \sigma$
implies
$|a_1-b_1|+
\cdots+
|a_q-b_q|\le 2\}$.
For 
$\sigma\in \Sigma_{q,e}$,
let $N_\sigma$
denote the submonoid
of $N_{q,e}$
generated by
$\sigma$.
Then,
$\Sigma_{q,e}$
defines
a regular proper subdivision
of $N_{q,e}$
and we obtain
a log blow-up
$X_{\Sigma_{q,e}}\to X$.

We show 
that the scheme
$X_{\Sigma_{q,e}}$
is semi-stable
over ${\cal O}_K$.
Let $f\colon
{\mathbb N}
\to P_{q,e}$
be the map sending
$1$ to 
$\langle \frac 1e,
\ldots, \frac 1e\rangle$.
Then, the dual 
map $f^*\colon
N_{q,e}\to
{\mathbb N}$
sends 
$a\in N_{q,e}$ to
$(a_1+\cdots+a_q)/e$
and hence
an arbitrary element
of $B_{q,e}$ to $1$.
A numbering
on $\sigma
\in \Sigma_{q,e}$
defines an isomorphism
${\mathbb N}^s
\to N_\sigma
\subset
N_{q,e}$.
Hence,
the composition
${\mathbb N}^s
\to {\mathbb N}$
of the restriction
$f^*|_{N_\sigma}\colon
N_\sigma\to {\mathbb N}$
with an isomorphism
${\mathbb N}^s
\to N_\sigma$
sends every member
of the canonical basis
of ${\mathbb N}^s$
to $1$.
From this, it follows
immediately 
that $X_\sigma$
is semi-stable
over ${\cal O}_K$
for every $\sigma
\in {\Sigma_{q,e}}$.

By Lemma \ref{lmsst},
the exponent $e$
and the divisors
defined by
$T_1,\ldots,T_q$
are well-defined
\'etale locally
up to permutation.
Since
the regular proper subdivision
$\Sigma_{q,e}$
is invariant
under permutations
of $q$ letters,
the \'etale locally constructed
log blow-ups
$X_{\Sigma_{q,e}}\to X$
patch each other
and define
a semi-stable
modification
$X'\to X$ globally.
\qed

Next, we
reformulate
\cite[Proposition 3.6]{dJI}
in our terminology.
This together with Lemma \ref{lmD1} will be used in the proof of Lemma \ref{lmcysd}
in the case (\ref{lmcysd}.1a).

\begin{lm}\label{lmlbS}
Let $S$ be a regular
noetherian scheme and
$D\subset S$ be a
divisor with simple normal
crossings.
Let
$f\colon
X\to S$ be a 
weakly strictly semi-stable
curve
such that 
the base change
$X_W=X\times_SW \to W
=S\setminus D$
is smooth.

Then, there exists
a proper modification
$X'\to X$
such that
$X'_W\to X_W$
is an isomorphism, that
$X'$ is regular
and weakly strictly semi-stable
over $S$
and that $X'\times_SD$
is a divisor with
simple normal crossings.
\end{lm}

{\it Proof.}
First, we consider
the case where the following
data are given:

Let ${\rm Spec}\ R
\subset S$ be
an affine open
subscheme,
$s_1,\ldots,s_n\in R$
be elements
defining irreducible
components
$D_1,\ldots,D_n$
of $D\cap {\rm Spec}\ R$,
$d_1,\ldots,d_n>0$
be integers
and let 
$X\to
{\rm Spec}\ R[T_1,T_2]/
(T_1T_2-
s_1^{d_1}\cdots s_n^{d_n})$
be an \'etale morphism
over $S$.
Let 
$c\colon 
\{1,\ldots,n\}
\to \{1,2\}$
be a function.

We define maps
${\mathbb N}\to
{\mathbb N}^2,
{\mathbb N}\to
{\mathbb N}^n$
of monoids by
$(1,1)$ and $(d_1,\ldots,d_n)$
and consider
the amalgamate
sum $P={\mathbb N}^2+
_{\mathbb N}
{\mathbb N}^n$.
The dual
$N=P^*$
is identified with
$\{(a,b)\in 
{\mathbb N}^2
\times
{\mathbb N}^n
\mid
a_1+a_2=d_1b_1+
\cdots +d_nb_n\}$.
Let
$e_1,\ldots,e_n
\in {\mathbb N}^n$
be the standard
basis.
For $i=1,\ldots,n$,
we put
$B_i=
\{(a,b)\in N\mid
b=e_i\}$.
We identify
$(1,i),(2,i)\in 
A=\{1,2\}
\times\{1,\ldots,n\}$
with
$((d_i,0),e_i),
((0,d_i),e_i)
\in B_i$
and regard $A$
as a subset of 
$B=\coprod_iB_i
\subset N$.
For each
$j\in \{1,\ldots n\}$,
let $\Sigma_j$ be
the finite set
consisting of
$\sigma\subset B_j\cup A$
satisfying 
the following
conditions:
\begin{itemize}
\item
If $(a,i)\in \sigma\cap A$
for $i<j$, we have
$a=c(i)$.
\item
If $((a_1,a_2),e_j),
((a'_1,a'_2),e_j)\in \sigma
\cap B_j$,
we have
$|a_1-a_1'|\le 1$.
\item
We have
$\{a\in \{1,2\}\mid
(a,i)\in \sigma\cap A,
i>j\}\subsetneqq\{1,2\}$.
\end{itemize}
We put $\Sigma=\bigcup_{j=1}^n
\Sigma_j$.
For each $\sigma
\in \Sigma$,
the submonoid $N_\sigma
\subset N$ generated by $\sigma$
is isomorphic to 
${\mathbb N}^s$
for $s={\rm Card}\ \sigma\ge 0$.
For $(a,b)\in N$,
if there exists 
an integer $1\le j\le n$
not satisfying the inequalities
$$
a_1\ge 
\sum_{i\le j,\ c(i)=1}b_id_i
\quad
\text{ and }
\quad
a_2\ge 
\sum_{i\le j,\ c(i)=2}b_id_i,$$
then, for the smallest such $j$,
there exists 
$\sigma\in \Sigma_j$
such that $(a,b)\in N_\sigma$.
If otherwise,
we have $(a,b)\in N_C$
for $C=\{(c(i),i)\mid
i\in \{1,\ldots,n\}\}
\in \Sigma_n$.
Hence,
$\Sigma$
defines a regular
proper
subdivision
of $N$.

The \'etale morphism
$X\to
{\rm Spec}\ R[T_1,T_2]/
(T_1T_2-
s_1^{d_1}\cdots s_n^{d_n})$
induces a morphism
$X\to 
{\rm Spec}\ {\mathbb Z}
[T_1,T_2,S_1,\ldots,S_n]/
(T_1T_2-
S_1^{d_1}\cdots S_n^{d_n})
=
{\rm Spec}\
{\mathbb Z}[P]$.
Hence the log blow-up
$X_\Sigma$
is defined by the 
regular proper
subdivision $\Sigma$.

We show that
the scheme $X_\Sigma$
satisfies the condition.
We consider
the dual
$\pi^*\colon
N\to {\mathbb N}^n$
of the canonical map
$\pi\colon
{\mathbb N}^n
\to P=
{\mathbb N}^2
+_{\mathbb N}
{\mathbb N}^n$.
Let $e_1,\ldots,e_s$
and $e'_1,\ldots,e'_n$
be standard bases
of ${\mathbb N}^s$
and 
of ${\mathbb N}^n$.
Then,
$\pi^*$
maps the elements
of $B_i$
to $e'_i
\in {\mathbb N}^n$.
Let $\sigma
\in \Sigma_i$
and 
take an isomorphism
${\mathbb N}^s
\to N_\sigma$
to the submonoid
generated by $\sigma$.
We consider 
the composition
$\varphi
\colon
{\mathbb N}^s
\to {\mathbb N}^n$
with the restriction
$N_\sigma\to {\mathbb N}^n$.
Then, there exists
a map
$g\colon 
\{1,\ldots,s\}
\to
\{1,\ldots,n\}$
such that 
$\varphi(e_j)=
e_{g(j)}$ for
$j=1,\ldots,s$.
Further,
for $i'=1,\ldots,n$,
we have
${\rm Card}(g^{-1}(i'))
\le 1$ for $i\neq i'$
and
${\rm Card}(g^{-1}(i))
\le 2$.
Thus,
we have either
an \'etale map
$X_\sigma
\to {\rm Spec}\ R
[T_1,T_2]/(T_1T_2-s_i)$
or an \'etale map
$X_\sigma
\to {\rm Spec}\ R[T]$.
Hence, 
the log blow-up
$X_\Sigma$
is weakly strictly
semi-stable
over $S$ and
regular.
Further 
$D\times_SX_\Sigma$
is a
divisor with
simple normal crossings.

We prove
the general case.
To patch
the local construction
above,
we fix a numbering
of irreducible components
of $E=X\times_SD$.
Let ${\rm Spec}\ R
\subset S$
be an affine open
and $V\to
{\rm Spec}\ R[T_1,T_2]/
(T_1T_2-
s_1^{d_1}\cdots s_n^{d_n})$
be an \'etale map
defined on an open 
subscheme $V$ of $X$.
We assume that
each $V\times_SD_i$
has two irreducible 
components
$E_{1,i}$ and $E_{2,i}$
defined by $(T_1,s_i)$
and $(T_2,s_i)$
respectively.
We define a function
$c\colon \{1,\ldots,n\}
\to \{1,2\}$
by requiring that
the index of the 
irreducible component
$E_{c(i),i}$
is the smaller
among 
$E_{1,i}$ and $E_{2,i}$
with respect to
the fixed numbering
of the irreducible components
of $E$.
By changing the numbering
of $D_1,\ldots,D_n$,
we may assume that
the indices
of the sequence
$E_{c(1),1},
\ldots,
E_{c(n),n}$
is increasing.
With this numbering
and the definition
of $c$,
it is easily seen
that the log blow-ups
$V_\Sigma$
patch globally
and define
a modification $X'\to X$.
\qed

The following Lemma
will be used
in the proof of
Corollary \ref{lmcysd}
the case  (\ref{lmcysd}.1b)
but {\em not} 
in the proof of
the conductor formula.

\begin{lm}\label{lmlbsn}
Let $S$ be a regular
noetherian scheme and
$D\subset S$ be a
divisor with simple normal
crossings.
Let $f\colon
X\to S$ be a 
weakly strictly semi-stable
scheme
such that 
the base change
$X_W=X\times_SW \to W
=S\setminus D$
is smooth.

For an irreducible component
$D_i$ of $D$,
let $I_i$
be the set of irreducible
components
of $X\times_SD_i$
and,
for $x\in X$
and $s=f(x)\in S$,
let $I_x$ be the set
of irreducible components
of the fiber $X_s$
containing $x$.
We assume
that the following condition
is satisfied:
\begin{itemize}
\item[{\rm (\ref{lmlbsn}.1)}]
There exist
a family of
functions
$\varphi_i\colon
I_i\to {\mathbb N}$
and a total order
on the finite set $I_x$ for
every $x\in X$
satisfying
the following condition:
If $s=f(x)\in D_i$
and if the map
$I_x\to I_i$
induced by
the inclusion
$X_s\to X_{D_i}$
is injective,
then the composition
$I_x\to {\mathbb N}$
with $\varphi_i$
is injective
and increasing.
\end{itemize}

Then, there exists
a proper modification
$X'\to X$
such that
$X'_W\to X_W$
is an isomorphism, that
$X'$ is regular
and weakly strictly semi-stable
over $S$
and that $X'\times_SD$
is a divisor with
simple normal crossings.
\end{lm}

{\it Proof.}
First, we consider
the case
where 
the irreducible
components
$D_1,\ldots,D_q$ of $D$
are defined by
$t_1,\ldots,t_q$,
there exists
a smooth map
$X\to S[T_1,\ldots,T_r]/
(T_1\cdots T_r-
t_1^{m_1}\cdots
t_q^{m_q})$
for integers
$m_1,\ldots,m_q\ge 0$
and
the total order
on $I_x$
is induced by
the natural order
on $\{1,\ldots,r\}$.
We define morphisms
${\mathbb N}\to {\mathbb N}^r$
and
$e\colon
{\mathbb N}\to {\mathbb N}^q$
by $1\mapsto (1,\ldots,1)$
and
by $1\mapsto (m_1,\ldots,m_q)$
respectively.
Let $P$ be the amalgamete sum
${\mathbb N}^r+_{\mathbb N}
{\mathbb N}^q$
with respect to the morphisms
above.
We consider the map
$P\to
\Gamma(X,{\cal O}_X)$
of monoids
defined by
$T_1,\ldots, T_r$
and
$t_1,\ldots, t_q$.

We define the dual morphisms
$|\ |\colon
{\mathbb N}^r\to {\mathbb N}$
and
$m^*\colon
{\mathbb N}^q\to {\mathbb N}$
by $|(a_1,\ldots,a_r)|
=a_1+\cdots+a_r$
and
$m^*(b_1,\ldots,b_q)
=m_1b_1+\cdots+m_qb_q$.
Then, the dual monoid
$N=P^*$
is identified with
$\{(a,b)\in {\mathbb N}^r
\times {\mathbb N}^q\mid
|a|=m^*(b)\}$.
We define a regular
proper subdivision $\Sigma$
of $N$.
Let $V$ be the finite set
$\{(a,j)
\in {\mathbb N}^r
\times \{1,\ldots,q\}
\mid |a|=m_j\}$.
We regard
$V$ as a subset of $N$
by identifying $(a,j)\in V$
with $f_{a,j}=(a,f_j)\in N$
where $f_1,\ldots,f_q$
denote 
the canonical basis
of ${\mathbb N}^q$.
For a vector
$a\in {\mathbb N}^r$,
we put ${\rm Supp}(a)
=\{i\in \{1,\ldots,r\}
\mid a_i>0\}$.
For elements
$(a,j), (a',j')\in V$,
we write
$(a,j)\le (a',j')$
if $\max {\rm Supp}(a)
\le \min {\rm Supp}(a')$
and $j\le j'$.
The relation
$\le$ satisfies
the anti-symmetry law
and the transitivity law
but not the reflexive law.
By abuse of terminology,
we say a subset
$\sigma\subset V$
is totally ordered
if  $(a,j),(a',j')
\in \sigma$
implies
either
$(a,j)\le (a',j')$,
$(a',j')\le (a,j)$
or
$(a,j)=(a',j')$.
We put $\Sigma=
\{\sigma\subset V\mid
\sigma$
is totally ordered
$\}$.
For $\sigma\in \Sigma$,
we consider 
the submonoid
$N_\sigma
\subset N$
generated by
$f_{a,j}$ for $(a,j)\in \sigma$.
For each $(a,b)\in N$,
one can easily find
the minimum totally
ordered subset
$\sigma\in \Sigma$
satisfying $(a,b)\in N_\sigma$.
Thus,
$\Sigma$ defines
a regular proper subdivision.
Hence the log blow-up
$X'=X_\Sigma$
is regular
and 
$X'\times_SD$
is a divisor with
simple normal crossings.

By the assumption
on the existence of
the functions
and the total orders,
the log blow-ups constructed
above patch globally
to give the required $X'$.
\qed

\subsection{Log products
and log blow-ups}
\label{sslb}
\setcounter{equation}0
We fix some terminology
and notation
on log products,
which will be constantly
used throughout this paper.
For the generality on
log schemes,
we refer to
\cite{kk-log},
\cite{IlV},
\cite[Section 4]{KSI}.
In this paper,
unless otherwise explicitly stated,
a log structure means an fs-log
structure defined Zariski locally.
In particular,
a log structure ${\mathcal M}_X$
is a sheaf of commutative monoids
on the Zariski site of a scheme $X$
endowed with a morphism of
sheaf of monoids
${\mathcal M}_X\to 
{\mathcal O}_X$
where ${\mathcal O}_X$
is regarded as a sheaf of monoids
with respect to the multiplication.
Further, Zariski locally on $X$,
the log structure
${\mathcal M}_X$ admits
a chart by an fs-monoid.
For a log structure
${\mathcal M}_X$,
let
$\bar {\mathcal M}_X$
denote the quotient
${\mathcal M}_X/
{\mathcal O}^\times_X$.

We recall some basic facts
on log schemes
from \cite{kk-log},
\cite[Section 4.3]{KSI}.
For morphisms $X\to S$
and $Y\to S$ of log schemes,
the fiber product
$X\times_S^{\log} Y$
is defined as a log scheme.
Note that $X,Y,S$ are
assumed to be fs-log schemes
and 
$X\times_S^{\log} Y$
is the fiber product
in the category of fs-log schemes.
We put $\log$ in the notation
to indicate that the underlying
scheme can be different
from 
$X\times_S Y$
in the category is schemes.
However, for example
if at least one of the morphisms
$X\to S$
and $Y\to S$
is strict,
the underlying scheme is
$X\times_S Y$.
In such a case,
we will drop $\log$ in the notation.

For a Cartier divisor $D$ of
a scheme $X$
defined by the ideal sheaf
${\mathcal I}_D
\subset 
{\mathcal O}_X$,
the associated log structure
is defined to be 
${\mathcal M}_X=
\coprod_{n\in {\mathbb N}}
{\cal I}som_{{\mathcal O}_X}
({\mathcal O}_X,
{\mathcal I}^n_D)$
endowed with the morphism
${\mathcal M}_X\to
{\mathcal O}_X$
induced by
the inclusions
${\mathcal I}^n_D\to
{\mathcal O}_X$.
For a finite family ${\cal D}=
(D_i)_{i\in I}$
of Cartier divisors
$D_i\subset X$,
the associated log structure
is similarly defined by
${\mathcal M}_X=
\coprod_{n\in {\mathbb N}^I}
{\cal I}som_{{\mathcal O}_X}
({\mathcal O}_X,
\prod_{i\in I}{\mathcal I}^{n_i}_{D_i})$.
It is the push-out 
of the log structures
defined by $D_i$ for $i\in I$.
We have a canonical map
${\mathbb N}^I\to 
\Gamma(X,\bar {\mathcal M}_X)$
that can be lifted to
a chart locally on $X$.

Let $P$ be an fs-monoid
and we consider two morphisms 
$P\to \Gamma(X,\bar 
{\mathcal M}_X)$ of monoids.
Then, by
applying  \cite[Proposition 4.2.3]{KSI}
to the surjection $P+P\to P$,
we conclude that 
the functor sending a log scheme $T$
to the set
$$\{
f\colon T\to X\mid
\text{the two compositions
$P\to \Gamma(X,\bar 
{\mathcal M}_X)
\overset{f^*}
\to \Gamma(T,\bar 
{\mathcal M}_T)$
are equal to each other}\}$$
of morphisms of log schemes
is representable
by a log \'etale scheme over
$X$,
that may be denoted by 
$X\times^{\log}_{X,P}X$.
Locally on $X$,
it is constructed as follows.
Let $\widetilde P$
be the inverse image of
$P$ by the sum
$P^{\rm gp}\oplus P^{\rm gp}
\to P^{\rm gp}$.
Locally on $X$,
we take liftings 
$P\to \Gamma(X,{\mathcal M}_X)$
of
$P\to \Gamma(X,\bar {\mathcal M}_X)$
and let 
$X\to
{\rm Spec}\ {\mathbb Z}[P+P]$
be the induced morphism
of log schemes.
Then, 
$X\times^{\log}_{X,P}X$
is constructed as
$X\times^{\log}_{
{\rm Spec}\ {\mathbb Z}[P+P]}
{\rm Spec}\ {\mathbb Z}[\widetilde P]$.

We apply the construction
in the following case.
Let $X\to S$ and $Y\to S$
be morphisms of log schemes,
$P$ be an fs-monoid
and 
$P\to \Gamma(X,\bar 
{\mathcal M}_X)$ and
$P\to \Gamma(Y,\bar 
{\mathcal M}_Y)$ 
be morphisms of monoids.
Then, 
they induces two morphisms
$P\to \Gamma(X\times^{\log}_SY,\bar 
{\mathcal M}_{X\times^{\log}_SY})$.
By applying the construction above,
we define the log product
$X\times^{\log}_{S,P}Y$.
It represents
the functor sending a log scheme $T$
over $S$ to the set
$$\left\{
(f\colon T\to X,
g\colon T\to Y)\left|
\text{the diagram
$\begin{CD}
P@>>> \Gamma(X,\bar 
{\mathcal M}_X)\\
@VVV @VV{f^*}V\\
\Gamma(Y,\bar 
{\mathcal M}_Y)
@>{g^*}>>
\Gamma(T,\bar 
{\mathcal M}_T)
\end{CD}$
is commutative}\right.\right\}$$
of pairs of morphisms of log schemes
over $S$.

\begin{lm}\label{lmlp}
If $P\to \Gamma(X,\bar 
{\mathcal M}_X)$ and
$P\to \Gamma(Y,\bar 
{\mathcal M}_Y)$ 
are locally lifted to charts,
then the projections
$X\times^{\log}_{S,P}Y\to X$ and
$X\times^{\log}_{S,P}Y
\to Y$ are strict morphisms
of log schemes.
\end{lm}

{\it Proof.}
By the construction above,
the monoid $P$
also defines charts on
$X\times^{\log}_{S,P}Y$.
\qed

We consider the following variant.
Let $X\to S$ and $Y\to S$
be morphisms of log schemes,
$P$ and $Q$ be fs-monoids
and 
$$
\begin{CD}
P@<<<Q@>>>P\\
@VVV@VVV@VVV\\
\Gamma(X,\bar 
{\mathcal M}_X)
@<<<
\Gamma(S,\bar 
{\mathcal M}_S)
@>>>
\Gamma(Y,\bar 
{\mathcal M}_Y)
\end{CD}$$ 
be a commutative diagram
of morphisms of monoids.
Then, we define
 the log product
$X\times^{\log}_{S,P/Q}Y$
by the cartesian diagram
\begin{equation}
\begin{CD}
X\times^{\log}_{S,P/Q}Y@>>>
X\times^{\log}_{S,P}Y
\\
@VVV @VVV\\
S@>>>S\times^{\log}_{S,Q}S
\end{CD}
\label{eqllp}
\end{equation}
where the bottom arrow
is the diagonal map
and the right vertical arrow
is defined by functoriality.

We make the construction 
explicit in the case where
the log structures of
$X$ and $Y$ are defined by 
finite families ${\cal D}=
(D_i)_{i\in I}$ and
${\cal E}=
(E_i)_{i\in I}$
of Cartier divisors
$D_i\subset X$ and of
$E_i\subset Y$
with the same index set
and the log structure of $S$
is trivial.
We define the log product
\begin{equation}
(X\times_SY)^\sim_
{{\cal D},{\cal E}}
\label{eqXYDE}
\end{equation}
to be
$X\times^{\log}_{S,
{\mathbb N}^I}Y$
defined by the canonical morphisms
${\mathbb N}^I
\to \Gamma(X,\bar 
{\mathcal M}_X),
{\mathbb N}^I\to \Gamma(Y,\bar 
{\mathcal M}_Y)$.
The canonical morphism
$(X\times_SY)^\sim_
{{\cal D},{\cal E}}
\to X\times_SY$ is log \'etale. 
If $X=Y$ and
${\cal D}={\cal E}$,
we let
$(X\times_SY)^\sim_
{{\cal D},{\cal E}}$
denoted by
$(X\times_SX)^\sim_
{\cal D}$.
Further if ${\cal D}$
is clear from the context,
we drop the subscript
${\cal D}$.

Locally, the log product
$(X\times_SY)^\sim_
{{\cal D},{\cal E}}$
is described as follows.
Assume that
$D_i$ and $E_i$
are defined by
$f_i\in\Gamma(X,{\cal O}_X)$
and
$g_i\in\Gamma(Y,{\cal O}_Y)$
respectively.
Then, $(f_i)_{i\in I}$
and $(g_i)_{i\in I}$
define
maps of monoids
${\mathbb N}^I
\to \Gamma(X,{\cal O}_X)$
and
${\mathbb N}^I
\to \Gamma(Y,{\cal O}_Y)$
and they further induce
a map
$P={\mathbb N}^I
\times {\mathbb N}^I
\to \Gamma(X\times_SY,
{\cal O}_{X\times_SY})$
from the direct sum.
We identify the
dual monoid $N=P^*$
with 
${\mathbb N}^I
\times {\mathbb N}^I$
and let
$N_\sigma=
{\mathbb N}^I
\subset {\mathbb N}^I
\times {\mathbb N}^I$
be the diagonal submonoid.
Then, the corresponding
submonoid
$P_{\sigma}
=\{p\in 
P^{\rm gp}\mid
f(p)\in {\mathbb N}
\text{ for }
f\in N_\sigma\}
\subset 
P^{\rm gp}
={\mathbb Z}^I
\times 
{\mathbb Z}^I$
is equal to
$\{(a,b)\in 
{\mathbb Z}^I
\times 
{\mathbb Z}^I\mid
a+b
\in {\mathbb N}^I\}$.
The log product
$(X\times_SY)^\sim_
{{\cal D},{\cal E}}$
is then equal to
\begin{align}
&(X\times_SY)_{\sigma}
\label{eqlp}
\\
&=
(X\times_SY)
\times_{
{\rm Spec}\ {\mathbb Z}[P]}
{\rm Spec}\ {\mathbb Z}[P_\sigma]
\nonumber\\&=
(X\times_SY)
\times_{
{\rm Spec}\ {\mathbb Z}[S_i,T_i;i\in I]}
{\rm Spec}\ {\mathbb Z}[S_i,T_i,U_i^{\pm1};i\in I]/
(S_i-U_iT_i;i\in I).
\nonumber
\end{align}

We have a global embedding
as follows.
For each $i\in I$,
let ${\cal I}_{D_i}
\subset {\cal O}_X$
and ${\cal I}_{E_i}
\subset {\cal O}_Y$
be the ideal sheaves.
We consider
the ${\mathbf P}^1$-bundle
${\mathbf P}(
{\rm pr}_1^*
{\cal I}_{D_i}\oplus
{\rm pr}_2^*
{\cal I}_{E_i})$
over $X\times_SY$.
The complement
$P_i\subset 
{\mathbf P}(
{\rm pr}_1^*
{\cal I}_{D_i}\oplus
{\rm pr}_2^*
{\cal I}_{E_i})$
of the two sections
defined by the surjections
${\rm pr}_1^*
{\cal I}_{D_i}\oplus
{\rm pr}_2^*
{\cal I}_{E_i}
\to
{\rm pr}_1^*
{\cal I}_{D_i}$
and
${\rm pr}_1^*
{\cal I}_{D_i}\oplus
{\rm pr}_2^*
{\cal I}_{E_i}
\to
{\rm pr}_2^*
{\cal I}_{E_i}$
is a
${\mathbf G}_m$-torsor
over $X\times_SY$.
The log product
$(X\times_SY)^\sim_
{{\cal D},{\cal E}}$
is a closed subscheme
of the fiber product
$\prod_{i\in I\ X\times_SY}P_i$
over $X\times_SY$.

We consider 
the variant of log product.
Further, let $B$ be a Cartier
divisor of $S$
and $(n_i)$
be a family of
integers $n_i\ge 0$
satisfying
$f^*B=
\sum_{i\in I}n_iD_i$
and
$g^*B=
\sum_{i\in I}n_iE_i$
for the same family $(n_i)_{i\in I}$
of integers $n_i\ge 1$.
We consider
the log structure of $S$
defined by $B$ 
and define the log product
\begin{equation}
(X\times_SY)^\sim_
{{\cal D},{\cal E}/B}
\label{eqXYDEB}
\end{equation}
to be
$X\times^{\log}_{S,
{\mathbb N}^I/{\mathbb N}}Y$
defined by the canonical morphisms
${\mathbb N}
\to \Gamma(S,\bar 
{\mathcal M}_S),
{\mathbb N}^I
\to \Gamma(X,\bar 
{\mathcal M}_X),
{\mathbb N}^I\to \Gamma(Y,\bar 
{\mathcal M}_Y)$.
It is a closed
subscheme of
$(X\times_SY)^\sim_
{{\cal D},{\cal E}}$.
When $B$ is clear from the context,
we let
$(X\times_SY)^\sim_
{{\cal D},{\cal E}/B}$
denoted by
$(X\times_{\mathbb S}Y)^\sim_
{{\cal D},{\cal E}}$
in order to distinguish
it from
$(X\times_SY)^\sim_
{{\cal D},{\cal E}}$.

The log product
$(X\times_SY)^\sim_
{{\cal D},{\cal E}/B}$
with respect to
${\cal D},{\cal E}$ and 
$B$
is locally described as follows.
Suppose that
$D_i,E_i$ and $B$
are defined by
$f_i\in\Gamma(X,{\cal O}_X),
g_i\in\Gamma(Y,{\cal O}_Y)$
and
$a\in \Gamma(S,{\cal O}_S)$
respectively.
We put $a=v\prod_if_i^{n_i}$
and $a=w\prod_ig_i^{n_i}$
for $v\in\Gamma(X,{\cal O}_X^\times)$
and
$w\in\Gamma(Y,{\cal O}_Y^\times)$.
Then, $((f_i)_{i\in I},v)$
and $((g_i)_{i\in I},w)$
define
maps of monoids
${\mathbb N}^I\times{\mathbb Z}
\to \Gamma(X,{\cal O}_X)$
and
${\mathbb N}^I\times{\mathbb Z}
\to \Gamma(Y,{\cal O}_Y)$.
Let $P$
be the amalgamate sum
$({\mathbb N}^I\times{\mathbb Z})
+_{\mathbb N}
({\mathbb N}^I\times{\mathbb Z})$
with respect to the
map
${\mathbb N}
\to
{\mathbb N}^I\times{\mathbb Z}$
sending 1 to
$((n_i),1)$.
Then, they further induce
a map
$P\to \Gamma(X\times_SY,
{\cal O}_{X\times_SY})$.
We identify the
dual monoid $N=P^*$
with 
$\{((a_i),(b_i))\in
{\mathbb N}^I
\times {\mathbb N}^I\mid
\sum_in_ia_i=
\sum_in_ib_i\}$
and let
$N_\sigma=
{\mathbb N}^I
\subset 
N\subset 
{\mathbb N}^I
\times {\mathbb N}^I$
be the diagonal submonoid.
Then, the corresponding
submonoid
$P_{\sigma}
\subset 
P^{\rm gp}
=({\mathbb Z}^I\times 
{\mathbb Z}
\oplus
{\mathbb Z}^I\times 
{\mathbb Z})/
\langle((n_i),1,
(-n_i),-1)\rangle$
is equal to
$\{\overline{(a,a',b,b')}
\in 
P^{\rm gp}\mid
a+b
\in {\mathbb N}^I\}$.
The log product
$(X\times_SY)^\sim_
{{\cal D},{\cal E}/B}$
is then equal to
\begin{align}
&(X\times_SY)_{\sigma}
\label{eqlpB}
\\&=
(X\times_SY)
\times_{
{\rm Spec}\ {\mathbb Z}[P]}
{\rm Spec}\ {\mathbb Z}[P_\sigma]
\nonumber\\
&=
(X\times_SY)
\times_{
{\rm Spec}\ {\mathbb Z}
[S_i,T_i;i\in I,
V^{\pm 1},W^{\pm 1}]}
\nonumber\\
&\qquad\qquad
{\rm Spec}\ {\mathbb Z}
[S_i,T_i,U_i^{\pm1};i\in I,
V^{\pm 1},W^{\pm 1}]/
(S_i-U_iT_i;i\in I,
W-V\textstyle{\prod_i}U_i^{n_i}).
\nonumber
\end{align}
In other words,
in the presentation
(\ref{eqlp}),
it is the closed
subscheme
defined by the relation
${\rm pr}_2^*w/
{\rm pr}_1^*v=
\prod_iU_i^{n_i}$.

We study the boundary
of log products.
Let $i\in I$
and put
$\check I_i=I
\setminus \{i\}$
and $\check{\cal D}_i
=(D_j)_{j\in \check I_i}$.
We define
$(D_i\times_SD_i)^\sim_{{\cal D}_i}$
to be
$(D_i\times_SD_i)\times
_{X\times_SX}
(X\times_SX)^\sim_{
\check {\cal D}_i}$.
If $D_i\cap D_j=
D_i\times_X D_j$
is a Cartier
divisor of $D_i$
for every 
$j\in \check I_i$,
the scheme
$(D_i\times_SD_i)^\sim_{{\cal D}_i}$
is the log product
with respect to
the family
$(D_i\times_X D_j)
_{j\in \check I_i}$
denoted by ${\cal D}_i$.

\begin{lm}\label{lmEi}
Let $X\to S$ and
${\cal D}=(D_i)_{i\in I}$
be as above.
Let $i\in I$
and assume that
$D_i\cap D_j=
D_i\times_X D_j$
is a Cartier
divisor of $D_i$
for every 
$j\in I'=I
\setminus \{i\}$.

{\rm 1.}
The scheme  $E_i=
(D_i\times_SD_i)
\times_{(X\times_SX)
}
(X\times_SX)^\sim_{\cal D}$
is equal to the
inverse images
${\rm pr}_1^{-1}(D_i)=
{\rm pr}_2^{-1}(D_i)$
of $D_i\subset X$
by the projections
$(X\times_SX)^\sim_{\cal D}
\to X$.
It is a ${\mathbf G}_m$-torsor
over
$(D_i\times_SD_i)^\sim
_{{\cal D}_i}$.
The restriction
of the log diagonal map
$D_i\to E_i$
defines a trivialization
of the restriction
of the ${\mathbf G}_m$-torsor
$E_i\to
(D_i\times_SD_i)^\sim
_{{\cal D}_i}$
to
$D_i\subset
(D_i\times_SD_i)^\sim
_{{\cal D}_i}$.

{\rm 2.}
Let $B$ be a Cartier divisor of $S$.
Assume that
$f^*B=\sum_jn_jD_j$
and that the coefficient
$n_i$ of $D_i$
in $f^*B$
is strictly positive
$n_i>0$.
Then,
the intersection
$E_i\cap 
(X\times_{\mathbb S}X)^\sim_{{\cal D}/B}$
is
a subscheme
of a $\mu_{n_i}$-torsor
over
$(D_i\times_SD_i)^\sim
_{{\cal D}_i}$.
The restriction
of the log diagonal map
$D_i\to E_i$
defines a trivialization
of the restriction
of the $\mu_{n_i}$-torsor
$E_i\cap 
(X\times_{\mathbb S}X)^\sim_{{\cal D}/B}
\to
(D_i\times_SD_i)^\sim
_{{\cal D}_i}$
to
$D_i\subset
(D_i\times_SD_i)^\sim
_{{\cal D}_i}$.
\end{lm}

{\it Proof.}
1. Clear from the 
inductive construction
$(X\times_SX)^\sim_{\cal D}=
(X\times_SX)^\sim_{
\check {\cal D}_i}
\times_{X\times_SX}
(X\times_SX)^\sim_{D_i}$
of the log product.

2. Clear from
the remark after
(\ref{eqlpB}).
\qed

We define a log blow-up
$(X\times_SY)'
_{{\cal D},{\cal E}}$
of $X\times_SY$
containing
the log product
$(X\times_SY)^\sim
_{{\cal D},{\cal E}}$
as an open subscheme.
For $i\in I$,
let $B_i\subset
N={\mathbb N}^I
\times
{\mathbb N}^I$ 
be the subset
$\{((a_k),
(b_k))
\in {\mathbb N}^I
\times
{\mathbb N}^I\mid
a_k=b_k=0
\text{ for }
k\neq i \text{ and }
(a_i,b_i)
\in \{(1,0),
(0,1),(1,1)\}\}$
consisting of three elements
and we put
$B=\bigcup_{i\in I}B_i$.
Then the set
$\Sigma=
\{\sigma\mid 
\sigma \subset B,\
{\rm Card}(\sigma
\cap B_i)\le 2$
for every $i\in I\}$
defines a regular proper
subdivision of
$N$.
We let the log blow-up
$(X\times_SY)_\Sigma$
denoted by
$(X\times_SY)'
_{{\cal D},{\cal E}}$.
Since the diagonal
submonoid
${\mathbb N}^I
\subset
{\mathbb N}^I
\times
{\mathbb N}^I$
is generated by the subset
$\sigma
=\{((a_i),(b_i))
\in B\mid
a_i=b_i$ for
every $i\in I\}$,
the log product
$(X\times_SY)^\sim
_{{\cal D},{\cal E}}$
is an open subscheme
of 
$(X\times_SY)'
_{{\cal D},{\cal E}}$.

We define a log blow-up
$(X\times_SY)'
_{{\cal D},{\cal E}/B}$
of $X\times_SY$
containing
the log product
$(X\times_SY)^\sim
_{{\cal D},{\cal E}/B}$
as an open subscheme,
assuming
$n_i\in \{0,1\}$
for every $i\in I$.
In order to define
the log blow-up,
we choose and fix
a total order of 
the subset
$I'=\{i\in I\mid
n_i=1\}$
of the index set $I$.

First, we consider
the case where
$I'=I$ namely
$n_i=1$ for every $i\in I$.
The dual $N=P^*$ of
$P={\mathbb N}^I
+_{\mathbb N}
{\mathbb N}^I$
is identified with
$\{(a,b)\in
{\mathbb N}^I
\times
{\mathbb N}^I
\mid
\sum_ia_i=\sum_ib_i\}$.
Let $(e_i)$
be the standard
basis of ${\mathbb N}^I$.
We identify an element
$(i,j)\in I\times I$
with 
$(e_i,e_j)\in N$
and regard
$I\times I$
as a subset of $N$.
We consider
the product order
on the product
$I\times I$.
Let $\Sigma$
be the set of
totally ordered subsets
$\sigma
\subset I\times I$.
For $\sigma\in \Sigma$,
let $N_\sigma
\subset N$
be the submonoid
generated by $\sigma$.
Then,
$\Sigma$
defines a regular proper
subdivision of
$N$.
We let the log blow-up
$(X\times_SY)_\Sigma$
denoted
$(X\times_SY)'
_{{\cal D},{\cal E}/B}$.
Since the diagonal
$\sigma
=\Delta_I
\subset I\times I$
corresponds
to the diagonal
submonoid
${\mathbb N}^I
\subset
N\subset
{\mathbb N}^I
\times
{\mathbb N}^I$,
the log product
$(X\times_SY)^\sim
_{{\cal D},{\cal E}/B}$
is an open subscheme
of 
$(X\times_SY)'
_{{\cal D},{\cal E}/B}$.

In the general case,
we put $I''=I\setminus I'$
and consider
the subfamilies
${\cal D}'=(D_i)_{i\in I'},
{\cal E}'=(E_i)_{i\in I'},
{\cal D}''=(D_i)_{i\in I''},
{\cal E}''=(E_i)_{i\in I''}$.
Then, we define the log product
$(X\times_SY)'_{
{\cal D},{\cal E}/B}$
as the fiber product
by
\setcounter{equation}0
\begin{equation}
(X\times_SY)'_{
{\cal D},{\cal E}/B}=
(X\times_SY)'_{
{\cal D}',{\cal E}'/B}
\times_{X\times_SY}
(X\times_SY)'_{
{\cal D}'',{\cal E}''}.
\label{eqlogbp}
\end{equation}

\subsection{Log Lefschetz trace formula
over a discrete valuation ring}\label{sslL}

We state and prove
a log Lefschetz trace formula
over a discrete valuation ring.
Let $L$ be a henselian discrete valuation
field.
We regard
$T={\rm Spec}\ {\cal O}_L$
as a log scheme
with the log structure
defined by the closed
point $t$
and also regard
$t$ as a log point.

Let $X$ be a 
weakly semi-stable scheme
over $T={\rm Spec}\ {\cal O}_L$
with smooth generic fiber $X_L$
and $D\subset X$
be a Cartier divisor
with normal crossings
relatively to $T$.
Let 
$j\colon X_L\to X$
be the open immersion.
In this section,
we regard the scheme $X$
as a log scheme
with the log structure
${\cal O}_X\cap 
j_*{\cal O}_{X_L}^\times$.
It is log smooth over $T$.
We consider
the fiber $X_t=X\times_Tt$
also as a log scheme
over a log point $t$.
We put
$U=X\setminus D$
and let
$j_U\colon U\to X$
be the open immersion.

If $X$ is proper,
we define the
log \'etale cohomology
with compact support by
$$H^q_{\log c}(U_{\bar t},
{\mathbb Q}_\ell)=
H^q_{\log}(X_{\bar t},
j_{U!}{\mathbb Q}_\ell)$$
where $j_{U!}$ is defined
on the log \'etale site.

\begin{lm}\label{lmcospH}
Let $X$ be a 
proper weakly 
semi-stable scheme
over $T={\rm Spec}\ {\cal O}_L$
and $D\subset X$
be a Cartier divisor
with normal crossings
relatively to $T$.
Then the cospecialization map
\setcounter{equation}0
\begin{equation}
H^*_{\log c}(U_{\bar t},{\mathbb Q}_\ell)
\longrightarrow
H^*_c(U_{\bar L},{\mathbb Q}_\ell)
\label{eqcosp}
\end{equation}
is an isomorphism.
\end{lm}

{\it Proof.}
In the case $D=\emptyset$,
it follows from
\cite[Theorem (3.2) (ii)]{Na2}.
We reduce the general
case to this case.
Let $\bar D$
be the normalization $D$
and let
$\pi\colon \bar D
\to X$ be the canonical
map.
Then, we have an exact sequence
$0\to j_{U!}{\mathbb Q}_{\ell,U}
\to {\mathbb Q}_{\ell,X}
\to \pi_*
{\mathbb Q}_{\ell,\bar D}
\to \Lambda^2\pi_*
{\mathbb Q}_{\ell,\bar D}
\to \cdots .$
Thus the assertion follows.
\qed

Let $X'$ be another
weakly semi-stable scheme
over $T$
with smooth generic fiber $X'_L$
and $D'\subset X'$
be a Cartier divisor
with normal crossings
relatively to $T$.
We also regard
$X'$ and $X'_t$ 
as log schemes
over $T$ and over $t$.
Let $\iota_t\colon X_t\to X'_t$
be an isomorphism of
log schemes over $t$
inducing an isomorphism
$D_t\to D'_t$.
Then, it induces an isomorphism
$\iota_{t*}\colon
H^q_{\log c}(U_{\bar t},{\mathbb Q}_\ell)
\to
H^q_{\log c}(U'_{\bar t},{\mathbb Q}_\ell)$.

Let $\Gamma\subset U_L\times U'_L$
be a closed subscheme of 
dimension $d=\dim U_L$.
We assume that
the second projection
$p_2\colon
\Gamma\to U'_L$
is proper.
Then, in \cite[Section 2.3]{KSA},
the map
\setcounter{equation}0
\addtocounter{thm}1
\begin{equation}
\Gamma^*\colon
H^q_c(U'_{\bar L},{\mathbb Q}_\ell)
\to
H^q_c(U_{\bar L},{\mathbb Q}_\ell)
\label{eqGm*}
\end{equation}
is defined as $p_{1*}p_2^*$.
We also
let $\Gamma^*$
denote the composition
\begin{equation}
\begin{CD}
H^q_c(U_{\bar L},{\mathbb Q}_\ell)
@.
H^q_c(U'_{\bar L},{\mathbb Q}_\ell)
@>{\Gamma^*}>>
H^q_c(U_{\bar L},{\mathbb Q}_\ell)\\
@A{\rm cosp.}A{\simeq}A
@A{\rm cosp.}A{\simeq}A
@.\\
H^q_{\log c}(U_{\bar t},{\mathbb Q}_\ell)
@>{\iota_{t*}}>>
H^q_{\log c}(U'_{\bar t},{\mathbb Q}_\ell)
@.
\end{CD}
\label{eqGm2*}
\end{equation}
by abuse of notation.
In this subsection,
we give a 
Lefschetz trace formula 
computing the alternating
sum 
\begin{equation}
{\rm Tr}(
\Gamma^* 
\colon
H^*_c(U_{\bar L},
{\mathbb Q}_\ell))=
\sum_{q=0}^{2d}
(-1)^q{\rm Tr}(
\Gamma^* 
\colon
H^q_c(U_{\bar L},
{\mathbb Q}_\ell))
\in {\mathbb Q}_\ell,
\label{eqTr*}
\end{equation}
assuming that
$X$ is weakly strictly
semi-stable.

Let $X, X',D$ and $D'$
be as above.
We assume further that
$X$ and $X'$ are
weakly strictly semi-stable
and that 
$D=D_1+\cdots+D_n$
and $D'=D'_1+\cdots+D'_n$
have simple normal crossings
with the same indices.
Let $\iota_t\colon X_t\to X'_t$
be an isomorphism of
log schemes over $t$
inducing isomorphisms
$D_{i,t}\to
D'_{i,t}$ for every
$1,\ldots,n$.

Let ${\mathbb T}$ denote
the log scheme $T$
endowed with the log structure
defined by the Cartier divisor $t$.
We consider an fs-monoid  $P$,
a morphism $ {\mathbb N}\to P$
of monoids and
a commutative diagram
\begin{equation}
\begin{CD}
P@<<< {\mathbb N}@>>>P\\
@VVV @VVV @VVV\\
\Gamma(X,
\bar {\mathcal M}_X)
@<<<
\Gamma({\mathbb T},
\bar {\mathcal M}_{\mathbb T})
@>>>
\Gamma(X',
\bar {\mathcal M}_{X'})
\end{CD}
\label{eqPT}
\end{equation}
of monoids
satisfying
the following condition:
\begin{itemize}
\item[{\rm (P)}]
The vertical
arrows are locally lifted to charts
and compatible with the isomorphism
$\iota_t\colon X_t\to X'_t$.
\end{itemize}

To define the log product
$(X\times_{\mathbb T}
X')^\sim$,
we define ${\mathbb X}$ 
to be the log scheme $X$
defined by the push-out
${\mathcal M}'_X$
of the log structure ${\mathcal M}_X$
and that defined by the family
${\cal D}=
(D_1,\ldots,D_n)$.
Similarly, we define
${\mathbb X}'$.
We consider
$$\begin{CD}
P\oplus {\mathbb N}^n
@<<< {\mathbb N}@>>>
P\oplus {\mathbb N}^n
\\
@VVV @VVV @VVV\\
\Gamma({\mathbb X},
\bar {\mathcal M}_{\mathbb X})
@<<<
\Gamma({\mathbb T},
\bar {\mathcal M}_{\mathbb T})
@>>>
\Gamma({\mathbb X}',
\bar {\mathcal M}_{{\mathbb X}'})
\end{CD}$$
and define the log product 
\begin{equation}
(X\times_{\mathbb T}
X')^\sim
\label{eqXX'T}
\end{equation}
to be
${\mathbb X}
\times_{{\mathbb T},P\oplus
{\mathbb N}^n/{\mathbb N}}
{\mathbb X}'$.
Since the canonical
map
$(X\times_{\mathbb T}
X')^\sim\to 
X\times_{\mathbb T}X'$
is log \'etale,
the projections
$(X\times_{\mathbb T}
X')^\sim\to X$
and
$(X\times_{\mathbb T}
X')^\sim\to X'$
are log smooth.
Similarly as 
Lemma \ref{lmlp},
the projections
are strict and hence
smooth.

By the universality of 
$(X\times_{\mathbb T}X')^\sim$,
the immersion
$X_t\to X$
and the composition
$\iota_t\colon
X_t\to X'_t \to X$
defines an immersion
$X_t\to
(X\times_{\mathbb T}X')^\sim$.
By identifying
$X'_t$ with $X_t$
by the isomorphism $\iota_t$,
let
$\delta_t\colon
X_t\to
(X\times_{\mathbb T}X')^\sim$
denote the immersion.
The generic fiber
$(X\times_{\mathbb T}X')^\sim
\times_T{\rm Spec}\ L$
is identified
with the log product
$(X_L\times_LX'_L)^\sim$
with respect to
the families of Cartier divisors
$D_{1,L},\ldots,D_{n,L}$
and
$D'_{1,L},\ldots,D'_{n,L}$.

\begin{df}\label{dfGK}
{\rm 1.}
For a scheme $S$,
let $K(S)$
denote the Grothendieck
group of
the exact category
of locally free ${\cal O}_S$-modules
of finite rank.
For a locally free ${\cal O}_S$-modules
${\mathcal E}$ of finite rank,
its class is denoted by
$[{\mathcal E}]\in K(S)$.

{\rm 2.}
For a noetherian scheme
$S$,
let $G(S)$
denote the Grothendieck
group of
the abelian category
of coherent ${\cal O}_S$-modules.
For a coherent ${\cal O}_S$-modules
${\mathcal F}$,
its class is denoted by
$[{\mathcal F}]\in G(S)$.
For an integral closed subscheme
$V$,
the class $[{\mathcal O}_V]$
is also denoted by $[V]$.
\end{df}

We define a map
$G((X_L\times_LX'_L)^\sim)
\to G(X_t)$ as follows.

\begin{lm}\label{lmGt}
Let $X$ and $X'$
be weakly strictly
semi-stable schemes
over
$T={\rm Spec}\ {\cal O}_L$
with smooth generic fibers
and
$D=D_1+\cdots+D_n
\subset X$
and
$D'=D'_1+\cdots+D'_n
\subset X'$
be
divisors
with simple normal
crossings relatively to $T$
with the same indices.
Let $\iota_t\colon
X_t\to X'_t$
be an isomorphism
of log schemes
compatible with
the numberings of $D$
and $D'$
and we consider
a commutative diagram
{\rm (\ref{eqPT})} of
monoids satisfying the condition {\rm (P)}.

{\rm 1.}
The pull-back
$G((X\times_{\mathbb T}X')^\sim)
\to 
G((X\times_{\mathbb T}X')_t^\sim)$
by the closed
immersion
$(X\times_{\mathbb T}X')_t^\sim
=(X\times_{\mathbb T}X')^\sim
\times_Tt
\to
(X\times_{\mathbb T}X')^\sim$
induces a map
\setcounter{equation}0
\begin{equation}
G((X_L\times_L X'_L)^\sim)
\to 
G((X\times_{\mathbb T}X')_t^\sim).
\label{eqGt1}
\end{equation}

{\rm 2.}
The map
$\delta_t\colon
X_t\to
(X\times_{\mathbb T}X')_t^\sim$
is a regular immersion
and it defines a pull-back
\begin{equation}
G((X\times_{\mathbb T}X')_t^\sim)
\to G(X_t).
\label{eqGt2}
\end{equation}
\end{lm}

{\it Proof.}
1.
Since
the projection
$(X\times_{\mathbb T}X')^\sim$
is smooth over $X$,
the scheme
$(X\times_{\mathbb T}X')^\sim$
is flat over $T$.
Hence the closed immersion
$(X\times_{\mathbb T}X')_t^\sim\to
(X\times_{\mathbb T}X')^\sim$
is a regular immersion
and is of finite tor-dimension.
Thus, it induces a map
$G((X\times_{\mathbb T}X')^\sim)
\to 
G((X\times_{\mathbb T}X')_t^\sim)$.

Since the sequence
$$ 
G((X\times_{\mathbb T}X')_t^\sim)
\to
G((X\times_{\mathbb T}X')^\sim)
\to 
G((X_L\times_L
X'_L)^\sim)\to 0$$
is exact
and since the composition
$$ 
G((X\times_{\mathbb T}X')_t^\sim)
\to
G((X\times_{\mathbb T}X')^\sim)
\to 
G((X\times_{\mathbb T}X')_t^\sim)$$
is the zero-map,
the assertion follows.

2.
Since the projection
$(X\times_{\mathbb T}X')^\sim
\to X$
is smooth,
the immersion
$\delta_t\colon
X_t\to
(X\times_{\mathbb T}X')^\sim_t$
is a section
of a smooth map
and is a regular immersion.
Hence the pull-back
on the Grothendieck groups
is defined.
\qed

For an element 
$\widetilde \Gamma\in 
G((X_L\times_L X'_L)^\sim)$,
by Lemma \ref{lmGt},
its reduction
$\widetilde \Gamma_t
\in
G((X\times_{\mathbb T}X')_t^\sim)$
and
the intersection product
\begin{equation}
(\widetilde \Gamma_t,\Delta_{X_t})
=\delta_t^*(\widetilde \Gamma_t)
\in G(X_t)
\label{eqdlgm}
\end{equation}
are defined.

Recall that
for a weakly semi-stable
scheme $X$ over $T$,
a semi-stable modification
$X_\Sigma$ is
constructed
in Lemma \ref{lmsOK}
by patching log blow-ups.
By the construction,
the pull-back $D_\Sigma
=D\times_XX_\Sigma$
is a divisor 
of $X_\Sigma$ with
simple normal crossings
relatively to $T$.
The canonical map
$X_\Sigma\to X$
induces an isomorphism
$X_{\Sigma, L}
\to X_L$
on the generic fiber.

\begin{cor}\label{corGt}
Let weakly strictly
semi-stable schemes
$X,X'$ over $T$
and an isomorphism
$\iota_t\colon
X_t\to X'_t$
and
a commutative diagram
{\rm (\ref{eqPT})} of monoids
be as in Lemma {\rm \ref{lmGt}}.
Let $f\colon X_\Sigma
\to X$ and $f'\colon
X'_\Sigma\to X'$
be the semi-stable
modification as above.
Then, the diagram
$$\begin{CD}
G((X_L\times_LX'_L)^\sim)
@>>> G((X_\Sigma)_t)\\
@| @VV{f_*}V\\
G((X_L\times_LX'_L)^\sim)
@>>> G(X_t)
\end{CD}$$
is commutative.
\end{cor}

{\it Proof.}
We show that the diagram
\setcounter{equation}0
\begin{equation}
\begin{CD}
(X_{\Sigma}\times_{\mathbb T}
X'_{\Sigma})^\sim
@>{{\rm pr}_1}>>
X_\Sigma
\\
@V{(f\times f')^\sim}VV
@VVfV\\
(X\times_{\mathbb T}X')^\sim
@>{{\rm pr}_1}>> X
\end{CD}
\label{eqSigt}
\end{equation}
is cartesian.
By the definition (\ref{eqllp}) of
$(X\times_{\mathbb T}X')^\sim$,
it suffices to show that
the diagram with $\times_{\mathbb T}$
replaced by $\times_T$
is cartesian.
Since the assertion is
local on $X\times_TX'$,
we may assume that we have
charts $P \to
\Gamma(X,{\mathcal O}_X)$
and $P \to
\Gamma(X',{\mathcal O}_{X'})$.
Let $Q$ be 
a sub-fs-monoid of
$P^{\rm gp}$
containing
$P$ as a submonoid.
Let $\widetilde P$
be the inverse image of
$P$ by the sum
$P^{\rm gp}\oplus P^{\rm gp}
\to P^{\rm gp}$
and define $\widetilde Q$
similarly.
Then, since
$P^{\rm gp}=
Q^{\rm gp}$,
the inclusions
$P\to \widetilde P,
Q\to \widetilde Q$
to the first factors
induce an isomorphism
${\mathbb Z}[\widetilde P]
\otimes_
{{\mathbb Z}[P]}
{\mathbb Z}[Q]
\to
{\mathbb Z}[\widetilde Q]$.
Hence the diagram
with $X_\Sigma$
and $X'_\Sigma$
replaced by
with $X\times_{
{\rm Spec}\ {\mathbb Z}[P]}
{\rm Spec}\ {\mathbb Z}[Q]$
and $X'\times_{
{\rm Spec}\ {\mathbb Z}[P]}
{\rm Spec}\ {\mathbb Z}[Q]$
is cartesian.
Since $X_\Sigma$
and $X'_\Sigma$
are defined by patching
them,
the diagram (\ref{eqSigt}) is
cartesian
by the local construction of
log product.

By the cartesian
diagram (\ref{eqSigt}), 
we obtain a commutative
diagram
$$\begin{CD}
G((X_{\Sigma}\times_{\mathbb T}
X'_{\Sigma})^\sim)
@>>> G((X_\Sigma)_t)\\
@V{(f\times f')^\sim_*}VV
@VV{f_*}V\\
G((X\times_{\mathbb T}X')^\sim)
@>>>
G(X_t).
\end{CD}$$
Since
the diagram
$$\begin{CD}
G((X_{\Sigma}\times_{\mathbb T}
X'_{\Sigma})^\sim)
@>>> G(X_L\times_LX'_L)\\
@V{(f\times f')^\sim_*}VV
@|\\
G((X\times_{\mathbb T}X')^\sim)
@>>>
G(X_L\times_LX'_L)
\end{CD}$$
is commutative,
the assertion follows.
\qed

Let $X,X',D,D'$
and $\iota_t\colon
X_t\to X'_t$ be
as above.
Assume $X$ and $X'$ are strictly
semi-stable
and let $E_1,\ldots,E_m$
be the irreducible components of $X_t$
and $E'_1,\ldots,E'_m$
be the irreducible components of $X'_t$
such that $\iota_t$ maps
$E_j$ to $E'_j$
for $j=1,\ldots,m$.
Then, the log product
$(X\times_{\mathbb T}X')^\sim$
is equal to
$(X\times_TX')_
{{\mathcal D}\cup {\mathcal E},
{\mathcal D}'\cup {\mathcal E}'/t}$
defined by the families
of Cartier divisors
${\mathcal D}=
(D_1,\ldots,D_n),
{\mathcal E}=
(E_1,\ldots,E_m)$ of $X$,
${\mathcal D}'=
(D'_1,\ldots,D'_n),
{\mathcal E}=
(E'_1,\ldots,E'_m)$ of $X'$
and $t$ of $T$.
Consequently, in this case,
the log blow-up
$(X\times_{\mathbb T}X')'$
is defined as
(\ref{eqlogbp}) and
contains
$(X\times_{\mathbb T}X')^\sim$
as an open subscheme.
It contains
the log product
$(X\times_{\mathbb T}
X')^\sim$
as an open subscheme.
The generic fiber
of the log blow-up
$(X\times_{\mathbb T}X')'\to
X\times_TX'$
is equal to the log
blow-up
$(X_L\times_LX'_L)'
\to 
X_L\times_LX'_L$
used in \cite{KSA}.
If $D=D'=\emptyset$,
the log blow-up
$(X\times_{\mathbb T}X')'$
is equal to the log
blow-up
$(X\times_TX')'$
used in \cite{KSI}.

\begin{lm}\label{lmres}
Let
$(X\times_{\mathbb T}X')'\to
X\times_TX'$
be the log blow-up
and let
$(D\times_{\mathbb T}X')',
(X\times_{\mathbb T}D')'
\subset
(X\times_{\mathbb T}X')'$
be the proper
transforms of 
$D\times_TX',
X\times_TD'
\subset
X\times_TX'$.
We consider
the open immersions
$$\begin{CD}
(X\times_{\mathbb T}X')'
\setminus 
((D\times_{\mathbb T}X')'
\cup
(X\times_{\mathbb T}D')')
@>{j_1}>>
(X\times_{\mathbb T}X')'
\setminus 
(X\times_{\mathbb T}D')'
\\
@AAA @AA{j}A\\
(X_L\times_LX'_L)^\sim
@>{j_{1L}}>>
(X_L\times_LX'_L)'
\setminus
(X_L\times_LD'_L)'.
\end{CD}$$
Then the canonical map
\setcounter{equation}0
\begin{equation}
j_{1!}{\mathbb Q}_\ell
\to 
Rj_*j_{1L!}{\mathbb Q}_\ell
\label{eqres'}
\end{equation}
is an isomorphism
on the log \'etale site.
\end{lm}

{\it Proof.}
For an irreducible 
component $D_i$ of $D$,
let $(D_i\times_{\mathbb T}X')'
\subset
(X\times_{\mathbb T}X')'$
denote the proper transform.
For a subset
$I\subset \{1,\ldots,n\}$,
we put
$\bigcap_{i\in I}
(D_i\times_{\mathbb T}X')'
=(D_I\times_{\mathbb T}X')'$
and 
let $i_I\colon
(D_I\times_{\mathbb T}X')'
\setminus 
((D_I\times_{\mathbb T}X')'
\cap
(X\times_{\mathbb T}D')')
\to 
(X\times_{\mathbb T}X')'
\setminus 
(X\times_{\mathbb T}D')'$
be the closed immersion.
Then, 
we have an exact sequence
$$0\to
j_{1!}{\mathbb Q}_\ell
\to
{\mathbb Q}_\ell
\to
\bigoplus_{i=1}^n
i_{i*}{\mathbb Q}_\ell
\to 
\bigoplus_{1\le i<j\le n}
i_{\{i,j\}*}{\mathbb Q}_\ell
\to 
\cdots.$$
Since $(D_I\times_{\mathbb T}X')'$
is log smooth over $T$,
the canonical map
$i_{I*}{\mathbb Q}_\ell\to
Rj_*i_{I L*}{\mathbb Q}_\ell$
is an isomorphism
by \cite[Theorem (0.2)]{Na2}.
\qed

The following theorem is
a key ingredient
in the proof of
a crucial step Proposition \ref{prcysd}
of the proof of the conductor formula.

\begin{thm}\label{thmlLTF}
Let $X$ and $X'$
be proper 
weakly strictly
semi-stable
schemes of relative
dimension $d$ over 
$T$ with smooth generic fibers and 
$D=D_1+\cdots+D_n
\subset X$ and 
$D'=D'_1+\cdots+D'_n
\subset X'$
be divisors with
simple
normal crossings
relatively to $T$.
Let
$\iota_t\colon
X_t\to X'_t$ be
an isomorphism
of proper log smooth schemes
over $t$
inducing isomorphisms
$(D_i)_t\to (D'_i)_t$
for $i=1,\ldots,n$
and we consider
a commutative diagram
{\rm (\ref{eqPT})} of
monoids satisfying the condition {\rm (P)}.
We put
$U=X\setminus D$
and
$U'=X'\setminus D'$.

Let $\Gamma'
\subset 
(X_L\times_LX'_L)'$
be a closed subscheme
of dimension $d$
satisfying 
\setcounter{equation}0
\begin{equation}
\Gamma'\cap
(D_L\times_LX'_L)'
\subset
\Gamma'\cap
(X_L\times_LD'_L)'
\label{eqG'}
\end{equation}
and we put
$\widetilde \Gamma=
\Gamma'\cap 
(X_L\times_LX'_L)^\sim$.
Then, for 
$\Gamma=\Gamma'
\cap (U_L\times_LU'_L)$,
the second projection
$p_2\colon
\Gamma\to U'_L$
is proper and,
for the composition $\Gamma^*$
{\rm (\ref{eqGm2*})}, 
we have
\begin{equation}
{\rm Tr}(\Gamma^*\colon
H^*_c(U_{\bar L},{\mathbb Q}_\ell))=
\deg (\widetilde \Gamma_t,\Delta_{X_t}).
\label{eqlLTF}
\end{equation}
\end{thm}

Note that
$(\widetilde \Gamma_t,\Delta_{X_t})$
in the right hand side
is defined in
(\ref{eqdlgm})
using $\iota_t\colon
X_t\to X'_t$.
The proof is a combination of
that of
\cite[Theorem 2.3.4]{KSA}
and that of \cite[Theorem 6.5.1]{KSI}.
In the proper case where $X=U$,
it is proved in \cite[Theorem 6.5.1]{KSI}.
The proof consists of verifying 
that the method
in the proof of
\cite[Theorem 2.3.4]{KSA}
to treat the open case also works
in this context.

{\it Proof.}
By the argument
in the beginning of
the proof of 
\cite[Theorem 2.3.4]{KSA},
the inclusion 
(\ref{eqG'})
implies that
the second projection
${\rm pr}_2\colon
\Gamma\to U'_L$
is proper.
Thus the endomorphism
$\Gamma^*$ of
$H^*_c(U_{\bar L},
{\mathbb Q}_\ell)$
is defined.

We prove
the equality
(\ref{eqlLTF})
first in the case where
$X$ and $X'$
are strictly
semi-stable.
We put
\begin{eqnarray*}
H^{2d}_{!*}(X_L
\times_L X'_L,
{\mathbb Q}_\ell(d))
&=&
H^{2d}(X_L
\times_L U'_L,
(j_U\times 1)_!
{\mathbb Q}_\ell(d)),
\\
H^{2d}_{\log,!*}
(X\times_TX',
{\mathbb Q}_\ell(d))
&=&
H^{2d}_{\log}
(X\times_TU',
(j_U\times 1)_!
{\mathbb Q}_\ell(d)),\\
H^{2d}_{\log,!*}
(X_t\times_tX'_t,
{\mathbb Q}_\ell(d))
&=&
H^{2d}_{\log}
(X_t\times_tU'_t,
(j_U\times 1)_!
{\mathbb Q}_\ell(d)).
\end{eqnarray*}
Then, 
in \cite[Lemma 2.3.2]{KSA},
the cycle class
$[\Gamma]
\in
H^{2d}_{!*}(X_L
\times_L X'_L,
{\mathbb Q}_\ell(d))$
is defined and the map
$\Gamma^*$
is described
in terms of the
cycle class $[\Gamma]$
as in the upper line
of the diagram
(\ref{eqcosmp}) below.
We consider the image
of $[\Gamma]$
by the composition
\begin{equation}
\begin{CD}
[\Gamma]\in&
H^{2d}_{!*}
(X_L\times_LX'_L,
{\mathbb Q}_\ell(d))\\
&@A{\rm restriction}A{\simeq}A\\
&H^{2d}_{\log,!*}
(X\times_TX',
{\mathbb Q}_\ell(d))\\
&@V{\rm restriction}VV\\
{[}\Gamma_t] \in
&H^{2d}_{\log,!*}
(X_t\times_tX'_t,
{\mathbb Q}_\ell(d))\\
&@V{\delta_t^*}VV\\
\delta_t^*[\Gamma_t]\in
&H^{2d}_{\log,c}
(U_t,{\mathbb Q}_\ell(d))\\
&@V{\rm Tr}VV\\
&{\mathbb Q}_\ell.
\end{CD}
\label{eqdt}
\end{equation}
Similarly as Lemma \ref{lmres},
the first arrow is an isomorphism
by \cite[Proposition (4.2)]{Na2}.

Since the cospecialization maps
are compatible
with the pull-back,
cup-product
and the trace maps,
we have a commutative diagram
\begin{equation}
\begin{CD}
H^*_c(U'_{\bar L},{\mathbb Q}_\ell)
@>{{\rm pr}_2^*}>>
H^*_{*!}(X_{\bar L}
\times X'_{\bar L},{\mathbb Q}_\ell)
@>{\cup[\Gamma]}>>
H^*_c(U_{\bar L}
\times U'_{\bar L},{\mathbb Q}_\ell)
@>{{\rm pr}_{1*}}>>
H^*_c(U_{\bar L},{\mathbb Q}_\ell)
\\
@V{\rm cosp.}VV
@V{\rm cosp.}VV
@VV{\rm cosp.}V
@VV{\rm cosp.}V
\\
H^*_{\log c}(U'_{\bar t},{\mathbb Q}_\ell)
@>{{\rm pr}_2^*}>>
H^*_{\log *!}(X_{\bar t}
\times X'_{\bar t},{\mathbb Q}_\ell)
@>{\cup[\Gamma_t]}>>
H^*_{\log c}(U_{\bar t}
\times U'_{\bar t},{\mathbb Q}_\ell)
@>{{\rm pr}_{1*}}>>
H^*_{\log c}(U_{\bar t},{\mathbb Q}_\ell).
\end{CD}
\label{eqcosmp}
\end{equation}
The composition of 
the arrows in the upper line
is the map $\Gamma^*$ (\ref{eqGm*}).
We define $\Gamma_t^*$
to be the composition of 
the arrows in the lower line.
Then, we have
\begin{equation}
{\rm Tr}(\Gamma^*\colon
H^*_c(U_{\bar L},{\mathbb Q}_\ell))=
{\rm Tr}(\Gamma_t^*\circ \iota_{t*}\colon
H^*_{\log c}(U_{\bar t},
{\mathbb Q}_\ell)).
\label{eqcomp}
\end{equation}

The standard argument
of the proof of
Lefschetz trace formula
using the K\"unneth formula
\cite[Theorem (6.2)]{Na1} 
and the Poincar\'e duality
\cite[Proposition (4.4)]{Na2}
for log \'etale cohomology,
shows that the diagram
$$\begin{CD}
\bigoplus_{q=0}^{2d}{\rm End}
(H^q_{\log c}(U_{\bar t},
{\mathbb Q}_\ell))
@>{\simeq}>>
H^{2d}_{\log,!*}(X_t\times X'_{\bar t},
{\mathbb Q}_\ell(d))\\
@V{\rm Tr}VV
@VV{\delta_t^*}V\\
{\mathbb Q}_\ell
@<{\rm Tr}<<
H^{2d}_{\log,c}(U_{\bar t},
{\mathbb Q}_\ell(d))
\end{CD}$$
is commutative.
Hence, the right hand side
of (\ref{eqcomp})
is equal to
${\rm Tr}(\delta_t^*[\Gamma_t])$.
Thus the proof of
(\ref{eqlLTF}) is reduced to showing
\begin{equation}
{\rm Tr}(\delta_t^*[\Gamma_t])=
\deg (\widetilde
\Gamma_t,\Delta_{X_t}).
\label{eqTd}
\end{equation}

In the definition
of the left hand side 
${\rm Tr}(\delta_t^*[\Gamma_t])$,
we modify the diagram
(\ref{eqdt})
using Lemma \ref{lmres}.
We consider the
log blow-up
$(X\times_{\mathbb T}X')'
\to X\times_TX'$
defined by (\ref{eqlogbp}).
We consider
the proper
transforms
$(D\times_{\mathbb T}X')',
(X\times_{\mathbb T}D')'
\subset
(X\times_{\mathbb T}X')'$
of 
$D\times_TX',
X\times_TD'
\subset
X\times_TX'$
and let $j_1\colon
(X\times_{\mathbb T}X')'
\setminus 
((D\times_{\mathbb T}X')'
\cup
(X\times_{\mathbb T}D')')
\to 
(X\times_{\mathbb T}X')'
\setminus 
(D\times_{\mathbb T}X')'$
be the open immersion
as in Lemma \ref{lmres}.
In the following equalities,
we define
the left hand sides
by the right hand sides
\begin{eqnarray*}
H^*_{\log !*}((X\times_{\mathbb T}X')',
{\mathbb Q}_\ell)
&=&
H^*_{\log}((X\times_{\mathbb T}X')'
\setminus 
(X\times_{\mathbb T}D')',
j_{1!}{\mathbb Q}_\ell)
\\
H^*_{!*}((X_L\times_LX'_L)',{\mathbb Q}_\ell)
&=&
H^*((X_L\times_LX'_L)'
\setminus 
(X_L\times_LD'_L)',j_{1L!}{\mathbb Q}_\ell).
\end{eqnarray*}

We consider a commutative
diagram
\begin{equation}
\begin{CD}
H^{2d}_{!*}
(X_L\times_LX'_L,
{\mathbb Q}_\ell(d))
@>>>
H^{2d}_{!*}
((X_L\times_LX'_L)',
{\mathbb Q}_\ell(d))
@>>>
H^{2d}
((X_L\times_LX'_L)^\sim,
{\mathbb Q}_\ell(d))
\\
@AAA@AAA@AAA\\
H^{2d}_{\log,!*}
(X\times_TX',
{\mathbb Q}_\ell(d))
@>>>
H^{2d}_{\log,!*}
((X\times_{\mathbb T}X')',
{\mathbb Q}_\ell(d))
@>>>
H^{2d}_{\log}
((X\times_{\mathbb T}X')^\sim,
{\mathbb Q}_\ell(d))\\
@V{\delta_t^*}VV @VV{\delta_t^*}V 
@VV{\delta_t^*}V\\
H^{2d}_{\log, c}(U_t,
{\mathbb Q}_\ell(d))
@>>>
H^{2d}_{\log}(X_t,
{\mathbb Q}_\ell(d))
@=
H^{2d}_{\log}(X_t,
{\mathbb Q}_\ell(d)).
\end{CD}
\end{equation}
The upper vertical
arrows are the
restrictions to the generic fiber
and are isomorphisms
by Lemma \ref{lmres}
and \cite[Proposition (4.2)]{Na2}.
The top left horizontal
arrow
is an isomorphism
by \cite[Corollary 2.2.2]{KSA}.
By the proof of 
\cite[Theorem 2.3.4]{KSA},
the cycle class
$[\Gamma']
\in 
H^{2d}_{!*}((X_L
\times_L X'_L)',
{\mathbb Q}_\ell(d))$
is the image of
$[\Gamma]
\in 
H^{2d}_{!*}(X_L
\times_L X'_L,
{\mathbb Q}_\ell(d))$.
Thus, in the diagram
(\ref{eqdt}),
we may replace
$X\times_TX'$ etc.\ by 
$(X\times_{\mathbb T}X')'$
etc.,
$H^{2d}_{\log, c}(U_t,
{\mathbb Q}_\ell(d))$
by
$H^{2d}_{\log}(X_t,
{\mathbb Q}_\ell(d))$
and $[\Gamma]$
by $[\Gamma']$.
Further,
we may replace
$(X\times_{\mathbb T}X')'$ etc.\ by 
$(X\times_{\mathbb T}X')^\sim$ etc.
and 
$[\Gamma']$
by $[\widetilde \Gamma]$
and drop $_{!*}$.
Therefore the proof of
(\ref{eqlLTF}) is further
reduced to showing
\begin{equation}
{\rm Tr}(\delta_t^*[\widetilde
\Gamma_t])=
\deg (\widetilde
\Gamma_t,\Delta_{X_t})
\label{eqTd2}
\end{equation}
where the left
hand side is defined
by 
\begin{equation}
\begin{CD}
[\widetilde \Gamma]\in&
H^{2d}
((X_L\times_LX'_L)^\sim,
{\mathbb Q}_\ell(d))\\
&@A{\rm restriction}A{\simeq}A\\
&H^{2d}_{\log}
((X\times_{\mathbb T}X')^\sim,
{\mathbb Q}_\ell(d))\\
&@V{\rm restriction}VV\\
[\widetilde \Gamma_t]\in
&H^{2d}_{\log}
((X\times_{\mathbb T}X')^\sim_t,
{\mathbb Q}_\ell(d))\\
&@V{\delta_t^*}VV\\
\delta_t^*[\widetilde \Gamma_t]\in
&H^{2d}_{\log}
(X_t,{\mathbb Q}_\ell(d))\\
&@V{\rm Tr}VV\\
&{\mathbb Q}_\ell.
\end{CD}
\label{eqdt2}
\end{equation}

We prove (\ref{eqTd2}).
The proof goes similarly
as that of \cite[Theorem 6.5.1]{KSI}.
We identify
$Gr_F^dK((X_L\times_LX'_L)^\sim)_{\mathbb Q}$
with
$Gr^F_0G((X_L\times_LX'_L)^\sim)_{\mathbb Q}$
by the canonical isomorphism,
cf.\ \cite[Lemma 2.1.4]{KSI}.
Since we are
assuming that $X$ is strictly semi-stable,
the log product
$(X\times_{\mathbb T}X')^\sim$
is regular.
Hence, the restriction map
$Gr_F^dK((X\times_{\mathbb T}X')^\sim)
\to
Gr_F^dK((X_L\times_LX'_L)^\sim)$
is a surjection
and the class 
$[\widetilde \Gamma]\in 
Gr_F^dK((X_L\times_LX'_L)^\sim)$
is lifted to
an element
$[\widetilde \Gamma]\in
Gr_F^dK((X\times_{\mathbb T}X')^\sim)$.
We define the class
$[\widetilde \Gamma]\in
H^{2d}
((X\times_{\mathbb T}X')^\sim,
{\mathbb Q}_\ell(d))$
as the Chern character.
Since the Chern character
is compatible
with the pull-back,
the class
$[\widetilde \Gamma]\in
H^{2d}
((X_L\times_LX'_L)^\sim,
{\mathbb Q}_\ell(d))$
on the top
is the restriction
of the class
$[\widetilde \Gamma]\in
H^{2d}
((X\times_{\mathbb T}X')^\sim,
{\mathbb Q}_\ell(d))$.
Further the trace map
$H^{2d}
(X_{\bar t},
{\mathbb Q}_\ell(d))
\to {\mathbb Q}_\ell$
is the composition
of the canonical map
$H^{2d}
(X_{\bar t},
{\mathbb Q}_\ell(d))
\to
H^{2d}_{\log}
(X_{\bar t},
{\mathbb Q}_\ell(d))$
with the trace map
$H^{2d}_{\log}
(X_{\bar t},
{\mathbb Q}_\ell(d))
\to {\mathbb Q}_\ell$.
Hence
the left hand side
of the equality
(\ref{eqTd2})
is the image of
$[\widetilde \Gamma]\in
H^{2d}
((X\times_{\mathbb T}X')^\sim,
{\mathbb Q}_\ell(d))$
in the second line
of the diagram (\ref{eqdt2})
with log removed everywhere.
Thus the equality
(\ref{eqTd2})
follows from the compatibility
of the trace map
with the degree map
\cite[Lemma 6.5.4]{KSI}.

We reduce the proof of
(\ref{eqlLTF})
to the case where
$X$ and $X'$
are strictly
semi-stable.
As in Corollary
\ref{corGt},
we consider
the semi-stable
modifications
$X_\Sigma\to X$ and
$X'_\Sigma\to X'$
constructed in Lemma \ref{lmsOK}.
The isomorphism
$\iota_t\colon X_t\to X'_t$
induces
an isomorphism
$\iota_{t \Sigma }\colon 
(X_\Sigma)_t\to (X'_\Sigma)_t$.
Thus,
by Corollary \ref{corGt},
it is reduced
to the strictly
semi-stable case.
\qed

\subsection{Log stalks}\label{sslst}

In the next subsection,
we prove an important
complement Proposition \ref{prcc}
to the log Lefschetz
trace formula Theorem \ref{thmlLTF}.
As a preliminary,
we briefly recall some
elementary terminology
on log points
and the stalks
of a tamely ramified
sheaf at a log geometric
point.
For more systematic
account,
we refer to
\cite[Section 4]{IlV}.
A reader familiar
with the generalities
on log schemes
may skip this subsection.

Let $t$ be the spectrum of
a field $F$.
If $t$ is endowed
with the log structure
defined by the chart
${\mathbb N}\to F$
sending $1$ to $0$,
we call $t$
a log point.
Let $\bar t$ be the spectrum of
a separably closed field $F$
of characteristic $p\ge 0$.
If $\bar t$ is endowed
with the log structure
defined by the chart
${\mathbb Z}_{(p)}\cap [0,\infty)
\to F$
sending any element
$a>0$ to $0$,
we call $\bar t$
a log geometric point.
For a log scheme $Y$,
we call a morphism
$t\to Y$ from a log point
a log point of $Y$.
Similarly,
we call a morphism
$\bar t\to Y$ from 
a log geometric point
a log geometric point of $Y$.

A typical example of log points
and log geometric points
are constructed as follows.
Let ${\cal O}_L$
be a discrete valuation ring
and regard 
$T={\rm Spec}\
{\cal O}_L$ as a log scheme
with the log structure
defined by the closed point $t$.
Then, the scheme $t$
endowed with
the pull-back log structure
is a log point.
Assume further
${\cal O}_L$ is henselian
and let $L^{\rm tr}$
denote
the maximal tamely ramified
extension of $L$.
Then, the limit
of the standard
log structures
on ${\rm Spec}\ {\mathcal O}_{L'}$
for finite extensions $L'$
of $L$ in $L^{\rm tr}$
defines a structure
of log geometric point
on the closed point
$\bar t$ of 
${\rm Spec}\ {\cal O}_{L^{\rm tr}}$.
A morphism
$T\to Y$ of log schemes
define a log point
$t\to Y$
and further
a log geometric point
$\bar t\to Y$.

For a log geometric point
$\bar t$
of a log scheme $Y$,
the log strict localization
$\widetilde Y_{\bar t}$
is defined in
\cite[4.5]{IlV}.
The definition
of the
log stalk ${\cal G}_{\bar t}$
of a sheaf ${\cal G}$
on the Kummer \'etale
site on a log scheme $Y$
at a log geometric point
$\bar t$ of $Y$
is given in
\cite[Definition 4.3]{IlV}.
We will make it explicit
in a special case.

Let $S$ be a
regular noetherian scheme
and $D$ be 
a divisor 
with normal crossings.
We put
$W=S\setminus D$
and 
$j\colon W\to S$
the open immersion.
Then, the log scheme $S$
with the log structure
${\cal M}_S
=
{\cal O}_S\cap
j_*{\cal O}_W^\times$
is log regular (\cite{KT}).
We consider
a locally constant
sheaf ${\cal F}$
on $W$
tamely ramified along $D$.
The direct image
$j_*{\cal F}$ 
on the Kummer \'etale site
of $S$ is
a locally constant sheaf
by Abhyankar's lemma.
Let $g\colon Y\to S$
be a morphism of
log schemes
and we consider
the pull-back
${\cal G}=g^*j_*{\cal F}$
to the Kummer \'etale site of $Y$.

Let $\bar t\to Y$
be a log geometric point
and let $\bar s$
denote the
geometric point of $S$
defined by the
composition
$\bar t\to Y\to S$.
Let $\tilde g\colon
\widetilde Y_{\bar t}
\to
\widetilde S_{\bar s}$
denote the map
of the log strict localizations induced
by $g$
and $\tilde j\colon
W\times_S\widetilde
S_{\bar s}\to
\widetilde S_{\bar s}$
denote
the open immersion.
The pull-back of
${\cal F}$
on $W\times_S\widetilde
S_{\bar s}$
is a constant sheaf
and hence
the direct image
$\widetilde{\cal F}=
\tilde j_*({\cal F}|_{
W\times_S\widetilde
S_{\bar s}})$ is
a constant sheaf
on the usual \'etale site
of $\widetilde S_{\bar s}$.
The log stalk
${\cal G}_{\bar t}$
is canonically identified with
the stalk 
$(\tilde g^*
\widetilde{\cal F})_{\bar t}$ 
at ${\bar t}$
of the pull-back of
the constant sheaf
$\widetilde{\cal F}$.
The map
$\tilde g\colon
\widetilde Y_{\bar t}\to
\widetilde S_{\bar s}$
induces an isomorphism
$(j_*{\cal F})_{\bar s}
\to
{\cal G}_{\bar t}$
of log stalks.

We consider the 
log cospecialization
map (\cite[(2.8) 6]{Na1}).
We will use it only in the following situation.
Let $g\colon Y\to S$ be a
morphism of log schemes
and ${\cal G}=g^*j_*{\cal F}$
be as above.
Let
$V\subset Y$
be an open subscheme
where the log structure
is trivial and
let $\bar \eta$ be 
a geometric point
of $V$.
Assume that the image of
the log geometric point
$\bar t$ in $Y$
lies in the closure
of the image of $\bar \eta$.
Then, by choosing
a lifting of $\bar \eta$
in $\widetilde Y_{\bar t}$,
a log
cospecialization map
${\cal G}_{\bar t}
\to 
{\cal G}_{\bar \eta}$
is defined
as the usual
cospecialization map
$(\tilde g^*\widetilde
{\cal F})_{\bar t}
\to 
(\tilde g^*\widetilde
{\cal F})_{\bar \eta}$.
Let $\bar \xi$
be an intermediate 
geometric point
of $V$
such that the image of
$\bar t$ lies in the closure
of the image of 
$\bar \xi$ and
that the image of 
$\bar \xi$
lies in the closure
of the image of 
$\bar \eta$.
Then, by choosing
liftings 
$\bar \xi\to
\widetilde Y_{\bar t}$
and
$\bar \eta
\to Y_{\bar \xi}\to
\widetilde Y_{\bar t}$
successively,
we obtain
the transitivity of
cospecialization maps
${\cal G}_{\bar t}
\to 
{\cal G}_{\bar \xi}
\to 
{\cal G}_{\bar \eta}$.

The following compatibility
of the cospecialization map
with the pull-back
will be used
in the proof
of Proposition \ref{prcc}.

\begin{lm}\label{lmcosp}
Let $S$
be a regular noetherian
scheme and
$W=S\setminus D
\subset S$ be the complement
of a divisor $D$ with
normal crossings.
We consider
a commutative diagram
$$
\xymatrix{
{\bar t}\ar[r]&
t\ar[r]\ar[dr]&
Y\ar[r]^g&
S.\\
&&Y\ar[u]_h\ar[ur]_{g'}&
}$$
of morphisms of
log schemes.
We assume $t$
is a log point
and 
$\bar t$
is a log geometric point.

Let
$\tilde h\colon
\widetilde Y_{\bar t}
\to
\widetilde Y_{\bar t}$
be the morphism on
the log strict localization
induced by $h$.
Let $\bar \eta$
be a usual geometric point
of an open subscheme
$V\subset Y$
where the log structure
is trivial.
We take a lifting
$\bar \eta
\to
\widetilde Y_{\bar t}$
and let $\bar h
\colon \bar \eta\to
\bar \eta$
be a morphism such that
the diagram
\setcounter{equation}0
\begin{equation}
\begin{CD}
\bar \eta@>>> 
\widetilde Y_{\bar t}\\
@V{\bar h}VV @VV{\tilde h}V\\
\bar \eta@>>> 
\widetilde Y_{\bar t}
\end{CD}
\label{eqcsp}
\end{equation}
is commutative.

Let ${\cal F}$
be a locally constant
sheaf on $W$
tamely ramified
along $D$
and we put
${\cal G}=
g^*(j_*{\cal F})$
and 
${\cal G}'=
g^{\prime*}(j_*{\cal F})
=h^*{\cal G}$.
Then,
for the isomorphism
$\bar h^*\colon
{\cal G}_{\bar \eta}
\to {\cal G}'_{\bar \eta}$,
we have a commutative diagram
$$\begin{CD}
{\cal G}_{\bar t}
@>{\rm cosp.}>>
{\cal G}_{\bar \eta}\\
@| @VV{\bar h^*}V\\
{\cal G}'_{\bar t}
@>{\rm cosp.}>>
{\cal G}'_{\bar \eta}
\end{CD}$$
\end{lm}

{\it Proof.}
By the commutative
diagram (\ref{eqcsp}),
we have a commutative diagram
$$\begin{CD}
\tilde g^*(\widetilde {\cal F})_{\bar t}
@>{\rm cosp.}>>
\tilde g^*(\widetilde {\cal F})_{\bar \eta}\\
@| @VV{\bar h^*}V\\
\tilde g^{\prime*}(\widetilde {\cal F})_{\bar t}
@>{\rm cosp.}>>
\tilde g^{\prime*}(\widetilde {\cal F})_{\bar \eta}
\end{CD}$$
Thus it
follows
from the descriptions
of the log stalks
and the cospecialization 
maps.
\qed

We consider a geometric situation.
Let $S$
be a regular noetherian
scheme and
$W=S\setminus D
\subset S$ be the complement
of a divisor $D$ with
normal crossings as above.
Let $f\colon X\to S$
be a proper 
weakly semi-stable scheme
such that
the base change
$X\times_SW\to W$
is smooth
and let
$f_U\colon U\to S$
be the restriction to
the complement
$U\subset X$
of a divisor $E\subset X$
with normal crossings
relatively to $S$
as in Lemma \ref{lmsstm}.
Then, for an integer
$n$ invertible on $S$,
by Lemma \ref{lmsstm},
the higher direct
image
${\cal F}=
R^qf_{U!}{\mathbb Z}/n{\mathbb Z}$
is locally constant on
$W$
and is tamely ramified 
along $D$.

Let $T={\rm Spec}\ {\cal O}_L$
be the spectrum of
a discrete valuation ring
with the log structure
defined by the closed
point $t\in T$
and $g\colon T\to S$
be a morphism of log schemes
such that the image
of the generic point is
in $W$.
Then, 
the pull-back
${\cal G}=g^*j_*{\cal F}$
is a locally constant sheaf
on the Kummer \'etale site of $T$
and the stalk ${\mathcal G}_{\bar L}$
at the geometric point
defined by an algebraic closure
${\bar L}$ is 
identified with
$H^q_c(U_{\bar L},
{\mathbb Z}/n{\mathbb Z})$
by the usual proper base change theorem.
By the proof of \cite[Proposition (4.3)]{Na2},
we obtain a commutative diagram
\begin{equation}
\begin{CD}
{\mathcal G}_{\bar t}
@>{\rm cosp.}>>
{\mathcal G}_{\bar L}\\
@VVV@VVV\\
H^q_{\log c}(U_{\bar t},
{\mathbb Z}/n{\mathbb Z})
@>>>
H^q_c(U_{\bar L},
{\mathbb Z}/n{\mathbb Z})
\end{CD}
\label{eqNak}
\end{equation}
of isomorphisms.

\subsection{Compatibility
with cospecializations}\label{sscc}

\addtocounter{thm}1
\setcounter{equation}0

We prove a
compatibility
with cospecialization
maps, that gives
an important complement to
the log Lefschetz
trace formula.
We consider the following data:
\begin{itemize}
\item[
{\rm (\ref{eqLgen}a)}]
Let $Y$ be a log scheme
and $V\subset Y$ be
an open subscheme
where the log structure
is trivial.
Let $h\colon Y\to Y$
be a morphism
of log schemes
satisfying 
$h(V)\subset V$.
\item[
{\rm (\ref{eqLgen}b)}]
Let $T$ be the spectrum
${\rm Spec}\
{\cal O}_L$ 
of a discrete valuation
ring ${\cal O}_L$
regarded as a log
scheme with
the log structure
defined by the closed
point $t$.
Let $T\to Y$ be a
morphism of log schemes
such that the image
of the generic point
${\rm Spec}\ L\in T$
is in $V$ and that
the map $t\to Y$ of 
log schemes
is the same as
the composition of
$t\to Y
\overset h\to Y$.
Let $\bar t$
be a log geometric point
above the log point $t$.
\item[
{\rm (\ref{eqLgen}c)}]
Let $\bar \eta$ be
a geometric point
of the strict henselization $V_{\bar \xi}$
of $V$ at a geometric point
$\bar \xi$ above
the image $\xi\in V$ of
${\rm Spec}\ L\to T$
and 
$\bar h\colon
\bar \eta\to
\bar \eta$
be an automorphism
compatible with $h$.
\item[
{\rm (\ref{eqLgen}d)}]
Let $f\colon X\to Y$ be a proper
and weakly 
strictly semi-stable
scheme 
of relative dimension $d$
over $Y$
such that
the base change
$X_V=X\times_YV\to V$
is smooth.
Let $D=D_1+\cdots+D_n$
be a divisor of $X$
with simple
normal crossings
relatively to $Y$.
\end{itemize}

The data above are
summarized in
the following diagram:
\begin{equation}
\xymatrix{
&&\qquad
X\supset D\ar[d]
\\
&
T\ar[r]&
Y&
Y\ar[l]_h
\\
t\ar[ur]\ar[urr]\ar[urrr]&
L\ar[r]\ar[u]&
V\ar[u]&
V\ar[l]^{h|_V}\ar[u]
\\
{\bar t}\ar[u]
&&
{V_{\bar \xi}}\ar[u]&
{V_{\bar \xi}}\ar[u]
\\
&&
{\bar \eta}\ar[u]&
{\bar \eta}\ar[l]^{\bar h}\ar[u]
}
\label{eqLgen}
\end{equation}

Let $U=X\setminus D$
be the complement
and $f_U\colon
U\to Y$
be the restriction of $f$.
Let $U_V\subset X_V$
and $U'_V\subset X'_V$ denote
the base changes of
$U\subset X$ 
by the inclusion
and the restriction
$h|_V$
respectively.
We consider
the log product
and the log blow-up
$(X_V\times_VX'_V)^\sim
\subset 
(X_V\times_VX'_V)'$
with respect to
the pull-backs
$(D_{1,V},\ldots,D_{n,V})$ 
and 
$(D'_{1,V},\ldots,D'_{n,V})$ 
by the inclusion
and by $h|_V$
of $(D_1,\ldots,D_n)$.
Let
$(D_V\times_VX'_V)',
(X_V\times_VD'_V)'
\subset
(X_V\times_VX'_V)'$
denote
the proper transforms
of
$D_V\times_VX'_V$
and of
$X_V\times_VD'_V$
respectively.

We consider
a closed subscheme
$\widetilde\Gamma\subset
(X_V\times_VX'_V)^\sim$
flat of relative
dimension $d$ over $V$.
Assume that
the second projection
${\rm pr}_2\colon
\Gamma
=\widetilde\Gamma\cap
(U_V\times_VU'_V)
\to U'_V$
is proper.
Then, the geometric fiber
$\Gamma_{\bar \eta}
\subset U_{\bar \eta}
\times_{\bar \eta}
U'_{\bar \eta}$
defines a linear map
$\Gamma^*\colon
H^q_c(U'_{\bar \eta},
{\mathbb Q}_\ell)
\to
H^q_c(U_{\bar \eta},
{\mathbb Q}_\ell)$
and the morphism
${\rm id}\times \bar h\colon
U_{\bar \eta}
\to 
U'_{\bar \eta}$
induces
an isomorphism
$\bar h^*\colon
H^q_c(U_{\bar \eta},
{\mathbb Q}_\ell)
\to
H^q_c(U'_{\bar \eta},
{\mathbb Q}_\ell)$.
Consequently,
the alternating sum
\begin{equation}
{\rm Tr}(
\Gamma^* \circ \bar h^*
\colon
H^*_c(U_{\bar \eta},
{\mathbb Q}_\ell))=
\sum_{q=0}^{2d}
(-1)^q{\rm Tr}(
\Gamma^* \circ \bar h^*
\colon
H^q_c(U_{\bar \eta},
{\mathbb Q}_\ell))
\in {\mathbb Q}_\ell
\label{eqsscc2}
\end{equation}
is defined.

Let $\iota_t\colon
X_t\to X'_t$ be
the isomorphism
defined by the assumption
that the map $t\to Y$ 
is the same as
the composition with 
$h\colon Y\to Y$.
It induces
an isomorphism
$\iota_{t*}\colon
H^q_{\log c}(U_{\bar t},{\mathbb Q}_\ell)
\to
H^q_{\log c}(U'_{\bar t},{\mathbb Q}_\ell)$
as in (\ref{eqGm2*}).
Since
$\Gamma_L
\subset U_L
\times_LU'_L$ induces
$\Gamma_L^*
\colon
H^*_c(U'_{\bar L},
{\mathbb Q}_\ell))
\to
H^*_c(U_{\bar L},
{\mathbb Q}_\ell)$,
the alternating sum
${\rm Tr}(
\Gamma^*
\colon
H^*_c(U_{\bar L},
{\mathbb Q}_\ell))$
is defined
by {\rm (\ref{eqTr*})}.

The following complement
to the Lefschetz trace formula
Theorem \ref{thmlLTF}
will be used in the proof
of a crucial step
Proposition \ref{prcfs}
of the proof of the conductor formula.

\begin{pr}\label{prcc}
Let the notation be
as in the diagram
{\rm (\ref{eqLgen})}.
We take a lifting
of the geometric point
$\bar \eta$ to
the log strict
localization
$\widetilde Y_{\bar t}$
and assume that
$\bar h\colon
\bar \eta\to
\bar \eta$
is compatible with
the morphism
$\tilde h\colon
\widetilde Y_{\bar t}
\to
\widetilde Y_{\bar t}$
induced by $h$.
Let $\Gamma$ be 
a closed subscheme
of $U_V\times_VU'_V$
flat
of relative dimension $d$ over $V$.
We assume that the second projection
${\rm pr}_2\colon
\Gamma\to U'_V$
is proper.

We assume the following condition:
\begin{itemize}
\item[{\rm (\ref{prcc}.2)}]
\addtocounter{equation}1
There exist a regular
noetherian scheme $S$,
a proper weakly semi-stable
scheme $\bar f\colon X_S\to S$,
a divisor $D_S\subset X_S$
with normal crossings
relatively to $S$
and a morphism
$g\colon Y_{\bar t}\to S$
from the usual strict
localization
satisfying the following
conditions.
The pull-back of
$D\subset X$ over $Y$
to $Y_{\bar t}$
is isomorphic to
that of $D_S\subset X_S$.
There exists a 
divisor $D_S$ of
$S$ with simple normal
crossings such that
the pull-back of $X_S$
to the complement
$W=S\setminus D_S$
is smooth.
\end{itemize}

Then, for the alternating sum
{\rm (\ref{eqsscc2})},
we have
\begin{equation}
{\rm Tr}(
\Gamma^*\circ \bar h^*
\colon
H^*_c(U_{\bar \eta},
{\mathbb Q}_\ell))=
{\rm Tr}(
\Gamma^*
\colon
H^*_c(U_{\bar L},
{\mathbb Q}_\ell)).
\label{eqlLTFe1}
\end{equation}
\end{pr}
{\it Proof.}
By the definition of
$\Gamma^*
\colon
H^*_c(U_{\bar L},
{\mathbb Q}_\ell)
\to
H^*_c(U_{\bar L},
{\mathbb Q}_\ell)$,
it suffices to show
the commutativity
of the diagram
$$\xymatrix{
H^*_c(U_{\bar \eta},{\mathbb Q}_\ell)
\ar[r]^{\bar h^*}&
H^*_c(U'_{\bar \eta},{\mathbb Q}_\ell)
\ar[r]^{\Gamma^*}&
H^*_c(U_{\bar \eta},{\mathbb Q}_\ell)
\\
H^*_c(U_{\bar L},{\mathbb Q}_\ell)
\ar[u]
&H^*_c(U'_{\bar L},{\mathbb Q}_\ell)
\ar[r]^{\Gamma^*}
\ar[u]
&H^*_c(U_{\bar L},{\mathbb Q}_\ell)
\ar[u]
\\
H^*_{\log,c}(U_{\bar t},{\mathbb Q}_\ell)
\ar[u]\ar[ur]
&&}$$
where the non-horizontal arrows
are the cospecialization maps.
For the right square,
it is a consequence
of the compatibility
of a correspondence
with usual cospecializations.

We show the
commutativity of
the left quadrangle.
By replacing $Y$ by
the strict localization,
we may assume $Y=Y_{\bar t}$.
We put $U_S=X_S\setminus D_S$
in (\ref{prcc}.2)
and let $f_W\colon 
U_W=U_S\times_SW
\to W$ be the restriction of
$f_S\colon X_S\to S$.
We consider the smooth
${\mathbb Q}_\ell$-sheaf
${\mathcal F}=
R^qf_{W!}{\mathbb Q}_\ell$
on $W$
tamely ramified along $D_S$.
Let $j_W\colon 
W\to S$ denote the open immersion
and define a smooth sheaf
${\mathcal G}
=g^*j_*{\mathcal F}$
on $Y_{\bar t}$.
We consider the diagram
$$\begin{CD}
{\mathcal G}_{\bar t}
@>{\rm cosp.}>>
{\mathcal G}_{\bar L}
@>{\rm cosp.}>>
{\mathcal G}_{\bar \eta}\\
@VVV @VVV @VVV\\
H^q_{\log,c}(U_{\bar t},{\mathbb Q}_\ell)
@>>>
H^q_c(U_{\bar L},{\mathbb Q}_\ell)
@>>>
H^q_c(U_{\bar \eta},{\mathbb Q}_\ell).
\end{CD}$$
The right square is
the usual commutative diagram
for \'etale cohomology
and the left is (\ref{eqNak}).
We have a similar commutative diagram for
$U'$.
By the transitivity
of cospecialization
maps,
the
commutativity of
the left quadrangle
follows from Lemma \ref{lmcosp}.
\qed

\newpage
\section{Tamely ramified coverings}
\label{stm}

In this section,
we give a definition
for an \'etale
morphism
of schemes
to be tamely ramified
along the boundary.
The purpose of studying
tame ramification
first is 
to define the
wild ramification
locus and to focus on it.

First, we formulate
the definition of
an unramified morphism
to be tamely ramified
along the boundary
using proper modifications and
log products
in Section \ref{sstmd}.
We give a tameness criterion,
Proposition \ref{prtY},
in terms of valuation rings
in Section \ref{sstvl}
after recalling tamely
ramified extensions
of valuation fields
and the limit of 
proper modifications
in Sections \ref{sstmv}
and \ref{sslim}
respectively.
We study 
the relation with Kummer coverings
in \ref{sstKm}.
Finally, we give
a criterion
for a Galois covering to
be tamely ramified
in terms of inertia
groups
in \ref{sstG}.

Although we don't need to
assume for schemes
to be separated
in a large part of this section,
we will assume it for simplicity.

\subsection{Tame ramification and
log products}\label{sstmd}

Recall that
a morphism of schemes
$V\to U$
of finite type
is said to be unramified
if the diagonal
map
$\delta_V\colon
V\to V\times_UV$
is an open immersion.
We consider
a separated scheme $Y$ containing
$V$ as an open subscheme
and
introduce a notion
that an unramified
morphism
$f\colon V\to U$
is tamely ramified
with respect to $Y$.

\begin{lm}\label{lmSig}
Let
$f\colon V\to U$
be an unramified
separated morphism of finite
type of schemes and $j\colon
V\to Y$
be an open immersion
of separated schemes.
Let 
${\cal D}=(D_i)_{i\in I}$ 
be a finite family
of Cartier divisors
of $Y$ such that
$V\cap D_i=\emptyset$
for every $i\in I$.

For a commutative diagram
\begin{equation}
\begin{CD}
Y@<j<< V\\
@VVV @VVfV\\
S@<<<U
\end{CD}
\label{eqYV}
\end{equation}
of separated schemes, let 
$(Y\times_SY)^\sim_{\cal D}$
be the log product
and define a closed subset
$\Sigma_{V/U}^{\cal D}Y\subset Y$
to be the intersection
$\Delta_Y^{\log} \cap \overline W$
of the log diagonal
with the closure of
the open and closed subscheme
$W=(V\times_UV)\setminus \Delta_V
\subset V\times_UV$
in the log product
$(Y\times_SY)^\sim_{\cal D}$.

Then,
the closed subset
$\Sigma_{V/U}^{\cal D}Y\subset Y$
is independent of the choice
of a diagram
{\rm (\ref{eqYV})}.
\end{lm}

For any morphisms
$Y\gets V\to U$ of
schemes,
we can complete them into a
diagram (\ref{eqYV})
by putting $S={\rm Spec}\ {\mathbb Z}$.

{\it Proof.}
We consider a
commutative
diagram (\ref{eqYV})
with $S$ replaced
by another separated
scheme $S'$
and show
that
$\Sigma_{V/U}^{\cal D}Y\subset Y$
are the same.
It suffices
to consider the case
where $S'={\rm Spec}\ {\mathbb Z}$.
Hence,
we may assume that
there exists a morphism
$S\to S'$
that makes the diagram
$$\xymatrix{
Y\ar[r]\ar[d] &
S\ar[dl] \\
S' &
U\ar[u]\ar[l]}$$
commutative.
The canonical map
$(Y\times_SY)^\sim_{\cal D}
\to
(Y\times_{S'}Y)^\sim_{\cal D}$
is a closed immersion
since it is a base change
of the diagonal
$S\to S\times_{S'}S$.
Hence, the assertion follows.
\qed

\begin{df}\label{dfSig}
Let
$f\colon V\to U$
be an unramified
separated morphism of finite
type of schemes
and $j\colon V\to Y$
be an open immersion
of separated schemes.

{\rm 1.}
For a finite family
of Cartier divisors
${\cal D}=(D_i)_{i\in I}$ 
of $Y$ such that
$V\cap D_i=\emptyset$
for every $i\in I$,
define a closed subset
$\Sigma_{V/U}^{\cal D}Y
\subset Y$
to be the intersection
$\Delta_Y^{\log} 
\cap \overline W$
of the log diagonal
with the closure of
the open and closed subscheme
$W=(V\times_UV)\setminus \Delta_V
\subset V\times_UV$
in the log product
$(Y\times_SY)^\sim_{\cal D}$
as in Lemma {\rm \ref{lmSig}}
by taking a commutative
diagram 
{\rm (\ref{eqYV})}.

Define a closed sunset
$\Sigma^+_{V/U}Y\subset Y$
to be the intersection
$\bigcap_{\cal D}
\Sigma_{V/U}^{\cal D}Y\subset Y$
where
${\cal D}=(D_i)_{i\in I}$ 
runs through finite families
of Cartier divisors
of $Y$ as above.

{\rm 2.}
We say $f\colon V\to U$
is {\rm tamely ramified
with respect to} $Y$,
if there exists
a proper scheme $Y'$
over $Y$ containing $V$
as an open subscheme
such that
$\Sigma^+_{V/U}Y'=
\emptyset$.
\end{df}

We will define 
the wild ramification locus
$\Sigma_{V/U}Y$
as a closed subset of
$\Sigma^+_{V/U}Y$
in Definition \ref{dftmT}.
Since 
$\Sigma_{V/U}^{\cal D}Y
\supset
\Sigma_{V/U}^{{\cal D}'}Y$
for ${\cal D}\subset {\cal D}'$,
there exists a finite family
${\cal D}$ of Cartier divisors 
of $Y$ such that
$\Sigma^+_{V/U}Y
=
\Sigma_{V/U}^{\cal D}Y$
if $Y$ is quasi-compact.
In particular,
the condition $\Sigma^+_{V/U}Y
=\emptyset$ is equivalent to the
existence of
a finite family
${\cal D}$ 
of Cartier divisors
of $Y$ as in Definition \ref{dfSig}.1
satisfying
$\Sigma_{V/U}^{\cal D}Y
=\emptyset$.

\begin{lm}\label{lmVYY'}
Let $S$ be a scheme
and let
$$\begin{CD}
Y'@<{j'}<<V'@>>> U'\\
@VgVV@VVV@VVV\\
Y@<j<< V@>>> U
\end{CD}$$
be a commutative
diagram of separated schemes over $S$.
Assume that 
$V\to U$ and $V'\to U'$
are unramified
and that the canonical map
$V'\to V\times_UU'$
is an immersion.
Assume also
that $j\colon V\to Y$ and 
$j'\colon V'\to Y'$ are
open immersions.
If 
${\mathcal O}_{Y'}
\to j'_*{\cal O}_{V'}$
is an injection,
then we have
$g(\Sigma^+_{V'/U'}Y')
\subset
\Sigma^+_{V/U}Y$.
\end{lm}

{\it Proof.}
Let ${\mathcal D}=(D_i)_{i\in I}$
be a finite family of
Cartier divisors
of $Y$ such that
$D_i\cap V=\emptyset$
for every $i\in I$.
By the assumption that
${\mathcal O}_{Y'}
\to j'_*{\cal O}_V$
is injective,
the pull-backs
$g^*D_i$
defines a family
${\mathcal D}'=(D'_i)_{i\in I}$
of Cartier divisors
of $Y'$ satisfying
$D'_i\cap V=\emptyset$
for every $i\in I$.
Hence, the morphism
$(g\times g)^\sim\colon
(Y'\times_SY')^\sim_{{\cal D}'}
\to
(Y\times_SY)^\sim_{\cal D}$
is defined.
By the assumption that
$V'\to V\times_UU'$
is an immersion,
the inverse image
of $\Delta_V$
by $V'\times_{U'}V'
\to V\times_UV$
is $\Delta_{V'}$.
Hence 
$V'\times_{U'}V'
\setminus \Delta_{V'}
\subset
(Y'\times_SY')^\sim_{{\cal D}'}$
is a subset of 
the inverse image of
$V\times_UV
\setminus \Delta_V$.
Thus,
we have
$\overline{V'\times_{U'}V'\setminus \Delta_{V'}}
\subset
(g\times g)^{\sim-1}
(\overline{
V\times_UV\setminus \Delta_V})$
and 
$\Sigma_{V'/U'}^{{\cal D}'}Y'
\subset
g^{-1}(\Sigma_{V/U}^{\cal D}Y)$.
By taking the intersection,
the assertion follows.
\qed

Whether 
$\Sigma^+_{V/U}Y$
is empty or not may depend
on $Y$
as the following example shows.

\begin{ex}
Let $A$ be a ring
where $2$ is invertible
and let $V=
{\rm Spec}\
A[T_1^{\pm1},T_2^{\pm1}]
\subset
Y={\mathbf A}^2_A=
{\rm Spec}\
A[T_1,T_2]$.
We define an action 
of a cyclic group $G$ of
order $4$ by
 $T_1\mapsto -T_2,
T_2\mapsto T_1$.
Then, $V$
is a $G$-torsor over $U=V/G$.
For the blow-up $Y'$ of $Y$
at the $0$-section,
an elementary computation
shows that $\Sigma^+_{V/U}Y'
=\emptyset$
while $\Sigma^+_{V/U}Y$
consists of the $0$-section of $Y$.
\end{ex}

\subsection{Tamely ramified extension of 
valuation fields}\label{sstmv}

We will study
tame ramification
defined in the previous
subsection in detail
in Section \ref{sstvl}.
As preliminaries,
we first
study tamely ramified
extensions of
valuation fields
in this subsection
and limit of
compactifications
in the next subsection.

We recall the
definition of 
tamely ramified
extensions of
valuation fields.
For generality on valuation rings,
we refer to \cite[Chapitre 6]{B} and
\cite[Chapter VI]{ZS}.
If $L$ is a finite separable
extension of a field $K$
and if $B$ is the integral closure
in $L$ of a valuation ring $A$ of $K$,
then the map
from the
finite set of maximal ideals
of $B$ to
the set of valuation rings
of $L$ dominating $A$
sending
a maximal ideal ${\mathfrak m}$
to the local ring
$B_{\mathfrak m}$
is a bijection
\cite[Chapitre 6,
Section 8,
n$^{\rm o}$ 3, Remarque]{B}.

\begin{df}\label{dfvt}
Let $L$ be a finite
separable extension
of a field $K$
and $B$ be a valuation
ring of $L$.
We put
$A=B\cap K$.
Let $A^{\rm sh}$
and $B^{\rm sh}$
be strict henselizations
and let $K^{\rm sh}$
and $L^{\rm sh}$
be the fraction fields.

We say that $L$
is tamely ramified
over $K$ with respect to
$B$ if the degree
$[L^{\rm sh}:K^{\rm sh}]$
is invertible in $B$.
\end{df}

If the residue field
of $B$ is of characteristic 0,
an arbitrary
finite separable
extension $L$ over $K$
is tamely ramified
with respect
to $B$.

We recall some
standard terminologies
on inertia subgroups.
Let $M$ be a
finite Galois extension
of a field $K$
of Galois group
$G={\rm Gal}(M/K)$ 
and $A$ be
a valuation ring of $K$.
Let $C\subset M$ 
be the integral closure
of $A$ and  
${\mathfrak m}$
be a maximal ideal
of $C$.
The subgroup
$D=\{\sigma \in G\mid
\sigma({\mathfrak m})
={\mathfrak m}\}$
is called the decomposition
group of ${\mathfrak m}$
and 
$I={\rm Ker}(D
\to {\rm Aut}(C/
{\mathfrak m}))$
the inertia group.
The local ring
$C_{\mathfrak m}$
is a valuation ring.
Let $C^{\rm sh}
_{\mathfrak m}$
be a strict henselization
and let $M^{\rm sh}$
be the fraction field
of $C^{\rm sh}
_{\mathfrak m}$.
Then, the $I$-fixed
part of $C^{\rm sh}
_{\mathfrak m}$
is a strict henselization
$A^{\rm sh}$ of $A$.

We regard
the value group
$\Gamma_A=K^\times/A^\times$
as a subgroup of
$\Gamma_{\mathfrak m}=
M^\times/C_{\mathfrak m}^\times$.
Then the map
$I\times M^\times
\to (C/{\mathfrak m})^\times$
defined by
$(\sigma,c)
\mapsto 
\overline{\sigma(c)/c}$
induces a pairing
$I\times \Gamma_{\mathfrak m}/
\Gamma_A
\to (C/{\mathfrak m})^\times$.

\begin{lm}\label{lmvalP}
Let $M$ be a
finite Galois extension
of a field $K$
of Galois group
$G={\rm Gal}(M/K)$ 
and $A$ be
a valuation ring of $K$.
Let $I\subset G$
be the inertia group
of a maximal ideal
${\mathfrak m}$
of the integral closure
$C\subset M$ 
of $A$.
Let $p$ be the characteristic
of the residue field
$A/{\mathfrak m}_A$.

{\rm 1.}
{\rm (\cite[Chapter VI \S12 Corollary of Theorem 24]{ZS} )}
If $p=0$,
the pairing
$I\times \Gamma_{\mathfrak m}
/\Gamma_A\to (C/{\mathfrak m})^\times$
is a perfect pairing
of finite abelian groups.

{\rm 2.}
{\rm ([loc.\ cit.\ Theorems 24 and 25])}
Assume $p>0$
and let
$(\Gamma_{\mathfrak m}/
\Gamma_A)'$
denote the prime-to-$p$
part of
$\Gamma_{\mathfrak m}/
\Gamma_A$.
Then, the kernel $P$ of
the induced map
$I\to 
{\rm Hom}(\Gamma_{\mathfrak m}
/\Gamma_A
,(C/{\mathfrak m})^\times)$
is the unique $p$-Sylow subgroup
of $I$
and the induced pairing
$I/P\times
(\Gamma_{\mathfrak m}/
\Gamma_A)'
\to 
(C/{\mathfrak m})^\times$
is a perfect
pairing 
of finite abelian groups.
\end{lm}

For the rest of this subsection,
in the case $p=0$,
we put
$(\Gamma_{\mathfrak m}/
\Gamma_A)'=
\Gamma_{\mathfrak m}/
\Gamma_A$
and $P=1$.

\begin{cor}\label{corvalt}
Let $K$ be a field and
$L$ be a finite 
separable extension of $K$.
Let $A$ be a valuation ring of $K$
and $B$ be
the integral closure of
$A$ in $L$.
Then the following conditions are equivalent:

{\rm (1)}
For every maximal ideal
${\mathfrak m}$ of $B$,
$L$ is tamely ramified over $K$
with respect to $B_
{\mathfrak m}$.

{\rm (2)}
There exist non-zero elements
$t_1,\ldots,t_n$ of 
the maximal ideal ${\mathfrak m}_A$
and integers
$m_1,\ldots,m_n$ 
invertible in $A$
such that
the normalization
$B'$ of $A$ in
$L[S_1,\ldots,S_n]
/(S_1^{m_1}-t_1,\ \ldots,
S_n^{m_n}-t_n)$
is finite \'etale over
the normalization
$A'$ of $A$ in
$K[S_1,\ldots,S_n]
/(S_1^{m_1}-t_1,\ \ldots,
S_n^{m_n}-t_n)$.
\end{cor}

{\it Proof.}
(1)$\Rightarrow$(2):
We may assume $L$
is a Galois extension of $K$.
Let ${\mathfrak m}$
be a maximal ideal of $B$
and let $L^{\rm sh}$
be the fraction field
of the strict henselization
of $B$ at a geometric
point above ${\mathfrak m}$
and define 
$K^{\rm sh}$ similarly.
By Lemma \ref{lmvalP},
$L^{\rm sh}$ is an abelian
extension of 
$K^{\rm sh}$
and the pairing
${\rm Gal}(L^{\rm sh}/K^{\rm sh})
\times 
(\Gamma_B/\Gamma_A)'
\to (B/{\mathfrak m})^\times$
is a perfect
pairing of
finite abelian groups
of order prime to $p$.
We take
an isomorphism
${\mathbb Z}/m_1{\mathbb Z}
\oplus \cdots \oplus
{\mathbb Z}/m_n{\mathbb Z}\to
(\Gamma_B/\Gamma_A)'$
and
its lifting
$\gamma\colon {\mathbb Z}^n\to \Gamma_B$.
Let $e_1,\ldots,e_n\in 
{\mathbb Z}^n$
be the standard basis
and we take
elements
$t_1,\ldots,t_n\in {\mathfrak m}_A$
satisfying
$v_A(t_i)=m_i \gamma(e_i)$
for $i=1,\ldots,n$.
Then, we have
$L^{\rm sh}=
K^{\rm sh}(t_1^{1/m_1},
\ldots, t_n^{1/m_n})$
and the assertion follows.

(2)$\Rightarrow$(1):
Since $[L^{\rm sh}:K^{\rm sh}]$
divides $m_1\cdots m_n$,
the extension $L$ is
tamely ramified.
\qed

We give a criterion
for a finite separable extension
of valuation field
to be tamely ramified.

\begin{pr}\label{prvalt}
Let $L$ be a finite
separable extension of
a field $K$
and $B\subset L$
be a valuation ring of $L$.
We put $A=B\cap K$,
$U={\rm Spec}\ K
\subset
S={\rm Spec}\ A$
and
$V={\rm Spec}\ L
\subset
Y={\rm Spec}\ B$.
Then, the
following conditions
{\rm (1)}--{\rm (4)}
are equivalent:

{\rm (1)}
$L$ is tamely ramified
over $K$ with respect to
$B$.

{\rm (2)}
There
exists a finite
family ${\cal D}=
(D_i)_{i\in I}$
of Cartier divisors
of $Y$
such that
the intersection
$\Sigma_{V/U}^{\cal D}Y=
\overline W\cap
\Delta_Y^{\log}$
with
the closure
of $W=
V\times_UV\setminus \Delta_V$
in the log
product $(Y\times_SY)^\sim_
{\cal D}$
is empty.

{\rm (3)}
Let $M$ be
an arbitrary finite separable
extension of $L$
and $\sigma\colon
L\to M$
be a morphism over $K$
different from the inclusion.
For an arbitrary
valuation ring $C$
of $M$
dominating $B$ and $\sigma(B)$,
there exists a non-zero element
$b\in B$
such that
$\sigma(b)/b\not\equiv 1
\bmod {\mathfrak m}_C$.

{\rm (4)}
Let $M$ be a
finite Galois extension
of $K$ 
of Galois group
$G={\rm Gal}(M/K)$ 
containing
$L$ as a subextension
and 
${\mathfrak m}$
be a maximal ideal
of the integral closure
$C\subset M$ of $B$ 
such that
$B\cap {\mathfrak m}
={\mathfrak m}_B$.
Then, 
the subgroup $H={\rm Gal}(M/L)
\subset G={\rm Gal}(M/K)$ 
contains the conjugates 
of the $p$-Sylow
subgroup $P$ of the 
inertia group
$I$ of ${\mathfrak m}$.
\end{pr}

In (2), we did not say
that $V\to U$
is tamely ramified
with respect to $Y$
because
the canonical map
$V\to Y$
may not be an open immersion.

{\it Proof.}
(2)$\Rightarrow$(3)
Let ${\cal D}=
(D_i)_{i\in I}$
be a finite family
of Cartier divisors
of $Y={\rm Spec}\ B$
such that
$\Sigma_{U/V}^{\cal D}Y
=\emptyset$.
Since $C$
dominates
$B$ and $\sigma(B)$,
the compositions
$L^\times
\to M^\times 
\to \Gamma_C$
and
$L^\times
\overset{\sigma}\to 
M^\times 
\to \Gamma_C$
are equal.
Hence,
the map
$\gamma=(i,\sigma)
\colon Z={\rm Spec}\ C
\to Y\times_SY$
induces
a map $\tilde
\gamma
\colon Z
\to (Y\times_SY)^\sim_{\cal D}$
to the log product.
Since $\sigma$
is different from
the inclusion,
the image
$\tilde \gamma(Z)$
is in the closure
$\overline W$
of $W=V\times_UV\setminus 
\Delta_V$.
Hence, we have
$\gamma(Z)\cap
\Delta_Y^{\log}
\subset
\Sigma_{U/V}^{\cal D}Y
=\emptyset$.

For $i\in I$,
let $b_i\in B$
an element defining
the divisor $D_i$.
Then,
the closed subscheme
$\Delta_Y^{\log}
\subset
(Y\times_SY)^{\sim}_{\cal D}$
is defined by
the ideal
$(b\otimes 1-1\otimes b; b\in B,\
(b_i\otimes 1)/(1\otimes b_i)-1; i\in I)$.
Hence the closed subscheme
of $Z$ defined by
$(\sigma(b)-b; b\in B,\
\sigma(b_i)/b_i-1; i\in I)$
is empty.
Namely the ideal of $C$
generated by
$\sigma(b)-b$ for $b\in B$
and $\sigma(b_i)/b_i-1$
for $i\in I$
contains a unit.
Thus the assertion is proved.

(3)$\Rightarrow$(2)
We put
$W=V\times_UV\setminus 
\Delta_V
=\coprod_{j\in J}
{\rm Spec}\ M_j$
where $M_j$ are fields.
We regard $M_j$
as an extension of
$L$ by the map
defined by the 
first projection
and let 
$\sigma_j\colon 
L\to M_j$
be the map
defined by the second projection.
For each $j\in J$,
the set
$\{C_{ij}\mid i\in I_j\}$
of valuation
rings of $M_j$
dominating 
both $B$ and $\sigma(B)$
is a finite set.
For each 
$C_{ij}$,
take a non-zero element 
$b_{ij}\in B$
such that
$\sigma_j(b_{ij})/
b_{ij}\not\equiv 1
\bmod {\mathfrak m}_{C_{ij}}$
and let
$D_{ij}$
be the Cartier divisor
of $Y$ defined by $b_{ij}$.
We put
$I=\coprod_{j\in J}I_j$
and
let ${\cal D}
=(D_{ij})_{j\in J,i\in I_j}$
be the family
of Cartier divisors.

We show that 
$\Sigma_{V/U}^{\cal D}Y$
is empty.
Let 
$Z_j$ be the closure
of 
${\rm Spec}\ M_j$
in the log product
$(Y\times_SY)^\sim_{\cal D}$.
Then, we have
$\overline W
=\bigcup_{j\in J}Z_j$
and 
$\Sigma_{V/U}^{\cal D}
Y=
\bigcup_{j\in J}
(Z_j\cap \Delta_Y^{\log})$.
Hence, if
$\Sigma_{V/U}^{\cal D}Y$
was not empty,
the intersection
$Z_j\cap \Delta_Y^{\log}$
would contain
the closed point $y$
of $Y$
for some
index $j\in J$.
Take a valuation ring $C$
of $M_j$
dominating the local ring
${\cal O}_{Z_j,y}$.
Then, $C$ dominates
$B$ and $\sigma_j(B)$.
Hence it is equal to
$C_{ij}$ for some 
$i\in I_j$
and $\sigma_j(b_{ij})
/b_{ij}-1$
is a unit of $C=C_{ij}$.
On the other hand,
since the image $y$ of
the closed
point of $C$ is in
$\Delta_Y^{\log}$, 
the ideal
$(\sigma_j(b)-b; b\in B,\
\sigma_j(b_{ik})/b_{ik}-1; 
k\in J,i\in I_k)$
is different from $C$,
as in the proof of
(2)$\Rightarrow$(3).
Thus, we obtain a contradiction.

(1)$\Rightarrow$(3)
By replacing $K$ by
an unramified extension,
we may assume that
the residue field of $L$
is a purely inseparable extension
of the residue field of $K$.
We may assume
$M$ is a Galois extension
and extend $\sigma
\colon L\to M$
to an element
of the Galois group
$G={\rm Gal}(M/K)$.
Then, since 
both $C$ and $\sigma(C)$
dominates $\sigma(B)$,
there exists $\tau\in G$
such that
$\tau(C)=\sigma(C)$
and that $\tau|_{\sigma(L)}
\colon \sigma(L)\to M$
is the inclusion.
Replacing $\sigma$
by $\tau^{-1}\sigma$
if necessary,
we may assume 
$\sigma(C)=C$.
Namely, by the assumption on the residue field,
we may assume that $\sigma$
is in the inertia
group $I\subset G$
of the maximal ideal
of ${\mathfrak m}_C$.
By the assumption
that $\sigma|_L
\neq {\rm id}_L$,
it is not an element of
the subgroup 
$H={\rm Gal}(M/L)\subset G$
corresponding to $L$.

We put
$\bar C=C/{\mathfrak m}_C$ and
we consider the
perfect pairings
$I/P
\times
(\Gamma_{{\mathfrak m}_C}
/\Gamma_A)'
\to \bar C^\times$
and
$(I\cap H)/(P\cap H)
\times
(\Gamma_{{\mathfrak m}_C}
/\Gamma_B)'
\to \bar C^\times$.
Since $[L^{\rm sh}:
K^{\rm sh}]
=[I:I\cap H]$
is invertible in $B$,
the induced
pairing
$I/(I\cap H)
\times 
(\Gamma_B/\Gamma_A)'
\to \bar C^\times$
is perfect.
Since $\sigma|_L$
is not the identity,
there exists an element
$b\in L^\times$
such that
$\sigma(b)/b
\not\equiv 1\bmod 
{\mathfrak m}_C$.

(3)$\Rightarrow$(4)
Let $\sigma$ be an element
of a conjugate $\tau P\tau^{-1}$ of $P$.
We regard
$L$ as a subfield of $M$
by $\tau|_L\colon L
\to M$.
Then, the maximal ideal
${\mathfrak m}'=\tau(
{\mathfrak m})$
satisfies
$\sigma({\mathfrak m}')=
{\mathfrak m}'$
and
$\tau|_L^{-1}(
{\mathfrak m}')
=
(\sigma\circ 
\tau|_L)^{-1}(
{\mathfrak m}')$
is equal to 
$B\cap {\mathfrak m}
=
{\mathfrak m}_B$.
If $\sigma$ was not
an element of $H$,
the condition (3) 
would imply
the existence of
an element $b\in L^\times$
such that
$\sigma(\tau(b))/
\tau(b)
\not\equiv 1
\bmod {\mathfrak m}'$.
This implies
that the
order of $\sigma$
is not a power
of $p$.
Thus, we get a contradiction.

(4)$\Rightarrow$(1)
Let $M$ be a finite
Galois extension of
$K$ containing $L$
as a subfield
and we put
$G={\rm Gal}(M/K)
\supset
H={\rm Gal}(M/L)$.
We take a maximal ideal
${\mathfrak m}$
of the integral closure
$C$ above 
the maximal ideal of
$B$
and let $I$ be
the inertia group of 
${\mathfrak m}$.
Then, 
the inertia group
$I$ is identified
with the Galois group
${\rm Gal}(M^{\rm sh}
/K^{\rm sh})$
and $H\cap I$
is the subgroup of
$I$ corresponding
to the field
$L^{\rm sh}$.
Hence
$P\subset H$
implies that
$[L^{\rm sh}:
K^{\rm sh}]=
[I:I\cap H]$
is prime to $p$.
\qed

\subsection{Limit of
compactifications and valuation rings}\label{sslim}

We study local rings
of the limit of compactifications.
Let $S$ be a separated
noetherian scheme
and $U$ be a 
separated scheme of
finite type over $S$.
We consider the
category ${\cal C}_{U/S}$
of compactifications
of $U$ over $S$.
Namely an object
of ${\cal C}_{U/S}$
is a pair
$(X,j)$ consisting of a
proper scheme
$X$ over $S$
and an open immersion
$j\colon U\to X$
over $S$.
A morphism
$f\colon (X',j')
\to (X,j)$
is a morphism
$f\colon X'\to X$
of schemes over $S$
such that $f\circ j'=j$.
In the following,
we omit $j$ from the notation 
and write simply $X$ for a 
compactification $(X, j)$.

\begin{lm}\label{lmCUS}
Let $S$ be a separated
noetherian scheme
and $U$ be a 
separated scheme of
finite type over $S$.

{\rm 1.} The
category ${\cal C}_{U/S}$
is cofiltered.
In particular,
it is non-empty.

{\rm 2.}
The objects 
containing $U$
as the complement
of a Cartier divisor
are cofinal in
${\cal C}_{U/S}$.
\end{lm}

{\it Proof.}
By Nagata's embedding theorem
\cite{nagata},
the category
${\cal C}_{U/S}$
is non-empty.
Since a blow-up
$X'\to X$ is proper,
the objects
containing $U$
as the complement
of a Cartier divisor
are cofinal in
${\cal C}_{U/S}$.
For objects
$(X,j)$ and $(X',j')$
of 
${\cal C}_{U/S}$,
a morphism
$X\to X'$ is unique if
${\mathcal O}_X
\to j_*{\cal O}_U$ is
injective.
If $X$ and $X'$
are objects of
${\cal C}_{U/S}$,
the schematic closure
of the diagonal map
$U\to X\times_SX'$
is an object
of ${\cal C}_{U/S}$.
Hence, the
category ${\cal C}_{U/S}$
is cofiltered.
\qed

We consider the projective limit
$\widetilde X
=\varprojlim_{X\in {\cal C}_{U/S}}X$
in the category
of locally ringed spaces.
The underlying
topological space
$\widetilde X$
is known to be quasi-compact
\cite[Theorem 5.14]{FK}.
For a point $\tilde x=(x_X)
\in \widetilde X$,
we have
${\cal O}_{\widetilde X,\tilde x}=
\varinjlim_{X\in {\cal C}^o_{U/S}}
{\cal O}_{X,x_X}$.
We will describe the limit
$\widetilde X$
and the local rings
${\cal O}_{\widetilde X,\tilde x}$
in terms of valuation rings.

\begin{df}\label{dfval}
Let $U$ be a scheme,
$u\in U$ be a point and
$A$ be a valuation
ring of the residue
field $\kappa(u)$.

{\rm 1.}
We say $A$ is $U$-external,
if $A\subsetneqq
\kappa(u)$ and if
there exists no intermediate
ring $A\subset A'\subsetneqq
\kappa(u)$
such that
the map $u\to U$ is extended
to ${\rm Spec}\ A'\to U$.

{\rm 2.}
Let $U\to S$
be a separated morphism of
schemes.
We say $A$ is $S$-integral,
if the composition
$u\to U\to S$
is extended to
a morphism
${\rm Spec}\ A\to S$.
\end{df}
\setcounter{equation}0

Let $u\in U$
and $A\subset
 \kappa(u)$
be an $S$-integral
valuation ring.
Then, for an object $X$
of ${\cal C}_{U/S}$,
the inclusion
$u\to U$
is uniquely extended to
a morphism ${\rm Spec}\ A
\to X$ over $S$
by the valuative criterion
of properness.
The images $x_A\in X\setminus U$ 
of the closed point
of ${\rm Spec}\ A$
define a point $\tilde x_A=
(x_A)$ of the projective limit
$\widetilde X=
\varprojlim_{X\in {\cal C}_{U/S}}X$.
Thus, we obtain a natural map
\begin{equation}
\coprod_{u\in U}
\{\text{\rm $S$-integral 
valuation ring of }\kappa(u)\}
\to
\widetilde X.
\label{eqvalX}
\end{equation}

An $S$-integral valuation ring
$A$ of $\kappa(u)$
is $U$-external if and only if
$\{u\}=U\times_X{\rm Spec}\ A
\subsetneqq
{\rm Spec}\ A$
for an object $X$ of
${\mathcal C}_{U/S}$.
Consequently, the map
(\ref{eqvalX}) induces
\begin{equation}
\coprod_{u\in U}
\{U\text{\rm-external and $S$-integral 
valuation ring of }\kappa(u)\}
\to
\widetilde X\setminus U.
\label{eqval}
\end{equation}
We show that the map (\ref{eqval})
is a bijection.

\begin{lm}\label{lmval}
{\rm (cf.\ \cite[5.4]{FK})}
Let $S$ be a separated
noetherian scheme
and $U$ be a separated
scheme of finite type
over $S$.
Let $\tilde x=(x_X)_X
\in \widetilde X\setminus U$
be a point in the complement and
put ${\cal O}_{\widetilde X,\tilde x}=
\varinjlim {\cal O}_{X,x_X}$.

Then, there exists a unique point
$u\in U$ such that
$U\times_X
{\rm Spec}\ {\cal O}_{\widetilde X,\tilde x}
= {\rm Spec}\ {\cal O}_{U,u}$
for every object
$X$ of ${\cal C}_{U/S}$.
The canonical map
${\cal O}_{\widetilde X,\tilde x}
\to 
{\cal O}_{U,u}$ is injective
and its image is
the inverse image of
a $U$-external and $S$-integral
valuation ring $A$ of $\kappa(u)$.
For each object $X$
of ${\cal C}_{U/S}$,
the point $x_X$
is the image of the closed point
of ${\rm Spec}\ A$
by the unique map
${\rm Spec}\ A
\to X$ over $S$
extending $u\to U$.
\end{lm}

{\it Proof.}
Let $\tilde x=(x_X)_X
\in \widetilde X\setminus U$
be a point in the complement.
For a morphism
$X'\to X$ of ${\mathcal C}_{U/S}$,
we have
$U\times_XX'=U$
if $U$ is dense in $X'$.
Hence the inverse image
$U\times_X
{\rm Spec}\ 
{\cal O}_{\widetilde X,\tilde x}$
is independent of $X$.
Thus, to show the 
existence of $u\in U$
such that 
$U\times_X
{\rm Spec}\ {\cal O}_{\widetilde X,\tilde x}
= {\rm Spec}\ {\cal O}_{U,u}$,
it suffices to show
the existence 
for one object $X
\in {\cal C}_{U/S}$.
Take an object $X
\in {\cal C}_{U/S}$.
We may assume
$U$ is the complement
of a Cartier divisor $D\subset X$.

By \cite[Proposition 5.12]{FK},
the local ring
${\cal O}_{\widetilde X,\tilde x}$
is $I_{\tilde x}$-valuative 
for $I_{\tilde x}={\cal I}_{D,x_X}
{\cal O}_{\widetilde X,\tilde x}$
in the terminology loc.\ cit.
Hence, 
$\Gamma(U\times_X{\rm Spec}\ 
{\cal O}_{\widetilde X,\tilde x}
,{\cal O}_{\widetilde X})=
\varinjlim_{X'}
\Gamma(U\times_{X'}
{\rm Spec}\ {\cal O}_{X',x_{X'}},
{\cal O}_{X'})$
is a local ring
further by
\cite[Proposition 5.11]{FK}.
Since
${\rm Spec}\ {\cal O}_{X',x_{X'}}
\to X'\times_X
{\rm Spec}\ {\cal O}_{X,x_X}$ is a limit of
open immersions,
its restriction
$U\times_{X'}
{\rm Spec}\ {\cal O}_{X',x_{X'}},
\to
U\times_X
{\rm Spec}\ {\cal O}_{X,x_X}$
is also a limit of
open immersions.
Hence
the local ring
$\Gamma(U\times_X{\rm Spec}\ 
{\cal O}_{\widetilde X,\tilde x}
,{\cal O}_{\widetilde X})$
is a localization of
$\Gamma(U\times_X{\rm Spec}\ 
{\cal O}_{X,x_{X}},
{\cal O}_X)$
and is equal to the local ring
${\cal O}_{U,u}$ at a point $u\in U$.
Further by 
\cite[Proposition 5.11]{FK},
the canonical map
${\cal O}_{\widetilde X,\tilde x}\to
{\cal O}_{U,u}$ is an injection
and its image is the inverse 
image of a valuation ring
$A$ of $\kappa(u)$
by the surjection
${\cal O}_{U,u}\to \kappa(u)$.

Since 
$\Gamma(U\times_X{\rm Spec}\ 
{\cal O}_{\widetilde X,\tilde x}
,{\cal O}_{\widetilde X})=
{\cal O}_{U,u}$, 
the valuation ring
$A$
is $U$-external.
Since $A$ is $X$-integral,
it is $S$-integral.
The image $x_X$
of the closed point
of ${\rm Spec}\ 
{\cal O}_{\widetilde X,\tilde x}$
is the same as
the image
of the closed point
of ${\rm Spec}\ A$
by the induced map
${\rm Spec}\ A\to X$.
\qed

\begin{cor}\label{corval}
The map {\rm (\ref{eqval})}
is a bijection.
The inverse is 
defined by
sending $\tilde x$
to the
valuation ring 
${\cal O}_{\widetilde X,\tilde x}
/{\mathfrak m}_u
={\rm Image}(
{\cal O}_{\widetilde X,\tilde x}
\to 
\kappa(u))\subset
\kappa(u)$,
in the notation of
Lemma {\rm \ref{lmval}}.
\end{cor}

{\it Proof.}
Let $u\in U$ be a point and
$A$ be a $U$-external and $S$-integral 
valuation ring of $\kappa(u)$.
Let $\tilde x\in \widetilde X$ 
be the point defined by
the images of the closed point
of ${\rm Spec}\ A$.
We consider
$A'=
{\cal O}_{\widetilde X,\tilde x}
/{\mathfrak m}_{u'}
\subset
\kappa(u')$
as in
Lemma {\rm \ref{lmval}}.
We have a natural local homomorphism
${\cal O}_{\widetilde X,\tilde x}\to A$.
By \cite[Proposition 5.11]{FK},
the ideal
${\mathfrak m}_{u'}
\subset {\cal O}_{\widetilde X,\tilde x}$
is the intersection
$\bigcap_nI_{\tilde x}^n$.
Since $A$ is $U$-external,
we have 
$\bigcap_nI_{\tilde x}^nA=0$.
Hence, it induces a
local homomorphism $A'\to A$.
Since $A'$ is also $U$-external,
we obtain $u'=u$.
Further,
since the valuation ring
$A\subset \kappa(u)$ dominates
$A'\subset \kappa(u)$,
we obtain $A'=A$.
\qed

\subsection{Tame ramification
and valuation rings}\label{sstvl}

We give a criterion for
tame ramification
in terms of valuation rings,
in
Proposition \ref{prtY}.
We slightly generalize
Definition
\ref{dfSig}.2.

\begin{df}\label{dftmT}
Let $f\colon V\to U$
be an unramified morphism
of finite type of 
separated schemes.

{\rm 1.}
Let $Y$ be a separated
scheme of finite type
containing $V$
as an open subschme.
We define 
{\rm the wild ramification
locus}
$\Sigma_{V/U}Y$
to be the closed subset
$$\Sigma_{V/U}Y
=
\bigcap_{Y'}
\pi_{YY'}(\Sigma_{V/U}^+Y')$$
where $\pi_{YY'}\colon Y'\to Y$
runs through objects of 
${\mathcal C}_{V/Y}$.

{\rm 2.}
Let $T$ be a separated
noetherian scheme
and assume that
$V$ is a separated
scheme of finite type over $T$.
We say $f\colon V\to U$
is {\rm tamely ramified
with respect to} $T$,
if there exists
an object $Y$
of ${\cal C}_{V/T}$
such that
$\Sigma_{V/U}Y$
is empty.
\end{df}

Since the category
${\cal C}_{V/T}$
is cofiltered,
the map
$f\colon V\to U$
is tamely ramified
with respect to 
$T$ if and only if
$\Sigma_{V/U}Y$
is empty for
every object
of ${\cal C}_{V/T}$.
If $U$ is a scheme over
$S$ and if finite \'etale morphisms
$V\to U$
and
$V'\to U$
are tamely ramified with respect
to $S$,
then the fiber product
$V\times_UV'\to U$
is also tamely ramified with respect
to $S$.
In particular, if a finite \'etale morphism
$V\to U$ of connected normal schemes
is tamely ramified with respect
to $S$,
its Galois closure $W\to U$
is also tamely ramified with respect
to $S$.

\begin{lm}\label{lmtYS}
Let $T$ be a separated
noetherian scheme,
$V$ be a separated
scheme of finite type over $T$
and
let $f\colon V\to U$
be an unramified morphism
of finite type of 
separated schemes.
Then, the objects $Y$ of
${\mathcal C}_{V/T}$
such that
there exists 
a finite family ${\cal D}$
of Cartier divisors
satisfying
$\Sigma_{V/U}Y=
\Sigma^{\cal D}_{V/U}Y$
are cofinal
in ${\mathcal C}_{V/T}$.
\end{lm}

{\it Proof.}
For an object $Y$ of
${\mathcal C}_{V/T}$,
let $\pi_Y\colon
\widetilde Y=
\varprojlim_{Y'\in
{\mathcal C}_{V/T}} Y'\to Y$
denote the projection.
Since ${\mathcal C}_{V/T}$ is cofiltered
and 
$\widetilde Y$ is quasi-compact
\cite[Theorem 5.14]{FK},
the objects $Y$ such that
$\pi_Y^{-1}(\Sigma^+_{V/U}Y)
=
\bigcap_{Y'}
\pi_{Y'}^{-1}(\Sigma^+_{V/U}Y')$
and $Y=\pi_Y(\widetilde Y)$
are cofinal in
${\mathcal C}_{V/T}$
by Lemma \ref{lmVYY'}.
For such $Y$,
we have
$\Sigma^+_{V/U}Y
=
\pi_Y(\pi_Y^{-1}(\Sigma^+_{V/U}Y))
=
\bigcap_{Y'}
\pi_Y(\pi_{Y'}^{-1}(\Sigma^+_{V/U}Y'))$.
Hence the assertion follows
from the quasi-compactness of $Y$.
\qed

\begin{lm}\label{lmtY}
Let $T$ be a separated
noetherian scheme
and $V$ be a separated
scheme of finite type over $T$.
Let
$f\colon V\to U$
be an unramified
separated morphism
of finite type of
schemes.
Let $v\in V$ be a point 
and $B$ be a $T$-integral valuation ring
of $\kappa(v)$.
We consider the conditions:

{\rm (1)}
There exists a proper
scheme $Y$ over $T$
containing $V$
as an open subscheme
such that
the closed subset
$\Sigma^+_{V/U}Y\subset Y$
does not meet 
the image of the map
${\rm Spec}\ B\to Y$.

{\rm (2)}
The finite separable
extension $\kappa(v)$
over $\kappa(f(v))$
is tamely ramified
with respect to $B$.

The condition
{\rm (1)} implies {\rm (2)}.
If $B$ is $V$-external  {\rm (Definition \ref{dfval}.1)},
then the conditions 
{\rm (1)} and {\rm (2)}
are equivalent.
\end{lm}

{\it Proof.}
(1)$\Rightarrow$(2) 
The condition (1) implies that
there exists
a finite family
${\cal D}=(D_i)_{i\in I}$
of Cartier divisors
of ${\rm Spec}\ B$
such that
$\Sigma_{v/f(v)}^{\cal D}
{\rm Spec}\ B$ is empty,
as the pull-back by 
${\rm Spec}\ B\to Y$.
Then, the condition (2)
is satisfied by
Proposition \ref{prvalt}
(2)$\Rightarrow$(1).

(2)$\Rightarrow$(1) 
Assume that $B$ is $V$-external.
By Proposition \ref{prvalt}
(1)$\Rightarrow$(2),
there exists
a finite family
${\cal D}=(D_i)_{i\in I}$
of Cartier divisors
of ${\rm Spec}\ B$
such that
$\Sigma_{v/f(v)}^{\cal D}
{\rm Spec}\ B$ is empty.
Let $\widetilde Y
=\varprojlim_{Y\in
{\cal C}_{V/T}}Y$
be the limit of
compactifications
and $\tilde y=(y_Y)_Y
\in \widetilde Y$
be the point
corresponding
to $B\subset \kappa(v)$.
We take a non-zero divisor
$f_i\in B$
defining $D_i$ 
and a lifting
$\tilde f_i\in {\cal O}_{
\widetilde Y,\tilde y}$
for each $i\in I$.
Since ${\cal O}_{
\widetilde Y,\tilde y}=
\varinjlim_Y{\cal O}_{Y,y_Y}$,
there exist an object
$Y$ of ${\cal C}_{V/T}$,
an open neighborhood 
$W$ of $y_Y$
and non-zero divisors
$g_i\in 
\Gamma(W,{\cal O}_Y)$
invertible on
$V\cap W$
sent to $f_i$ for $i\in I$.
By replacing
$Y$ by the blow-up
of the closure of the divisors
of $g_i$,
we may assume that
there exists
a finite family
${\cal D}_Y=(D_{Y,i})_{i\in I}$
of Cartier divisors
of $Y$
such that
$D_{Y,i}\cap V=\emptyset$
and
$D_{Y,i}\cap W$
is defined by $g_i$
for each $i\in I$.
Then the inverse image
of $\Sigma_{V/U}^{{\cal D}_Y}Y$
by the map
${\rm Spec}\ B\to Y$
is equal to
$\Sigma_{v/f(v)}^{\cal D}
{\rm Spec}\ B$
and hence is empty.
Thus the assertion is proved.
\qed

\begin{pr}\label{prtY}
Let $T$ be a separated
noetherian scheme
and $V$ be a separated
scheme of finite type over $T$.
For an unramified
separated morphism
$f\colon V\to U$
of finite type of schemes,
the following conditions
are equivalent:

{\rm (1)}
$f\colon V\to U$
is tamely ramified with
respect to $T$.

{\rm (2)}
For every point $v\in V$
and for every
$T$-integral and $V$-external
valuation ring 
{\rm (Definition \ref{dfval})}
$B$ of $\kappa(v)$,
the extension $\kappa(v)$
over $\kappa(f(v))$
is tamely ramified
with respect to $B$.

{\rm (3)}
For every point $v$ of $V$ and 
for every $T$-integral valuation ring $B$ 
of $\kappa(v)$, the
extension $\kappa(v)$ over 
$\kappa(f(v))$ is tamely ramified with respect to $B$.
\end{pr}

{\it Proof.}
(3)$\Rightarrow$(2) 
Clear.

(1)$\Rightarrow$(3) 
It follows from Lemma \ref{lmtY}
(1)$\Rightarrow$(2).

(2)$\Rightarrow$(1) 
Let $\widetilde Y
=\varinjlim_{Y\in
{\cal C}_{V/T}}Y$
be the limit of
compactifications
and let
$\pi_Y\colon
\widetilde Y
\to Y$
denote the projection.
By Lemma \ref{lmtY}
(2)$\Rightarrow$(1),
for every point $\tilde y
\in \widetilde Y$
of the boundary,
there exists a proper
scheme $Y$ over $T$
containing $V$
as an open subscheme
such that the inverse image
of $\Sigma^+_{V/U}Y\subset Y$
does not contain $\tilde y$.
In other words,
the intersection
$\bigcap_{Y\in {\cal C}_{V/T}}
\pi_Y^{-1}
(\Sigma^+_{V/U} Y)\subset 
\widetilde Y$ is empty.
Since $\widetilde Y$
is quasi-compact
and since
the category ${\cal C}_{V/T}$
is cofiltered,
there exists
an object $Y$ of
${\cal C}_{V/T}$
such that
$\pi_Y^{-1}
(\Sigma^+_{V/U} Y)$
is empty.
Thus the assertion follows.
\qed

\begin{cor}\label{corptp}
Let $T$ be a separated
noetherian scheme
and $V$ be a separated
scheme of finite type over $T$.
Let 
$f\colon V\to U$
be a separated
unramified
morphism of
finite type.
If one of the
following conditions
is satisfied,
then 
$f\colon V\to U$
is tamely ramified 
with respect to $T$.

{\rm (1)} 
$T$ is a scheme
over ${\mathbb Q}$.

{\rm (2)} 
$V$ is a
$G$-torsor over $U$
for a finite group $G$
of order
invertible on $T$.
\end{cor}

{\it Proof.}
Let $v\in V$
be a point
and $B$ be a
$T$-integral and $V$-external
valuation ring 
of $\kappa(v)$.
Then, 
either of the
conditions (1) and (2)
implies
that
the extension $\kappa(v)$
over $\kappa(f(v))$
is tamely ramified
with respect to $B$.
Hence the assertion follows
from Proposition
\ref{prtY}
(2)$\Rightarrow$(1).
\qed

\begin{cor}\label{cortgen}
Let $T$ be a separated
noetherian scheme
and $V$ be a separated
scheme of finite type over $T$.
Let
$f\colon V\to U$
be an unramified
and separated 
dominant morphism
of schemes of finite type.
Assume that $U$ and $V$
are integral and let
$\xi\in U$ and $\eta\in V$
be the generic point respectively.
We consider
the following conditions:

{\rm (1)}
$f\colon V\to U$
is tamely ramified with
respect to $T$.

{\rm (2)}
The extension
$\kappa(\eta)$
of $\kappa(\xi)$
is tamely ramified
with respect to
an arbitrary
$T$-integral
valuation ring
of $\kappa(\eta)$.

Then, the implication
{\rm (1)}$\Rightarrow${\rm (2)}
always holds.
The other implication
{\rm (2)}$\Rightarrow${\rm (1)}
holds if $V$ is regular.
\end{cor}

{\it Proof.}
(1)$\Rightarrow$(2)
It follows from Proposition 
\ref{prtY}
(1)$\Rightarrow$(3).

(2)$\Rightarrow$(1)
Assume $V$ is regular.
Let $v\in V$ and
let $B\subset \kappa(v)$
be a $T$-integral valuation ring.
Since $V$ is regular,
there exists a valuation ring
$B_0\subset \kappa(\eta)$ 
dominating ${\cal O}_{V,v}$
such that the residue field
of $B_0$ is $\kappa(v)$.
Let $B_1$ be the inverse image
of $B\subset \kappa(v)$ by the surjection
$B_0\to \kappa(v)$.
Then $B_1\subset \kappa(\eta)$
is a valuation ring and $T$-integral.
Thus,
by Proposition
\ref{prtY}, 
it suffices to show that
$\kappa(v)$
is tamely ramified over $\kappa(u)$
with respect to $B$
assuming $\kappa(\eta)$
is tamely ramified over $\kappa(\xi)$
with respect to $B_1$.

We put $A=B\cap \kappa(u)$
and $A_1=B_1\cap \kappa(\xi)$.
Then, the map $A_1\to A$
is a surjection and
we have
$\kappa(v)=\kappa(u)\otimes_{A_1}B_1$.
Hence we have $B=A\otimes_{A_1}B_1$.
Let $A_1^{\rm sh}$
and $B_1^{\rm sh}$
be the strict henselizations
and $\kappa(\xi)_1^{\rm sh}$
and $\kappa(\eta)_1^{\rm sh}$
be their fraction fields.
By the assumption,
the degree
$[\kappa(\eta)_1^{\rm sh}:
\kappa(\xi)_1^{\rm sh}]$
is invertible in $B$.
Let $A^{\rm sh}$
and $B^{\rm sh}$
be the strict henselizations
and $\kappa(u)^{\rm sh}$
and $\kappa(v)^{\rm sh}$
be their fraction fields.
Then,
$A^{\rm sh}$
and $B^{\rm sh}$
are quotients
of $A_1^{\rm sh}$
and of $B_1^{\rm sh}$
and the canonical map
$\kappa(u)^{\rm sh}\otimes_
{A_1^{\rm sh}}B_1^{\rm sh}
\to \kappa(v)^{\rm sh}$
is an isomorphism.
Hence,
we have
$[\kappa(v)^{\rm sh}:
\kappa(u)^{\rm sh}]
=[\kappa(\eta)_1^{\rm sh}:
\kappa(\xi)_1^{\rm sh}]$
and the assertion follows.
\qed

The following example
shows that the condition (2)
need not imply (1)
if we replace ``regular'' by ``normal''.

\begin{ex}
{\rm
Let $k$ be an algebraically closed
field of characteristic $p>0$,
$E$ be an ordinary elliptic curve
over $k$
and ${\cal L}$ be
a very ample invertible ${\cal O}_E$-module
on $E$ e.\ g.\ 
${\cal O}(3\cdot[0])$.
Let $X_0={\rm Spec}\ 
\bigoplus \Gamma(E,{\cal L}^{\otimes n})$
be the affine cone.
The blow-up $X_1$ of $X_0$ at
the origin is the line bundle
over $E$ associated to ${\cal L}$.
Let $Y_1\to X_1$ be
the base change of
the map $V:E^{(p)}\to E$
and $Y_0\to X_0$ be
the Stein factorization of
the composition $Y_1\to X_1\to X_0$.

Let $C\to C'$ be a 
finite \'etale cyclic 
covering of affine curves
of degree $p$.
We assume that the map
$\overline C\to \overline C'$
of the compactifications
is wildly ramified.
We put $V=Y_0\times C$,
$Y=Y_0\times \overline C$
and we consider the action of
$E[p](k)\times {\rm Gal}(C/C')\simeq
({\mathbb Z}/p{\mathbb Z})^2$
on $Y$.
Let $G\subset 
E[p](k)\times {\rm Gal}(C/C')$
be a diagonal subgroup
and $X=Y/G$ be the quotient.
Since the action of $G$
on $V$ is free,
the map $f\colon V\to U=V/G$
is finite and \'etale.

We show that $Y\to X$
satisfies (2).
Since $V\colon E^{(p)}\to E$
is finite \'etale,
the blow-up $Y_1\to X_1$
of $Y_0\to X_0$ is finite \'etale.
Hence, the action of
$G$ on $Y_1\times \overline C$
is free and the map
$Y_1\times \overline C\to
(Y_1\times \overline C)/G$
is finite \'etale.

We show that $Y\to X$
does not satisfy (1).
The inclusion
$\overline C\to Y$
at the origin of $Y_0$
induces
$\overline C'\to X$.
By the assumption,
the covering
$\overline C\to
\overline C'$
is widely
ramified at the boundary
$c\in \overline C
\setminus C$.
Since the valuation ring
${\cal O}_{\overline C,c}$ is
$Y$-integral,
the assertion follows.}
\end{ex}

\subsection{Tame ramification
and Kummer coverings}\label{sstKm}

We consider a finite
\'etale morphism
$f\colon V\to U$
of separated
schemes of finite
type over a separated
noetherian scheme $S$.
We study the condition
for $f\colon V\to U$
to be tamely ramified
with respect to $S$
in terms of Kummer coverings.

\begin{df}\label{dflur}
Let $X$ be a scheme,
$U\subset X$
be an open subscheme
and
$f\colon V\to U$
be a finite \'etale
morphism.

{\rm 1.}
Let $x\in X\setminus U$
be a point of boundary.
We say 
$f\colon V\to U$
is of Kummer type
at $x$ if
there exist
an open neighborhood $W$ of $x$,
functions
$t_1,\ldots,t_n\in \Gamma(W,{\cal O}_X)$
invertible on $U_W=U\cap W$,
integers
$m_1,\ldots,m_n\ge 1$
invertible on $W$
and
an \'etale surjective morphism
$W'\to W$,
such that
the base change
of $f\colon V\to U$
by the \'etale map
$U\times_X(W
[S_1,\ldots,S_n]/
(S_1^{m_1}-t_1,
\ldots,S_n^{m_n}-t_n)
\times_WW')
\to U$
is a constant
\'etale covering.

{\rm 2.}
We say 
$f\colon V\to U$
is of Kummer type
with respect to $X$
if it is
of Kummer type
at every point
of $x\in X\setminus U$.
\end{df}

\begin{lm}\label{lmtX}
Let $S$ be a separated
noetherian scheme
and $f\colon V\to U$ be a
finite \'etale morphism
of separated
schemes of finite type over $S$.
Let $u\in U$ be a point 
let $A$ be an $S$-integral valuation ring of $\kappa(u)$.
We consider the conditions:

{\rm (1)}
For every point
$v\in f^{-1}(u)$
and for every 
valuation ring $B$ of
$\kappa(v)$
dominating $A$,
the extension $\kappa(v)$
over $\kappa(u)$
is tamely ramified
with respect to $B$.

{\rm (2)}
There exists a proper
scheme $X$ over $S$
containing $U$
as an open subscheme
such that
$V\to U$ is of Kummer type
at  the image of
the closed point 
of ${\rm Spec}\ A$
by 
${\rm Spec}\ A\to X$.

The condition
{\rm (1)} implies {\rm (2)}.
If $A$ is $U$-external,
then the conditions 
{\rm (1)} and {\rm (2)}
are equivalent.
\end{lm}

{\it Proof.}
(2)$\Rightarrow$(1) 
Since the assertion is 
\'etale local on ${\rm Spec}\ A$,
it follows from
Corollary \ref{corvalt}
(2)$\Rightarrow$(1).

(1)$\Rightarrow$(2)
Assume that
$A$ is $U$-external.
By Corollary \ref{corvalt}
(1)$\Rightarrow$(2),
there exists non-zero
elements
$t_1,\ldots,t_n\in 
{\mathfrak m}_A$
and integers
$m_1,\ldots,m_n\ge 1$
invertible in $A$
such that,
for every $v\in f^{-1}(u)$,
the normalization
of $A$ in
$\kappa(v)[S_1,\ldots,S_n]
/(S_1^{m_1}-t_1,\ldots,
S_n^{m_n}-t_n)$
is finite \'etale over
the normalization
of $A$ in
$\kappa(u)[S_1,\ldots,S_n]
/(S_1^{m_1}-t_1,\ldots,
S_n^{m_n}-t_n)$.
We put
$A_1=
A[S_1,\ldots,S_n]
/(S_1^{m_1}-t_1,\ldots,
S_n^{m_n}-t_n)$
and let $A^{\rm sh}$
denote the strict henselization of $A$.
Then, 
the \'etale covering
$V_u=V\times_Uu$ of $u$
is trivialized by
the base change
$A\to A^{\rm sh}\otimes_AA_1$.
Consequently,
there exist
an \'etale $A$-algebra $A'$
and a maximal ideal 
${\mathfrak m}'$ of $A'$
above the maximal ideal 
${\mathfrak m}_A$
such that
the \'etale covering
$V_u=V\times_Uu$ of $u$
is trivialized by
the base change
$A\to A'\otimes_AA_1$.

Let $\widetilde X
=\varprojlim_{X\in
{\cal C}_{U/S}}X$
be the limit of
compactifications
and $\tilde x
\in \widetilde X$
be the point
corresponding
to $A\subset \kappa(u)$.
We write
$\widetilde A=
{\cal O}_{
\widetilde X,\tilde x}$
and let ${\mathfrak p}
={\mathfrak m}_u
\subset \widetilde A$
denote the kernel
of the surjection
$\widetilde A\to A$.
We take liftings
$\tilde t_1,\ldots,
\tilde t_n\in 
\widetilde A$
of
$t_1,\ldots,t_n\in A$
and put
$\widetilde A_1=
\widetilde A[S_1,\ldots,S_n]
/(S_1^{m_1}-
\tilde t_1,\ldots,
S_n^{m_n}-\tilde t_n)$.
We also take an \'etale 
$\widetilde A$-algebra
$\widetilde A'$
such that
$\widetilde A'
\otimes_{\widetilde A}
A=A'$
and put 
$\widetilde A'_1=
\widetilde A'\otimes_
{\widetilde A}
\widetilde A_1$.
Then, we have an isomorphism
$\widetilde A'_1
\otimes_
{\widetilde A}A
\to A'\otimes_AA_1$.
Let
$\tilde {\mathfrak m}'$
be the maximal ideal of 
$\widetilde A'$
above ${\mathfrak m}'$
and 
$\tilde {\mathfrak m}'_1$
be a maximal ideal 
$\widetilde A'_1$
above $\tilde {\mathfrak m}'$.
We will apply
the following lemma
to the localization
$\widetilde A'_{1,
\tilde {\mathfrak m}'_1}$.

\begin{lm}\label{lmAp}
Let $A$ be a local ring
and ${\mathfrak p}$
be a prime ideal of $A$
such that
$A$ is canonically
isomorphic
to the inverse image
of $A/{\mathfrak p}$
by the surjection
$A_{\mathfrak p}\to
\kappa({\mathfrak p})$.
Let $B$ be a finite
$A$-algebra
such that
$B\otimes_AA_{\mathfrak p}$
is flat over
$A_{\mathfrak p}$
and that
$B\otimes_A\kappa({\mathfrak p})$
is isomorphic to
the product 
$\kappa({\mathfrak p})^n$.
We define an $A$-subalgebra
$B'$ of $B\otimes_A
A_{\mathfrak p}$
to be the inverse image
of $(A/{\mathfrak p})^n$
by the surjection
$B\otimes_A A_{\mathfrak p}\to
B\otimes_A \kappa({\mathfrak p})$.
Then, $B'$ is finite \'etale
over $A$
and the canonical maps
$B\otimes_AA_{\mathfrak p}\to
B'\otimes_AA_{\mathfrak p}$
and 
$B'/{\mathfrak p}B'
\to (A/{\mathfrak p})^n$
are isomorphisms.
\end{lm}

{\it Proof.}
We may assume
$B$ is ${\mathfrak p}$-torsion
free
and identify
$B\subset
B\otimes_AA_{\mathfrak p}$.
By the assumption on $A$,
the prime ideal
${\mathfrak p}
={\rm Ker}(A\to
A/{\mathfrak p})$
is equal to
the maximal ideal
${\mathfrak p}
A_{\mathfrak p}
={\rm Ker}(
A_{\mathfrak p}\to
\kappa({\mathfrak p}))$.
Hence, we
have 
${\mathfrak p}B'\subset
{\mathfrak p}
A_{\mathfrak p}
(B\otimes_A
A_{\mathfrak p})
=
{\mathfrak p}
A_{\mathfrak p}\cdot B
={\mathfrak p}B\subset
{\mathfrak p}B'$.
Thus, 
we have an equality
${\mathfrak p}B'
={\mathfrak p}
A_{\mathfrak p}
(B\otimes_A
A_{\mathfrak p})$
and an isomorphism
$B'/{\mathfrak p}B'
\to (A/{\mathfrak p})^n$.
Since the $A$-module
$B'/B$
is isomorphic to
$(A/{\mathfrak p})^n/
(B/{\mathfrak p}B)$,
the $A$-module $B'$
is of finite type.
Hence, by Nakayama's lemma,
$B'$ is finite flat over $A$
and hence is \'etale over $A$.
The canonical map
$B\otimes_A\kappa({\mathfrak p})\to
B'\otimes_A\kappa({\mathfrak p})$
is an isomorphism and hence
$B\otimes_AA_{\mathfrak p}\to
B'\otimes_AA_{\mathfrak p}$
is also an isomorphism.
\qed

We go back to the proof of
Lemma \ref{lmtX}.
By \cite[Proposition 5.11]{FK},
the local ring
$\widetilde A$
is canonically isomorphic
to the inverse image
of $\widetilde A/{\mathfrak p}$
by the surjection
$\widetilde A_{\mathfrak p}\to
\kappa({\mathfrak p})$.
We show that
the local ring
$\widetilde A'_{1,
\tilde {\mathfrak m}'_1}$
also satisfies
the condition of Lemma \ref{lmAp}.
Since $\widetilde A'_1$
is flat 
over $\widetilde A$,
it follows that
$\widetilde A'_1$
is canonically isomorphic
to the inverse image
of $\widetilde A'_1/
{\mathfrak p}\widetilde A'_1
=A'\otimes_AA_1$
by the surjection
$\widetilde A'_1
\otimes_{\widetilde A}
\widetilde A_{\mathfrak p}
\to
\widetilde A'_1
\otimes_{\widetilde A}
\kappa({\mathfrak p})$.
Thus the claim follows
by localization.

By Zariski's main theorem,
there exists
a finite
$\widetilde A$-algebra
$\widetilde B$
such that
$\widetilde B
\otimes_{\widetilde A}
{\cal O}_{U,u}$
is isomorphic
to 
$\Gamma(V\times_U
{\rm Spec}\ {\cal O}_{U,u},
{\cal O}_V)$.
Hence, by Lemma \ref{lmAp},
the base change
of $V\times_U
{\rm Spec}\ {\cal O}_{U,u}
\to
{\rm Spec}\ {\cal O}_{U,u}$
by
$\widetilde A\to
\widetilde A'_1$
is extended to
a finite \'etale covering
on a neighborhood
of $\tilde {\mathfrak m}'_1$.
By replacing
$A'$ by an \'etale algebra
contained in $A^{\rm sh}$,
we may assume that
the base change
of $V\times_U
{\rm Spec}\ {\cal O}_{U,u}
\to
{\rm Spec}\ {\cal O}_{U,u}$
by
$\widetilde A\to
\widetilde A'_1$
is a constant
finite \'etale covering.

Since
$\widetilde A=
\varinjlim_X{\cal O}_{X,x_X}$,
there exist an object
$X$ of ${\cal C}_{U/S}$,
an open neighborhood 
$W$ of $x_X$
and non-zero divisors
$f_1,\ldots,f_n
\in \Gamma(W,{\cal O}_X)$
invertible on
$U\cap W$
sent to $\tilde t_i$ for $i=1,\ldots,n$
and an \'etale morphism
$W'\to W$ such that
$\widetilde A'$ is the pull-back
of
$W'\to W$ by
${\rm Spec}\ \widetilde A\to W$.
We put
$W_1=W[S_1,\ldots,S_n]
/(S_1^{m_1}-
f_1,\ldots,
S_n^{m_n}-f_n)$.
Then further by $\widetilde A=
\varinjlim_X{\cal O}_{X,x_X}$,
replacing $X$ if
necessary,
the base change of
$V\to U$
by
$W_1\times_WW'\to X$
is a constant
finite \'etale covering.
\qed

\begin{pr}\label{prtX}
Let $S$ be a separated
noetherian scheme
and $f\colon V\to U$ be a
finite \'etale morphism
of separated
schemes of finite type over $S$.
Then, 
the following conditions
are equivalent:

{\rm (1)}
$f\colon V\to U$
is tamely ramified with
respect to $S$.

{\rm (2)}
There exists a proper
scheme $X$ over $S$
containing $U$
as an open subscheme
such that
$V\to U$ is of Kummer type
with respect to $X$.
\end{pr}

{\it Proof.}
(2)$\Rightarrow$(1) 
It follows from Lemma \ref{lmtX}
(2)$\Rightarrow$(1)
and Proposition \ref{prtY}
(2)$\Rightarrow$(1).

(1)$\Rightarrow$(2) 
By Proposition \ref{prtY}
(1)$\Rightarrow$(2)
and Lemma \ref{lmtX}
(1)$\Rightarrow$(2),
for every point $\tilde x
\in \widetilde X
=\varprojlim_{X\in
{\cal C}_{U/S}}X$
of the boundary,
there exists a proper
scheme $X$ over $S$
containing $U$
such that
the maximum
open subscheme
$W_X\subset X$
where $V\to U$
is of Kummer type
contains the image $\pi_X(\tilde x)\in X$
of $\tilde x\in \widetilde X$
by the projection $\pi_X
\colon \widetilde X\to X$.
In other words,
the family
$(\pi_X^{-1}(W_X)
)_{X\in {\cal C}_{U/S}}$
is an open covering of
$\widetilde X$.
Since $\widetilde X$
is quasi-compact
and since
the category ${\cal C}_{U/S}$
is cofiltered,
there exists
an object $X$ of
${\cal C}_{U/S}$
such that
$\pi_X^{-1}(W_X)=
\widetilde X$.
Thus the assertion follows.
\qed

\subsection{Tame ramification
of Galois coverings}\label{sstG}

The following proposition
shows that the definition 
of tame ramification
here
is equivalent to 
that given by Gabber
\cite[Section 2.1]{Vidal}
for Galois coverings.
Let $X$ be a normal scheme,
$U\subset X$,
be a dense open
subscheme and $V\to U$ 
be a finite $G$-torsor.
Then, for a geometric point
$\bar x$ of $X$,
the inertia subgroup $I_{\bar x}$
of $G$ is defined up to conjugacy.

\begin{pr}\label{prtmin}
Let $S$ be a separated
noetherian scheme
and $U$ be a separated 
normal integral scheme of
finite type over $S$.
For 
a $G$-torsor
$f\colon V\to U$
for a finite group $G$,
the following conditions
are equivalent:

{\rm (1)}
$f\colon V\to U$
is tamely ramified
with respect to $S$.

{\rm (2)}
There exists a proper
normal scheme $X$ over $S$
and an open immersion
$U\to X$ over $S$
such that
for every geometric
point $\bar x$ of $X$,
the order of the inertia
subgroup $I_{\bar x}\subset G$
is invertible at $\bar x$.
\end{pr}

{\it Proof.}
(1)$\Rightarrow$(2)
By Proposition \ref{prtX}
(1)$\Rightarrow$(2),
there exists a proper
normal scheme $X$ over $S$
and an open immersion
$U\to X$ over $S$
such that
$V\to U$
is of Kummer type
with respect to $X$.
Hence the assertion follows.

(2)$\Rightarrow$(1)
For every point
$u\in U$,
for every $S$-integral
and $U$-external
valuation ring $A$
of $\kappa(u)$,
for every $v\in f^{-1}(u)$
and for every 
valuation ring $B$ of
$\kappa(v)$
dominating $A$,
the order of
the inertia group
$I_{B/A}$ is invertible
in $A$.
Hence, by 
Proposition \ref{prtY}
(2)$\Rightarrow$(1),
the map $V\to U$
is tamely ramified
with respect to $S$.
\qed

\begin{cor}\label{corunrX}
Let $f\colon V\to U$ be a finite \'etale 
morphism and let $X$ be a normal scheme
containing $U$ as a dense open subscheme.
We assume $f\colon V\to U$
is a $G$-torsor
for a finite group $G$.
We consider the following conditions:

{\rm (1)}
$f\colon V\to U$
is tamely ramified with respect to $X$.

{\rm (2)}
For every geometric
point $\bar x$ of $X$,
the order of the inertia
subgroup $I_{\bar x}\subset G$
is invertible at $\bar x$.

{\rm (3)}
Let $x$ be an arbitrary point of $X$
such that the local ring
${\cal O}_{X,x}$ is a discrete valuation ring.
Then, $f\colon V\to U$
is tamely ramified
with respect to
${\cal O}_{X,x}$.

Then, we have implications 
{\rm (2)$\Rightarrow$(1)$\Rightarrow$(3)}.
If $X$ is a regular separated
noetherian scheme
and if $U$ is the complement
of a divisor with normal crossings,
then {\rm (3)} implies {\rm (2)}.
\end{cor}

{\it Proof.}
(2)$\Rightarrow$(1)
It follows from
Proposition \ref{prtmin}
(2)$\Rightarrow$(1).

(1)$\Rightarrow$(3)
It follows from
Proposition \ref{prtY}
(1)$\Rightarrow$(2).

(3)$\Rightarrow$(2)
By \cite[Proposition 5.2]{SGA1}
 (Lemme d'Abhyankhar absolu),
the condition (3) implies
that $V\to U$ is
of Kummer type
with respect to $X$.
Hence (3) implies (2).
\qed

If we drop the assumption
that $U$ is the complement
of a divisor with normal crossings, 
the implication
(1)$\Rightarrow$(2)
nor
(3)$\Rightarrow$(1)
need not hold
even if $X$ is regular,
as the following examples show.
The authors thank
M.\ Raynaud 
for the help to find
Example \ref{extm}.2.

\begin{ex}\label{extm}
{\rm 1. 
Let $k$ be an algebraically closed
field of characteristic $p>0$
and $V_0\to U_0=
{\rm Spec}\ k[t,t^{-1},(t-1)^{-1}]=
{\mathbf P}_k^1\setminus 
\{0,1,\infty\}$
be a finite \'etale connected
Galois covering
of degree divisible by $p$
tamely ramified
at $0,1,\infty$.
We put $X={\mathbf A}_k^2=
{\rm Spec}\ k[x,y]
\supset
U={\rm Spec}\ k[x,y,(xy(x-y))^{-1}]$
and define a map
$U\to U_0$
by sending $t$ to $x/y$.
Let $V=V_0\times_{U_0}U$
be the pull-back by the map
$U\to U_0$.
Then, since $V_0\to U_0$
is assumed tamely ramified,
the covering $V\to U$
is tamely ramified with respect to $X$ and
satisfies the condition (1).
Since the inertia group
at the origin $0\in X$
is equal to the Galois group
${\rm Gal}(U_0/V_0)$,
the condition (2) is not satisfied.

2.
Let $k$ be a
field of characteristic $p\ge 3$
and $m\ge 1$ be an integer.
We consider the cyclic covering
$$\begin{CD}
Z&={\rm Spec}\ k[x,y,z]/
(z^{p-1}&-(x^{2m(p-1)-1}+
x^{m(p-1)}y^{m(p-1)}+
y^{2m(p-1)}))\\
@V{\pi}VV&\\
X&={\mathbf A}^2_k
={\rm Spec}\ k[x,y]&
\end{CD}$$
of degree $p-1$
ramified at the divisor
$D=(x^{2m(p-1)-1}+
x^{m(p-1)}y^{m(p-1)}+
y^{2m(p-1)})$.
We put $U=X\setminus D$
and $W=Z\times_XU$.
We consider the 
Artin-Schreier covering
of $Z\times_X
({\mathbf G}_m\times 
{\mathbf A}^1)$
defined by
$T^p-T=y^{m(p-2)}z/x^{mp}$.
Since 
\begin{eqnarray*}
&&
\frac{y^{m(p-2)}z}{x^{mp}}-
\left(
\Bigl(\frac z{x^my^m}\Bigr)^p-
\frac z{x^my^m}\right)\\
&=&
\frac{(y^{2m(p-1)}-
(x^{2m(p-1)-1}+
x^{m(p-1)}y^{m(p-1)}+
y^{2m(p-1)})+
x^{m(p-1)}y^{m(p-1)})z}
{x^{mp}y^{mp}}
\\&=&
-\frac{x^{m(p-2)-1}z}{y^{mp}},
\end{eqnarray*}
it is extended to a finite
\'etale covering
on $Z\setminus \pi^{-1}(0)$.
Hence it defines
a finite \'etale
Galois covering
$V\to U$ of Galois 
group ${\mathbb F}_p^\times
\ltimes {\mathbb F}_p$,
tamely ramified
with respect
to $X\setminus \{0\}$.
Thus the Galois covering
$V\to U\subset X$
satisfies the condition
(3) in Corollary \ref{corunrX}.

Let $X'\to X$ be the blow-up
at the origin
and $Z'$ be the normalization
of $X'$ in $W$.
We put $t=x/y$.
Then, since
$x^{2m(p-1)-1}+
x^{m(p-1)}y^{m(p-1)}+
y^{2m(p-1)}
=y^{2m(p-1)-1}(
t^{2m(p-1)-1}+
t^{m(p-1)}y+y)$,
the cyclic covering
$Z'\to X'$
is totally ramified
along the exceptional divisor
$E\subset X'$
and the valuation of
$y^{m(p-2)}z/x^{mp}=
z/(x^{2m}t^{m(p-2)})$
at the generic point
of the inverse image
$E'=E\times_{X'}Z'$
is $
(2m(p-1)-1)-
2m(p-1)=-1$.
Hence
the Artin-Schreier covering
$V\to W$ is totally ramified
along $E'$
and 
$V\to U$ is not
tamely ramified with respect to $X$.
Thus the Galois covering
$V\to U\subset X$
does not 
satisfy the condition
(1) in Corollary \ref{corunrX}.}
\end{ex}

\newpage 
\section{Complements
on localized intersection products}\label{slp}

We compute certain tor-sheaves
in Section \ref{sstor}.
This computation
plays a crucial role
in the proof of
the excision formula
and of the
blow-up formula 
in Section \ref{ssexc}.
We recall the definition
of the localized intersection 
product and some useful formulas
in Section \ref{sslip}.

The results
in Sections \ref{sslcc}
and \ref{ssrex}
are used only
in an explicit computation
of the logarithmic different in 
Section \ref{ssdiff}.
In Section \ref{ssrex},
we prove a refinement Proposition \ref{prrex}
of the excess intersection formula
\cite[Proposition 3.4.2]{KSI},
which relates
the classes of certain tor-sheaves
defining the localized
intersection product
with the mapping
cone of exterior derived power
complexes.
In Section \ref{sslcc},
we show
that the class of
the mapping
cone of exterior derived power
complexes is given by the
localized Chern class.

\subsection{A computation of Tor sheaves}
\label{sstor}

We compute certain Tor sheaves
related to blow-up.
We recall some terminology
on tor-sheaves.
For a cartesian diagram
\begin{equation}
\begin{CD}
Y@<<< Z\\
@VgVV @VVV\\
S@<f<< X
\end{CD}
\label{eqtch}
\end{equation} of
schemes,
a quasi-coherent
${\cal O}_X$-module ${\cal F}$,
a quasi-coherent
${\cal O}_Y$-module ${\cal G}$
and an integer $q\ge 0$,
a quasi-coherent
${\cal O}_Z$-module 
${\cal T}or^{{\cal O}_S}_q(
{\cal F},{\cal G})$
is defined
in \cite[(6.5.3)]{EGA3}.
If $S={\rm Spec}\ A,
X={\rm Spec}\ B,
Y={\rm Spec}\ C$
are affine
and if
${\cal F}$ and ${\cal G}$
are the quasi-coherent sheaves
associated to an
$B$-module $M$
and an $C$-module $N$
respectively,
then
${\cal T}or^{{\cal O}_S}_q(
{\cal F},{\cal G})$
is associated to the
$B\otimes_AC$-module 
$Tor^A_q(M,N)$.
If ${\cal F}={\cal O}_X$,
we put
$L_qf^*{\cal G}=
{\cal T}or^{{\cal O}_S}_q(
{\cal O}_X,{\cal G})$.

\begin{df}\label{dftor}
{\rm 1.
(\cite[Definition 1.5]{sga6})}
Let $X$ and $Y$ be schemes over
a scheme $S$.
We say that
$X$ and $Y$ are {\rm tor-independent
over} $S$
if ${\cal T}or^{{\cal O}_S}_q(
{\cal O}_X,{\cal O}_Y)=0$
for every $q>0$.

{\rm 2. (\cite[Definition 3.1]{sga6})}
Let $f\colon X\to S$ be
a morphism of schemes.
We say that
$f$ is {\rm of finite tor-dimension},
if there exists an integer $n\ge 0$
such that, for every quasi-coherent
${\cal O}_S$-module ${\cal F}$
and every integer $q>n$,
we have
$L_qf^*{\cal F}=0$.
\end{df}
If $X$ or $Y$ is flat over $S$,
then $X$ and $Y$ are tor-independent
over $S$.

\begin{lm}
\label{lmtidp}
We consider morphisms
$$\begin{CD}
Y@>g>> S@<f<<X@<{f'}<< X'\end{CD}$$
of schemes.
Assume that
$X$ and $Y$ are tor-independent
over $S$.
Then,
$X'$ and $Y$ are tor-independent
over $S$
if and only if
$X'$ and $X\times_SY$ are tor-independent
over $X$.
\end{lm}

{\it Proof.}
By the assumption that
${\cal T}or^{{\cal O}_S}_q(
{\cal O}_X,{\cal O}_Y)=0$
for every $q>0$,
we obtain an isomorphism
${\cal T}or^{{\cal O}_S}_q(
{\cal O}_{X'},{\cal O}_Y)
\to 
{\cal T}or^{{\cal O}_X}_q(
{\cal O}_{X'},{\cal O}_{X\times_SY})$
for every $q>0$.
\qed

\begin{lm}\label{lmtch}
Assume that the
schemes $S$ and $X$
are noetherian
and $Y$ is of
finite type over $S$ 
in the diagram
{\rm (\ref{eqtch})}.
Then,
for a coherent
${\cal O}_X$-module ${\cal F}$
and
for a coherent
${\cal O}_Y$-module ${\cal G}$,
the
${\cal O}_Z$-modules
${\cal T}or_q^{{\cal O}_S}
({\cal F},{\cal G})$
are coherent.
\end{lm}

{\it Proof.}
Since the question is local
on $Z$, we may assume
that schemes $S,X,Y$ and $Z$
are affine.
We take a closed
immersion $Y\to P={\mathbf A}^n_S$
to an affine space.
Since an ${\mathcal O}_Z$-module
is coherent if it is
coherent as
an ${\mathcal O}_{X\times_SP}$-module,
we may replace
$Y$ by $P$.
Hence, we may assume 
further that $Y$ is flat over
$S$.
Then,
a resolution ${\cal L}$ of
${\cal G}$
by free ${\cal O}_Y$-modules
of finite rank
is a resolution by flat
${\cal O}_S$-modules.
Since 
${\cal T}or_q^{{\cal O}_S}
({\cal F},{\cal G})$
is a cohomology sheaf
of the complex
${\mathcal F}
\otimes_{{\mathcal O}_S}
{\mathcal L}$,
it is a coherent
${\mathcal O}_Z$-module.
\qed

Let $S$ be
a regular noetherian
scheme of finite
dimension.
Then, for 
a scheme $f\colon X\to S$
of finite type over $S$,
the dimension function
$X\to {\mathbb N}$
is defined 
as in \cite[Section 2.1]{KSI}.
Namely,
for a point $x\in X$ and $s=f(x)$,
we put
$\dim x=
{\rm tr. deg}
(\kappa(x)/\kappa(s))+
\dim S-\dim {\cal O}_{S,s}$.
Using this dimension function,
the topological filtration
$F_\bullet G(X)$ and
the lower numbering
Chow groups
$CH_\bullet (X)$
are defined.
We have a canonical map
$CH_\bullet (X)
\to 
{\rm Gr}^F_\bullet G(X)$
sending the class
$[V]$ 
of an integral closed subscheme
$V$ of $X$ to the class 
$[{\mathcal O}_V]$
also denoted by $[V]$.

By Lemma \ref{lmtch},
for a morphism $f\colon X\to Y$
of noetherian schemes
of finite tor-dimension
and for a scheme $A$
of finite type
over $Y$,
the pull-back map
$$f^*\colon
G(A)\to G(A\times_YX)$$
is defined by
$f^*([{\cal F}])=
\sum_{q\ge 0}
(-1)^q
[{\cal T}or_q^{{\cal O}_Y}
({\cal F},{\cal O}_X)]$
for a coherent
${\cal O}_A$-module
${\cal F}$.

\begin{lm}\label{lmtfil}
\setcounter{equation}0
Let $S$ be
a regular noetherian
scheme of finite
dimension
and $f\colon
X\to Y$
be a quasi-projective
morphism locally
of complete intersection
of relative
virtual dimension $r$
of schemes 
of finite type over $S$.

Then, for
a scheme $A$ of finite
type over $Y$,
the pull-back $f^*\colon
G(A)\to G(A\times_YX)$
preserves
topological
filtration
in the sense that
$f^*$ maps
$F_\bullet G(A)$
to $F_{\bullet+r} 
G(A\times_YX)$.
Further, for
an integer $q\ge0$,
we have
a commutative diagram
\begin{equation}
\begin{CD}
CH_q(A)@>{f^!}>>
CH_{q+r}(A\times_YX)
\\
@VVV @VVV\\
{\rm Gr}^F_q G(A)
@>{f^*}>>
{\rm Gr}^F_{q+r} 
G(A\times_YX).
\end{CD}
\end{equation}
\end{lm}

{\it Proof.}
By the assumption on $f$,
it is the composition
$X\to P\to Y$
of a regular immersion
$X\to P$ and
a smooth morphism
$P\to Y$.
Since, it is clear
for a smooth map,
it is reduced to
the case where
$f$ is a regular immersion.
Then, it follows from
\cite[Proposition 2.2.2]{KSI}.
\qed

We compute
${\cal T}or_r^{{\cal O}_X}
({\cal O}_{X'},{\cal O}_Y)$
for morphisms
$Y\to X'\to X$
under certain conditions.
Corollaries \ref{corbp}
and \ref{cortlp} of
the following proposition
are crucial in the proof of
Proposition \ref{prblup}
and Theorem \ref{thmexc}
respectively.

\begin{pr}\label{prspsq}
\setcounter{equation}0
Let $X$ be a scheme and
${\cal N}$ be a locally free 
${\cal O}_X$-module 
of finite rank $c$.
Let $\alpha\colon
{\cal N}\to {\cal O}_X$
be an ${\cal O}_X$-linear map
and
$C\subset X$ be the closed
subscheme defined by
${\cal I}_C={\rm Im}(\alpha\colon
{\cal N}\to {\cal O}_X)$.
Let $P\subset {\mathbf P}({\cal N})$
be an open subscheme
of the associated
${\mathbf P}^{c-1}$-bundle
${\mathbf P}({\cal N})=
{\cal P}roj(S^\bullet{\cal N})$
and $p\colon
P\to X$ be the canonical map.

Let ${\cal E}={\rm Ker}(
p^*{\cal N}\to {\cal O}_P(1))
=\Omega^1_{P/X}(1)$
be the kernel
of the canonical surjection.
Let $X'\subset P$
be the closed subscheme
defined by the image
${\cal I}_{X'}=
{\rm Im}(\alpha'\colon
{\cal E}\to {\cal O}_P)$
of the restriction
$\alpha'=
p^*\alpha|_{\cal E}\colon
{\cal E}\to {\cal O}_P$
to ${\cal E}\subset p^*{\cal N}$
and $q\colon X'\to X$
denote the composition.

We consider
a cartesian diagram
\begin{equation}
\begin{CD}
E_Y@>>> E@>>>
P_C@>>>C\\
@VVV @VVV 
@VVV@VVV\\
Y@>g>>X'
@>{\subset}>> P@>p>> X
\end{CD}
\label{eqXXtor}
\end{equation}
of schemes.
Let $f\colon Y\to X$
denote the composition
of the bottom arrows.
We assume that
$E_Y=E\times_{X'}Y$
is a Cartier divisor
of $Y$.

{\rm 1.}
The upper middle arrow
$E\to P_C=P
\times_XC$
is an isomorphism.
The restriction 
$X'\setminus E
\to
X\setminus C$
of
$q\colon X'\to X$
is an isomorphism.
The composition
of the immersions
$X'\setminus E
\to
P\setminus P_C
\to
{\mathbf P}({\cal N})$
is the composition
of the isomorphism
$X'\setminus E
\to
X\setminus C$
and the section
$X\setminus C
\to {\mathbf P}
({\cal N})$
defined by the surjection
$\alpha|_{
X\setminus C}
\colon
{\cal N}|_{
X\setminus C}
\to
{\cal O}_{
X\setminus C}$.

The map $\alpha\colon
{\mathcal N}\to {\cal O}_X$
induces a surjection
$f^*{\mathcal N}\to {\cal I}_{E_Y}
\subset {\cal O}_Y$.
The composition
$Y\to X'\to {\mathbf P}({\cal N})$
is the section defined by
the surjection
$f^*{\mathcal N}\to {\cal I}_{E_Y}$.

{\rm 2.}
We assume that
the immersion
$X'\to P$ is a regular
immersion of codimension $c-1$
{\rm \cite[Definition 1.4]{SGA6B}}.
Let $\gamma\colon Y
\to P_Y=P\times_XY$
be the section
defined by $g$ and 
$\Gamma\subset P_Y$
be the image of $\gamma$
regarded as a closed subscheme of $P_Y$.
Let ${\rm pr}_1\colon P_Y\to P$
and ${\rm pr}_2\colon P_Y\to Y$
denote the projections.

Then, 
the composition
$q\colon
X'\to X$
is of finite tor-dimension
and
there exists a spectral sequence
$E^1_{p,q}
\Rightarrow
E_r$ of
${\cal O}_{X'\times_XY}$-modules
such that
\begin{equation}
E_r=
\begin{cases}
{\rm Ker}({\cal O}_{X'\times_XY}
\to {\cal O}_\Gamma)
&\qquad \text{ if }r=0,\\
{\cal T}or_r^{{\cal O}_X}
({\cal O}_{X'},{\cal O}_Y)
&\qquad \text{ if }r\neq 0,
\end{cases}
\label{eqspsq}
\end{equation}
$$
E^1_{p,q}=
{\rm pr}_1^*
\Omega_{P/X}^{p+q+1}(p+q+1)
\otimes
{\rm pr}_2^*{\cal N}_{E_Y/Y}^{-(q+1)}
\qquad \text{ if } p\ge 0,q\ge 0,
p+q\le c-2
$$
and 
$E^1_{p,q}=0$
otherwise.

{\rm 3.}
Assume that
$Y$ is noetherian.
The surjection
${\cal O}_{X'\times_XY}
\to {\cal O}_\Gamma$
is an isomorphism
outside $E\times_CE_Y$
and we have equalities
\begin{eqnarray}
&&
[{\rm Ker}({\cal O}_{X'\times_XY}
\to {\cal O}_\Gamma)]
+\sum_{r>0}(-1)^r
[{\cal T}or_r^{{\cal O}_X}
({\cal O}_{X'},{\cal O}_Y)]
\nonumber
\\
&=&\sum_{p=1}^{c-1}
(-1)^{p-1}
\sum_{q=1}^p
[{\rm pr}_1^*
\Omega^p_{E/C}(p)
\otimes 
{\rm pr}_2^*{\cal N}_{E_Y/Y}^{-q}]
\label{eqKos}
\\
&=&
\sum_{s=2}^c
(-1)^s
[\Lambda^s\pi^*{\cal N}]
\cdot 
\sum_{q\ge 1,r\ge 1,q+r\le s}
[{\rm pr}_1^*{\cal O}(-r)]
\cdot 
[{\rm pr}_2^*{\cal N}_{E_Y/Y}^{-q}]
\label{eqKos2}
\end{eqnarray}
in $G(E\times_CE_Y)$
where $\pi\colon 
E\times_CE_Y\to X$
denotes the canonical map.
\end{pr}

{\it Proof.}
1.
The first paragraph is
clear from the definition
of $X'$ and of $C$.
The map $\alpha\colon
{\cal N}\to {\cal O}_X$
defines a surjection
${\cal N}\to {\cal I}_C
\subset {\cal O}_X$
and hence induces a surjection
$f^*{\cal N}\to 
{\cal I}_{E_Y}
\subset {\cal O}_Y$.
By the first paragraph,
the kernel of
the surjection
$f^*{\cal N}\to 
g^*({\cal O}_{X'}(1))$
is equal to the
kernel of
$f^*{\cal N}\to 
{\cal I}_{E_Y}$
on the complement
$Y\setminus E_Y$.
Since a section
$Y\to 
{\mathbf P}({\cal N})$
is uniquely determined
by its restriction
to the complement
of a divisor,
the assertion follows.

2.
By the definition
of regular immersion
\cite[Definition 1.4]{SGA6B},
the Koszul complex
$${\cal K}={\rm Kos}(\alpha')
=[\Lambda^{c-1}{\cal E}
\to \cdots \to 
{\cal E}
\overset {\alpha'}
\to {\cal O}_P]$$
is a resolution of an 
${\cal O}_P$-module
${\cal O}_{X'}$.
Since $P$ is flat over $X$
and $X'\to P$ is 
assumed to be
a regular immersion, 
the composition 
$q\colon X'\to X$
is of finite tor-dimension.
Further,
we obtain
an isomorphism
$${\cal H}_r
({\cal K}
\otimes_
{{\cal O}_X}
{\cal O}_Y)
\to 
{\cal T}or_r^{{\cal O}_X}
({\cal O}_{X'},{\cal O}_Y)$$
from the resolution
${\cal K}\to {\cal O}_{X'}$.

We construct a resolution
of the ${\cal O}_{P_Y}$-module
${\cal O}_\Gamma$.
Let $\beta\colon
{\rm pr}_1^*{\cal E}
\to 
{\rm pr}_2^*
{\cal I}_{E_Y}$
be the restriction
of the map
${\rm pr}_1^*p^*{\cal N}
={\rm pr}_2^*f^*{\cal N}
\to {\rm pr}_2^*
{\cal I}_{E_Y}$
to the kernel
${\rm pr}_1^*{\cal E}
={\rm Ker}
({\rm pr}_1^*p^*{\cal N}
\to
{\rm pr}_1^*{\cal O}_P(1))$.
The map $\beta$
induces the pull-back
${\rm pr}_1^*{\cal E}
\to {\cal O}_{P_Y}$
of $\alpha'$.
Since the section
$\gamma\colon
Y\to P_Y$
is defined by the surjection
$f^*{\cal N}\to 
{\cal I}_{E_Y}$ by 1.,
the closed subscheme
$\Gamma\subset
P_Y$ is characterized
by the condition that
the cokernel 
${\rm Coker}(\beta\colon
{\rm pr}_1^*{\cal E}
\to {\rm pr}_2^*
{\cal I}_{E_Y})$
is an invertible
${\cal O}_\Gamma$-module.
Hence
the Koszul complex
${\cal K}'={\rm Kos}(\beta')$
defined by the twist
$\beta'\colon
{\rm pr}_1^*{\cal E}
\otimes
{\rm pr}_2^*{\cal I}_{E_Y}^{-1}
\to {\cal O}_{P_Y}$
of $\beta$
is a resolution of
the ${\cal O}_{P_Y}$-module
${\cal O}_\Gamma$.

We consider 
the morphism
of complexes
${\cal K}\otimes_{{\cal O}_X}
{\cal O}_Y
\to 
{\cal K}'$
induced by the inclusion
${\rm pr}_1^*{\cal E}
\to
{\rm pr}_1^*{\cal E}
\otimes
{\rm pr}_2^*{\cal I}_{E_Y}^{-1}$.
Then, it
induces the canonical surjection
${\cal O}_{X'\times_XY}
\to
{\cal O}_\Gamma$.
Hence, for
the complex
${\cal M}=
({\cal K}'/({\cal K}\otimes_{{\cal O}_X}
{\cal O}_Y))[-1]$,
we have an isomorphism
${\cal H}_r{\cal M}
\to E_r$.

The $p$-th component
${\cal M}_p$
of the complex
${\cal M}$ is given by
$${\cal M}_p=
(\Lambda^{p+1}
{\rm pr}_1^*{\cal E}
\otimes
{\rm pr}_2^*{\cal I}_{E_Y}^{-(p+1)})/
\Lambda^{p+1}
{\rm pr}_1^*{\cal E}=
\Lambda^{p+1}
{\rm pr}_1^*{\cal E}
\otimes
({\rm pr}_2^*{\cal I}_{E_Y}^{-(p+1)}/
{\cal O}_{P_Y}).$$
We define an increasing filtration
$F_\bullet$ on ${\cal M}$
by 
$F_q{\cal M}_p=
\Lambda^{p+1}
{\rm pr}_1^*{\cal E}
\otimes
({\rm pr}_2^*
{\cal I}_{E_Y}^{-(q+1)}/
{\cal O}_{P_Y})$.
Then, we obtain a spectral sequence
$$E^1_{p,q}=
{\rm Gr}^F_{-q}{\cal M}_{p+q}
\Rightarrow
{\cal H}_r{\cal M}.$$
Since
${\rm Gr}^F_{-q}{\cal M}_p=
\Lambda^{p+1}
{\rm pr}_1^*{\cal E}
\otimes
{\rm pr}_2^*{\cal N}_{E_Y/Y}^{-(q+1)}$
and ${\cal E}=
\Omega^1_{P/X}(1)$,
the assertion follows.

3.
By 2.,
we have the equality
(\ref{eqKos}).
By the exact sequence
$0\to {\cal E}
\to p^*{\cal N}
\to {\cal O}(1)\to 0$,
we have
an exact sequence
$0
\to
\Lambda^cp^*{\cal N}(-(c-p))
\to
\cdots
\to
\Lambda^{p+1}p^*{\cal N}(-1)
\to 
\Lambda^p{\cal E}\to 0$
and an equality
$[\Lambda^p{\cal E}]
=\sum_{p+r=p+1}^c
(-1)^{r-1}
[\Lambda^{p+r}p^*{\cal N}(-r)]$.
Substituting this
and putting $p+r=s$,
we obtain the second equality
(\ref{eqKos2}).
\qed

\begin{cor}\label{corbp}
\setcounter{equation}0
Let $X$ be a noetherian scheme
and $C\subset X$ be
a closed subscheme
such that the immersion
$C\to X$
is a regular immersion
of codimension $c$.
Let $q\colon X'\to X$
be the blow-up at $C$
and let $E=X'\times_XC$
denote the exceptional divisor.

Let $Y$ be a noetherian scheme
over $X'$
such that
$E_Y=E\times_{X'}Y\subset Y$
is a Cartier divisor
and let $\Gamma$
denote the image of 
the section
$Y\to X'\times_XY$.
Let
${\rm pr}_1\colon
E\times_CE_Y\to E$
and
${\rm pr}_2\colon
E\times_CE_Y\to E_Y$
denote the projections.

Assume that
$q\colon X'\to X$
is locally of complete intersection.
Then,
we have an equality
\begin{eqnarray}
&&
[{\rm Ker}({\cal O}_{X'\times_XY}
\to {\cal O}_\Gamma)]
+\sum_{r>0}(-1)^r
[{\cal T}or_r^{{\cal O}_X}
({\cal O}_{X'},{\cal O}_Y)]
\nonumber
\\
&=&\sum_{p=1}^{c-1}
(-1)^{p-1}
\sum_{q=1}^p
[{\rm pr}_1^*
\Omega^p_{E/C}(p)
\otimes 
{\rm pr}_2^*{\cal N}_{E_Y/Y}^{-q}]
\label{eqbpY}
\end{eqnarray}
in $G(E\times_CE_Y)$.
\end{cor}

{\it Proof.}
Locally on $X$,
there exists a surjection
${\cal N}={\cal O}_X^{\oplus c}
\to {\cal I}_C$ of
${\cal O}_X$-modules.
Hence, applying
Proposition \ref{prspsq}.3
to $P={\mathbf P}({\cal N})$
and the closed
immersion $Y=X'\to P$,
we obtain
a spectral sequence (\ref{eqspsq}),
locally on $X$.
By the proof of 
Proposition \ref{prspsq}.3.,
it suffices to construct the
spectral sequence (\ref{eqspsq})
globally.

If we have another
locally free
${\cal O}_X$-module ${\cal N}'$
of rank $c$
and a surjection
${\cal N}'\to {\cal I}_C$,
then locally
we have an isomorphism
${\cal N}\to {\cal N}'$
compatible with
the surjections to 
${\cal I}_C$.
It induces an isomorphism
of spectral sequences
and the assertion follows.
\qed

The authors do not know
how to construct globally the
spectral sequence (\ref{eqspsq})
under the assumption
of Corollary \ref{corbp}
without using patching.

\begin{cor}\label{cortlp}
\setcounter{equation}0
Let $U$ be a scheme
of finite type 
over a noetherian scheme $S$
and $D\subset U$ be a Cartier divisor.
Let $q\colon
(U\times_SU)^\sim
\to 
U\times_SU$
be the log product
with respect to $D$
and assume that
$q\colon
(U\times_SU)^\sim
\to 
U\times_SU$
is locally of complete intersection 
of relative dimension $0$.
Let $\Delta_U\subset U\times_SU$
and $\Delta^{\log}_U\subset 
(U\times_SU)^\sim$
denote the diagonal
and the log diagonal
respectively and
identify the inverse image
of $\Delta_D\subset U\times_SU$
by $(U\times_SU)^\sim
\to U\times_SU$
with ${\mathbf G}_{m,D}$

Then,
the kernel of the surjection
$q^*{\mathcal O}_{\Delta_U}
\to {\mathcal O}_{\Delta^{\log}_U}$
and 
$L_rq^*{\mathcal O}_{\Delta_U}$
for $r>0$
are coherent
${\mathcal O}_{{\mathbf G}_{m,D}}$-modules
and 
we have
$$
[{\rm Ker}(q^*{\mathcal O}_{\Delta_U}
\to {\mathcal O}_{\Delta^{\log}_U})]
+\sum_{r>0}
(-1)^r[L_rq^*{\mathcal O}_{\Delta_U}]
=
[{\mathbf G}_{m,D}]$$
in $G({\mathbf G}_{m,D})$.
\end{cor}

{\it Proof.}
We define 
a locally free
${\mathcal O}_{U\times_SU}$-module 
of rank 2 by
${\mathcal N}=
{\rm pr}_1^*{\cal I}_D
\oplus
{\rm pr}_2^*{\cal I}_D$ and
define a ${\mathbf P}^1$-bundle
${\mathbf P}(
{\rm pr}_1^*{\cal I}_D
\oplus
{\rm pr}_2^*{\cal I}_D)$
over $U\times_SU$.
The complement
$P
\subset
{\mathbf P}(
{\rm pr}_1^*{\cal I}_D
\oplus
{\rm pr}_2^*{\cal I}_D)$
of the sections
defined by the surjections
${\mathcal N}=
{\rm pr}_1^*{\cal I}_D
\oplus
{\rm pr}_2^*{\cal I}_D
\to
{\rm pr}_i^*{\cal I}_D$
for $i=1,2$
is a ${\mathbf G}_m$-bundle
on $U\times_SU$.
We regard
the log product
$(U\times_SU)^\sim$
as a subscheme
of $P$.
By the assumption that
$(U\times_SU)^\sim
\to 
U\times_SU$
is locally of complete intersection 
of relative dimension $0$,
the immersion
$(U\times_SU)^\sim
\to P$ is a regular
immersion of codimension 1
by Lemma below.

We apply Proposition
\ref{prspsq} to the diagram
\begin{equation}
\begin{CD}
D
@>>>
P_{D\times_SD}
@=
P_{D\times_SD}
@>>>
D\times_SD
\\
@VVV@VVV @VVV@VVV\\
U
@>>>
(U\times_SU)^\sim
@>{\subset}>>P@>>> 
U\times_SU
\end{CD}
\end{equation}
where the image of
the first arrow
$U\to (U\times_SU)^\sim$
in the bottom
is $\Delta_U^{\log}$
and that of the composition
$U\to U\times_SU$
is $\Delta_U$.
In the notation there,
we have $C=D\times_SD,
Y=U,E_Y=D$
and $E$ is
the ${\mathbf G}_m$-bundle
$P_{D\times_SD}$.
Hence
$E\times_CE_Y$
is $P_D={\mathbf G}_{m,D}$.

In the right hand
side of (\ref{eqKos2}),
the pull back of
${\cal N}=
{\rm pr}_1^*{\cal I}_D
\oplus
{\rm pr}_2^*{\cal I}_D$
to $D$ is
${\cal N}_{D/U}^
{\oplus 2}$
and 
the second exterior power
$\Lambda^2\pi^*{\mathcal N}$
is ${\cal N}_{D/U}^{\otimes 2}$.
Since
${\cal N}_{E_Y/Y}={\cal N}_{D/U}$
and the pull-back of
${\mathcal O}(1)$ to $D$
is also ${\cal N}_{D/U}$
by Proposition
\ref{prspsq}.1,
the assertion follows.
\qed

\begin{lm}\label{lmciri}
Let $S$ be a noetherian scheme,
$X\to S$ be a scheme
locally of complete intersection
of relative dimension $d$
and
$P\to S$ be a smooth scheme
of relative dimension $n$.
Then, an immersion
$X\to P$ over $S$
is a regular immersion
of codimension $n-d$.
\end{lm}

{\it Proof.}
Since the assertion is
local on $X$, 
we may take a regular
immersion $X\to Q$
of codimension $c$
over $S$
to a smooth scheme
$Q$ of relative dimension
$d+c$ over $S$.
The immersion 
$X\to P\times_SQ$
is the composition
of the section
$X\to P\times_SX$
of a smooth morphism
of relative dimension $n$
and 
the smooth base change
$P\times_SX
\to P\times_SQ$
of the regular immersion
$X\to Q$
of codimension $c$
and is 
a regular immersion 
of codimension
$n+c$.
It is also
the composition
of the section
$X\to X\times_SQ$
of a smooth morphism
of relative dimension $d+c$
and an immersion
$X\times_SQ
\to P\times_SQ$.
Hence, the immersion
$X\times_SQ
\to P\times_SQ$
is a regular immersion
of codimension
$(n+c)-(d+c)=n-d$
on the image
of $X$
by \cite[Proposition 19.1.5]{EGA4}.
Since
$X\times_SQ
\to P\times_SQ$
is a smooth base change
of 
$X\to P$,
the assertion follows.
\qed

\subsection{Derived exterior
power and localized Chern classes}
\label{sslcc}

We study the relation between
derived exterior
power complexes
and localized Chern classes.
Let $X$ be a 
noetherian scheme,
${\cal E}$
and ${\cal E}'$
be locally free
${\cal O}_X$-modules
of the same finite rank $n$
and $e\colon 
{\cal E}\to {\cal E}'$
be a morphism.
For an integer $k\ge0$,
we consider the complex
$[\Lambda^k{\cal E}\to 
\Lambda^k{\cal E}']$
where
$\Lambda^k{\cal E}'$
is put on degree $0$.
Let $D$ be a closed
subset of $X$
such that ${\cal K}=
[{\cal E}\to {\cal E}']$
is acyclic
on the complement $X\setminus D$.
We define a morphism
\setcounter{equation}0
\begin{equation}
[\Lambda^k{\cal E}\to 
\Lambda^k{\cal E}']_D
\colon
G(X)\longrightarrow G(D)
\label{eqlcc}
\end{equation}
by sending the class
$[{\cal F}]$
of a coherent ${\cal O}_X$-module
${\cal F}$ to
$[{\rm Coker}(
{\cal F}\otimes
\Lambda^k{\cal E}
\to{\cal F}\otimes
\Lambda^k{\cal E}')]
-[{\rm Ker}(
{\cal F}\otimes
\Lambda^k{\cal E}
\to{\cal F}\otimes
\Lambda^k{\cal E}')]$.

Recall that
homomorphisms
$\lambda_t,\gamma_t\colon
K(X)\to
1+tK(X)[[t]]
\subset
K(X)[[t]]^\times$
are defined by
$\lambda_t([{\cal E}])
=
\sum_{q=0}^n[\Lambda^q{\cal E}]t^q$
and
$\gamma_t([{\cal E}])=
\lambda_{\frac t{1-t}}([{\cal E}])$.
We define
operators
$\gamma_k({\cal K})_D
\colon
G(X)\to G(D)$ 
for $k\ge 1$ by
requiring that
$\gamma_t({\cal K})_D=
\sum_{k=1}^n
\gamma_k({\cal K})_D\cdot t^k$
is given by
\begin{equation}
\gamma_t({\cal K})_D=
\left(\sum_{k=1}^n
[\Lambda^k{\cal E}\to 
\Lambda^k{\cal E}']_D\cdot
\Bigl(\frac t{1-t}\Bigr)^k
\right)
\cdot
\gamma_t([{\cal E}])^{-1}.
\label{eqdlt}
\end{equation}
If $\gamma_k({\cal K})$
without the suffix $D$
denotes the composition
with $G(D)\to G(X)$,
we have
$1+\gamma_t({\cal K})
=1+\bigl(
\gamma_t([{\cal E}'])
-\gamma_t([{\cal E}])
\bigr)
\cdot
\gamma_t([{\cal E}])^{-1}
=
\gamma_t([{\cal E}']-[{\cal E}])$.

Let
$c_k({\cal K})_D\colon
CH_{\bullet}(X)
\to
CH_{\bullet-k}(D)$
be the localized
Chern class map
defined 
by using the graph
construction
in \cite[Section 18.1]{fulton}.

\begin{pr}\label{prlcc}
\setcounter{equation}0
Let $S$ be
a regular noetherian
scheme of finite dimension.
Let $X$ be a scheme
of finite type over $S$,
${\cal E}$
and ${\cal E}'$
locally free
${\cal O}_X$-modules
of the same finite rank
and $e\colon 
{\cal E}\to {\cal E}'$
be a morphism
such that the restriction
on the complement of $D$
is an isomorphism.
For the complex
${\cal K}=[e\colon
{\cal E}\to {\cal E}']$
and an integer $k>0$,
we have the following.

{\rm 1.}
The map
$\gamma_k({\cal K})_D\colon
G(X)\to G(D)$
sends the topological
filtration $F_jG(X)$ to 
$F_{j-k}G(D)$.

{\rm 2.}
The diagram
\begin{equation}
\begin{CD}
CH_{j}(X)@>{c_k({\cal K})_D}>>
CH_{j-k}(D)\\
@VVV @VVV\\
{\rm Gr}^F_{j}G(X)
@>{\gamma_k({\cal K})_D}>>
{\rm Gr}^F_{j-k}G(D)
\end{CD}
\end{equation}
is commutative.
\end{pr}

{\it Proof.}
1. 
It suffices to show
$\gamma_k({\cal K})_D([X])
\in F_{d-k}G(D)$
assuming that
$X$ is integral
of dimension $d$
and $D\neq X$,
by a standard argument.
We show this by using
the most elementary case
of MacPherson's graph
construction 
cf. \cite[Section 18.1]{fulton}.

Let $n$ be the rank of
${\cal E}$
and $p\colon
G\to X$ be the Grassmann
scheme ${\rm Grass}_n(
{\cal E}\oplus {\cal E}')$
classifying
subbundles of rank $n$.
The second factor
${\cal E}'
\subset 
{\cal E}\oplus {\cal E}'$
defines
a section $s_0\colon X\to G$.
Let $t$ denote
the coordinate of
${\mathbf G}_{m,X}$.
Then, the graph
of $t^{-1}\cdot e\colon
{\cal E}_{{\mathbf G}_{m,X}}
\to 
{\cal E}'_{{\mathbf G}_{m,X}}$
defines a section
$\tilde s\colon
{\mathbf G}_{m,X}
\to 
{\mathbf G}_{m,G}$.
At $t=1$,
the restriction
$\tilde s|_1\colon
X\to G$
is the section
defined by the graph of $e$.

On the complement
$U=X\setminus D$,
the restriction
$e|_U\colon 
{\cal E}_U\to 
{\cal E}'_U$ 
is an isomorphism.
The transpose
of the graph of
$t\cdot e|_U^{-1}\colon
{\cal E}'_{{\mathbf A}^1_U}
\to 
{\cal E}_{{\mathbf A}^1_U}$
defines
a section
$\bar s\colon
{\mathbf A}^1_U
\to 
{\mathbf A}^1_{G_U}$.
At $t=0$,
the restriction
$\bar s|_0\colon
U\to G_U$
is the restriction
$s_0|_U$.
The restrictions
of 
$\tilde s$ and $\bar s$
on 
${\mathbf G}_{m,U}
=
{\mathbf G}_{m,X}
\cap
{\mathbf A}^1_U$
are the same.
Let $\widetilde X
\subset
{\mathbf A}^1_G$
denote the schematic
closure of
$\tilde s({\mathbf G}_{m,X})
\cup
\bar s({\mathbf A}^1_U)$
and let
$\pi:\widetilde X
\to {\mathbf A}^1_X$
be the projection.
Since $\bar s|_0\colon
U\to G_U$
is the restriction
of $s_0$,
the fiber
$\widetilde X_0
=\widetilde X
\times_{{\mathbf A}^1_X}X$
at $t=0$
contains $s_0(X)$
as a closed subscheme.

Let $\widetilde{\cal E}
\subset 
{\cal E}_{\widetilde X}
\oplus 
{\cal E}'_{\widetilde X}$
be the restriction
to $\widetilde X
\subset {\mathbf A}^1_G$
of the
tautological subbundle
and $\tilde e\colon
\widetilde{\cal E}
\to
{\cal E}'_{\widetilde X}$
be the restriction
of the second projection.
We consider
the complex
$\widetilde{\cal K}=
[\widetilde{\cal E}
\to
{\cal E}'_{\widetilde X}]$.
The restriction
of $\widetilde
{\cal K}$
to $\widetilde U=
\bar s({\mathbf A}^1_U)$
is acyclic.
We put
$\widetilde D=
{\mathbf A}^1_D\times
_{{\mathbf A}^1_X}
\widetilde X$ and let
$\pi_D\colon
\widetilde D
\to
{\mathbf A}^1_D$ be the projection.
We consider the composition
$$\begin{CD}
G(\widetilde X)
@>{\gamma_k(\widetilde{\cal K})_
{\widetilde D}}>>
G(\widetilde D)
@>{\pi_{D*}}>>
G({\mathbf A}^1_D).
\end{CD}$$
The pull-backs $i_1^*,i_0^*\colon
G({\mathbf A}^1_D)
\to G(D)$ 
by the
sections $D\to {\mathbf A}^1_D$
at $t=1,0$
are the same isomorphisms.
Since
the fiber $\widetilde{\cal K}_1$
at $t=1$
recovers the original complex
${\cal K}$ on $X$,
we have
$\gamma_k({\cal K})_D([X])=
i_1^*\pi_{D*}(\gamma_k(\widetilde{\cal K})_
{\widetilde D}
([\widetilde X]))$.
Let $\widetilde{\cal K}_0$
denote 
the pull-back
of $\widetilde{\cal K}$
to $\widetilde X_0$
and let
$\pi_0\colon 
\widetilde X_0\to X_0$
be the projection.
We put $\widetilde D_0=
\widetilde D \cap
\widetilde X_0$.
We have further
$\gamma_k({\cal K})_D([X])=
i_0^*\pi_{D*}(\gamma_k(\widetilde{\cal K})_
{\widetilde D}
([\widetilde X]))
=
\pi_{0*}(\gamma_k(\widetilde{\cal K}_0)_
{\widetilde D_0}
([\widetilde X_0])).$

Let $[\overline X_0]$
denote the class
of the kernel
${\rm Ker}(
{\cal O}_{\widetilde X_0}
\to 
{\cal O}_{s_0(X)})$
of the surjection.
Since the restriction
$\widetilde {\cal K}|_{s_0(X)}$
is acyclic
and $\overline X_0
\subset \widetilde D_0$,
we have
$\gamma_k(\widetilde{\cal K}_0)_
{\widetilde D_0}
([\widetilde X_0])
=
\gamma_k(\widetilde{\cal K}_0)_
{\widetilde D_0}
([\overline X_0])
=
\gamma_k([{\cal E}']-
[\widetilde {\cal E}])
([\overline X_0])$.
Hence,
\begin{equation}
\gamma_k({\cal K})_D
([X])
=
\pi_{0*}(\gamma_k([{\cal E}']-
[\widetilde {\cal E}])
([\overline X_0]))
\label{eqX0}
\end{equation}
is an element of
$F_{d-k}G(D)$
as required.

2.
In the notation above,
the localized Chern class
$c_k({\cal K})_D\cap [X]\in
CH_{d-k}(D)$
is defined as
$\pi_{0*}\bigl(c({\cal E}')
c({\cal E})^{-1}
([\overline X_0])\bigl)_{\dim =d-k}$
\cite[Section 18.1]{fulton}.
Hence, the assertion follows.
\qed

The proof of
Proposition
\ref{prlcc} shows that
the system of
maps $\gamma_k({\cal K})_D$
is characterized
by the compatibility
with the Gysin maps
for regular immersions
and the normalization property
that the composition
with the natural map
$G(D)\to G(X)$
is equal to the map
$\gamma_k([{\cal E}']
-[{\cal E}])$.
Similarly as
(\ref{eqX0}),
we have
\begin{equation}
[\Lambda^k({\cal E})
\to 
\Lambda^k({\cal E}')]
_D([X])
=
\pi_{0*}(
([\Lambda^k({\cal E}')]
-
[\Lambda^k(\widetilde{\cal E})])
([\overline X_0])).
\label{eqX0E}
\end{equation}

Similarly as
\cite[Section 18.1]{fulton},
we have the following
properties.

\begin{cor}\label{corlcc}
Let $X$ be a 
scheme of finite type
over a regular noetherian
scheme $S$
of finite dimension
and $D\subset X$
be a closed subscheme.

{\rm 1.}
Let ${\cal E}\to {\cal E}'$
be a morphism
of locally free
${\cal O}_{{\mathbf A}^1_X}$-modules
of finite rank
such that the complex of ${\cal K}=
[{\cal E}\to {\cal E}']$
is acyclic outside 
${\mathbf A}^1_D$.
Let ${\cal K}_0$
and ${\cal K}_1$
be the pull-back
of ${\cal K}$
by the $0$-section
and the $1$-section
respectively.
Then, we have
$$\gamma_k({\cal K}_0)_D=
\gamma_k({\cal K}_1)_D.$$

{\rm 2.}
Let ${\cal E}_1\to {\cal E}'_1$
and ${\cal E}_2\to {\cal E}'_2$
be morphisms
of locally free
${\cal O}_X$-modules
of finite rank
such that the complexes
${\cal K}_1=
[{\cal E}_1\to {\cal E}'_1]$
and
${\cal K}_2=
[{\cal E}_2\to {\cal E}'_2]$
are acyclic outside $D$.
If 
${\cal K}_1\to
{\cal K}_2$
is a quasi-isomorphism, we have
$$\gamma_k({\cal K}_1)_D=
\gamma_k({\cal K}_2)_D.$$
\end{cor}

Let ${\cal K}$
be a complex of ${\cal O}_X$-modules
such that
there exist
a morphism
${\cal E}\to {\cal E}'$
of locally free
${\cal O}_X$-modules
of finite rank
and a quasi-isomorphism
$[{\cal E}\to {\cal E}']
\to {\cal K}$.
Then, we define
the map
$\gamma_k({\cal K})_D
\colon G(X)\to G(D)$
to be 
$\gamma_k([{\cal E}\to
{\cal E}'])_D$.
This is well-defined
by Corollary \ref{corlcc}.2.
The localized Chern class
$c_k({\cal K})_D
\colon CH_j(X)\to CH_{j-k}(D)$
is defined 
similarly as
$c_k([{\cal E}\to {\cal E}'])_D$
if $X$ is of finite type
over a regular noetherian
scheme $S$ of finite dimension.

We consider coherent 
${\cal O}_X$-modules
${\cal F},
{\cal F}'$
and
a morphism 
$f\colon
{\cal F}\to {\cal F}'$ of
${\cal O}_X$-modules
satisfying the following condition:
\begin{itemize}
\item[(\ref{sslcc}.3.1)]
There exists
a locally free ${\cal O}_X$-module
${\cal E}'$ of finite rank and
a surjection
${\cal E}'\to {\cal F}'$.
The kernel
${\cal E}=
{\rm Ker}({\cal F}\oplus
{\cal E}'\to {\cal F}')$
is a locally free ${\cal O}_X$-module
of the same finite rank
as ${\cal E}'$.
\end{itemize}
Let $D\subset X$
be a closed subscheme
such that $f$ is an isomorphism
outside $D$.
Then, since the map
$[{\cal E}\to {\cal E}']
\to
{\mathcal K}= 
[{\cal F}\to {\cal F}']$
is a quasi-isomorphism,
the map
$\gamma_k({\cal F}\to
{\cal F}')_D
\colon G(X)\to G(D)$
and 
the localized Chern class
$c_k({\cal F}\to
{\cal F}')_D
\colon CH_j(X)\to CH_{j-k}(D)$
are defined.

We further assume
that ${\cal F}$ is of
tor-dimension $\le 1$
and let $r$ be the virtual rank.
In other words,
for a surjection
${\cal E}\to {\cal F}$
as in (\ref{sslcc}.3.1),
the kernel
${\cal L}=
{\rm Ker}
({\cal E}\to {\cal F})$
is locally free of rank
${\rm rank}\ {\cal E}-r$.
Since the canonical map
${\rm Ker}
({\cal E}\to {\cal F})\to
{\rm Ker}
({\cal E}'\to {\cal F}')$
is an isomorphism,
the sheaf ${\cal F}'$ is also of
tor-dimension $\le 1$
and of virtual rank $r$.

We also define
a map
$\delta_k({\cal F}\to
{\cal F}')_D\colon
G(X)\to G(D)$
for $k>0$ by requiring that
$\delta_t({\cal F}\to
{\cal F}')_D=
\sum_{k=1}^\infty
\delta_k({\cal F}\to
{\cal F}')_D\cdot t^k$
is given by
$$\delta_t({\cal F}\to
{\cal F}')_D=
\gamma_t({\cal F}\to
{\cal F}')_D
\cdot
\gamma_t([{\cal F}]-r).$$
We recall that the localized Chern class
$c_k({\cal F}'- {\cal F})_D
\colon CH_j(X)\to CH_{j-k}(D)$
is defined by requiring
\addtocounter{thm}1
\setcounter{equation}1
\begin{equation}
\sum_{k>0}c_k({\cal F}'- {\cal F})_D
\cdot t^k=
\left(\sum_{k>0}
c_k({\cal F}\to {\cal F}')_D
\cdot t^k\right)
\cdot c_t({\cal F})
\label{eqck}
\end{equation}
in \cite[(3.24)]{KSA}.

For the definition and
properties
of the derived exterior power
$L\Lambda^k{\cal F}$,
we refer to \cite[Section 1.2]{KSI}.
We recall that
for a locally free resolution
$[{\cal L}\to {\cal E}]\to {\cal F}$
as above,
we have a quasi-isomorphism
$[\Gamma^k{\cal L}\to
\Gamma^{k-1}{\cal L}\otimes
{\cal E}
\to
\cdots
\to
{\cal L}\otimes
\Lambda^{k-1}{\cal E}
\to
\Lambda^k{\cal E}]\to
L\Lambda^k{\cal F}$,
where $\Gamma^\bullet$
denotes the divided power.
For an integer $k>0$, 
the mapping cone
$[L\Lambda^k{\cal F}
\to L\Lambda^k{\cal F}']$
of the derived exterior powers
is defined.
We define a map
$$\begin{CD}
[L\Lambda^k{\cal F}
\to L\Lambda^k{\cal F}']_D\colon
G(X)\to G(D)
\end{CD}$$
sending
$[{\cal G}]$
to
$\sum_q(-1)^q{\cal T}or^{{\cal O}_X}_q
([L\Lambda^k{\cal F}
\to L\Lambda^k{\cal F}'],
{\cal G})$.
We describe
the map
$[L\Lambda^k{\cal F}
\to L\Lambda^k{\cal F}']_D$
using the operators
$\gamma_k({\cal F}\to {\cal F}')_D
\colon G(X)\to G(D)$.

\begin{pr}\label{prLL}
\setcounter{equation}0
Let $S$ be a regular
noetherian scheme
of finite dimension
and $X$ be a scheme
of finite type over $S$.
Let 
${\cal F},{\cal F}'$
be coherent ${\cal O}_X$-modules
of tor-dimension $\le 1$
and of virtual rank $r$
and $f\colon 
{\cal F}\to
{\cal F}'$ be a morphism
of ${\cal O}_X$-modules
satisfying the condition
{\rm (\ref{sslcc}.3.1)}.
Let $D\subset X$
be a closed subscheme
such that $f$ is an isomorphism
outside $D$.

{\rm 1.}
We put formally
$\beta_t=
\sum_{k=1}^\infty
[L\Lambda^k{\cal F}
\to L\Lambda^k{\cal F}']_D
(\frac t{1-t})^k
\colon G(X)
\to G(D)[[t]]$.
Then, we have
$\delta_t({\cal F}
\to {\cal F}')_D
=(1-t)^r\beta_t$.

{\rm 2.}
For $n=r+1$,
we have
\begin{equation}
[L\Lambda^n{\cal F}
\to L\Lambda^n{\cal F}']_D
=\delta_n({\cal F}
\to {\cal F}')_D
\label{eqLL}
\end{equation}
and it
sends the topological filtration
$F_jG(X)$ to
$F_{j-n}G(D)$. Further
the diagram
\begin{equation}
\begin{CD}
CH_j(X)
@>{c_n({\cal F}'-
{\cal F})_D}>>
CH_{j-n}(D)\\
@VVV @VVV\\
{\rm Gr}^F_jG(X)
@>{[L\Lambda^n{\cal F}
\to L\Lambda^n{\cal F}']_D}>>
{\rm Gr}^F_{j-n}G(D)
\end{CD}
\label{eqLLc}
\end{equation}
is commutative.
\end{pr}

In the case where
${\mathcal F}=0$
and $D=X$,
Proposition \ref{prLL}.2
is proved in \cite[Proposition 2.4.4]{KSI}.

{\it Proof.}
1.
The equality $\delta_t({\cal F}
\to {\cal F}')_D
=(1-t)^r\beta_t$
is reduced
to the equality
\begin{equation}
\delta_t({\cal F}
\to {\cal F}')_D([X])
=(1-t)^r\beta_t([X])
\label{eqbt}
\end{equation}
in $G(D)[[t]]$,
by a standard argument.
In fact, since
$G(X)$ is generated
by the classes of integral
closed subschemes $V$,
it suffices to prove the
formula for $[V]$
and to take the push-forward.
By the assumption
that
${\cal F}$ and ${\cal F}'$
are of tor-dimension $\le 1$,
the kernel
${\cal L}=
{\rm Ker}({\cal E}\to {\cal F})=
{\rm Ker}({\cal E}'\to {\cal F}')$
is locally free
and we have quasi-isomorphisms
$[{\cal L}\to {\cal E}]
\to {\cal F}$ and
$[{\cal L}\to {\cal E}']
\to {\cal F}'$.
Hence, we obtain
$[L\Lambda^k{\cal F}
\to L\Lambda^k{\cal F}']_D
([X])
=
\sum_{q=0}^k(-1)^q
[\Lambda^{k-q}{\cal E}
\to \Lambda^{k-q}{\cal E}']_D
[\Gamma^q{\cal L}]([X])$.
We apply the graph construction
to ${\cal E}\to {\cal E}'$
and use the notation
in the proof of Proposition
\ref{prlcc}.
Then by
(\ref{eqX0E}), we have
$$
[L\Lambda^k{\cal F}
\to L\Lambda^k{\cal F}']_D([X])
=\pi_{0*}
\left(
\sum_{q=0}^k(-1)^q
([\Lambda^{k-q}{\cal E}']-
[\Lambda^{k-q}\widetilde{\cal E}])
[\Gamma^q{\cal L}]([\overline X_0])
\right).$$
Thus, 
we obtain
\begin{eqnarray*}
\beta_t([X])
&=&\pi_{0*}\Bigl((
\gamma_t([{\cal E}'])
-\gamma_t([\widetilde{\cal E}]))
\gamma_t([{\cal L}])^{-1}
([\overline X_0])\Bigr)\\
&=&
\gamma_t([{\cal F}'])
\pi_{0*}\left(\Bigl(1
-\frac1{1+\gamma_t(
[{\cal E}\to \widetilde{\cal E}])}\Bigr)
([\overline X_0])\right)
\\
&=&
\gamma_t([{\cal F}'])
\left(1
-\frac1{1+\gamma_t([{\cal E}\to 
{\cal E}'])_D}\right)
([X])
=\gamma_t({\cal F})
\gamma_t([{\cal E}\to {\cal E}'])_D([X])
.\end{eqnarray*}
Since
$\gamma_t({\cal F})
=
\gamma_t({\cal F}-r)\gamma_t(1)^r
=
\gamma_t({\cal F}-r)(1-t)^{-r}$,
we obtain
the equality (\ref{eqbt}).

2.
For $n=r+1$,
we have
\begin{eqnarray*}
&&\delta_t({\cal F}
\to {\cal F}')_D=
(1-t)^r\beta_t\\
&\equiv&
\sum_{k=1}^r
[L\Lambda^k{\cal F}
\to L\Lambda^k{\cal F}']_D
\cdot
t^k(1-t)^{r-k}
+
[L\Lambda^n{\cal F}
\to L\Lambda^n{\cal F}']_D
\cdot
t^n(1-t)^{-1}\bmod t^{n+1}.
\end{eqnarray*}
Comparing the coefficients
of $t^n$, we obtain
$\delta_n({\cal F}
\to {\cal F}')_D=
[L\Lambda^n{\cal F}
\to L\Lambda^n{\cal F}']_D$.
The remaining
assertion follows
from this and Proposition \ref{prlcc}.
\qed

\begin{lm}\label{lmkos}
\setcounter{equation}0
Let ${\cal E}$ and
${\cal E}'$ be locally free
${\cal O}_X$-modules
of the same rank $n$
and 
let ${\cal E}'\to {\cal E}
\to {\cal O}_X$ be
morphisms of ${\cal O}_X$-modules.
We consider the Koszul
complexes
${\cal K}={\rm Kos}(
{\cal E}\to {\cal O}_X)$
and
${\cal K}'={\rm Kos}(
{\cal E}'\to {\cal O}_X)$
and the induced morphism
${\cal K}'\to {\cal K}$
of complexes.

Let ${\cal M}=
[{\cal K}^*\to {\cal K}
^{\prime *}]$
be the mapping cone
of the morphism
of the dual complexes.
Then, the homology sheaves
${\cal H}_q({\cal M})$
are 
${\cal H}_0({\cal K}')
={\cal O}_X/({\rm Image}\ {\cal E}')
$-modules.
\end{lm}

{\it Proof.}
The product
${\cal K}\otimes {\cal K}
\to
{\cal K}$ of Koszul complex
induces a canonical map
${\cal K}\otimes {\cal K}^*
\to
{\cal K}^*$.
The canonical maps
${\cal K}\otimes {\cal K}^*
\to
{\cal K}^*$
and
${\cal K}'\otimes {\cal K}^{\prime *}
\to
{\cal K}^{\prime *}$
induces
${\cal K}'\otimes {\cal M}
\to
{\cal M}$.
This defines a multiplication
${\mathcal H}_0({\cal K}')
\otimes {\mathcal H}_q({\cal M})
\to
{\mathcal H}_q({\cal M})$
compatible with
the ${\cal O}_X$-module
structure.
Thus the assertion follows.
\qed

Let $X$ be a scheme
and ${\cal K}\to {\cal K}'$
be a morphism of chain
complexes of flat 
${\cal O}_X$-modules.
For an integer $q\ge 0$,
we consider the mapping
cone 
$[L\Lambda^q{\cal K}
\to L\Lambda^q{\cal K}']$.
The canonical maps
$L\Lambda^{q+1}{\cal K}\to
L\Lambda^q{\cal K}
\otimes {\cal K}$
and
$L\Lambda^{q+1}{\cal K}'\to
L\Lambda^q{\cal K}'
\otimes {\cal K}'$
\cite[(1.2.1.4)]{KSI}
induce a map
\begin{equation}
[L\Lambda^{q+1}{\cal K}\to
L\Lambda^{q+1}{\cal K}']
\longrightarrow
[L\Lambda^q{\cal K}
\to
L\Lambda^q{\cal K}']
\otimes {\cal K}'.
\label{eqLq}
\end{equation}

We consider the following condition
on a complex ${\cal K}$
of ${\cal O}_X$-modules.
\begin{itemize}
\item[(L($n$))]
For every $x\in X$,
there exist an open neighborhood $U$ of $x$,
locally free ${\cal O}_U$-modules
${\cal E}$ and ${\cal L}$
of rank $n$ and $1$ respectively
and a quasi-isomorphism
$[{\cal L}\to {\cal E}]\to {\cal K}|_U$.
\end{itemize}
If a complex ${\cal K}$
of ${\cal O}_X$-modules
satisfies
the condition (L($n$)),
it is a perfect complex
of tor-dimension $\le 1$.
For a perfect complex ${\cal K}$
of ${\cal O}_X$-modules
satisfying the condition (L($n$)),
let $Z$ denote the closed subscheme
of $X$ defined by the
annihilator ideal
${\cal I}_Z=
{\rm Ann}\Lambda^n{\cal H}_0({\cal K})$
and $i\colon Z\to X$
denote
the immersion.
Then, for
a quasi-isomorphism
$[{\cal L}\to {\cal E}]\to {\cal K}|_U$
on an open subscheme $U$
of $X$ as in (L($n$)),
the intersection $Z\cap U\subset U$
is the largest closed
subscheme where the map
${\cal L}\to {\cal E}$
is the zero-map.
Further, a quasi-isomorphism
$[{\cal L}\to {\cal E}]\to {\cal K}|_U$
induces an isomorphism
$i^*{\cal L}
\to{\cal L}_Z|_{Z\cap U}=
L_1i^*{\cal K}|_{Z\cap U}$.
Hence, 
${\cal L}_Z=L_1i^*{\cal K}$
is an invertible ${\cal O}_Z$-module
\cite[Lemma 2.4.1.1]{KSI}.

\begin{cor}\label{corkos}
Let ${\cal K}$
and ${\cal K}'$
be complexes
of ${\cal O}_X$-modules
satisfying the condition
{\rm (L($n$))}
and ${\cal K}\to
{\cal K}'$
be a morphism
such that the mapping
cone $[{\cal K}\to {\cal K}']$
is of tor-dimension $\le 1$.
For $q>0$,
let $[L\Lambda^q{\cal K}
\to L\Lambda^q{\cal K}']$
denote the mapping cone.

Let $Z$ and $Z'$
denote the closed subschemes
of $X$ defined by the
annihilator ideals
${\cal I}_Z=
{\rm Ann}\Lambda^n{\cal H}_0({\cal K})$
and
${\cal I}_{Z'}=
{\rm Ann}\Lambda^n{\cal H}_0({\cal K}')$.
Let $i\colon Z\to X$
and $i'\colon Z'\to X$
denote
the immersions
and let
${\cal L}_Z=L_1i^*{\cal K}$
and 
${\cal L}'_{Z'}=L_1i^{\prime*}{\cal K}'$
be the invertible ${\cal O}_Z$-module
and 
${\cal O}_{Z'}$-module respectively.

{\rm 1.}
The scheme $Z$ is a closed subscheme of $Z'$
and the canonical map
${\cal L}'_{Z'}
\otimes_{{\cal O}_{Z'}}
{\cal O}_Z\to 
{\cal L}_Z$ is an isomorphism.

{\rm 2.}
The homology sheaf
${\cal H}_p([L\Lambda^q{\cal K}
\to L\Lambda^q{\cal K}'])$
is an ${\cal O}_{Z'}$-module
if $p>0$ or $q\ge n$.

{\rm 3.}
The map
{\rm (\ref{eqLq})}
induces an isomorphism
\setcounter{equation}0
\begin{equation}
{\cal H}_{p+1}([L\Lambda^{q+1}{\cal K}
\to L\Lambda^{q+1}{\cal K}'])
\longrightarrow 
{\cal H}_p([L\Lambda^q{\cal K}
\to L\Lambda^q{\cal K}'])
\otimes {\cal L}'_{Z'}
\label{eqkos}
\end{equation}
for $p>0$ or $q\ge n$.
\end{cor}

In the case where
${\mathcal K}=0$,
Corollary \ref{corkos}.2 (resp.\ .3)
is proved in \cite[Lemma 2.4.2.1 (resp.\ .2)]{KSI}.

{\it Proof.}
1.
Since the
question is local,
we may assume
that there is a quasi-isomorphism
$[{\cal L}\to {\cal E}']
\to {\cal K}'$
for locally free sheaves
${\cal E}'$
and ${\cal L}$ of
rank $n$ and 1.
Then, the mapping cone
$[{\cal E}'\oplus {\cal K}
\to {\cal K}']$
is a perfect complex of
tor-amplitude $[1,1]$
and hence
is quasi-isomorphic
to ${\cal E}[1]$
for a locally free
sheaf ${\cal E}$
of rank $n$.
Thus,
we obtain 
a quasi-isomorphism
$[{\cal L}\to {\cal E}]
\to {\cal K}$.
Since the pull-back
of ${\cal L}\to {\cal E}$
to $Z$ 
is the zero-map,
the pull-back
of ${\cal L}\to {\cal E}'$
to $Z$ is also the zero-map.
This shows that $Z$
is a subscheme of $Z'$.
Further,
the canonical map
${\cal L}'_{Z'}
\otimes_{{\cal O}_{Z'}}
{\cal O}_Z\to 
{\cal L}_Z$ is induced
by the identity of
${\cal L}$ and is an isomorphism.

2.
Since the assertion is local, 
we may assume, as in the proof of 1.,  
that ${\cal K}=
[{\cal O}_X\overset e\to {\cal E}]$ 
and ${\cal K}'=
[{\cal O}\overset {e'}\to {\cal E}']$ 
where  ${\cal E}$ and ${\cal E}'$ 
are locally free of rank $n$, 
and ${\cal K}\to {\cal K}'$ is 
induced by
${\cal E}\to {\cal E}'$
compatible with $e$ and $e'$.
Then for $q\geq 0$, 
the derived exterior power complex
$L\Lambda^q {\cal K}$ 
is identified with the complex
$[{\cal O}_X\overset e\to {\cal E} 
\overset {\wedge e}\to \Lambda^2 {\cal E}
\overset {\wedge e}\to\cdots 
\overset {\wedge e}
\to \Lambda^q {\cal E}]$.
Consequently,
if $q\ge n$, its shift 
$(L\Lambda^q {\cal K})[-q]$ is isomorphic to the dual of the Koszul complex
${\rm Kos}(e^*\colon  {\cal E}^*\to {\cal O}_X)$ 
associated to the dual
of $e\colon{\cal O}_X\to {\cal E}$
and similarly for the complex
$(L\Lambda^q {\cal K}')[-q]$.
Since the ideal sheaf ${\mathcal I}_{Z'}
\subset {\mathcal O}_X$ of $Z'$ is 
the image of
$e^{\prime*}\colon  
{\cal E}^{\prime*}\to {\cal O}_X$, 
the assertion follows from Lemma \ref{lmkos}.

3.
We show that
(\ref{eqLq}) induces
(\ref{eqkos}).
The map  (\ref{eqLq}) induces a homomorphism 
\begin{align*}
{\cal H}_{p+1}([L\Lambda^{q+1} {\cal K}\to L\Lambda^{q+1}{\cal K}']) 
\to \
&{\cal H}_{p+1}([L\Lambda^q{\cal K} \to L\Lambda^q{\cal K}']\otimes {\cal K}')
\\&=
{\cal T}or_{p+1}^{{\cal O}_X}([L\Lambda^q{\cal K} \to L\Lambda^q{\cal K}'],{\cal K}').
\end{align*}
Since ${\cal K}'$ is of tor-dimension $\leq 1$, the spectral sequence 
$$E^2_{s,t}=
{\cal T}or_s^{{\cal O}_X}
({\cal H}_t([L\Lambda^q {\cal K}\to L\Lambda^q {\cal K}']), {\cal K}') 
\Rightarrow 
{\cal T}or_{s+t}^{{\cal O}_X}
([L\Lambda^q{\cal K}\to L\Lambda^q {\cal K}'], {\cal K}')$$ 
induces a homomorphism
$${\cal T}or_{p+1}^{{\cal O}_X}
([L\Lambda^q {\cal K} \to L\Lambda^q {\cal K}'], {\cal K}')
\to {\cal T}or_1^{{\cal O}_X}
({\cal H}_p([L\Lambda^q{\cal K}\to L\Lambda^q {\cal K}']), {\cal K}').$$
Since ${\cal H}_p([L\Lambda^q{\cal K}\to L\Lambda^q {\cal K}'])$ 
is an ${\cal O}_{Z'}$-module, 
we have an isomorphism
\begin{align*}
{\cal T}or_1^{{\cal O}_X}
({\cal H}_p([L\Lambda^q{\cal K}\to L\Lambda^q {\cal K}']), {\cal K}')
\to\ &
{\cal H}_p([L\Lambda^q{\cal K}\to L\Lambda^q {\cal K}'])
\otimes
L_1i^{\prime*}{\cal K}'\\
&=
{\cal H}_p([L\Lambda^q{\cal K}\to L\Lambda^q {\cal K}'])
\otimes_{{\cal O}_{Z'}}
{\cal L}'_{Z'}.
\end{align*}
(cf.\ \cite[Lemma 2.4.1.3]{KSI}).
They define a canonical map
(\ref{eqkos}).

By \cite[(2.4.2.2)]{KSI}, 
the canonical maps
$$L_{p+1}\Lambda^{q+1}{\cal K}\to L_p\Lambda^q{\cal K}\otimes_{{\cal O}_Z} {\cal L}_Z, \quad L_{p+1}\Lambda^{q+1}{\cal K}'\to L_p\Lambda^q{\cal K}'\otimes_{{\cal O}_{Z'}} {\cal L}'_{Z'}$$
are isomorphisms.
Further by 1., 
the canonical map
$L_p\Lambda^q{\cal K}\otimes_{{\cal O}_Z} {\cal L}_Z
\to 
L_p\Lambda^q{\cal K}\otimes_{{\cal O}_{Z'}} {\cal L}_{Z'}$
is an isomorphism.  
Thus the map (\ref{eqkos}) is an isomorphism.
\qed

\subsection{Localized intersection 
product}\label{sslip}
\setcounter{equation}0

We briefly recall terminologies
and properties
on cotangent complexes
\cite[Chapitre II]{ctcx}
and excess conormal complexes
\cite[Definition 1.6.3]{KSI}
which will play the central roles in this
and the next subsections.

For a morphism of schemes $X\to S$, 
the cotangent complex $L_{X/S}$ is defined  in \cite[Chapitre II]{ctcx} 
as an object of the derived category 
$D^-_{qcoh}
({\mathcal O}_X)$
of the category 
of quasi-coherent ${\mathcal O}_X$-modules.
We have ${\mathcal H}_0(L_{X/S})=\Omega_{X/S}^1$ and 
${\mathcal H}_i(L_{X/S})=0$ for $i<0$. 
For a morphism
$f\colon X\to Y$ of
schemes over $S$, we have a distinguished triangle 
\begin{equation}
Lf^*L_{Y/S}\to L_{X/S}\to L_{X/Y}\to.
\label{eqLtr}
\end{equation}
If $X\to P$ is a regular
closed immersion over $S$
and if $P$ is smooth over $S$, 
then there exists a canonical
isomorphism 
\begin{equation}
L_{X/S}\to[N_{X/P}
\overset d\to \Omega_{P/S}^1\otimes_{{\mathcal O}_P} {\mathcal O}_X]
\label{eqLXS}
\end{equation} 
to the complex
where the conormal sheaf
$N_{X/P}={\mathcal I}_X/
{\mathcal I}_X^2$
is put on degree 1
and $\Omega^1_{P/S}\otimes_{{\mathcal O}_P}{\mathcal O}_X$ is put on degree $0$.
If $T\to S$ is flat, 
then the canonical map
$Lp_1^*L_{X/S}\to
L_{X\times_S T/T}$ 
is an isomorphism
where 
$p_1\colon 
X\times_S T\to X$
denotes the first projection.

For an immersion $V\to X$ of schemes, 
the conormal complex $M_{V/X}$ 
is defined in \cite[Definition 1.6.3]{KSI} to be the shift $L_{V/X}[-1]$  of the cotangent complex $L_{V/X}$. We have a canonical isomorphism
${\cal H}_0(M_{V/X})\to N_{V/X}$
to the conormal sheaf.

\begin{lm}\label{lma08} 
If $X\to S$ is flat, then 
there exists
a canonical isomorphism
$L_{X/S}\to 
M_{X/(X\times_S X)}$. 
\end{lm}

{\it Proof.} 
We consider the distinguished
triangle (\ref{eqLtr})
for $X\to X\times_S X\to X$, where the first arrow is the diagonal
$\delta\colon
X\to X\times_S X$
and the second arrow is the second projection.
Since $L_{X/X}=0$, we obtain
an isomorphism
$M_{X/(X\times_S X)}=
L_{X/(X\times_S X)}[-1] \to 
L \delta^*L_{(X\times_S X)/X}$. 
Since $X\to S$ is assumed to be flat, 
the canonical map
$p_1^*L_{X/S}
\to
L_{(X\times_S X)/X}$
is an isomorphism 
where $p_1$ denotes the first projection. 
It induces
an isomorphism
$L_{X/S}\to L\delta^*Lp^*_1L_{X/S}\to 
L\delta^*L_{(X\times_S X)/X}$.  
\qed

For an immersion
$V\to X$ and for a morphism $W\to X$ such that $T= V\times_X W \to W$ is a regular immersion, the excess conormal complex $M'_{V/X,W}$ is defined 
in \cite[Definition 1.6.3]{KSI} as an object of the derived category 
$D^-_{qcoh}({\mathcal O}_V)$
of the category of quasi-coherent ${\mathcal O}_V$-modules. 
We have a distinguished triangle 
$$M'_{V/X,W}\to Lg^*M_{V/X}\to N_{T/W}\to $$
where $g$ is the morphism $T\to V$.

We recall the definition of
locally a hypersurface \cite[Definition 3.1.1]{KSI}.

\begin{df}\label{dflhs}
Let $S$ be a scheme.
A scheme $X$ of finite presentation over $S$
is called locally a hypersurface
of relative virtual dimension $n-1$
if, locally on $X$,
it is a Cartier divisor of
a smooth scheme 
of relative dimension $n$
over $S$.
\end{df}

For such $X$,
if $i\colon
X\to P$
is a regular immersion
to a smooth scheme over $S$,
the cotangent complex
$L_{X/S}$
is canonically 
quasi-isomorphic to
the complex
$[N_{X/P}\to
i^*\Omega^1_{P/S}]$
inducing a canonical isomorphism
${\cal H}_0(L_{X/S})
\to \Omega^1_{X/S}$.
Consequently,
the complex $L_{X/S}$
satisfies the condition
(L($n$)) in Section \ref{sslcc}.

We recall the definition
of the localized intersection product
\cite[Definition 3.2.2]{KSI}.
Let $S$ be a regular noetherian
scheme of finite dimension
and let $X$ be
a scheme of finite type over $S$
that is locally a hypersurface
of relative virtual dimension $n-1$.
Let $Z$ be the closed subscheme
of $X$ defined by the annihilator
of $\Omega_{X/S}^n$,
$i\colon Z\to X$
be the closed immersion
and let
${\cal L}_Z=L_1i^*L_{X/S}$
be the invertible ${\cal O}_Z$-module.
The underlying set of $Z$
is the complement
of the largest open 
subscheme of $X$
smooth over $S$.

Let $V$ be a closed subscheme of $X$
and $W$ be a noetherian scheme over $X$.
We put $T=V\times_XW$.
Then, 
for a coherent
${\cal O}_V$-module ${\cal F}$,
a complex ${\cal G}
\in D^b_{\text{coh}}({\cal O}_W)$
of ${\cal O}_W$-modules
with bounded coherent cohomology sheaves
and for 
a sufficiently large integer $q$,
the ${\cal O}_T$-module
${\cal T}or_q^{{\cal O}_X}(
{\cal F},{\cal G})$
is supported on the
inverse image $Z_T=Z\times_XT$
and the class
$[{\cal T}or_q^{{\cal O}_X}(
{\cal F},{\cal G})]$ 
in $G(Z_T)_{/{\cal L}_Z}=
{\rm Coker}({\cal L}_Z-1\colon 
G(Z_T)\to G(Z_T))$
depends only on
the parity of $q$ by
\cite[Theorem 3.2.1]{KSI}.

\begin{df}\label{dflip1}
\setcounter{equation}0
Let $S$ be a regular noetherian
scheme of finite dimension
and let $X$ be
a scheme of finite type over $S$
that is locally a hypersurface
of relative virtual dimension $n-1$.
Let $Z$ be the closed subscheme
of $X$ defined by the annihilator
of $\Omega_{X/S}^n$,
$i\colon Z\to X$
be the closed immersion
and let
${\cal L}_Z=L_1i^*L_{X/S}$
be the invertible ${\cal O}_Z$-module.
Let $V$ be a closed subscheme of $X$
and $W$ be a noetherian scheme over $X$.
We put $T=V\times_XW$, $Z_T=Z\times_XT$
and $G(Z_T)_{/{\cal L}_Z}=
{\rm Coker}({\cal L}_Z-1\colon 
G(Z_T)\to G(Z_T))$.

Then, the localized intersection product
\begin{equation}
\begin{CD}
((\ , \ ))_X\colon
G(V)\times G(W)
@>>>
G(Z_T)_{/{\cal L}_Z}
\end{CD}
\label{eqlip}
\end{equation}
is a biadditive pairing defined by
\begin{equation}
(([{\cal F}] ,[{\cal G}]))_X
=
(-1)^q
([{\cal T}or_q^{{\cal O}_X}(
{\cal F},{\cal G})]
-
[{\cal T}or_{q+1}^{{\cal O}_X}(
{\cal F},{\cal G})])
\label{eqlip2}
\end{equation}
for a coherent
${\cal O}_V$-module ${\cal F}$,
a coherent
${\cal O}_W$-module ${\cal G}$
and for 
a sufficiently large integer $q$.
\end{df}

In \cite{KSI},
it is denoted $[[\ ,\ ]]$.
We have changed the notation
to emphasize the similarity with
the usual intersection pairing $(\ ,\ )$
of algebraic cycles.
It is proved in \cite[Theorem 3.2.1]{KSI},
that the right hand side of
(\ref{eqlip2})
is independent of $q$ sufficiently large
and defines a pairing.

For ${\cal G}={\cal O}_W$,
we put
\begin{equation}
(([{\cal F}] ,[{\cal O}_W]))_X
=
(([{\cal F}] ,W))_X.
\label{eqlip3}
\end{equation}
The localized intersection product
with $W$
defines a map
\begin{equation}
\begin{CD}
((\ ,W))_X\colon
G(V)@>>> G(Z_T)_{/{\cal L}_Z}.
\end{CD}
\label{eqlip4}
\end{equation}
Similarly,
we define a map
$((V,\ ))_X\colon
G(W)\to G(Z_T)_{/{\cal L}_Z}$.
If $\Sigma\subset W$ 
is a closed subset
such that
the restriction 
$W\setminus \Sigma\to X$
is of finite tor-dimension,
the map (\ref{eqlip4})
is lifted to a map
\begin{equation}
\begin{CD}
((\ ,W))_X\colon
G(V)@>>> G(Z_T\times_W\Sigma)_{/{\cal L}_Z}.
\end{CD}
\label{eqlip5}
\end{equation}

We recall some formulas on
localized intersection product.

\begin{lm}\label{lmfm1}
Let $X$ and $X'$ be locally
hypersurfaces of relative
virtual dimension $n-1$ and $n'-1$
over a
regular noetherian scheme $S$
of finite dimension
and 
$$\begin{CD}
V@<<< T\\
@V{\cap}VV @VV{\cap}V\\
X'@<<< W'\\
@VfVV @VVV\\
X@<<< W
\end{CD}$$
be a cartesian diagram
of noetherian schemes
over $S$
where the vertical arrows
are closed immersions.

We assume that
the immersion $f\colon
X'\to X$ is of finite tor-dimension.
We also assume that 
$Z\times_XX'$
is a closed subscheme of
the closed subscheme $Z'$
of $X'$ defined by the annihilator
of $\Omega_{X'/S}^{n'}$
and that
the map $G(Z_T) \to G(Z'_T)$ induced by the inclusion $Z_T\to Z'_T$
induces a map
$G(Z_T)_{/{\cal L}_Z}
\to
G(Z'_T)_{/{\cal L}_{Z'}}.$

Then, the composition
$$\begin{CD}
G(V)@>{((\ ,W))_X}
>>
G(Z_T)_{/{\cal L}_Z}
@>>>
G(Z'_T)_{/{\cal L}_{Z'}}
\end{CD}$$
is equal to the localized
intersection product
$$\begin{CD}
((\ ,[Lf^*{\cal O}_W]))_{X'}\colon
G(V)
@>>>
G(Z'_T)_{/{\cal L}_{Z'}}
\end{CD}$$
with $[Lf^*{\cal O}_W]
\in G(W')$.
\end{lm}

{\it Proof.}
It follows from the 
canonical isomorphism
${\cal T}or^{{\cal O}_X}_q
({\cal F},{\cal O}_W)
\to
{\cal T}or^{{\cal O}_{X'}}_q
({\cal F},Lf^*{\cal O}_W)$
for a coherent
${\mathcal O}_V$-module
${\mathcal F}$
\cite[Lemma 1.5.1]{KSI}.
\qed

\begin{lm}\label{lmfm3}
Let $X$ be locally a
hypersurface of
relative dimension $n-1$ over a
regular noetherian scheme $S$
of finite dimension
and 
$$\begin{CD}
V@<<< T\\
@V{\cap}VV @VV{\cap}V\\
X@<<< W
\end{CD}$$
be a cartesian diagram
of noetherian schemes
over $S$
where the vertical arrows
are closed immersions.
We assume $V$ is regular.

Then, the map
$((\ ,W))_X\colon
G(V)
\to 
G(Z_T)_{/{\cal L}_Z}$
is equal to the usual intersection
product
$$\begin{CD}
(\ ,((V,W))_X)_V\colon
G(V)
@>>>
G(Z_T)_{/{\cal L}_Z}
\end{CD}$$
with the
localized intersection
product
$((V,W))_X\in 
G(Z_T)_{/{\cal L}_Z}$.

Assume in addition
that $V$ 
is locally of complete intersection
of relative dimension $n-c$
over $S$, 
that $W$ is of dimension $p$
and that
the immersion
$T\to W$ is a regular
immersion of codimension $c'$.
Let $M'_{V/X,W}$
denote the excess conormal complex.
Then,
we have
\setcounter{equation}0
\begin{equation}
((V,W))_X
=
(-1)^{c-c'}
{c_{c-c'}}^T_{Z_T}
(M'_{V/X,W})\cap [T]
\label{eqfm3}
\end{equation}
in ${\rm Gr}^F_{p-c}(G(Z_T)_{/{\mathcal L}_Z})$.
\end{lm}

{\it Proof.}
We apply
\cite[Lemma 3.3.1]{KSI}
to the spectral sequence
$E^2_{p,q}=
{\cal T}or_p^{{\cal O}_V}
({\cal F},{\cal T}or_q^{{\cal O}_X}
({\cal O}_V,
{\cal O}_W))
\Rightarrow
E_{p+q}=
{\cal T}or_{p+q}^{{\cal O}_X}
({\cal F},{\cal O}_W)$.
Then, 
similarly as
in the proof of
\cite[Proposition 3.3.2.1]{KSI},
we obtain an equality
$(({\cal F} ,W))_X
=
({\cal F} ,((V,W))_X)_V$
for
a coherent ${\cal O}_V$-module
${\cal F}$.

The equality (\ref{eqfm3})
is 
\cite[(3.4.4.1)]{KSI}.
\qed

\begin{lm}\label{lmfm2}
Let $X$ be locally
a hypersurface over a
regular noetherian scheme $S$
of finite dimension
and 
$$\begin{CD}
V@<<< T@<<< T'\\
@V{\cap}VV @VV{\cap}V@VV{\cap}V\\
X@<<< W@<g<< W'
\end{CD}$$
be a cartesian diagram
of noetherian schemes
over $S$
where the vertical arrows
are closed immersions.
We assume that
$g\colon W'\to W$ 
is of finite tor-dimension.

Then, the map
$((\ ,W'))_X\colon
G(V)
\to 
G(Z_{T'})_{/{\cal L}_Z}$
is equal to the composition of
$$\begin{CD}
G(V)
@>{((\ ,W))_X}>>
G(Z_T)_{/{\cal L}_Z}
@>{g^*}>>
G(Z_{T'})_{/{\cal L}_Z}.
\end{CD}$$
\end{lm}

{\it Proof.}
It suffices to put
${\cal G}={\cal O}_W$
and
${\cal H}={\cal O}_{W'}$
in \cite[Proposition 3.3.2.1]{KSI}.
\qed

\begin{lm}\label{lmfm4}
\setcounter{equation}0
Let $X$ be locally
a hypersurface over a
regular noetherian scheme $S$
of finite dimension
and let $X'$ be locally
a hypersurface over a
regular noetherian scheme $S'$
of finite dimension.
We consider
a cartesian diagram
$$\begin{CD}
V@<<< V'@<<< T\\
@V{\cap}VV @V{\cap}VV@VV{\cap}V\\
X@<f<< X'@<<< W
\end{CD}$$
of noetherian schemes
over $S$
where the vertical arrows
are closed immersions.
We assume that
$f\colon X'\to X$ 
is of finite tor-dimension.
We also assume
that 
$Z_T=Z\times_XT$
is a subset of
$Z'_T=Z'\times_{X'}T$
set-theoretically
and that
the canonical morphism
$G(Z_T)\to G(Z'_T)$
induces
a morphism
$G(Z_T)_{/{\cal L}_Z}
\to 
G(Z'_T)_{/{\cal L}_{Z'}}.$

Then, the composition
$((\ ,W))_X\colon
G(V)
\to 
G(Z_{T})_{/{\cal L}_Z}
\to 
G(Z'_T)_{/{\cal L}_{Z'}}$
is equal to the composition of
$$\begin{CD}
G(V)
@>{f^*}>>
G(V')
@>{((\ ,W))_{X'}}>>
G(Z'_T)_{/{\cal L}_{Z'}}.
\end{CD}$$
\end{lm}

{\it Proof.}
Similarly as
\cite[Corollary 3.3.4.3]{KSI}
applied by taking $W=X'$
and $V'$ to be $W$,
it follows from
\cite[Proposition 3.3.3]{KSI}.
\qed

\subsection{Relative excess
intersection formula}\label{ssrex}

We establish a refinement
of the excess intersection formula
\cite[Proposition 3.4.2]{KSI}.
The result in this
subsection will be
used in the proof of
an explicit computation of
the logarithmic different in
Proposition \ref{prrexX}.
Most of the proofs in this subsection
are immediate variations
of those in \cite[Sections 1.6, 1.7]{KSI}.
The reader is recommended to
read them before following the proof here.
The authors apologize for
possible inconvenience.

We consider a commutative diagram
\setcounter{equation}0
\begin{equation}
\begin{CD}
T@>\subset>> W\\
@VVV @VV{g'}V\\
V'@>\subset>> X'\\
@VVV @VVfV\\
V@>\subset>> X
\end{CD}
\label{eqgl}
\end{equation}
of schemes
satisfying the following condition:
\begin{itemize}
\item[(\ref{ssrex}.0.2)]
\addtocounter{equation}1
The horizontal
arrows are closed immersions
and the immersion $T\to W$
is a regular immersion.
The upper square 
and the tall rectangle
are cartesian.
\end{itemize}
Let $g\colon
W\to X$ denote the composition
of the right vertical arrows.
Recall that a simplicial algebra
$A_{V/X,W}$ on $T$
such that 
the normal complex
$N(A_{V/X,W})$
computes
$Lg^*{\cal O}_V$
is defined in
\cite[1.6.3, 1.6.4]{KSI}.
Further,
an ideal
$I_{V/X,W}\subset
A_{V/X,W}$ is defined as the kernel
of the surjection
$A_{V/X,W}\to {\cal O}_T$
and the excess conormal complex
$M'_{V/X,W}=N(I_{V/X,W}/I^2_{V/X,W})[-1]$
is defined as a
complex of ${\cal O}_T$-modules
loc.\ cit.

We construct a variant of
the spectral sequence 
$E^1_{p,q}
={\cal H}_{2p+q}
L\Lambda^{-p}
M'_{V/X,W}
\Rightarrow
{\cal T}or^{{\cal O}_{X'}}_{p+q}
({\cal O}_V,{\cal O}_W)$
\cite[(1.6.4.3)]{KSI}.
We further assume that
the map $f\colon X'\to X$
is of finite tor-dimension.
We define an object
${\cal C}$ of 
$D^-_{\rm coh}({\cal O}_{X'})$
fitting in the distinguished
triangle
\begin{equation}
\to {\cal C}
\to Lf^*{\cal O}_V
\to {\cal O}_{V'}\to
\label{eqCV}
\end{equation}
as follows.
Since $Lf^*{\cal O}_V$
is acyclic in degree $>0$,
the $0$-th subcomplex
$\tau_{\le 0}Lf^*{\cal O}_V$
in the canonical filtration
is quasi-isomorphic to 
$Lf^*{\cal O}_V$.
Namely, for a complex
${\cal C}'$ defining
$Lf^*{\cal O}_V$
in the derived category
and the subcomplex
$\tau_{\le 0}{\cal C}'$
defined by
replacing the components
${\cal C}'_q$ by $0$ for $q>0$
and
${\cal C}'_0$ by 
${\rm Ker}({\cal C}'_0
\to {\cal C}'_{-1})$, the inclusion
$\tau_{\le 0}{\cal C}'
\to 
{\cal C}'$ is a quasi-isomorphism.
Let ${\cal C}$
be the subcomplex of ${\cal C}'$
obtained by further
replacing
${\rm Ker}({\cal C}'_0
\to {\cal C}'_{-1})$ on degree 0
by the kernel of the surjection
${\rm Ker}({\cal C}'_0
\to {\cal C}'_{-1})
\to {\cal O}_{V'}$.
Then, the complex
${\cal C}$ fits in the
distinguished triangle
(\ref{eqCV})
and is independent
of the choice of ${\cal C}'$
up to a canonical quasi-isomorphism.

We use the notation
in \cite[Proposition 1.6.4]{KSI}.
A simplicial algebra
$A_{V/X,W}$ 
and an increasing filtration
$F^\bullet A_{V/X,W}
=I^\bullet_{V/X,W}$ 
are defined 
in \cite[Definition 1.6.3]{KSI}.
Similarly
$A_{V'/X',W}$
and its filtration $F^\bullet$ are defined.
They define a filtered
chain complex
$\widetilde {\cal C}
=[A_{V/X,W}\to A_{V'/X',W}]$
of simplicial modules.
Here $A_{V'/X',W}$ is 
put on degree 0
and $A_{V/X,W}$ is on degree 1.
Then, 
since $P_X({\cal O}_V)$
and $P_{X'}({\cal O}_{V'})$
are resolutions of ${\cal O}_V$
and of ${\cal O}_{V'}$
by free simplicial algebras
\cite[1.5.5.6]{ctcx},
we have
a quasi-isomorphism
$\int N\widetilde {\cal C}
\to Lg^{\prime*}{\cal C}$
from
the associated simple complex.
Thus, we obtain a spectral sequence
\begin{equation}
E^1_{p,q}=
{\cal H}_{p+q}
{\rm Gr}_F^{-p}
\textstyle{\int} N
\widetilde {\cal C}
\Rightarrow
{\cal T}or^{{\cal O}_{X'}}_{p+q}
({\cal C},{\cal O}_W).
\label{eqspsqC}
\end{equation}

\begin{lm}\label{lmspsqA}
\setcounter{equation}0
Let the notation be as above.
Then, for the $E^1$-term
of {\rm (\ref{eqspsqC})},
there exists a canonical 
isomorphism
\begin{equation}
E^1_{p,q}
\to {\cal H}_{2p+q}[
L\Lambda^{-p}
M'_{V/X,W}
\to
L\Lambda^{-p}
M'_{V'/X',W}]
\end{equation}
of ${\cal O}_T$-modules.
\end{lm}

By Lemma \ref{lmspsqA},
we obtain a spectral sequence
\begin{equation}
E^1_{p,q}
={\cal H}_{2p+q}[
L\Lambda^{-p}
M'_{V/X,W}
\to
L\Lambda^{-p}
M'_{V'/X',W}]
\Rightarrow
{\cal T}or^{{\cal O}_{X'}}_{p+q}
({\cal C},{\cal O}_W)
\label{eqspsqA}
\end{equation}
of ${\cal O}_T$-modules.

{\it Proof.}
The proof is similar to
\cite[Proposition 1.6.4]{KSI}.
Similarly as loc.\ cit.,
the canonical map
$S^p({\rm Gr}_F^1
\widetilde {\cal C})
\to {\rm Gr}_F^p
\widetilde {\cal C}$
is an isomorphism.
The normal complex
of the graded piece
${\rm Gr}_F^1
\widetilde {\cal C}
=[{\rm Gr}_F^1A_{V/X,W}
\to
{\rm Gr}_F^1A_{V'/X',W}]$
is defined
by the canonical map
$M'_{V/X,W}[1]
\to
M'_{V'/X',W}[1]$
of conormal complexes.
Hence the normal complex
$N{\rm Gr}_F^p
\widetilde {\cal C}
=N[{\rm Gr}_F^pA_{V/X,W}
\to
{\rm Gr}_F^pA_{V'/X',W}]$
is canonically quasi-isomorphic to
the mapping cone of
$NS^p(M'_{V/X,W}[1])
\to
NS^p(M'_{V'/X',W}[1])$.
Thus by \cite[Proposition 1.2.8]{KSI}, 
the $E^1$-term
$E^1_{p,q}=
{\cal H}_{p+q}
{\rm Gr}_F^{-p}\int N
\widetilde {\cal C}$
is given by
${\cal H}_{2p+q}[
L\Lambda^{-p}
M'_{V/X,W}
\to
L\Lambda^{-p}
M'_{V'/X',W}]$
as required.
\qed

\addtocounter{thm}1
\setcounter{equation}0

We consider a commutative diagram
\begin{equation}
\begin{CD}
T@>\subset>> W@.\\
@VVV @VVV@.\\
V'@>\subset>> X'@>\subset>>P'\\
@VVV @VVV @VVV\\
V@>\subset>> X@>\subset>>P
\end{CD}
\label{eqlc}
\end{equation}
of schemes
satisfying the following conditions:
\begin{itemize}
\item[(\ref{lmspmap}.2)]
\addtocounter{equation}1
The horizontal
arrows are closed immersions
and the immersions
$X\to P$ and
$X'\to P'$ 
are regular immersions
of the same codimension.
The immersion $T\to W$
is also a regular immersion.
The upper square,
the left tall rectangle
and the right square
are cartesian.
\end{itemize}

In \cite[(1.7.2.1)]{KSI},
a map $\lambda_{V/X/P,W}
\colon
L\Lambda^pM'_{V/X,W}
\to 
N_{X/P}\otimes
L\Lambda^{p-1}M'_{V/X,W}[1]$
is defined as the composition
of the maps
\begin{equation}
L\Lambda^pM'_{V/X,W}
\to
M'_{V/X,W}
\otimes
L\Lambda^{p-1}M'_{V/X,W}
\to 
N_{X/P}\otimes
L\Lambda^{p-1}
M'_{V/X,W}[1]
\label{eqlamb0}
\end{equation}
induced by the canonical map
$M'_{V/X,W}
\to Lg_T^*M_{V/X}
\to Lg_T^*M_{X/P}[1]=
N_{X/P}\otimes
{\cal O}_T[1]$.
We identify
$N_{X'/P'}=f^*N_{X/P}$
and define
$\lambda_{V'/X'/P',W}
\colon
L\Lambda^pM'_{V'/X',W}
\to 
N_{X/P}\otimes
L\Lambda^{p-1}
M'_{V'/X',W}[1]$
similarly.
We will construct a map
(\ref{eqlamb}) below compatible
with $\lambda_{V/X/P,W}$
and
$\lambda_{V'/X'/P',W}$.

We use the notation
in the proof of
\cite[Proposition 1.7.2]{KSI}.
Then,
we have a commutative diagram of
exact sequences
\begin{equation}
\begin{CD}
0@>>>
J_{B'}/J_{B'}^2
@>>>
B'/J_{B'}^2
@>>>
A'@>>>0\\
@.@AAA 
@AAA 
@AAA 
@.\\
0@>>>
J_B/J_B^2
@>>>
B/J_B^2
@>>>
A@>>>0
\end{CD}
\end{equation}
of filtered simplicial modules.
The lower exact sequence
is constructed for $V\to X\gets W$
and the upper one is
for $V'\to X'\gets W'$.
Further we have a quasi-isomorphism
$[A\to A']\to 
[A_{V/X,W}\to A_{V'/X',W}]$.
By the assumption that the right square
in (\ref{eqlc}) is cartesian,
we have an isomorphism
$N_{X/P}
\otimes[A\to A']
\to 
[J_B/J_B^2\to
J_{B'}/J_{B'}^2]$
by iii.\ in the step 1 of the proof of
\cite[Proposition 1.7.2]{KSI}.
They induce
a canonical map
\begin{equation}
\lambda\colon
[L\Lambda^pM'_{V/X,W}\to
L\Lambda^pM'_{V'/X',W}]
\to
N_{X/P}
\otimes
[L\Lambda^{p-1}M'_{V/X,W}\to
L\Lambda^{p-1}M'_{V'/X',W}][1].
\label{eqlamb}
\end{equation}
By the step 3 
in the proof of
\cite[Proposition 1.7.2]{KSI},
we see that it 
fits in a commutative
diagram of distinguished
triangles
\begin{equation}
\begin{CD}
L\Lambda^pM'_{V/X,W}
@>{\lambda_{V/X/P,W}}>>
N_{X/P}
\otimes
L\Lambda^{p-1}M'_{V/X,W}[1]\\
@VVV @VVV\\
L\Lambda^pM'_{V'/X',W}
@>{\lambda_{V'/X'/P',W}}>>
N_{X/P}
\otimes
L\Lambda^{p-1}M'_{V'/X',W}[1]\\
@VVV @VVV\\
[L\Lambda^pM'_{V/X,W}\to
L\Lambda^pM'_{V'/X',W}]
@>{\lambda}>>
N_{X/P}
\otimes
[L\Lambda^{p-1}M'_{V/X,W}\to
L\Lambda^{p-1}M'_{V'/X',W}][1].\\
@VVV @VVV
\end{CD}
\label{eqlmcp}
\end{equation}

\addtocounter{thm}{-1}
\begin{lm}\label{lmspmap}
Let $E_{\cal T}$
denote the spectral sequence
$E^1_{p,q}\Rightarrow E_{p+q}$
{\rm (\ref{eqspsqA})}
and let
$E_{\cal T}[-1,3]$
denote the spectral sequence
$E^1_{p+1,q-3}
\Rightarrow E_{p+q-2}$.

Then, there exists a map
\begin{equation}
\begin{CD}
\alpha\colon
E_{\cal T}@>>>
N_{X/P}\otimes
E_{\cal T}[-1,3]
\end{CD}
\label{eqal}
\end{equation}
of spectral sequences
where the maps of $E^1$-terms
are induced by the map
$\lambda$ in 
{\rm (\ref{eqlamb})}.
\end{lm}

{\it Proof.}
This is a variant
of the construction
in \cite[Proposition 1.7.2]{KSI}
where the corresponding map for
the spectral sequence
$E^1_{p,q}
={\cal H}_{2p+q}
L\Lambda^{-p}
M'_{V/X,W}
\Rightarrow
{\cal T}or^{{\cal O}_{X}}_{p+q}
({\cal O}_V,{\cal O}_W)$
is defined.
Similarly as 
the step 1 in the proof
loc.\ cit.,
we obtain a map of spectral sequences.
The assertion on
the map of $E^1$-terms
is clear from the definition
of $\lambda$.
\qed

We prove a relative version
of the excess intersection
formula
\cite[Proposition 3.4.2]{KSI}.

\begin{pr}\label{prrex}
Let $S$ be a regular
noetherian scheme of finite
dimension
and $X,X'$
be locally hypersurfaces
over $S$.
We consider a commutative diagram
\setcounter{equation}0
\begin{equation}
\begin{CD}
T@>\subset>> W\\
@VVV @VVV\\
V'@>\subset>> X'\\
@VVV @VVfV\\
V@>\subset>> X
\end{CD}
\label{eqrex}
\end{equation}
of schemes over $S$
satisfying the condition
{\rm (\ref{ssrex}.0.2)}.
Let $c'$ be the codimension
of the regular immersion $T\to W$.
Let $M'_{V/X,W}$
and $M'_{V'/X',W}$
be the excess
conormal complexes.

We assume
that, locally on $X$,
there exists 
a commutative diagram
\begin{equation}
\begin{CD}
V'@>\subset>> X'@>\subset>>P'\\
@VVV @VVfV @VV{\tilde f}V\\
V@>\subset>> X@>\subset>>P
\end{CD}
\label{eqrexP}
\end{equation}
of schemes over $S$
satisfying the following condition:
\begin{itemize}
\item[{\rm(\ref{prrex}.3)}]
\addtocounter{equation}1
The maps
$X\to P$
and $X'\to P'$
are regular immersion
of codimension $1$
and that 
$V\to P$
and $V'\to P'$
are regular immersion
of codimension $n$.
The right square is cartesian.
The schemes $P$ and $P'$
are smooth of
relative dimension $n$
over $S$.
\end{itemize}

Let $U'\subset X'$
be an open subscheme
such that
$V'\times_{X'}U'
\to
V\times_XU'$
is an open immersion
and put $\Sigma=
X'\setminus U'$.
Let $Z'\subset X'$
be the closed subscheme
defined by the annihilator
of $\Omega_{X'/S}^n$
and we put
${\cal L}_{Z'}=L_1i^{\prime *}L_{X'/S}$
where
$i'\colon Z'\to X'$
is the closed immersion.
Then,
we have the following.

{\rm 1.}
The map $f\colon
X'\to X$
is of finite tor-dimension.

{\rm 2.}
The canonical map
$L\Lambda^{n-c'}M'_{V/X,W}\to
L\Lambda^{n-c'}M'_{V'/X',W}$
is a quasi-isomorphism
on the complement
$T\setminus 
(Z'_T\times_{X'}\Sigma)$.

{\rm 3.}
We define
$f^![V]-[V']
\in G(\Sigma\cap 
(V\times_XX'))$
by
$$f^![V]-[V']
=[{\rm Ker}(
f^*{\cal O}_V\to {\cal O}_{V'})]
+\sum_{q>0}(-1)^q
[L_qf^*{\cal O}_V]$$
and put $d=\dim W$.
Then,
we have
\begin{equation}
(((f^![V]-[V']),W))_{X'}
=
[L\Lambda^{n-c'}(M'_{V/X,W}
\to M'_{V'/X',W})]
\label{eqrexf}
\end{equation}
in $F_{d-n+c'}G(Z'_T
\times_{X'}\Sigma)_{/{\cal L}_{Z'}}$.
\end{pr}

{\it Proof.}
1.
In the diagram
(\ref{eqrexP}),
the schemes $X$ and $P'$
are tor-independent
over $P$.
Since $\tilde f$ is of finite
tor-dimension,
the map $f$ is
also of finite
tor-dimension.

2.
The assertion is local on $X$.
Recall that
the diagram (\ref{eqrexP})
defines 
a quasi-isomorphism
$M'_{V/X,W}
\to
[N_{X/P}
\otimes {\cal O}_T
\to 
{\rm Ker}
(N_{V/P}\otimes {\cal O}_T
\to N_{T/W})]$
\cite[Lemma 1.7.1]{KSI}
and similarly for
$M'_{V'/X',W}$.
Hence the excess conormal
complexes
$M'_{V/X,W}$
and 
$M'_{V'/X',W}$
satisfy the condition
(L($n-c'$)) in Section \ref{sslcc}.

The map
$V'\to V\times_XX'$
defined by the diagram
(\ref{eqrex})
is an open immersion
outside $\Sigma$.
Hence the canonical map
$M'_{V/X,W}\to
M'_{V'/X',W}$
is a quasi-isomorphism
on the complement
$T\setminus T\times_{X'}\Sigma$
by the description
of the excess conormal
complex recalled above.

Let $Z'_1\subset T$
be the closed subscheme
defined by the annihilators
of $\Lambda^{n-c'}
{\mathcal H}_0(M'_{V'/X',W})$.
By the assumption
that the right square 
in (\ref{eqrexP}) is cartesian,
the canonical map
$N_{X/P}\otimes
_{{\mathcal O}_X}{\mathcal O}_{X'}
\to N_{X'/P'}$
is an isomorphism.
Hence the assumption
that the mapping cone
$[M'_{V/X,W}\to
M'_{V'/X',W}]$
is of tor-dimension $\le 1$
is satisfied.
Thus by Lemma \ref{lmkos}, 
the homology sheaves
${\mathcal H}_q([L\Lambda^{n-c'}M'_{V/X,W}
\to
L\Lambda^{n-c'}M'_{V'/X',W})]$
are ${\mathcal O}_{Z'_1}$-module
for $q\ge 0$.

Since the diagram (\ref{eqrexP})
defines an isomorphism
$[N_{X'/P'}\to 
\Omega^1_{P'/S}
\otimes
_{{\mathcal O}_{P'}}{\mathcal O}_{X'}]
\to L_{X'/S}$,
we have an inclusion
$Z'_1\subset Z'_T$
of closed subschemes of $T$.
Hence the assertion follows.

3.
The proof is similar to
\cite[Proposition 3.4.2]{KSI}
using Lemmas \ref{lmspsqA}
and \ref{lmspmap}.

By 2., the right hand
side of (\ref{eqrexf})
is defined as an element of
$F_{d-n+c'}G(Z'_{T}
\times_{X'}\Sigma)$.
Since $X$ and $P'$
are tor-independent
over $P$,
we see that the canonical map
$L\tilde f^*{\cal O}_V
\to Lf^*{\cal O}_V$
is a quasi-isomorphism.
Hence,
the complex ${\cal C}$
of ${\cal O}_{X'}$-modules
in Lemma \ref{lmspsqA}
is acyclic outside $\Sigma$.
We consider the spectral sequence
(\ref{eqspsqA})
and will show the equality
\begin{equation}
\sum_{p+q=r,r+1}
(-1)^{p+q}
[E^1_{p,q}]=
(-1)^r[E_r]
+(-1)^{r+1}[E_{r+1}]
\label{eqEr}
\end{equation}
for sufficiently large $r$.
The right hand side
of (\ref{eqEr})
is independent of
sufficiently large $r$
by \cite[Theorem 3.2.1]{KSI}
and defines 
the left hand side
of (\ref{eqrexf}).

We show that the left hand side
of (\ref{eqEr})
equals the right hand side
of (\ref{eqrexf}).
By \cite[Lemma 3.4.1]{KSI},
the map $\lambda_{V/X/P,W}$
induces an isomorphism
$L_{p+1}\Lambda^{q+1}
M'_{V/X,W}
\to
N_{X/P}
\otimes
L_p\Lambda^q
M'_{V/X,W}$
for $q\ge n-c'$
and similarly for
$\lambda_{V'/X'/P',W}$.
By the commutative
diagram (\ref{eqlmcp}),
it induces an isomorphism
$E^1_{p,q}
\to
{\mathcal L}_{Z'}
\otimes
E^1_{p+1,q-3}$
if $-(p+1)\ge n-c'$.
Hence,
the left hand side
of (\ref{eqEr}) is
equal to
$\sum_{q}
(-1)^{-(n-c')+q}
[E^1_{-(n-c'),q}]$
and to
the right hand side of
(\ref{eqrexf}).

We show
the equality (\ref{eqEr})
applying \cite[Lemma 3.3.1]{KSI}.
Note that, in 
\cite[Lemma 3.3.1]{KSI},
the map
$\alpha_r$
is used only to show
that
the right hand side
is independent of
$r\ge r_0$
and that
we can drop the compatibility
assumption
with the restriction
$\alpha_r|_{T\cap U}$
since it is not used
in the proof.
We apply
\cite[Lemma 3.3.1]{KSI},
to the maps
(\ref{eqlamb}) and
\cite[(3.1.3.1)]{KSI}
and the map (\ref{eqal}).
Then, by Lemma \ref{lmspmap},
the maps of 
$E^1$-terms
of (\ref{eqal})
are isomorphisms 
for sufficiently large $p+q$
and
the assumption
of \cite[Lemma 3.3.1]{KSI}
is satisfied.
Thus,
we obtain
the equality (\ref{eqEr}).
\qed

\newpage 
\section{Intersection
product with the
log diagonal}\label{siw}

From this section on,
we fix a complete discrete
valuation field $K$
with perfect residue
field $F$ of characteristic $p>0$.
We put $S={\rm Spec}\ 
{\cal O}_K$
and $s={\rm Spec}\ F$.
Both $0$ and $p$
are allowed as the characteristic of $K$.
A morphism of schemes over $S$
is always 
a morphism over $S$.

We introduce in Section \ref{sslipd}
the localized intersection
product with the log diagonal
by applying Definition \ref{dflip1}.
We establish an important 
property that 
the localized intersection
product with the log diagonal
is independent of the boundary
in Proposition \ref{prprod}.
In preliminary subsections
\ref{sslpX}
and \ref{sslctcx},
we study local structure
of log products
and the logarithmic cotangent complex
respectively.

\subsection{Log products}
\label{sslpX}

In this subsection,
we study the log self-product
of 
a regular flat scheme
$X$ of finite type
over $S={\rm Spec}\ {\cal O}_K$
with respect to a divisor
$D\subset X$
with simple normal crossings.
We do not assume
inclusion
between $D$
and the closed fiber $X_F$.
The case where
$D=X_F$ and $X_K$ is smooth
over $K$ is treated in
\cite[Sections 5.1, 5.2]{KSI}.
First, 
we study the local structure
of $X$.

\begin{lm}\label{lmXP}
Let $X$ be a regular
flat scheme of finite type
over $S$
and $D$ be a divisor
with simple normal
crossings.
Then, for every point
$x$ of $X$,
there exist an open 
neighborhood $U$ of $x$,
a smooth scheme
$P$ over $S$,
a divisor 
$\widetilde D$ of $P$ with 
simple normal crossings
relatively to $S$ and a
regular immersion $U\to P$
of codimension $1$
such that
$D\cap U=\widetilde D
\times_PU$.
\end{lm}

{\it Proof.}
It suffices
to prove the
cases where
$x$ is a closed point
of $X_F$
and of $X_K$
respectively.
First, we show the case
where
$x$ is a closed point
of $X_F$.
Let $D_1,\ldots,D_m$
be the irreducible components
of $D$ containing $x$
and take $t_1,\ldots,t_m
\in {\mathfrak m}_x$
defining $D_1,\ldots,D_m$
on a neighborhood of $x$.
We extend it
to a minimal system
$t_1,\ldots,t_n
\in {\mathfrak m}_x$
of generators.
Then, 
the map
$U\to {\mathbf A}^n_S$
defined by $t_1,\ldots,t_n$
on an open neighborhood $U$ of $x$
is unramified.
Hence, after shrinking $U$
if necessary,
there is an \'etale scheme
$P\to {\mathbf A}^n_S$
and a regular closed immersion
$U\to P$ of codimension 1
such that
$D\cap U$ is the sum
of the pull-back
of the first $m$
coordinate hyperplanes
by \cite[Corollaire (18.4.7)]{EGA4}.

Next, we show the case
where
$x$ is a closed point
of $X_K$.
We take a minimal system
$t_1,\ldots,t_n
\in {\mathfrak m}_x$
of generators as above.
There exists an
element $t_0\in 
{\mathcal O}_{X,x}$
such that
the residue field
$\kappa(x)$ is
a finite separable
extension of $K(\bar t_0)$
and that
$K(\bar t_0)$ 
is purely inseparable
over $K$ by
Lemma \ref{lmKp} below.
Then, 
the map
$U\to {\mathbf A}^{n+1}_S$
defined by $t_0,t_1,\ldots,t_n$
on an open neighborhood 
$U$ of $x$
is unramified.
Thus, we conclude
similarly as above.
\qed

\begin{lm}
\label{lmKp}
Let $k$ be a field of
characteristic $p>0$
such that $[k: k^p]=p$.
Then, for a finite extension
$L$ of $k$,
there exists an integer $e\ge 0$
such that $L$ is a separable
extension of
$k^{1/p^e}$.
\end{lm}

{\it Proof.}
Let $L_1$ be the separable closure of
$k$ in $L$ and put $[L:L_1]=p^e$.
Then since $[L_1:L_1^p]=[k:k^p]=p$,
it follows that
$L$ is a unique purely
inseparable extension 
$L_1^{1/p^e}$ of $L$
of degree $p^e$
and is a separable
extension of
$k^{1/p^e}$.
\qed

We show a relative version of Lemma \ref{lmXP}.

\begin{lm}\label{lmYQ}
\setcounter{equation}0
Let $X$ and 
$Y$ be regular
flat schemes
of finite type over $S$ 
and $f\colon
Y\to X$ be a morphism
over $S$.
Let $D\subset X$
and $E\subset Y$
be divisors
with simple normal
crossings
such that
$Y\setminus E\subset
f^{-1}(X\setminus D)$.
Let $y$ be a point of $Y$
and we put $x=f(y)\in X$.

Then, there exists 
open neighborhoods
$U$ and $V$ of
$x$ and $y$ respectively
satisfying $f(V)\subset U$
and a cartesian diagram
\setcounter{equation}0
\begin{equation}
\begin{CD}
V@>>> Q&\ \supset \widetilde E\\
@V{f|_V}VV @VV{\tilde f}V\\
U@>>> P&\ \supset \widetilde D
\end{CD}
\label{eqYQ}
\end{equation}
of schemes
over $S$ satisfying
the following
conditions:
\begin{itemize}
\item
The schemes $P$ and $Q$
are smooth over $S$ 
and $\widetilde D=
\sum_i\widetilde D_i\subset P$
and $\widetilde E=
\sum_j\widetilde E_j\subset Q$
are divisors with relative simple 
normal crossings
relatively to $S$
respectively.
For each $i$,
we have
$\tilde f^{-1}
(\widetilde D_i)
=\sum_je_{ij}
\widetilde E_j$
for some integers $e_{ij}\ge 0$.
The horizontal arrows
are regular immersions
of codimension $1$
and $D\cap U=
\widetilde D\times_PU$
and
$E\cap V=
\widetilde E\times_QV$
as in Lemma {\rm \ref{lmXP}}.
\end{itemize}
\end{lm}

{\it Proof.}
It suffices to prove
the cases where
$y$ is a closed point of $Y_F$
and of $Y_K$
respectively.
First we prove the
case where
$y$ is a closed point of $Y_F$
and hence
$x$ is a closed point of $X_F$.
By the proof of
Lemma \ref{lmXP},
we obtain a diagram
(\ref{eqYQ}) without $\tilde f$
together with
\'etale morphisms
$P\to {\mathbf A}^n_S
={\rm Spec}\
{\cal O}_K[T_1,\ldots,T_n]$
and
$Q\to {\mathbf A}^{n'}_S
={\rm Spec}\
{\cal O}_K[S_1,\ldots,S_{n'}]$
such that
$\widetilde D$ and 
$\widetilde E$
are the pull-back
of the first
$m$ and $m'$ coordinate hyperplanes
respectively
and that
the maximal ideals
${\mathfrak m}_x
\subset {\cal O}_{X,x}$
and
${\mathfrak m}_y
\subset {\cal O}_{Y,y}$
are generated by
$t_i=T_i|_U$
for $i=1,\ldots,n$
and
$s_j=S_j|_V$
for $j=1,\ldots,n'$
respectively.

For $i=1,\ldots,m$,
we put
$f^*t_i=v_i\prod_js_j^{e_{ij}}$
for some units
$v_i$ on $V$.
For $i=m+1,\ldots,n$,
we also put
$f^*t_i=\sum_ja_{ij}s_j$
for some functions
$a_{ij}$ on $V$.
After shrinking $Q$,
we take units 
$\tilde v_i$ on $Q$
lifting $v_i$
and functions
$\tilde a_{ij}$ on $Q$
lifting $a_{ij}$
and define a map
$g\colon Q\to 
{\mathbf A}^n_S$ by
sending $T_i$ to
$\tilde v_i\prod_j
S_j^{e_{ij}}$
for $i=1,\ldots,m$
and to
$\sum_j\tilde a_{ij}S_j$
for $i=m+1,\ldots,n$.
By replacing $Q$ by 
an \'etale neighborhood
$Q\times_{{\mathbf A}^n_S}P$,
we obtain a
map $\tilde f\colon
Q\to P$
that makes
(\ref{eqYQ}) 
a commutative diagram.

The equalities
$\tilde f^{-1}
(\widetilde D_i)
=\sum_je_{ij}
\widetilde E_j$
follow from 
the definition of $g$.
We show that the diagram
(\ref{eqYQ})
thus obtained is cartesian
on a neighborhood of $y$.
Let
$\tilde {\mathfrak m}_x
\subset {\cal O}_{P,x}$
and
$\tilde {\mathfrak m}_y
\subset {\cal O}_{Q,y}$
be the maximal ideals.
Since the horizontal
arrows are regular immersions
of codimension 1,
it suffices to show
that the canonical map
$\tilde f^*\colon
\tilde {\mathfrak m}_x
/\tilde {\mathfrak m}_x^2
\to
\tilde {\mathfrak m}_y
/\tilde {\mathfrak m}_y^2$
induces an isomorphism
${\rm Ker}(
\tilde {\mathfrak m}_x
/\tilde {\mathfrak m}_x^2
\to
{\mathfrak m}_x
/{\mathfrak m}_x^2)
\to
{\rm Ker}(
\tilde {\mathfrak m}_y
/\tilde {\mathfrak m}_y^2
\to
{\mathfrak m}_y
/{\mathfrak m}_y^2)$
on the subspaces
of dimension 1.
Since $
{\mathfrak m}_x
/{\mathfrak m}_x^2
=
\langle t_1,\ldots,t_n\rangle$
and
${\mathfrak m}_y
/{\mathfrak m}_y^2
=\langle s_1,\ldots,s_{n'}
\rangle$,
we have
$\tilde {\mathfrak m}_x
/\tilde {\mathfrak m}_x^2
=
{\rm Ker}(
\tilde {\mathfrak m}_x
/\tilde {\mathfrak m}_x^2
\to
{\mathfrak m}_x
/{\mathfrak m}_x^2)
\oplus
\langle T_1,\ldots,T_n\rangle$
and
$\tilde {\mathfrak m}_y
/\tilde {\mathfrak m}_y^2
=
{\rm Ker}(
\tilde {\mathfrak m}_y
/\tilde {\mathfrak m}_y^2
\to
{\mathfrak m}_y
/{\mathfrak m}_y^2)
\oplus
\langle S_1,\ldots,S_{n'}
\rangle$.
By the definition of $g$,
the map $\tilde f^*$
sends the subspace
$\langle T_1,\ldots,T_n\rangle
\subset
\tilde {\mathfrak m}_x
/\tilde {\mathfrak m}_x^2$
into
$\langle S_1,\ldots,S_{n'}
\rangle
\subset
\tilde {\mathfrak m}_y
/\tilde {\mathfrak m}_y^2$.
Since
$\pi,T_1,\ldots,T_n$
and
$\pi,S_1,\ldots,S_{n'}$
are bases of
$\tilde {\mathfrak m}_x
/\tilde {\mathfrak m}_x^2
=
{\rm Ker}(
\tilde {\mathfrak m}_x
/\tilde {\mathfrak m}_x^2
\to
{\mathfrak m}_x
/{\mathfrak m}_x^2)
\oplus
\langle T_1,\ldots,T_n\rangle$
and of
$\tilde {\mathfrak m}_y
/\tilde {\mathfrak m}_y^2$
respectively,
the image 
of ${\rm Ker}(
\tilde {\mathfrak m}_x
/\tilde {\mathfrak m}_x^2
\to
{\mathfrak m}_x
/{\mathfrak m}_x^2)$
is not in 
$\langle S_1,\ldots,S_{n'}
\rangle$.
Hence the map
${\rm Ker}(
\tilde {\mathfrak m}_x
/\tilde {\mathfrak m}_x^2
\to
{\mathfrak m}_x
/{\mathfrak m}_x^2)
\to
{\rm Ker}(
\tilde {\mathfrak m}_y
/\tilde {\mathfrak m}_y^2
\to
{\mathfrak m}_y
/{\mathfrak m}_y^2)$
induced by
$\tilde f^*$
is an isomorphism 
as required.

Next,
we prove the case
where $y$ is a closed point
of $Y_K$ and hence
$x$ is a closed point of $X_K$.
Let $\pi$ be a
prime element of
${\mathcal O}_K$.
We may take
$t_0\in {\cal O}_{X,x}$
and
$s_0\in {\cal O}_{Y,y}$
such that
$t_0\equiv \pi
\bmod {\mathfrak m}_x$
and
$s_0\equiv \pi
\bmod {\mathfrak m}_y$
if $K$ is of characterstic $0$
and
$t_0^{p^a}\equiv \pi
\bmod {\mathfrak m}_x$
and
$s_0^{p^{a+b}}\equiv \pi
\bmod {\mathfrak m}_y$
for some integers $a\ge 0,b\ge0$
and 
$\kappa(x)$ and
$\kappa(y)$ are
separable extentions
over $K(\bar t_0)$
and over $K(\bar s_0)$
respectively
if $K$ is of characterstic $p>0$
by Lemma \ref{lmKp}.
Then, we obtain a diagram 
(\ref{eqYQ}) without $\tilde f$
together with
\'etale morphisms
$P\to {\mathbf A}^{n+1}_S
={\rm Spec}\
{\cal O}_K[T_0,T_1,\ldots,T_n]$
and
$Q\to {\mathbf A}^{n'+1}_S
={\rm Spec}\
{\cal O}_K[S_0,S_1,\ldots,S_{n'}]$
similarly as above.
We put
$f^*t_0=s_0^{p^b}
+\sum_{j>0}a_{0j}s_j$
and take liftings
$\tilde a_{0j}$
as above.
By the same procedure
for $i>0$ as above,
we define a map
$\tilde f\colon
Q\to P$
that makes
(\ref{eqYQ}) 
a commutative diagram
satisfying
the equalities
$\tilde f^{-1}
(\widetilde D_i)
=\sum_je_{ij}
\widetilde E_j$.

We show that the diagram
(\ref{eqYQ})
thus obtained is cartesian
on a neighborhood of $y$.
By the definition of $\tilde f$,
the map $\tilde f^*$
sends the subspace
$\langle T_1,\ldots,T_n\rangle
\subset
\tilde {\mathfrak m}_x
/\tilde {\mathfrak m}_x^2$
into
$\langle S_1,\ldots,S_{n'}
\rangle
\subset
\tilde {\mathfrak m}_y
/\tilde {\mathfrak m}_y^2$.
Since
$T_0^{p^a}-\pi,T_1,\ldots,T_n$
and
$S_0^{p^{a+b}}-\pi,S_1,\ldots,S_{n'}$
are bases of
$\tilde {\mathfrak m}_x
/\tilde {\mathfrak m}_x^2
=
{\rm Ker}(
\tilde {\mathfrak m}_x
/\tilde {\mathfrak m}_x^2
\to
{\mathfrak m}_x
/{\mathfrak m}_x^2)
\oplus
\langle T_1,\ldots,T_n\rangle$
and of
$\tilde {\mathfrak m}_y
/\tilde {\mathfrak m}_y^2$
respectively,
the image 
of ${\rm Ker}(
\tilde {\mathfrak m}_x
/\tilde {\mathfrak m}_x^2
\to
{\mathfrak m}_x
/{\mathfrak m}_x^2)$
is not in 
$\langle S_1,\ldots,S_{n'}
\rangle$.
This implies the required
assertion as above.
\qed

Let $X$ be a regular flat
separated scheme 
of finite type over $S$
and $D\subset X$
be a divisor with simple normal
crossings.
We consider the log product
$(X\times_SX)^\sim$
defined as 
$(X\times_SX)^\sim_{\cal D}$
with respect to 
the family ${\cal D}
=(D_i)_{i\in I}$
of irreducible
components of $D$.

\begin{lm}\label{lmX}
\setcounter{equation}0
Let $X$ be a regular flat
separated scheme 
of finite type over $S$
and $D\subset X$
be a divisor with simple normal
crossings.

Then, the log product
$(X\times_SX)^\sim$
is locally a hypersurface
(Definition {\rm \ref{dflhs}})
over $X$ 
with respect
to either of the projections.
\end{lm}

\noindent{\it Proof.}
Let $x\in X$
be a point,
$U\subset X$ be an 
open neighborhood of $x$
and 
$U\to P$
be a regular immersion
of codimension 1
to a smooth scheme
$P$ over $S$
satisfying the condition
in Lemma \ref{lmXP}.
We define the log product
$(P\times_SX)^\sim$ similarly.
Since the second projection
$(P\times_SX)^\sim
\to X$
is log smooth and strict 
by Lemma \ref{lmlp},
it is smooth.
Hence
the log product
$(P\times_SX)^\sim$
is regular.

By the universality
of log product,
we have a cartesian 
diagram
\begin{equation}
\begin{CD}
(U\times_SX)^\sim
@>{{\rm pr}_1}>>U\\
@VVV @VVV\\
(P\times_SX)^\sim
@>{{\rm pr}_1}>>P.
\end{CD}
\label{eqPUX}
\end{equation}
Hence the ideal defining 
the immersion
$(U\times_SX)^\sim
\to
(P\times_SX)^\sim$
is locally monogenic.
Thus, to conclude
that
$(U\times_SX)^\sim$
is a divisor of
$(P\times_SX)^\sim$,
it is sufficient to show
that the immersion
$(U\times_SX)^\sim
\to
(P\times_SX)^\sim$
is nowhere dominant.
Thus, it 
is reduced to the
case where $D$ is empty.
In this case,
$U\times_SX$
is a divisor of
$P\times_SX$
since $X$
is flat over $S$.
\qed

\setcounter{equation}0
\begin{cor}\label{cortdf}
Let $X$ and $Y$
be regular flat separated schemes
of finite type
and $D\subset X$
and $E\subset Y$
be divisors with
simple normal crossings.
Let $f\colon Y\to X$
be a morphism over $S$
satisfying
$f^{-1}(D)\subset E$
set-theoretically.

Then, the map
$(f\times f)^\sim\colon
(Y\times_S Y)^\sim
\to
(X\times_S X)^\sim$
is locally of 
complete intersection
and hence is
of finite tor-dimension.
\end{cor}

{\it Proof.}
Since the assertion 
is local,
we may take a cartesian diagram
(\ref{eqYQ}).
By the cartesian diagram
(\ref{eqPUX}) and
the corresponding one for $Y$,
the diagram
\begin{equation}
\begin{CD}
(V\times_S Y)^\sim
@>{(f_V\times f)^\sim}>>
(U\times_S X)^\sim\\
@VVV @VVV\\
(Q\times_S Y)^\sim
@>{(\tilde f\times f)^\sim}>>
(P\times_S X)^\sim
\end{CD}
\label{eqXYPQ}
\end{equation}
is cartesian.
Since 
$(P\times_S X)^\sim$
is smooth over $X$
and
$(Q\times_S Y)^\sim$
is smooth over $Y$,
they are regular.
Hence the bottom
horizontal arrow
$(\tilde f\times f)^\sim$
is locally
of complete intersection.
Since the vertical arrows
are regular immersion
of codimension 1
and the diagram is
cartesian,
$(U\times_S X)^\sim$
and
$(Q\times_S Y)^\sim$
are tor-independent
over
$(P\times_S X)^\sim$.
Hence
the top arrow
$(f_V\times f)^\sim$
is locally 
of complete intersection
of the same virtual
relative dimension.
\qed

We study the local structure
of the log product
$(X\times_SX)^\sim$
inductively on
the number of irreducible
components of $D$.
Let $X$, $D$ 
and the log product
$(X\times_SX)^\sim$
be as in the beginning
of this subsection.
Let $X_1$ be
a regular divisor
of $X$ such that
$D_1=D\cap X_1$ is
a divisor of $X_1$
with simple normal
crossings.
Let ${\cal D}=
(D_i)_{i\in I}$
be the family of
irreducible components
of $D$
and we consider
the family
${\cal D}_1=
(D_i\cap X_1)_{i\in I}$
of smooth divisors of $X_1$.
Then,
the log product
$(X_1\times_SX_1)^\sim$
with respect to
${\cal D}_1$
is identified
with the inverse image
of $X_1\times_SX_1$
by the canonical map
$(X\times_SX)^\sim
\to X\times_SX$.

The sum
$D'=D\cup X_1$
is a divisor of $X$
with simple normal crossings.
We consider the log product
$(X\times_SX)^\approx$
with respect to $D'$.
By the inductive
construction of
the log product,
we have a canonical
isomorphism
$(X\times_SX)^\approx
\to
(X\times_SX)^\sim
\times
_{X\times_SX}
(X\times_SX)^\sim_{X_1}$.
The inverse image $E$
of $(X_1\times_SX_1)^\sim$
by the canonical map
$(X\times_SX)^\approx
\to
(X\times_SX)^\sim$
is a ${\mathbf G}_m$-torsor
over
$(X_1\times_SX_1)^\sim$
by Lemma \ref{lmEi}.
They are summarized in
the cartesian diagram
\begin{equation}
\begin{CD}
E@>>> 
(X_1\times_SX_1)^\sim
@>>>
X_1\times_SX_1\\
@VVV @VVV @VVV\\
(X\times_SX)^\approx
@>>> 
(X\times_SX)^\sim
@>>>
X\times_SX
\end{CD}
\end{equation}
where the vertical
arrows are
closed immersions.
The morphism
$(X\times_SX)^\approx
\to (X\times_SX)^\sim$
is an isomorphism
on the complements of
$E$ and $(X_1\times_SX_1)^\sim$.
The subscheme
$E\subset 
(X\times_SX)^\approx$
is the inverse
image of
$X_1\subset X$
by the composition
of the canonical map
$(X\times_SX)^\approx
\to X\times_SX$
with either of
the projections
$X\times_SX\to X$.

To understand 
the local structure
of the log product,
it suffices to study
it on a neighborhood
of $E$ by the inductive
construction 
of the log product.

\begin{lm}\label{lmX2}
\setcounter{equation}0
Let the notations
$X,D,X_1,E$ etc.\
be as above.

{\rm 1.}
Assume $X_1$ is flat over $S$.
Then 
the immersion
$(X_1\times_SX_1)^\sim
\to
(X\times_SX)^\sim$
is a regular immersion
of codimension $2$
and 
$E$ is a Cartier
divisor of $(X\times_SX)^\approx$.

{\rm 2.}
Assume $X_1$ is 
a subscheme of
the closed fiber $X_F$
and put $d=\dim X_F$.
Then the scheme
$(X_1\times_FX_1)^\sim$ 
is smooth of 
dimension $2d$
over $F$.
\end{lm}

\noindent{\it Proof.}
1.
The immersion
$(X_1\times_SX_1)^\sim
\to
(X\times_SX)^\sim$
is locally of complete
intersection
by Corollary \ref{cortdf}.
Hence, the assertion
follows from Lemma \ref{lmciri}.

2.
Since the projections
$(X_1\times_F X_1)^\sim
\to X_1$
are smooth
of relative dimension $d$,
the assertion follows.
\qed

We show some tor-independences
(Definition \ref{dftor}.1).
Its consequences
Corollary \ref{cortorE}.1
and \ref{cortorE}.2
will be used in the proof of
Propositions \ref{prU1} and
\ref{prU2} respectively.

\begin{lm}\label{lmtorE}
\setcounter{equation}0
Let $X$ and $Y$
be regular flat separated
schemes over $S$
of finite type
and $f\colon 
Y\to X$ be a morphism over $S$.
Let $D\subset X$
be a regular divisor
such that
$D_Y=D\times_XY$ is a divisor of $Y$
and
let $D'$ be a divisor 
of $Y$ with simple normal crossings.
We assume that
either both $D$ and
$D_Y$ are 
flat over $S$ or they
are schemes over $F$.

{\rm 1.}
The fiber products
$D\times_SD$
and 
$Y\times_SY$
are tor-independent over
$X\times_SX$.

{\rm 2.}
Let $(X\times_SX)^\sim$
and $(Y\times_SY)^\sim$
be the log product
with respect to $D$
and $D'$ respectively.
Assume that
$D_Y=D\times_XY$
is a subset of $D'$
set-theoretically.
Further assume that
either $D$ and $D'$ are
flat over $S$ 
or $D$ and $D'$ are schemes over $F$.

Let $E\subset
(X\times_SX)^\sim$
be the pull-back
of
$D\subset X$
by either of the two projections.
Then $E$ and 
$(Y\times_SY)^\sim$
are tor-independent over
$(X\times_SX)^\sim$.
\end{lm}

{\it Proof.}
1.
By Lemma \ref{lmtidp},
it suffices to show
that
$D\times_SD$ and
$X\times_SY$ 
are tor-independent
over
$X\times_SX$ 
and that
$D\times_SD_Y$ and
$Y\times_SY$ 
are tor-independent
over
$X\times_SY$.

By the assumption that
$D_Y=D\times_XY$ 
is a divisor of $Y$,
it follows that
$D$ and $Y$ are tor-independent
over $X$.
Either if $D$ is flat over $S$
or if $D$ is a scheme over $F$,
the fiber product
$D\times_SD$ is flat over $D$.
Hence
$D\times_SD$
and $Y$ are tor-independent
over $X$ with respect to
the second projection
$D\times_SD\to X$ by
Lemma \ref{lmtidp}.
Thus, by applying
Lemma \ref{lmtidp}
to 
$Y\to X\gets X\times_SX
\gets D\times_SD$,
we conclude that
that
$D\times_SD$ and
$X\times_SY$ 
are tor-independent
over
$X\times_SX$. 

Similarly, either if 
$D_Y$ is flat over $S$
or if $D$ and $D_Y$
are both schemes over $F$,
the fiber product
$D\times_SD_Y$ is flat over $D$.
Hence
$D\times_SD_Y$
and $Y$ are tor-independent
over $X$ with respect to
the first projection
$D\times_SD_Y\to X$ by
Lemma \ref{lmtidp}.
Thus, by applying
Lemma \ref{lmtidp}
to 
$Y\to X\gets X\times_SY
\gets D\times_SD_Y$,
we conclude that
that
$D\times_SD_Y$ and
$Y\times_SY$ 
are tor-independent
over
$X\times_SY$
as required. 

2.
First, we show the case
where $D$ and $D'$
are flat over $S$.
By Lemma \ref{lmX2}.1,
$E$ is a divisor
of $(X\times_SX)^\sim$.
Since every $D'_i$ is flat over $S$
in this case,
the pull-back 
$E'_i\subset (Y\times_SY)^\sim$
of an irreducible component
$D'_i$ of $D'$
by either of the two projections
$(Y\times_SY)^\sim\to Y$
is also a divisor.
Hence, 
if $D_Y=\sum_i e_iD'_i$,
the pull-back
$(f\times f)^{\sim*}E
=\sum_ie_iE'_i$ is
also a divisor
and the assertion follows
in this case.

Assume $D=D_F$.
Since the assertion is local,
we may take a cartesian diagram
(\ref{eqYQ}) satisfying
the condition in
Lemma \ref{lmYQ}
and we consider the diagram
(\ref{eqXYPQ}).
By applying
Lemma \ref{lmtidp}
to
$U\to P
\gets 
(X\times_S P)^\sim
\gets
(Y\times_S Q)^\sim$,
we conclude that
the schemes
$(X\times_S U)^\sim$
and
$(Y\times_S Q)^\sim$
are tor-independent
over
$(X\times_S P)^\sim$.
Since
$E\cap (X\times_S U)^\sim$
is a divisor
of $(X\times_S P)^\sim$
and
$(f\times f)^{\sim*}E\cap (Y\times_S V)^\sim$
is a divisor
of $(Y\times_S Q)^\sim$,
the schemes
$E\cap (X\times_S U)^\sim$
and
$(Y\times_S Q)^\sim$
are tor-independent
over
$(X\times_S P)^\sim$.
Then, applying Lemma \ref{lmtidp}
to
$(Y\times_S Q)^\sim
\to
(X\times_S P)^\sim
\gets
(X\times_S U)^\sim
\gets
E\cap (X\times_S U)^\sim$,
we conclude that
$E\cap (X\times_S U)^\sim$
and
$(Y\times_S V)^\sim$
are tor-independent over
$(X\times_S U)^\sim$.
\qed

\begin{cor}\label{cortorE}
Let the notation be
as in Lemma {\rm \ref{lmtorE}}
and we put
$D_Y=D\times_XY
=\sum_ie_iD'_i$ and
$c=\dim Y_K-\dim X_K$.

{\rm 1.}
Let $f_i\colon D'_i\to D$
be the restriction of $f
\colon Y\to X$
and $(f_i\times f_j)^*
\colon
G(D\times_SD)\to G(D'_i\times_SD'_j)$
be the pull-back.
Then, the map $(f\times f)^*\colon
{\rm Gr}^F_{\bullet}
G(D\times_SD)
\to
{\rm Gr}^F_{\bullet+2c}
G(D_Y\times_SD_Y)$
defined by 
$f\times f\colon
Y\times_SY\to X\times_SX$
is the composition of
$$\begin{CD}
{\rm Gr}^F_{\bullet}
G(D\times_SD)
@>{((f_i\times f_j)^*)_{i,j}}>>
\bigoplus_{i,j}
{\rm Gr}^F_{\bullet+2c}
G(D'_i\times_SD'_j)
@>{\sum_{i,j}e_i\cdot e_j\cdot
{\rm can}}>>
{\rm Gr}^F_{\bullet+2c}
G(D_Y\times_SD_Y)
\end{CD}.$$

{\rm 2.}
We put 
$E'=E\times_{(X\times_SX)^\sim}
(Y\times_SY)^\sim$ and
let $E'_i\subset (Y\times_SY)^\sim$
be the pull-back of
$D'_i\subset Y$
by either of the projections.
Then, 
the restriction
$g_i\colon E'_i\to E$
of $(f\times f)^\sim\colon
(Y\times_SY)^\sim\to
(X\times_SX)^\sim$
is of finite tor-dimension.
Further, 
the map $(f\times f)^{\sim *}
\colon
{\rm Gr}^F_{\bullet}
G(E)\to 
{\rm Gr}^F_{\bullet+2c}
G(E')$
defined by 
$(f\times f)^\sim\colon
(Y\times_SY)^\sim\to
(X\times_SX)^\sim$
is the composition of
$$\begin{CD}
{\rm Gr}^F_{\bullet}
G(E)
@>{(g_i^*)_i}>>
\bigoplus_i
{\rm Gr}^F_{\bullet+2c}
G(E'_i)
@>{\sum_ie_i\cdot
{\rm can}}>>
{\rm Gr}^F_{\bullet+2c}
G(E')
\end{CD}.$$
\end{cor}

{\it Proof.}
1.
By Corollary \ref{cortdf}
applied to
$(\amalg_i D'_i)\times_S
(\amalg_i D'_i)
\to D\times_SD$
with the trivial divisor,
the map $f_i\times f_j
\colon
D'_i\times_SD'_j
\to D\times_SD$
is of finite tor-dimension.
Hence
the map $(f_i\times f_j)^*
\colon
G(D\times_SD)\to G(D'_i\times_SD'_j)$
is defined.
Let $f_D\colon
D_Y\to D$ be the base change of $f$.
Then, by Lemma \ref{lmtorE}.1,
the map
$(f\times f)^*\colon
G(D\times_SD)
\to
G(D_Y\times_SD_Y)$
is the same as the pull-back by
$f_D\times f_D\colon
D_Y\times_SD_Y
\to D\times_SD$.
By the assumption,
either every $D'_i$ is flat
over $S$ or is a scheme over $F$.
Hence, there exists a filtration
on ${\cal O}_{D_Y\times_SD_Y}$
such that graded pieces are
invertible 
${\cal O}_{D'_i\times_SD'_j}$-modules
with multiplicities $e_ie_j$.
Thus the assertion follows.

2.
Both $E'_i$
and $E\times_{D\times_SD}
(D'_i\times_SD'_i)^\sim$
are ${\mathbf G}_m$-torsors
over
$(D'_i\times_SD'_i)^\sim$.
The map
$E'_i\to E\times_{D\times_SD}
(D'_i\times_SD'_i)^\sim$
 induced by
$g_i\colon E'_i\to E$
is compatible with
the $e_i$-th power map
of ${\mathbf G}_m$
and is flat.
Hence, the map
$g_i$ is of finite tor-dimension
by Corollary \ref{cortdf}.
By Lemma \ref{lmtorE}.2,
the map
$(f\times f)^{\sim*}\colon
G(E)\to G(E')$
is the same as the pull-back by
the restriction $g_E\colon E'
\to E$.
Since there exists a filtration
on ${\cal O}_{E'}$
such that graded pieces are
invertible
${\cal O}_{E'_i}$-modules
with multiplicities $e_i$,
the assertion follows.
\qed

\subsection{Logarithmic cotangent complex}\label{sslctcx}

We define a logarithmic
version of the cotangent complex.

\begin{df}\label{dfLXSl}
Let $X$ be a regular flat
separated scheme of finite
type over $S={\rm Spec}\
{\mathcal O}_K$ and
$D$ be a divisor of $X$
with simple normal crossings.
We put $n=\dim X_K+1$.
Let $(X\times_SX)^\sim$
denote the log product with
respect to the family $(D_i)$ of
Cartier divisors
consisting of the irreducible
components of the divisor $D$.
We regard $X$ as
a closed subscheme of 
$(X\times_SX)^\sim$ by
the log diagonal map
$\delta\colon
X\to (X\times_SX)^\sim$.

Define
{\rm the logarithmic 
cotangent complex}
$L_{X/S}(\log D)$
to be the conormal complex
$M_{X/(X\times_SX)^\sim}
=L_{X/(X\times_SX)^\sim}[1]$
and
a coherent ${\cal O}_X$-module
$\Omega^1_{X/S}(\log D)$
to be the conormal sheaf
$N_{X/(X\times_SX)^\sim}
={\mathcal H}_0
(L_{X/S}(\log D))$.
Define a 
closed subscheme
$\Sigma_{X/S}$ of $X$
to be that
defined by the annihilator
of the $n$-th exterior power
$\Omega^n_{X/S}(\log D)
=\Lambda^n
\Omega^1_{X/S}(\log D)$.
\end{df}

If the characteristic of $K$
is 0,
the coherent sheaf
$\Omega^1_{X/S}(\log D)$
is locally free of rank $n-1$
on the generic fiber
and hence
$\Sigma_{X/S}$
is supported on the
closed fiber set-theoretically.

Since the log diagonal map
$\delta\colon
X\to (X\times_SX)^\sim$
is a section of
the projection
$(X\times_SX)^\sim\to X$,
the pull-back
$L\delta^*L_{(X\times_SX)^\sim/X}$
of the cotangent complex
\cite{ctcx}
is canonically identified
with the conormal complex
$M_{X/(X\times_SX)^\sim}$
\cite[Definition 1.6.3.1]{KSI}
and hence with the
logarithmic cotangent complex
$L_{X/S}(\log D)$.

\setcounter{equation}0
\begin{lm}\label{lmOmU}
Let $X$ be a regular flat
separated scheme of finite
type over $S$
and $D$ be
a divisor with simple normal crossings.

{\rm 1.}
Let $U$ be an open subscheme
of $X$ and 
$U\to P$
be a regular immersion
of codimension $1$ into
a smooth scheme $P$ over $S$ 
satisfying $D\cap U=
\widetilde D\times_PU$
as in Lemma {\rm \ref{lmXP}}.
Then, we have
a quasi-isomorphism
\begin{equation}
\begin{CD}
[N_{U/P}
\to
\Omega^1_{P/S}(\log \widetilde D)
\otimes_{{\cal O}_P}{\cal O}_U]
\to
L_{X/S}(\log D)|_U.
\end{CD}
\label{eqOmU}
\end{equation}
Consequently,
the logarithmic cotangent
complex
$L_{X/S}(\log D)$
satisfies the condition {\rm (L($n$))} in
Section {\rm \ref{sslcc}}
and we have an exact sequence
\begin{equation}
N_{U/P}
\to
\Omega^1_{P/S}(\log \widetilde D)
\otimes_{{\cal O}_P}{\cal O}_U
\to
\Omega^1_{U/S}(\log D)\to 0.
\label{eqOmUS}
\end{equation}

{\rm 2.}
Let $D_1,\ldots,D_m$
be the irreducible
components of $D$.
Then, we have
a distinguished triangle
\begin{equation}
\begin{CD}
@>>>
L_{X/S}
@>>>
L_{X/S}(\log D)
@>>>
\bigoplus_{i=1}^m
{\cal O}_{D_i}
@>>>
\end{CD}
\label{eqOmL}
\end{equation}
and consequently
 an exact sequence
\begin{equation}
\begin{CD}
0@>>>
\Omega^1_{X/S}
@>>>
\Omega^1_{X/S}(\log D)
@>>>
\bigoplus_{i=1}^m
{\cal O}_{D_i}
@>>>0.
\end{CD}
\label{eqOm}
\end{equation}
\end{lm}
{\it Proof.}
1.
The distinguished triangle 
(\ref{eqLtr}) for
the immersions
$U\to (U\times_SX)^\sim
\to (U\times_SP)^\sim$
defines a distinguished
triangle
$$\delta|_U^*N_{(U\times_SX)^\sim
/(U\times_SP)^\sim}
\to
N_{U/(U\times_SP)^\sim}
\to M_{U/(U\times_SX)^\sim}\to.$$
Since 
$U\to (U\times_SP)^\sim$
is a section of
the smooth morphism
$(U\times_SP)^\sim\to U$,
the isomorphism
${\rm pr}_2^*
\Omega^1_{P/S}(\log \widetilde D)
\to
\Omega^1_{(U\times_SP)^\sim/U}$
induces a canonical isomorphism
$N_{U/(U\times_SP)^\sim}
\to 
\Omega^1_{P/S}(\log \widetilde D)
\otimes_{{\cal O}_P}{\cal O}_U$.
By the cartesian diagram
(\ref{eqPUX}), 
the canonical map
$N_{U/P}\to
\delta|_U^*N_{(U\times_SX)^\sim
/(U\times_SP)^\sim}$
is an isomorphism.
Thus the assertion follows.

2.
The distinguished triangle 
(\ref{eqLtr}) for
$X\to (X\times_SX)^\sim
\to X\times_SX$
defines a distinguished
triangle
\begin{equation}
L_{X/S}\to
L_{X/S}(\log D)
\to L\delta^*
L_{(X\times_SX)^\sim/
X\times_SX}\to. 
\label{eqLXStr}
\end{equation}
Let $E_1,\ldots,E_m$
be the inverse images
by either of the two projections
$(X\times_SX)^\sim\to X$.
Then, we have
a canonical isomorphism
$\Omega^1_{(X\times_SX)^\sim/
X\times_SX}
\to 
\bigoplus_{i=1}^m
{\mathcal O}_{E_i}$.
It suffices to show that
this induces an isomorphism
$L\delta^*L_{(X\times_SX)^\sim/
X\times_SX}
\to 
\bigoplus_{i=1}^m
{\mathcal O}_{E_i}$.
The assertion is local on $X$.
By comparing the isomorphism
(\ref{eqOmU})
with the isomorphism
$[N_{U/P}
\to
\Omega^1_{P/S}
\otimes_{{\cal O}_P}{\cal O}_U]
\to
L_{U/S}$,
the
distinguished triangle
(\ref{eqLXStr})
implies that
$L\delta^*L_{(X\times_SX)^\sim/
X\times_SX}
\to 
\bigoplus_{i=1}^m
{\mathcal O}_{E_i}$
is an isomorphism.
\qed

If there
exists a dense open subscheme
of $X$ smooth over $S$,
the coherent ${\mathcal O}_X$-module
$\Omega^1_{X/S}(\log D)$
is locally free of rank $n-1$
on it
and the first map
$N_{U/P}
\to
\Omega^1_{P/S}(\log \widetilde D)
\otimes_{{\cal O}_P}{\cal O}_U$
in (\ref{eqOmUS}) is an injection.

On an open
subscheme $U\subset X$
with a regular immersion
$U\to P$
as in Lemma \ref{lmXP},
if $e_1,\ldots,e_n$
is a basis of
$\Omega^1_{P/S}(\log \widetilde D)
\otimes_{{\cal O}_P}{\cal O}_U$
and if
$a_1e_1+
\cdots+a_ne_n$
is the image
of a basis of
$N_{U/P}$,
then the restriction
$\Omega^n_{X/S}(\log D)|_U$
is isomorphic
to ${\cal O}_U/(a_1,\ldots,a_n)$
and hence
the annihilator ideal
${\rm Ann}\
\Omega^n_{X/S}(\log D)|_U
\subset
{\cal O}_U$
is generated by
$a_1,\ldots,a_n$.

We study the logarithmic cotangent
complex $L_{X/S}(\log D)$ more in detail.
First, we consider the case where
there exists a dense open subscheme
of $X$ smooth over $S$.

\begin{lm}\label{lmZ}
\setcounter{equation}0
We put
$X_F=\sum_i
l_iD_i$
as a divisor of $X$
and define
a Cartier divisor
$D'$ of $X$
by 
$D'=\sum_{D_i\subset X_F,
D_i\subset D}
l_iD_i$.
We put $Z=\Sigma_{X/S}$
and let $i:Z\to X$
be the closed immersion.
We also put $n=\dim X_K+1$.
Assume that there
exists a dense open subscheme
of $X$ smooth over $S$.
Then, we have the following.

{\rm 1.}
There
exists a unique 
${\cal O}_X$-linear map
$$\cdot d\log \pi\colon
{\cal O}_{D'}(D'-X_F)
={\cal O}_{D'}
\otimes_{{\cal O}_X}
{\cal I}_{X_F-D'}
\to 
\Omega^1_{X/S}(\log D)$$
sending
a local generator
$g$ of the ideal
${\cal I}_{X_F-D'}\subset
{\cal O}_X$
to 
$dg +g\cdot d\log (\pi/g)$
for a prime element $\pi$
of $K$.
The map
$\cdot d\log \pi$
is 
independent of 
the choice of 
a uniformizer
$\pi$.

{\rm 2.}
The map 
$\cdot d\log \pi\colon
{\cal O}_{D'}(D'-X_F)
\to 
\Omega^1_{X/S}(\log D)$
is an injection and the
cokernel
$\Omega^1_{X/S}(\log D)/
{\cal O}_{D'}(D'-X_F)$
is an ${\cal O}_X$-module
of tor-dimension $\le 1$.

{\rm 3.}
The ${\cal O}_Z$-module
$L_1i^*\Omega^1_{X/S}(\log D)$
is invertible.
For a normal scheme $W$
over $F$ and a morphism
$\varphi\colon
W\to Z$
over $S$ and for
the pull-back
$\varphi^*
L_1i^*\Omega^1_{X/S}(\log D)$,
there exists
a canonical isomorphism
$$N_{s/S}\otimes
{\cal O}_W\to
\varphi^*
L_1i^*\Omega^1_{X/S}(\log D)$$
of trivial
invertible ${\cal O}_W$-modules,
where $N_{s/S}$
denotes the conormal sheaf
${\mathfrak m}_K/
{\mathfrak m}_K^2$
of the closed point
$s$ of $S$.
\end{lm}

The complement $X\setminus Z$
is the largest open subscheme
of $X$ smooth over $S$,
which is assumed to be dense in $X$.
In the case
$X\setminus D\subset X_K$,
we have
$D'=X_F$ and
the cokernel
$\Omega^1_{X/S}(\log D)/
{\cal O}_{X_F}$
will be denoted by
$\Omega^1_{X/S}(\log D/\log F)$.

{\it Proof.}
1.
The local section $d\log (\pi/g)$
of $\Omega^1_{X/S}(\log D)$
is independent of the choice
of a prime element $\pi$.
We have
$d(ug) +(ug)\cdot d\log (\pi/ug)
=
u(dg +g\cdot d\log (\pi/g))
+g(du-ud\log u)$
for a unit $u$
and the last term is $0$.
Hence
the ${\cal O}_X$-linear map
$d\log \pi\cdot \colon
{\cal O}_X(D'-X_F)
\to 
\Omega^1_{X/S}(\log D)$
is well-defined.

Since
$(\pi/g)(dg +g\cdot d\log (\pi/g))
=(\pi/g)dg +g d(\pi/g)=
d\pi=0$,
it induces
${\cal O}_{D'}(D'-X_F)
\to 
\Omega^1_{X/S}(\log D)$.

2.
For the injectivity,
it suffices to
show it at
the generic point
of each irreducible
component $D_i$ of
$D'$.
Hence,
we may assume $X_K=X\setminus D$.
Then, it follows
from \cite[Lemma 5.3.4.2]{KSI}.

We show that
$\Omega^1_{X/S}(\log D)/
{\cal O}_{D'}(D'-X_F)$
is of tor-dimension $\le 1$.
Since the question is local,
we take an immersion
$U\to P$ as in Lemma \ref{lmXP}.
Let $\tilde g$
be a function on
$P$ lifting 
$g=\pi/\prod_it_i^{l_i}$.
Then, 
on a neighborhood of $U$,
the divisor
$U\subset P$
is defined by
an equation
$\pi=\tilde g\cdot\prod_iT_i^{l_i}$.
Hence,
the image of
$N_{U/P}\to
\Omega^1_{P/S}(\log 
\widetilde D)
\otimes_{{\cal O}_P}{\cal O}_U$
is generated
by $d(\tilde g\cdot\prod_iT_i^{l_i})$.
Since the image of
the section
$(\prod_it_i^{l_i})^{-1}
d(\tilde g\cdot\prod_iT_i^{l_i})
=d\tilde g+
\tilde g\cdot
\sum_il_id\log T_i$
of $\Omega^1_{P/S}
(\log \widetilde D)
\otimes_{{\cal O}_P}{\cal O}_U$
in $\Omega^1_{X/S}(\log D)|_U$
is $dg+g\cdot d\log(\pi/g)$,
we have
a locally free resolution
$0\to
N_{U/P}(D')\to
\Omega^1_{P/S}(\log 
\widetilde D)
\otimes_{{\cal O}_P}{\cal O}_U
\to
(\Omega^1_{X/S}(\log D)/
{\cal O}_{D'}(D'-X_F))|_U
\to0$.

3.
Since $\Omega^1_{X/S}(\log D)$
satisfies the condition L($n$)
in Section \ref{sslcc},
the ${\cal O}_Z$-module
$L_1i^*\Omega^1_{X/S}(\log D)$
is invertible.

We may assume $W$
is integral.
Let $i_W=i\circ \varphi
\colon W\to X$ denote
the composition.
Then,
since 
$\Omega^1_{X/S}(\log D)$
is of tor-dimension $\le 1$
and the ${\cal O}_Z$-modules
$i^*\Omega^1_{X/S}(\log D)$
and
$L_1i^*\Omega^1_{X/S}(\log D)$
are locally free, 
the canonical map
$\varphi^*
L_1i^*\Omega^1_{X/S}(\log D)
\to
L_1i_W^*
\Omega^1_{X/S}(\log D)$
is an isomorphism.

First we consider the case
$i_W(W)\not\subset D$.
By $i_W(W)\not\subset D$,
we have
$L_1i_W^*
{\mathcal O}_{D_i}=0$.
Hence the exact sequence
(\ref{eqOm})
induces an isomorphism
$L_1i_W^*
\Omega^1_{X/S}
\to
L_1i_W^*
\Omega^1_{X/S}(\log D)$
of invertible
${\cal O}_W$-modules.
Let $Z'$ be the
closed subscheme of $X$ 
defined by
the annihilator
ideal of
$\Omega^n_{X/S}=
\Lambda^n\Omega^1_{X/S}$.
For a morphism $g\colon T\to X$
of schemes,
the ${\mathcal O}_T$-module $L_1g^*\Omega^1_{X/S}$
is invertible if and only
if $g$ factors through $Z'$,
since $\Omega^1_{X/S}$
satisfies the condition L($n$)
in Section \ref{sslcc}.
Since
$L_1i_W^*
\Omega^1_{X/S}$
is invertible,
the map
$i_W\colon
W\to X$ 
factors through
the closed subscheme $Z'$.
Hence, the assertion follows from
\cite[Lemma 5.1.3.1]{KSI}.

Next, we consider the case
$i_W(W)\subset D$.
Then the exact sequence
$0\to
{\cal O}_X(-X_F)
\to
{\cal O}_X(D'-X_F)
\to
{\cal O}_{D'}(D'-X_F)\to 0$
defines an isomorphism
$i_W^*{\cal O}_X(-X_F)
=N_{s/S}\otimes {\cal O}_W\gets
L_1i_W^*{\cal O}_{D'}(D'-X_F)$.
Further the map
$d\log \pi\cdot\colon
{\cal O}_{D'}(D'-X_F)
\to 
\Omega^1_{X/S}(\log D)$
induces a map
$L_1i_W^*{\cal O}_{D'}(D'-X_F)
\to 
L_1i_W^*\Omega^1_{X/S}(\log D)$.
We show that it is an isomorphism.

Since the question is local,
we take an immersion
$U\to P$ as in Lemma \ref{lmXP}.
Then, by the proof of 2.,
we have a commutative
diagram of exact sequences
\begin{equation}
\begin{CD}
0@>>>
N_{U/P}
@>>>
N_{U/P}(D')
@>>>
{\cal O}_{D'}(D'-X_F)|_U
@>>>0\\
@.@|@VVV @VV{d\log \pi\cdot|_U}V@.\\
0@>>>
N_{U/P}
@>>>
\Omega^1_{P/S}(\log \widetilde D)
\otimes_{{\cal O}_P}{\cal O}_U
@>>>
\Omega^1_{X/S}(\log D)|_U
@>>>0.
\end{CD}
\end{equation}
Then both
$L_1i_W^*{\cal O}_{D'}(D'-X_F)$
and 
$L_1i_W^*\Omega^1_{X/S}(\log D)$
are identified with
$i_W^*N_{U/P}$
and
the assertion follows.
\qed

Next, we consider
the case where
$\Sigma_{X/S}=X$.
This occurs only if the characteristic of
$K$ is $p>0$.

\setcounter{equation}0
\begin{lm}\label{lmLXS}
Assume $\Sigma_{X/S}=X$.
Then, 
there exists a canonical isomorphism
$ \Omega^1_S
\otimes{\mathcal O}_X
\to
{\mathcal H}_1(L_{X/S}(\log D))$ of
invertible ${\mathcal O}_X$-modules.
\end{lm}

{\it Proof.}
Let $K_0=K^p\subset K$
and put $S_0={\rm Spec}\ {\mathcal O}_{K_0}$.
The composition of closed immersions
$X\to (X\times_SX)^\sim
\to (X\times_{S_0}X)^\sim$
defines a distinguished triangle
\begin{equation}
\begin{CD}
\to
L\delta^*
L_{(X\times_SX)^\sim/(X\times_{S_0}X)^\sim}[1]
@>>>
L_{X/S_0}(\log D)
@>>>
L_{X/S}(\log D)\to
\end{CD}
\label{eqLS0}
\end{equation}
of cotangent complexes.
By the cartesian diagram
$$\begin{CD}
(X\times_SX)^\sim
@>>>
(X\times_{S_0}X)^\sim\\
@VVV @VVV\\
S@>>> S\times_{S_0}S,
\end{CD}$$
we obtain a surjection
$\Omega^1_{S/S_0}
\otimes_{{\mathcal O}_S}
{\mathcal O}_X
\to
N_{(X\times_SX)^\sim/
(X\times_{S_0}X)^\sim}$.
We claim that this map
and the map
${\mathcal H}_1L_{X/S}(\log D)
\to
N_{(X\times_SX)^\sim/
(X\times_{S_0}X)^\sim}$
defined by (\ref{eqLS0}) 
are isomorphisms.
This will complete the proof
since $\Omega^1_{S/S_0}=
\Omega^1_S$
is an invertible ${\mathcal O}_S$-module.

To show the claim,
it suffices to show that
the canonical map
${\mathcal H}_1L_{X/S_0}(\log D)
\to
{\mathcal H}_1L_{X/S}(\log D)$
is the 0-map since
${\mathcal H}_0L_{X/S_0}(\log D)
\to
{\mathcal H}_0L_{X/S}(\log D)$
is an isomorphism
and ${\mathcal H}_1L_{X/S}(\log D)$
is invertible.
We show that
${\mathcal H}_1L_{X/S_0}(\log D)
\to
{\mathcal H}_1L_{X/S}(\log D)$
is the 0-map.
Since the assertion is local on $X$,
we may assume that
there is a regular immersion $X\to P_0$
of codimension 1
to a smooth scheme $P_0$ over
$S_0$ as in Lemma \ref{lmXP}.
It induces an immersion $X\to 
P=P_0\times_{S_0}S$
to a smooth scheme $P$ over $S$.
Then, the map
${\mathcal H}_1L_{X/S_0}(\log D)
\to
{\mathcal H}_1L_{X/S}(\log D)$
is identified with
the canonical map
$N_{X/P_0}\to N_{X/P}$
of the conormal sheaves
induced by the projection 
$P=P_0\times_{S_0}S\to P_0$.
Since $X\times_{S_0}S
=X\times_{P_0}P$
regarded as a Cartier divisor
of $P$ is
$p$-times the Cartier
divisor $X$ of $P$,
the assertion follows.
\qed

We study a consequence of 
Lemma {\rm \ref{lmYQ}}.

\begin{lm}\label{corYQ}
\setcounter{equation}0
Let $f\colon Y\to X$
be as in Lemma {\rm \ref{lmYQ}}
and put $n=\dim X_K+1$
and $n'=\dim Y_K+1$ respectively.
Then, for the closed subschemes
$\Sigma_{X/S}\subset X$
and 
$\Sigma_{Y/S}\subset Y$
defined by
the annihilators
$\Lambda^n\Omega^1_{X/S}(\log D)$
and
$\Lambda^{n'}\Omega^1_{Y/S}(\log E)$,
the pull-back
$\Sigma_{X/S}\times_XY$
is a subscheme of
$\Sigma_{Y/S}$.
\end{lm}

{\it Proof.}
The assertion is local on $Y$.
We consider a cartesian diagram
(\ref{eqYQ}) satisfying
the condition in
Lemma \ref{lmYQ}.
Then, we have a commutative diagram
of exact sequences
\begin{equation}
\begin{CD}
0@>>>
f|_V^*N_{U/P}
@>>>
f|_V^*(\Omega^1_{P/S}(\log \widetilde D)
\otimes_{{\cal O}_P}{\cal O}_U)
@>>>
f|_V^*\Omega^1_{X/S}(\log D)|_U
@>>>0\\
@.@V{\simeq}VV@VVV@VVV@.\\
0@>>>
N_{V/Q}
@>>>
\Omega^1_{Q/S}(\log \widetilde E)
\otimes_{{\cal O}_Q}{\cal O}_V
@>>>
\Omega^1_{Y/S}(\log E)|_V
@>>>0
\end{CD}
\label{eqOmUV}
\end{equation}
The closed subscheme
$(\Sigma_{X/S}\times_XY)\cap V$
is the largest closed subscheme
where the pull-back
$f|_V^*N_{U/P}
\to
f|_V^*(\Omega^1_{P/S}(\log \widetilde D)
\otimes_{{\cal O}_P}{\cal O}_U)$
is the zero-map.
Similarly 
$\Sigma_{Y/S}\cap V$
is the largest closed subscheme
where the pull-back
$N_{V/Q}
\to\Omega^1_{Q/S}(\log \widetilde E)
\otimes_{{\cal O}_Q}{\cal O}_V$
is the zero-map.
Hence
the assertion follows
from (\ref{eqOmUV}).
\qed

Next, we study consequences of 
Lemma {\rm \ref{lmX2}}.
As loc.\ cit.,
we assume that $X_1$ is a regular divisor
of $X$
and that $D_1=D\cap X_1$
is a divisor with simple normal crossings.
Let $(X\times_SX)^\sim$
and $(X\times_SX)^\approx$
be the log product with respect
to $D$ and $D'=D\cup X_1$
respectively
and $E$ be the inverse image of
$X_1$ by either of the two projections
$(X\times_SX)^\approx\to X$.
We have cartesian diagrams
\addtocounter{thm}1
\setcounter{equation}0
\begin{equation}
\begin{CD}
X_1@>{\subset}>>E\\
@VVV @VVV\\
X@>{\subset}>>
(X\times_SX)^\approx,
\end{CD}
\quad
\begin{CD}
X_1@>{\subset}>>
(X_1\times_SX_1)^\sim\\
@VVV @VVV\\
X@>{\subset}>>
(X\times_SX)^\sim
\end{CD}
\label{eqcartE}
\end{equation}
by the definition of $E$
and by the universality of log product,
respectively.

First, we consider
the case where $X_1$ is flat over $S$.

\addtocounter{thm}{-1}
\begin{lm}\label{lmXX1}
Let the notation be as above
and assume that
$X_1$ is flat over $S$.
Let $i\colon
X_1\to X$ denote
the closed immersion.
Then, 
the left cartesian diagram
{\rm (\ref{eqcartE})}
defines a distinguished
triangle
$$
L_{X_1/S}(\log D_1)
\to
Li^*L_{X/S}(\log D')
\to 
{\mathcal O}_{X_1}
\to$$
on $X_1$.
\end{lm}

{\it Proof.}
By Lemma  {\rm \ref{lmX2}.1},
the vertical arrows in
the left cartesian diagram
{\rm (\ref{eqcartE})}
are regular immersion of
codimension 1.
Hence the canonical map
$Li^*L_{X/S}(\log D')
\to M_{X_1/E}$
to the conormal complex
is a quasi-isomorphism.
Since $E\to
(X_1\times_SX_1)^\sim$
is smooth,
the morphisms $X_1\to E\to
(X_1\times_SX_1)^\sim$
define a distinguished triangle
$L_{X_1/S}(\log D_1)
\to M_{X_1/E}\to
\Omega^1_{E/
(X_1\times_SX_1)^\sim}
\otimes{\mathcal O}_{X_1}$.
Since $E$ is a ${\mathbf G}_m$-torsor
over $(X_1\times_SX_1)^\sim$
splitting on $X_1$,
we have a canonical trivialization
$\Omega^1_{E/
(X_1\times_SX_1)^\sim}
\otimes{\mathcal O}_{X_1}
\to
{\mathcal O}_{X_1}$.
Thus, the assertion follows.
\qed

Next, we assume that $X_1$ is
contained
in the closed fiber $X_F$.
Since the immersions
$X_1\to E$
and $X_1\to
(X_1\times_FX_1)^\sim$
are reguler immersions
of smooth schemes over $F$,
the cartesian diagram (\ref{eqcartE})
define the excess conormal
complexes
$M'_{X/(X\times_SX)^\approx,E}$ 
and
$M'_{X/(X\times_SX)^\sim,
(X_1\times_FX_1)^\sim}$
as complexes of 
${\cal O}_{X_1}$-modules
as in \cite[Definition 1.6.3.2]{KSI}.
By the definitions
$L_{X/S}(\log D)
=M_{X/(X\times_SX)^\sim}$
and
$L_{X/S}(\log D')
=M_{X/(X\times_SX)^\approx}$,
they fit in the distinguished 
triangles
\addtocounter{thm}1
\setcounter{equation}0
\begin{align}
\to
M'_{X/(X\times_SX)^\approx,E}
\to & Li_1^*
L_{X/S}(\log D')
\to
N_{X_1/E}
\to,
\label{eqME}
\\
\to
M'_{
X/(X\times_SX)^\sim,
(X_1\times_FX_1)^\sim}
\to & Li_1^*
L_{X/S}(\log D)
\to
\Omega^1_{X_1/F}(\log D_1)
\to
\label{eqMX1}
\end{align}
by \cite[Proposition 1.6.4.2]{KSI}.
More concretely,
they are described as follows.

\addtocounter{thm}{-1}
\begin{lm}\label{lmX1}
Let the notation be as above
and assume that
$X_1$ is
a smooth scheme
of dimension $d$
over $F$.
Let $N_{s/S}$
be the conormal sheaf
${\mathfrak m}_K/
{\mathfrak m}_K^2$
of the closed point
$s$ of $S$.

{\rm 1.}
There exists 
a canonical isomorphism
$M'_{X/(X\times_SX)^\approx,
E}
\to N_{s/S}
\otimes_F{\cal O}_{X_1}[1]$.

{\rm 2.}
There exists a canonical isomorphism
$M'_{X/(X\times_SX)^\sim,
(X_1\times_FX_1)^\sim}
\to [N_{s/S}
\otimes_F{\cal O}_{X_1}
\to N_{X_1/X}]$.
The complex
$M'_{X/(X\times_SX)^\sim,
(X_1\times_FX_1)^\sim}$
on $X_1$
is acyclic outside the intersection
$X_1\cap \Sigma_{X/S}$
where $\Sigma_{X/S}$
is defined in Definition
{\rm \ref{dfLXSl}}.
\end{lm}

{\it Proof.}
1.
Let $i_1\colon X_1\to X$
be the closed immersion.
Since the immersion
$X_1\to E$ is a
regular immersion
of codimension $d+1$, 
the canonical surjection
$i_1^*
\Omega^1_{X/S}(\log D')
\to
N_{X_1/E}$
is an isomorphism.
Hence, by the distinguished triangle
(\ref{eqME}),
we obtain an isomorphism
$M'_{X/(X\times_SX)^\approx,E}
\to L_1i_1^*
L_{X/S}(\log D')[1]$.
Thus the assertion follows
from
Lemma \ref{lmZ}.3
and Lemma \ref{lmLXS}.

2. 
Similarly as Lemma 
\ref{lmOmU}.2,
we have a distinguished triangle
$\to
L_{X/S}(\log D)
\to
L_{X/S}(\log D')
\to 
{\cal O}_{X_1}\to $
and hence
a distinguished triangle
$\to
Li_1^*L_{X/S}(\log D)
\to
Li_1^*L_{X/S}(\log D')
\to 
Li_1^*{\cal O}_{X_1}\to $.
By (\ref{eqME}), (\ref{eqMX1})
and by the exact sequence
$0\to
\Omega^1_{X_1/F}(\log D_1)
\to
N_{X_1/E}
\to 
{\cal O}_{X_1}\to 0$,
we obtain a distinguished
triangle
$\to 
M'_{X/(X\times_SX)^\sim,
(X_1\times_FX_1)^\sim}
\to 
M'_{
X/(X\times_SX)^\approx,E_1}
\to 
L_1i_1^*{\cal O}_{X_1}[1]
\to.$
Hence the assertion
follows from 1.\
and the isomorphism
$L_1i_1^*{\cal O}_{X_1}
\to N_{X_1/X}.$

On the complement
of $\Sigma_{X/S}$,
the immersion
$X_1\to X_F$
is an open immersion
and 
the complex
$M'_{X/(X\times_SX)^\sim,
(X_1\times_FX_1)^\sim}$
is acyclic.
\qed

Assume that the underlying
set of the closed fiber
$X_F$ is a subset of $D$
and put $n=\dim X_K+1$.
We give a variant of Lemma \ref{lmX}
and Lemma \ref{lmZ}.
We define a variant
$(X\times_{\mathbb S}X)^\sim$
of the log product
by the cartesian diagram
$$\begin{CD}
(X\times_{\mathbb S}X)^\sim
@>>> (X\times_SX)^\sim\\
@VVV @VVV\\
S@>>> (S\times_SS)^\sim
\end{CD}$$
where
$(S\times_SS)^\sim$
is the log product defined
defined with respect to the Cartier
divisor $s$ of $S$.

\begin{lm}\label{lmXX0}
Let $X$ be a regular flat
scheme of finite type
over $S$
and $D\subset X$
be a divisor with simple normal
crossings.
Assume that the underlying
set of the
closed fiber
$X_F$ is a subset of $D$.

Then the scheme
$(X\times_{\mathbb S}X)^\sim$
is flat and locally a hypersurface
over $X$ with respect
to either of the projections.
\end{lm}

\noindent{\it Proof.}
The log scheme $X$ with the log
structure defined by $D$
is log flat 
(\cite[Section 4.3]{KSI})
and log locally of complete
intersection
(\cite[Definition 4.4.2]{KSI})
over 
$S$ with the log
structure defined by 
the closed point,
similarly as \cite[Lemma 5.2.1]{KSI}.
Since the projection
$(X\times_{\mathbb S}X)^\sim
\to X$ is strict,
it is flat.
Since the assertion is local,
we take a regular
immersion $U\to P$ of codimension 1
as in Lemma \ref{lmXP}.
Then, we have a closed immersion
$(U\times_{\mathbb S}X)^\sim
\to
(P\times_SX)^\sim$.
Then,
since $(P\times_SX)^\sim$
is smooth over $X$
and 
$(U\times_{\mathbb S}X)^\sim$
is locally of complete intersection,
the immersion
$(U\times_{\mathbb S}X)^\sim
\to
(P\times_SX)^\sim$
is a regular immersion
by Lemma \ref{lmciri}.
We verify that
it is a regular immersion of 
codimension 1
by reducing to the case where $D$
is empty.
\qed

We also define a variant
$L_{X/S}(\log D/\log F)$ of
the logarithmic cotangent complex
to be the conormal complex
$M_{X/(X\times_{\mathbb S}X)^\sim}$.
The coherent ${\mathcal O}_X$-module
$\Omega^1_{X/S}(\log D/\log F)$
defined as
${\mathcal H}_0(L_{X/S}(\log D/\log F))$
is the conormal sheaf
$N_{X/(X\times_{\mathbb S}X)^\sim}$.
For a regular immersion
$U\to P$ as in Lemma \ref{lmXP},
we have shown in the proof
of Lemma \ref{lmXX0}
that the immersion
$(U\times_{\mathbb S}X)^\sim
\to (P\times_SX)^\sim$
is a regular immersion
of codimension 1.
The immersions
$U\to (U\times_{\mathbb S}X)^\sim
\to (P\times_SX)^\sim$ define a distinguished
triangle
$\delta|_U^*N_{(U\times_{\mathbb S}X)^\sim
/(P\times_SX)^\sim}
\to
N_{U/(P\times_SX)^\sim}
\to M_{U/(U\times_{\mathbb S}X)^\sim}\to$.
Similarly as for $L_{X/S}(\log D)$,
it defines a quasi-isomorphism
\setcounter{equation}0
\addtocounter{thm}{1}
\begin{equation}
\begin{CD}
[N_{(U\times_{\mathbb S}X)^\sim
/ (P\times_SX)^\sim}
\otimes{\cal O}_U
\to
\Omega^1_{P/S}(\log \widetilde D)
\otimes_{{\cal O}_P}{\cal O}_U]
\to
L_{X/S}(\log D/\log F)|_U.
\end{CD}
\label{eqOmUl}
\end{equation}
This shows that
the logarithmic cotangent
complex
$L_{X/S}(\log D/\log F)$
satisfies the condition {\rm (L($n$))} in
Section {\rm \ref{sslcc}}.
If there
exists a dense open subscheme
of $X$ smooth over $S$,
the coherent ${\mathcal O}_X$-module
$\Omega^1_{X/S}(\log D/\log F)$
is locally free of rank $n-1$
on a dense open subscheme
and the map
$N_{(U\times_{\mathbb S}X)^\sim
/ (P\times_SX)^\sim}
\otimes{\cal O}_U
\to
\Omega^1_{P/S}(\log \widetilde D)
\otimes_{{\cal O}_P}{\cal O}_U$
in (\ref{eqOmUl}) is an injection.

The immersions
$X\to (X\times_{\mathbb S}X)^\sim
\to (X\times_SX)^\sim$ define a
distinguished triangle
\begin{equation}
\begin{CD}
L\delta^*
L_{(X\times_{\mathbb S}X)^\sim/
(X\times_SX)^\sim}
@>>>
L_{X/S}(\log D)
@>>>
L_{X/S}(\log D/\log F)
@>>>.
\end{CD}
\label{eqLlog}
\end{equation}
Since $(S\times_SS)^\sim
={\rm Spec}\
{\mathcal O}_K[U^{\pm1}]/
((U-1)\pi)$
for a prime element $\pi$,
the conormal sheaf
$N_{S/(S\times_SS)^\sim}$
is isomorphic to $F$
and is generated by $d\log \pi$.
Hence,
the distinguished triangle
(\ref{eqLlog})
gives an exact sequence
\begin{equation}
\begin{CD}
{\mathcal O}_{X_F}
@>>>
\Omega^1_{X/S}(\log D)
@>>>
\Omega^1_{X/S}(\log D/\log F)
@>>>0.
\end{CD}
\label{eqomom}
\end{equation}
If there exists a dense open
subscheme of $X$
smooth over $S$,
the first arrow
${\mathcal O}_{X_F}
\to
\Omega^1_{X/S}(\log D)$
is an injection 
by Lemma \ref{lmZ}.2.
If $X$ is nowhere smooth over $S$,
the second arrow
$\Omega^1_{X/S}(\log D)
\to
\Omega^1_{X/S}(\log D/\log F)$
is an isomorphism
of locally free ${\mathcal O}_X$-modules
of rank $n+1$.
We put
$\Omega^n_{X/S}(\log D/\log F)
=\Lambda^n
\Omega^1_{X/S}(\log D/\log F)$.

\addtocounter{thm}{-1}
\begin{lm}\label{lmXX}
Let $Z$ be the closed subscheme
of $X$ defined by
the ideal ${\rm Ann}\
\Omega^n_{X/S}(\log D/\log F)$
and 
$i:Z\to X$
be the closed immersion.
Then the restriction
of the invertible ${\cal O}_Z$-module
$L_1i^*\Omega^1_{X/S}
(\log D/\log F)$
to the reduced closed fiber
$Z_{F,{\rm red}}$
is trivial.
\end{lm}

\noindent{\it Proof.}
It suffices to consider
the cases where
the complement $X\setminus Z$
is dense
and $Z=X$ respectively.
First, we prove
the case where
the complement $X\setminus Z$
is dense.
In this case,
the proof is similar
to \cite[Lemma 5.3.5.1]{KSI}.
Let $i'\colon Z_{F,{\rm red}}
\to X$ be the immersion.
Similarly as in 
Lemma \ref{lmZ}.3.,
the restriction
of the invertible ${\cal O}_Z$-module
$L_1i^*\Omega^1_{X/S}
(\log D/\log F)$
to 
$Z_{F,{\rm red}}$
is isomorphic to
$L_1i^{\prime *}
\Omega^1_{X/S}
(\log D/\log F)$.

By the exact sequence
(\ref{eqomom})
together with the
injectivity of
${\mathcal O}_{X_F}
\to
\Omega^1_{X/S}(\log D)$,
we obtain an exact sequence
\begin{align*}
0&\to L_1i^{\prime *}{\cal O}_{X_F}
\to 
L_1i^{\prime *}\Omega^1_{X/S}
(\log D)
\to 
L_1i^{\prime *}\Omega^1_{X/S}
(\log D/\log F)\\
&
\to i^{\prime *}{\cal O}_{X_F}
\to 
i^{\prime *}\Omega^1_{X/S}
(\log D)
\to 
i^{\prime *}\Omega^1_{X/S}
(\log D/\log F)
\to 0.
\end{align*}
By a local description
of $Z$ similar to
that given before 
Lemma \ref{lmZ},
it follows that
the last map
$i^{\prime *}\Omega^1_{X/S}
(\log D)
\to 
i^{\prime *}\Omega^1_{X/S}
(\log D/\log F)$
is an isomorphism
of locally free
${\cal O}_{Z_{F,{\rm red}}}$-modules.
Hence the map
$L_1i^{\prime *}\Omega^1_{X/S}
(\log D/\log F)
\to i^{\prime *}{\cal O}_{X_F}$
is a surjection 
of invertible 
${\cal O}_{Z_{F,{\rm red}}}$-modules
and is an isomorphism.
Therefore,
the map
$L_1i^{\prime *}{\cal O}_{X_F}
\to 
L_1i^{\prime *}\Omega^1_{X/S}
(\log D)$
is also an isomorphism
of invertible
${\cal O}_{Z_{F,{\rm red}}}$-modules.
Since
$L_1i^{\prime *}{\cal O}_{X_F}$
is isomorphic to
${\cal O}_{Z_{F,{\rm red}}}$,
the assertion follows.

We show the case where
$Z=X$.
The proof in this case
is similar to that of Lemma \ref{lmLXS}.
We show that
there exists a canonical isomorphism
$ \Omega^1_S
\otimes{\mathcal O}_X
\to
{\mathcal H}_1L_{X/S}(\log D/\log F)$ of
invertible ${\mathcal O}_X$-modules.
Let $K_0=K^p\subset K$
and put $S_0={\rm Spec}\ {\mathcal O}_{K_0}$
as in the proof of
Lemma \ref{lmLXS}.
The composition of closed immersions
$X\to (X\times_{\mathbb S}X)^\sim
\to (X\times_{S_0}X)^\sim$
defines a distinguished triangle
\setcounter{equation}0
\begin{equation}
\begin{CD}
\to
L\delta^*
L_{(X\times_{\mathbb S}X)^\sim/(X\times_{S_0}X)^\sim}[1]
@>>>
L_{X/S_0}(\log D)
@>>>
L_{X/S}(\log D/\log F)\to
\end{CD}
\label{eqLS0l}
\end{equation}
of cotangent complexes.
By the cartesian diagram
$$\begin{CD}
(X\times_{\mathbb S}X)^\sim
@>>>
(X\times_{S_0}X)^\sim\\
@VVV @VVV\\
S@>>> (S\times_{S_0}S)^\sim,
\end{CD}$$
we obtain a surjection
$\Omega^1_{S/S_0}(\log F)
\otimes_{{\mathcal O}_S}
{\mathcal O}_X
\to
N_{(X\times_{\mathbb S}X)^\sim/
(X\times_{S_0}X)^\sim}$.
Similarly as in the proof of
Lemma \ref{lmLXS},
this map
and the map
${\mathcal H}_1L_{X/S}(\log D/\log F)
\to
N_{(X\times_{\mathbb S}X)^\sim/
(X\times_{S_0}X)^\sim}$
defined by (\ref{eqLS0l}) 
are isomorphisms.
This complete the proof
since $\Omega^1_{S/S_0}(\log F)$
is isomorphic to ${\mathcal O}_S$.
\qed

\subsection{Intersection 
product with the log diagonal}
\label{sslipd}

\setcounter{equation}0
Let $X$ be a regular flat 
separated scheme
of finite type
over $S={\rm Spec}\ {\cal O}_K$ 
and $D\subset X$
be a divisor with simple normal crossings.
We put
$n=\dim X_K+1$.
We define the localized
intersection product
with the log diagonal as follows.

We recall the notation
from the previous subsection.
The log product
$(X\times_SX)^\sim$
is defined
with respect to
the family
${\cal D}=
(D_i)_{i\in I}$
of irreducible components of $D$.
The logarithmic cotangent
complex $L_{X/S}(\log D)$
is defined as
the conormal complex
$M_{X/(X\times_SX)^\sim}$
of the log diagonal
$\delta\colon
X=\Delta_X^{\log}
\to (X\times_SX)^\sim$
and we have a canonical isomorphism
${\cal H}_0
(L_{X/S}(\log D))
\to \Omega^1_{X/S}(\log D)$.
We define a closed subscheme
$\Sigma_{X/S}$ of $X$ to be that 
defined by the annihilator
ideal of
$\Omega^n_{X/S}(\log D)=
\Lambda^n( \Omega^1_{X/S}(\log D))$.
Let
${\cal L}_{\Sigma_{X/S}}$
denote the invertible 
${\mathcal O}_{\Sigma_{X/S}}$-module
$L_1i^*L_{X/S}(\log D)$
where
$i\colon \Sigma_{X/S}
\to X$ denote
closed immersion.

We consider the following
special case of
Definition \ref{dflip1}.
We consider
$X$ and the log product
$(X\times_SX)^\sim$
as $S$ and $X$ in
Definition \ref{dflip1}.
The log product
$(X\times_SX)^\sim$
is locally a hypersurface
of relative dimension $n-1$
over $X$
by either of the two projections
by Lemma \ref{lmX}.
As $V$ in 
Definition \ref{dflip1},
we take $X$ regarded
as a closed subscheme of
$(X\times_SX)^\sim$
by the log diagonal.
Since the canonical map 
$\delta^*\Omega^1_{(X\times_SX)^\sim/X}\to \Omega^1_{X/S}(\log D)$
is an isomorphism,
the intersection
$Z\times_XV\subset V$  in
Definition \ref{dflip1}
is
$\Sigma_{X/S}\subset X$
in our setting.

We consider
a scheme $W$
of finite type over $S$
and a morphism 
$g\colon W\to 
(X\times_SX)^\sim$ over $S$
as $W\to X$ in
Definition \ref{dflip1}.
Then, $Z_T\subset T\subset W$ in
Definition \ref{dflip1}
are the inverse images
$g^{-1}(\Sigma_{X/S})
\subset
g^{-1}(\Delta_X^{\log})
\subset W$
where $\Delta_X^{\log}$
denotes $X$ regarded
as a closed subscheme
of 
$(X\times_SX)^\sim$
by the log diagonal.
Since the canonical map 
$L\delta^*L_{(X\times_SX)^\sim/X}\to L_{X/S}(\log D)$
is an isomorphism,
the pull-back of
${\mathcal L}_Z$ to $Z\times_XV$
in Definition \ref{dflip1}
is the invertible sheaf
${\cal L}_{\Sigma_{X/S}}
=L_1i^*L_{X/S}(\log D)$
on $\Sigma_{X/S}$
in our setting.
By Lemma below,
under the assumption (A)
there,
the group $G(Z_T)_{/{\cal L}_Z}$
in Definition \ref{dflip1}
is
$G(g^{-1}(\Sigma_{X/S}))$
in our setting.

\begin{lm}\label{lmW}
Let the notation be
as above.
In particular, 
let $W$ be a scheme
of finite type over $S$
and 
$g\colon W\to 
(X\times_SX)^\sim$ 
a morphism over $S$.
Assume that
the morphism 
$g\colon W\to 
(X\times_SX)^\sim$ 
satisfies the
condition:
\begin{itemize}
\item[{\rm (A)}]
The inverse image
$g^{-1}(\Sigma_{X/S})$
is supported on 
the closed fiber $W_F$
set-theoretically.
\end{itemize}
Then the multiplication
of the pull-back of
the invertible 
${\mathcal O}_{\Sigma_{X/S}}$-module
${\cal L}_{\Sigma_{X/S}}$
on
the Grothendieck group
$G(g^{-1}(\Sigma_{X/S}))$
is the identity.
\end{lm}

{\it Proof.}
Since the Grothendieck group
$G(g^{-1}(\Sigma_{X/S}))$
is generated by the 
push forward of the classes
of the normalizations
of integral closed subschemes,
it follows from
Lemma \ref{lmZ}.3.\
and Lemma \ref{lmLXS}.
\qed

Thus, we make the following definition
as in \cite[Definition 5.1.5]{KSI}.

\begin{df}\label{dflip}
\setcounter{equation}0
Let $X$ be a regular flat 
separated scheme
of finite type
over $S={\rm Spec}\ {\cal O}_K$ and $D\subset X$
be a divisor with simple normal crossings.
Let $W$ be a 
scheme of finite type
over $S$
and $g\colon W
\to (X\times_SX)^\sim$
be a morphism over $S$
satisfying the condition
{\rm (A)}
in Lemma {\rm \ref{lmW}}.
Then, 
we define {\rm the localized intersection
product
with the log diagonal}
\begin{equation}
\begin{CD}
((\Delta_X^{\log},\ \ ))_{
(X\times_SX)^\sim}
\colon
G(W)@>>>
G(g^{-1}(\Sigma_{X/S}))
\end{CD}
\label{eqprod1}
\end{equation}
as the product {\rm(\ref{eqlip})}
with the class of
${\mathcal F}=
{\mathcal O}_X$.
\end{df}

In \cite[Definition 5.1.5]{KSI},
we defined the localized
intersection product
with the log diagonal
under the assumption
that the generic fiber
is smooth and $D=X_F$.
Here, we replace the assumption
by the condition (A) in Lemma \ref{lmW}.

The logarithmic
localized intersection
product
$((\Delta_X^{\log},\ \ ))_{
(X\times_SX)^\sim}
\colon
G(W)\to
G(g^{-1}(\Sigma_{X/S}))$
preserves the topological
filtration in the sense
that
it induces a map
$$F_q
G(W)\to
F_{q-n}
G(g^{-1}(\Sigma_{X/S}))$$
\cite[Theorem 3.4.3.1]{KSI}.

If we take
a closed immersion
$A\to(X\times_SX)^\sim$
as $W\to
(X\times_SX)^\sim$,
the condition (A) in Lemma \ref{lmW}
can be written as
\begin{itemize}
\item[{\rm (A$'$)}]
The intersection
$\delta^{-1}(A)\cap \Sigma_{X/S}$
is supported on 
the closed fiber $X_F$
set-theoretically.
\end{itemize}
Under this assumption,
the localized intersection
product
with the log diagonal
\begin{equation}
\begin{CD}
((\ \ ,\Delta_X^{\log}))_{
(X\times_SX)^\sim}
\colon
G(A)@>>>
G(\delta^{-1}(A)\cap \Sigma_{X/S})
\end{CD}
\label{eqprod12}
\end{equation}
is defined as
the product {\rm(\ref{eqlip})}
with the class of
${\mathcal G}=
{\mathcal O}_X$
by taking
$X\to (X\times_SX)^\sim$
as $W\to X$ in Definition \ref{dflip1}.
By the symmetry of ${\mathcal T}or$,
we have
$$((\ \ ,\Delta_X^{\log}))_{
(X\times_SX)^\sim}
=
((\Delta_X^{\log},\ \ ))_{
(X\times_SX)^\sim}.$$

The localized intersection
product
with log diagonal
has the following
functoriality.

\begin{lm}\label{corprod}
Let $Y$ be another regular
flat separated scheme
of finite type over
$S$ and $E\subset Y$
be a divisor with simple
normal crossings.
Let $(Y\times_SY)^\sim$
be the log product
with respect to $E$
and we put $V=Y\setminus E$.
Let $f\colon Y\to X$
be a morphism over $S$
such that $f(V)\subset U
=X\setminus D$
and we consider the map
$(f\times f)^\sim\colon
(Y\times_S Y)^\sim\to
(X\times_S X)^\sim$
of log products.

Let  $A$ be 
a closed subscheme
of $(X\times_SX)^\sim$
satisfying the condition
{\rm (A$'$)} after Definition
{\rm \ref{dflip}}
and assume that
$A_Y=
(f\times f)^{\sim-1}(A)
\subset (Y\times_SY)^\sim$
also satisfies the corresponding condition
that $\delta_Y^{-1}(A_Y)\cap \Sigma_{Y/S}$
is supported on the
closed fiber $Y_F$
set-theoretically.
Let
$f^*\colon
G(\delta_X^{-1}(A)\cap
\Sigma_{X/S})\to
G(\delta_Y^{-1}(A_Y)\cap
\Sigma_{Y/S})$ be the
the pull-back by
$f\colon Y\to X$.

Then the pull-back
$(f\times f)^{\sim*}
\colon G(A)
\to G(A_Y)$
by 
$(f\times f)^\sim\colon
(Y\times_SY)^\sim\to
(X\times_S X)^\sim$
is defined. 
Further,
the diagram
$$\begin{CD}
G(A)@>{((\ ,\Delta_X^{\log}))}>>
G(\delta_X^{-1}(A)\cap
\Sigma_{X/S})\\
@V{(f\times f)^{\sim*}}VV 
@VV{f^*}V\\
G(A_Y)
@>{((\ ,\Delta_Y^{\log}))}>>
G(\delta_Y^{-1}(A_Y)\cap
\Sigma_{Y/S})
\end{CD}$$
is commutative.
\end{lm}

{\it Proof.}
By Lemma \ref{corYQ},
we have
$f^{-1}(\Sigma_{X/S})\subset
\Sigma_{Y/S}$.
Since
$f^{-1}(\delta_X^{-1}(A))
=\delta_Y^{-1}(A_Y)$,
the map
$G(\delta_X^{-1}(A)\cap
\Sigma_{X/S})\to
G(\delta_Y^{-1}(A_Y)\cap
\Sigma_{Y/S})$ 
is defined 
by the assumption (A$'$).

By Corollary \ref{cortdf},
the map
$(f\times f)^\sim\colon
(Y\times_S Y)^\sim\to
(X\times_S X)^\sim$
is of finite tor-dimension.
Hence,
the pull-back
$(f\times f)^{\sim*}
\colon G(A)
\to G(A_Y)$
by 
$(f\times f)^\sim\colon
(Y\times_S Y)^\sim\to
(X\times_S X)^\sim$
is defined. 

We apply the associativity,
Lemma \ref{lmfm2},
by taking 
$A\to (X\times_S X)^\sim
\gets X\gets Y$
as
$V\to X\gets W\gets W'$.
Then the composition
via upper right
is equal to
the map
$((\ ,\Delta_Y^{\log}))_
{(X\times_SX)^\sim}$.
We also apply
the associativity,
Lemma \ref{lmfm4},
by taking 
$A\to (X\times_S X)^\sim
\gets (Y\times_S Y)^\sim\gets Y$
as
$V\to X\gets X'\gets W'$.
Then,
the composition
via lower left
is also equal to
the same map.
\qed

We establish an important 
property that 
the localized intersection
product with the log diagonal
is independent of the boundary,
in Proposition \ref{prprod} below.
We begin with preliminary
computations.
Let $X_1$ be a regular
divisor of $X$
such that
the intersection $D_1=
X_1\cap D$
is a divisor of $X_1$
with simple normal crossings.
Let ${\cal D}=
(D_i)_{i\in I}$
be the family of
irreducible components
of $D$
and we consider
the family
${\cal D}_1=
(D_i\cap X_1)_{i\in I}$
of smooth divisors of $X_1$.
We identify the log product
$(X_1\times_SX_1)^\sim$
with respect to
${\cal D}_1$
with the inverse image
of $X_1\times_SX_1$
by the canonical map
$(X\times_SX)^\sim
\to X\times_SX$.
The sum
$D'=D\cup X_1$
is a divisor of $X$
with simple normal crossings.

We consider the log product
$(X\times_SX)^\approx$
with respect to $D'$,
the log diagonal map
$ \delta
\colon X\to
(X\times_SX)^\approx$
and the canonical map
$(X\times_SX)^\approx
\to
(X\times_SX)^\sim$.
The inverse image $E$
of $(X_1\times_SX_1)^\sim
\subset
(X\times_SX)^\sim$
by
$(X\times_SX)^\approx
\to
(X\times_SX)^\sim$
is a ${\mathbf G}_m$-torsor
over
$(X_1\times_SX_1)^\sim$
by Lemma \ref{lmEi}.
The pull-back of 
the ${\mathbf G}_m$-torsor $E$
by the log diagonal
$X_1\to 
(X_1\times_SX_1)^\sim$
is trivialized 
by
the restiction to $X_1$
of
the log diagonal 
$X\to (X\times_SX)^\approx$.
We identify $E\times_
{(X_1\times_SX_1)^\sim}X_1$
with ${\mathbf G}_{m,X_1}$
and the 
restriction of the log
diagonal
$X_1\to E\times_
{(X_1\times_SX_1)^\sim}X_1$
with the $1$-section
$1_{X_1}$.
They are summarized
in the cartesian diagram
\setcounter{equation}0
\addtocounter{thm}1
\begin{equation}
\begin{CD}
{\mathbf G}_{m,X_1}
@>>>
E
@>>>
(X\times_SX)^\approx\\
@VVV @VVV @VVV\\
X_1
@>>>
(X_1\times_SX_1)^\sim
@>>>
(X\times_SX)^\sim.
\end{CD}
\label{eqdX1}
\end{equation}

\addtocounter{thm}{-1}

\begin{lm}\label{lmX1p}
Let $X_1$ be a regular
irreducible divisor of $X$
such that
$D_1=X_1\cap D$
is a divisor of $X_1$
with simple normal crossings.
Let $(X\times_SX)^\approx$
denote the log product
with respect to
$D'=D\cup X_1$
and let $\Sigma'_{X/S}$
be the closed subscheme of $X$
defined by the annihilator
ideal ${\rm Ann}\
\Omega^n_{X/S}(\log D')$.

Let $A$ be 
a closed subscheme
of $E\subset (X\times_SX)^\approx$
satisfying the 
condition {\rm (A$'$)}
after Definition {\rm \ref{dflip}}.
We identify
$E\times_{
(X_1\times_SX_1)^\sim}
X_1$
with ${\mathbf G}_{m,X_1}$
and the section
$X_1\to
E\times_{
(X_1\times_SX_1)^\sim}X_1$
defined by the restriction
of the log diagonal
with $1_X\colon
X_1\to {\mathbf G}_{m,X_1}$.
Then, the intersection product
\begin{equation}
((\ ,\Delta^{\log}_X))_
{(X\times_SX)^\approx}
\colon
G(A)\to G(\delta^{-1}(A)\cap 
\Sigma'_{X/S})
\label{eqcpE}
\end{equation}
with the log diagonal
satisfies the following.

{\rm 1.}
Assume that $X_1$ is flat over $S$.
Then the composition
$A\to E\to 
(X_1\times_SX_1)^\sim$
satisfies the condition 
{\rm (A)} in Lemma {\rm \ref{lmW}}
and the localized intersection
product 
\begin{equation}
((\Delta^{\log}_{X_1},\ \ ))_
{(X_1\times_SX_1)^\sim}\colon
G(A)\to
G(A\cap 
{\mathbf G}_{m,\Sigma_{X_1/S}})
\label{eqA1}
\end{equation}
with the log diagonal
$X_1\to
(X_1\times_SX_1)^\sim$
is defined.
Further, 
we have
$\Sigma_{X_1/S}
=
\Sigma'_{X/S}\cap X_1$
and 
the map {\rm (\ref{eqcpE})}
is induced by the 
composition
\begin{equation}
\begin{CD}
G(A)@>{((\Delta^{\log}_{X_1},\ \ ))_
{(X_1\times_SX_1)^\sim}}>>
G(A\cap 
{\mathbf G}_{m,\Sigma_{X_1/S}})
@>{(\ ,1_{X_1})
_{{\mathbf G}_{m,X_1}}
}>>G(\delta^{-1}(A)\cap 
\Sigma_{X_1/S}).
\end{CD}
\label{eqX1p1}
\end{equation}

{\rm 2.}
Assume that
$X_1$ is a subscheme
of the closed fiber
$X_F$.
Then, we have $X_1\subset
\Sigma'_{X/S}$
and 
the map 
on the graded quotients
$((\ ,\Delta^{\log}_X))_
{(X\times_SX)^\approx}
\colon
{\rm Gr}^F_{\bullet}
G(A)\to 
{\rm Gr}^F_{\bullet-n}
G(\delta^{-1}(A)\cap 
\Sigma'_{X/S})$
induced by
{\rm (\ref{eqcpE})}
is induced by
the usual intersection product 
\begin{equation}
\begin{CD}
G(A)@>{(\ ,1_{X_1})_E}>>
G( \delta^{-1}(A)\cap X_1).
\end{CD}
\label{eqX1p2}
\end{equation}
\end{lm}

{\it Proof.}
1. 
By Lemma \ref{lmXX1},
we have an 
exact sequence
$0\to
\Omega^1_{X_1/S}(\log D_1)
\to \Omega^1_{X/S}(\log D')
\otimes {\mathcal O}_{X_1}
\to {\mathcal O}_{X_1}
\to 0$.
This defines an isomorphism
$\Omega^{n-1}_{X_1/S}(\log D_1)
\to \Omega^n_{X/S}(\log D')
\otimes {\mathcal O}_{X_1}$.
Hence, we have
$\Sigma_{X_1/S}
=
\Sigma'_{X/S}\cap X_1$.
Thus the condition 
(A$'$) for $A\to (X\times_SX)^\approx$
implies 
the condition 
(A) for $A\to (X_1\times_SX_1)^\sim$.

By Lemma \ref{lmX2}.1,
$E$ is a Cartier
divisor of $(X\times_SX)^\approx$
and we have
$(E,\Delta^{\log}_X)
_{(X\times_SX)^\approx}
=[\Delta^{\log}_{X_1}]$.
We apply
Lemma \ref{lmfm1}
to the diagram
$$\begin{CD}
A@<<< A\cap X_1 \\
@VVV @VVV\\
E@<<<X_1\\
@VVV @VVV\\
(X\times_SX)^\approx 
@<<<X
\end{CD}$$
Then, the map
(\ref{eqcpE})
is the same as
$((\ ,1_{X_1}))_E
\colon
G(A)\to G(\delta^{-1}(A)\cap 
\Sigma_{X_1/S})$.
Further 
we apply
Lemma \ref{lmfm2}
by taking the upper line
in the cartesian diagram
$$\begin{CD}
E@<<<
{\mathbf G}_{m,X_1}@<<< X_1
\\
@VVV @VVV@.\\
(X_1\times_SX_1)^\sim
@<<<X_1@.
\end{CD}$$
as $X\gets W\gets W'$
on the lower line
in the diagram of
Lemma \ref{lmfm2}.
Then the map
$((\ ,1_{X_1}))_E
\colon
G(A)\to G(\delta^{-1}(A)\cap 
\Sigma_{X_1/S})$
is equal to the composition
of
$((\ ,[{\mathbf G}_{m,X_1}]))_E
\colon
G(A)\to G(A\cap 
{\mathbf G}_{m,\Sigma_{X_1/S}})$
with the usual intersection
product
$(\ ,X_1)_{{\mathbf G}_{m,X_1}}\colon
G(A\cap 
{\mathbf G}_{m,\Sigma_{X_1/S}})
\to
G(\delta^{-1}(A)\cap 
\Sigma_{X_1/S})$.
Since $E$ is flat over
$(X_1\times_SX_1)^\sim$,
the first map
$((\ ,[{\mathbf G}_{m,X_1}]))_E
\colon
G(A)\to G(A\cap 
{\mathbf G}_{m,\Sigma_{X_1/S}})$
is the same as (\ref{eqA1})

2. 
We show $X_1\subset
\Sigma'_{X/S}$.
Since $X_1$ is assumed
irreducible,
it suffices to show
that the generic point $\xi_1$ of $X_1$
is contained in $\Sigma'_{X/S}$.
It suffices to consider
the case where
the complement
$X\setminus 
\Sigma'_{X/S}$ is dense.
Then, we have
an exact sequence
$0\to 
{\mathcal O}_{X_F}
\to 
\Omega^1_{X/S}(\log D')
\to 
\Omega^1_{X/S}(\log D'/\log F)
\to 0$
on a neighborhood $\xi_1$.
By the assumption
that the complement
$X\setminus 
\Sigma'_{X/S}$ is dense,
the free part
of the module
$\Omega^1_{X/S}(\log D')_{\xi_1}$
over the discrete valuation ring
${\mathcal O}_{X,\xi_1}$
has rank $n-1$.
Since it is not torsion free,
the smallest number
of generators is $n$.
Hence $\xi_1$
is contained in
$\Sigma'_{X/S}$ as required.

We apply Lemma \ref{lmfm3}
by taking
$E\to (X\times_SX)^\approx
\gets \Delta_X^{\log}$
as 
$V\to X\gets W$.
Then, the localized
intersection product
$((\ ,\Delta^{\log}_X))_
{(X\times_SX)^\approx}$
is equal to the
usual intersection product
$(\ ,((E,\Delta^{\log}_X))_
{(X\times_SX)^\approx})_E$.
By the excess intersection formula
(\ref{eqfm3}),
we have
$((E,\Delta^{\log}_X))_
{(X\times_SX)^\approx}
=c_0(M'_{X/(X\times_SX)^\approx,E})
\cap[X_1]$.
By Lemma \ref{lmX1}.1,
the right hand side
is equal to
$1_{X_1}$.
Hence, the assertion follows.
\qed

\begin{pr}\label{prprod}
\setcounter{equation}0
Let $X$ be a regular flat
separated scheme 
of finite type over $S$
and $D\subset X$
be a divisor with simple normal
crossings.
We put $U=X\setminus D$.
Let $A$ be a closed subscheme
of $(X\times_SX)^\sim$
satisfying the condition
{\rm (A$'$)} after Definition
{\rm \ref{dflip}}
and the following condition:
\begin{itemize}
\item[{\rm (B)}]
For each irreducible
component $D_i$ of $D$,
we regard ${\mathbf G}_{m,D_i}$
as a closed subscheme of
$(X\times_SX)^\sim$
as in {\rm (\ref{eqdX1})}.
Then,
there exists an integer $l_i\ge 1$
such that the intersection
$A\cap {\mathbf G}_{m,D_i}$
is supported on 
the subscheme
$\mu_{l_i,D_i}
\subset
{\mathbf G}_{m,D_i}$.
\end{itemize}
We put $A^\circ =
A\cap (U\times_SU)$.
Then, there exists a unique map
${\rm Gr}^F_{\bullet}
G(A^\circ)\to 
{\rm Gr}^F_{\bullet-n}
G(\delta^{-1}(A)\cap 
\Sigma_{X/S})$,
also denoted by
$((\ ,\Delta_X^{\log}))$,
that makes the diagram
\begin{equation}
\xymatrix{
{\rm Gr}^F_{\bullet}
G(A)\ar[rr]^{\!\!\!\!\!
\!\!\!\!\!\!\!\!\!\!\!\!\!\!\!\!\!\!\!\!\!\!\!\!\!((\ ,\Delta_X^{\log}))}
\ar[d]_{\rm restriction}&&
{\rm Gr}^F_{\bullet-n}
G(\delta^{-1}(A)\cap 
\Sigma_{X/S})\\
{\rm Gr}^F_{\bullet}
G(A^\circ)\ar[rru]
}
\label{eqprpr}
\end{equation}
commutative.
\end{pr}

\noindent{\it Proof.}
For each irreducible component
$D_i$ of $D$,
let $E_i$ be its inverse image
by either of the two projections
$(X\times_SX)^\sim \to X$ and
we put $A_i=A\cap E_i$.
By the exact sequence
$\bigoplus_i
{\rm Gr}^F_{\bullet}
G(A_i)
\to 
{\rm Gr}^F_{\bullet}
G(A)
\to 
{\rm Gr}^F_{\bullet}
G(A^\circ)\to 0,$
it suffices to show
that the composition of
\begin{equation}
\begin{CD}
{\rm Gr}^F_{\bullet}
G(A_i)
@>>> 
{\rm Gr}^F_{\bullet}
G(A)
@>{((\ ,\Delta_X^{\log}))}>> 
{\rm Gr}^F_{\bullet-n}
G(\delta^{-1}(A)\cap 
\Sigma_{X/S})
\end{CD}
\label{eqAi}
\end{equation}
is the zero map
for each $i$.

First, we consider the case 
where $D_i$ is flat over $S$.
By Lemma \ref{lmX1p}.1,
the composition of (\ref{eqAi})
is induced by the composition
of 
$$\begin{CD}
G(A_i)
@>{((\ ,\Delta^{\log}_{D_i}))_{
(D_i\times_SD_i)^\sim}}>> 
G(A_i\cap {\mathbf G}_{m,
\Sigma_{D_i/S}})
@>{(\ ,1_{D_i})_{{\mathbf G}_{m,D_i}}}>> 
G(\delta^{-1}(A_i)
\cap \Sigma_{D_i/S}).
\end{CD}$$
Hence, it suffices
to show the second map
is the zero-map.
By the assumption (B),
the intersection
$A_i\cap {\mathbf G}_{m,D_i}$
is a closed subscheme of
$\mu_{l_i,D_i}$.
Hence,
if the characteristic of $K$
is $p>0$,
the assertion follows from
Lemma \ref{lmvan} below.
If the characteristic of $K$
is 0,
the generic fiber
$\Sigma_{D_i/S}\times_S
{\rm Spec}\ K$ is empty.
Hence it
also follows from
Lemma \ref{lmvan} below.

Next, we consider the case 
where $D_i$ is a subscheme of the
closed fiber $X_F$.
In this case,
by Lemma \ref{lmX1p}.2,
the composition of (\ref{eqAi})
is induced by the composition
of 
$$\begin{CD}
G(A_i)
@>{(\ ,\Delta^{\log}_{D_i})_{
(D_i\times_FD_i)^\sim}}>> 
G(A_i\cap {\mathbf G}_{m,D_i})
@>{(\ ,1_{D_i})_{{\mathbf G}_{m,D_i}}}>> 
G(A_i\cap \Delta_{D_i}^{\log})
\end{CD}$$
Hence, it suffices
to show the second map
is the zero-map.
By the assumption (B),
the assertion 
in this case
is also reduced to the
following Lemma \ref{lmvan}.

\begin{lm}\label{lmvan}
Let $D$ be a 
noetherian scheme
over ${\mathbb F}_p$
and $l\ge 1$ be an integer.
Let $A$ be a closed subscheme
of $\mu_{l,D}
\subset {\mathbf G}_{m,D}$.
Then, the intersection product
$$\begin{CD}
(\ ,1_D)_{{\mathbf G}_{m,D}}\colon
G(A)@>>> 
G(A\cap 1_D)\end{CD}$$
with the unit section
$1_D\subset \mu_{l,D}$
is the zero map.
\end{lm}

\noindent{\it Proof.}
By replacing $l$ by its $p$-part $l'$,
we may assume that $l$
is a power of the characteristic $p>0$
of $F$
since $\mu_{l',D}$
is a closed and open subscheme of
$\mu_{l,D}$ and
has the same intersection
with the $1$-section.
Further, since the 
closed immersion
$A\cap 1_D\to A$
defined by a nilpotent ideal induces an isomorphism
$G(A\cap 1_D)\to G(A)$,
we may assume $A$ is a closed
subscheme of $1_D$.
For
a coherent
${\cal O}_{1_D}$-module ${\cal F}$,
we have
$$([{\cal F}],1_D)_{{\mathbf G}_{m,D}}
=[{\cal T}or^{{\cal O}_{{\mathbf G}_{m,D}}}_0
({\cal F},{\cal O}_{1_D})]
-[{\cal T}or^{{\cal O}_{{\mathbf G}_{m,D}}}_1
({\cal F},{\cal O}_{1_D})]
=[{\cal F}]-[{\cal F}]=0.$$
Hence the assertion follows.
\qed

The following Lemma,
analogous to
Lemma \ref{lmX1p},
will be used in
the proof of Proposition
\ref{prU1},
which in turn will be
used in the proof
of a blow-up formula
Proposition
\ref{prblup}.

\begin{lm}\label{lmX1pp}
Let $X_1$ be a regular
divisor of $X$
such that
$D_1=X_1\cap D$
is a divisor of $X_1$
with simple normal crossings.
Let $(X\times_SX)^\sim$
denote the log product
with respect to
$D$
and let $\Sigma_{X/S}$
be the closed subscheme of $X$
defined by the annihilator
ideal ${\rm Ann}\
\Omega^n_{X/S}(\log D)$.
We regard
$(X_1\times_SX_1)^\sim$
as a closed subscheme
of $(X\times_SX)^\sim$
as in {\rm (\ref{eqdX1})}.

Let $A$ be 
a closed subscheme
of $(X_1\times_SX_1)^\sim$
satisfying the 
condition {\rm (A$'$)}
after Definition {\rm \ref{dflip}}
and let
$i\colon
(X_1\times_SX_1)^\sim\to
(X\times_SX)^\sim$
be the closed immersion.
Then, 
the map on the graded pieces
\setcounter{equation}0
\begin{equation}
((\ ,\Delta^{\log}_X))_
{(X\times_SX)^\sim}
\colon
{\rm Gr}^F_\bullet
G(A)\to 
{\rm Gr}^F_{\bullet-n}
G(\delta^{-1}(A)\cap 
\Sigma_{X/S})
\label{eqcp1}
\end{equation}
induced by
the intersection product
with the log diagonal
is computed as follows.

{\rm 1.}
Assume $X_1$ is flat over $S$.
Then,
we have
$\Sigma_{X_1/S}
=
\Sigma_{X/S}\cap X_1$
and 
the map
{\rm (\ref{eqcp1})}
is the composition of
\begin{equation}
\begin{CD}
{\rm Gr}^F_\bullet
G(A)@>{((\ ,\Delta_{X_1}^{\log}))_
{(X_1\times_SX_1)^\sim}}>>
{\rm Gr}^F_{\bullet-(n-1)}
G(\delta^{-1}(A)\cap 
\Sigma_{X_1/S})\\
@>{-i_*\circ c_1(N_{X_1/X})}>>
{\rm Gr}^F_{\bullet-n}
G(\delta^{-1}(A)\cap 
\Sigma_{X_1/S}).
\end{CD}
\label{eqX1p3}
\end{equation}

{\rm 2.}
Assume that
$X_1$ is a subscheme
of the closed fiber
$X_F$.
Then, the map
on the graded pieces
induced by
{\rm (\ref{eqcp1})}
is the composition of
\begin{equation}
\begin{CD}
{\rm Gr}^F_\bullet
G(A)
@>{(\ ,\Delta_{X_1}^{\log})_
{(X_1\times_FX_1)^\sim}}>>
{\rm Gr}^F_{\bullet-(n-1)}
G(\delta^{-1}(A))\\
@>{-c_1([N_{s/S}
\otimes {\cal O}_{X_1}
\to N_{X_1/X}])_{
X_1\cap \Sigma_{X/S}}}>>
{\rm Gr}^F_{\bullet-n}
G(\delta^{-1}(A)\cap \Sigma_{X/S}).
\end{CD}
\label{eqX1p4}
\end{equation}
\end{lm}

{\it Proof.}
1.
We have shown the equality
$\Sigma_{X_1/S}
=
\Sigma_{X/S}\cap X_1$
at the beginning of
the proof of Lemma \ref{lmX1p}.1.
By Lemma \ref{lmX2}.1,
the immersion
$(X_1\times_SX_1)^\sim
\to
(X\times_SX)^\sim$ is a 
regular immersion of
codimension 2.
Since the excess conormal
sheaf
${\rm Ker}(
({\rm pr}_1^*N_{X_1/X}
\oplus
{\rm pr}_2^*N_{X_1/X})|_{X_1}
\to
N_{X_1/X})$
is isomorphic to
$N_{X_1/X}$,
we have
$((X_1\times_SX_1)^\sim,
\Delta_X^{\log})_
{(X\times_SX)^\sim}=
-c_1(N_{X_1/X})\cap
\Delta_{X_1}^{\log}$.
We apply
Lemma \ref{lmfm1}
by taking
$(X\times_SX)^\sim 
\supset (X_1\times_SX_1)^\sim,
X\supset X_1$
as $X\supset X',W\supset W'$.
Then, 
the map
{\rm (\ref{eqcp1})}
is induced by
$((\ ,-c_1(N_{X_1/X})
\cap \Delta_{X_1}^{\log}))_
{(X_1\times_SX_1)^\sim}
\colon
G(A)\to G(\delta^{-1}(A)\cap 
\Sigma_{X_1/S})$.
Further, it is equal
to the composition
of (\ref{eqX1p3}).

2. 
We apply Lemma \ref{lmfm3}
by taking
$(X_1\times_SX_1)^\sim
\to (X\times_SX)^\sim
\gets \Delta_X^{\log}$
as 
$V\to X\gets W$.
Then, 
the map
{\rm (\ref{eqcp1})}
is equal to the
usual intersection product
$$(\ ,(((X_1\times_SX_1)^\sim,
\Delta_X^{\log}))_
{(X\times_SX)^\sim})_{
(X_1\times_SX_1)^\sim}$$
with 
$(((X_1\times_SX_1)^\sim,
\Delta_X^{\log}))_
{(X\times_SX)^\sim}
=-c_1(M'_{X/(X\times_SX)^\sim,
(X_1\times_SX_1)^\sim})\cap [X_1]$.
By Lemma \ref{lmX1}.2,
the right hand side
is equal to
$-c_1([N_{s/S}\otimes{\cal O}_{X_1}
\to N_{X_1/X}])_{X_1\cap \Sigma_{X/S}}
\cap 
\Delta_{X_1}^{\log}$.
Hence, the assertion follows.
\qed

If $(X\setminus D)_F=\emptyset$,
we have an alternative
construction.
Let $A\subset 
(X\times_{\mathbb S}X)^\sim
\subset 
(X\times_SX)^\sim$
be a closed subscheme
satisfying the condition 
(A).
Then, a localized intersection
product
\begin{equation}
\begin{CD}
((\ ,\Delta_X^{\log}))_{
(X\times_{\mathbb S}X)^\sim}
\colon
G(A)@>>>
G(\delta^{-1}(A)\cap 
\Sigma_{X/S}).
\end{CD}
\label{eqprod2}
\end{equation}
is defined similarly
as (\ref{eqprod1}),
by Lemma \ref{lmXX}.
We show that it gives
the same invariants.

\begin{pr}\label{prbfS}
Let $X$ be a regular flat 
separated scheme
of finite type
over $S$ and $D\subset X$
be a divisor with simple normal crossings.
We assume 
$(X\setminus D)_F=\emptyset$.

Then,
for a closed subscheme
$A$ of $(X\times_{\mathbb S}X)^\sim$
satisfying the condition
{\rm (A)},
we have an equality
$$((\ ,\Delta^{\log}_X))
_{(X\times_SX)^\sim}=
((\ ,\Delta^{\log}_X))
_{(X\times_{\mathbb S}X)^\sim}$$
of maps
${\rm Gr}^F_\bullet
G(A)\to 
{\rm Gr}^F_{\bullet-n}
G(\delta^{-1}(A)\cap 
\Sigma_{X/S})$.
\end{pr}

{\it Proof.}
Let $W\subset A$
be an integral closed subscheme.
If $W$ is a subscheme
of $\Delta_X^{\log}$,
let $\pi\colon W'\to W$
denote the identity of $W$.
If not,
let $\pi\colon W'\to W$
be the blow-up at
$W'\cap \Delta_X^{\log}$.
Since $G(A)$
is generated by 
the classes $\pi_*[W']$
for integral
closed subschemes $W$ of $A$,
it suffices to show the equality
$((W' ,\Delta^{\log}_X))
_{(X\times_SX)^\sim}=
((W' ,\Delta^{\log}_X))
_{(X\times_{\mathbb S}X)^\sim}.$

We put $T'=W'\times_
{(X\times_{\mathbb S}X)^\sim}
\Delta_X^{\log},
Z=\Sigma_{X/S}$,
$d=\dim W_K+1$
and let $\varphi\colon
T'\to X$ denote the canonical map.
Then, by the excess intersection
formula \cite[Theorem 3.4.3]{KSI},
we have 
$$((W' ,\Delta^{\log}_X))
_{(X\times_SX)^\sim}=
\pi_*((-1)^d{c_d}^{T'}_{Z_{T'}}
(M'_{X/(X\times_SX)^\sim,W'})
\cap [T']),$$
$$((W' ,\Delta^{\log}_X))
_{(X\times_{\mathbb S}X)^\sim}=
\pi_*((-1)^d{c_d}^{T'}_{Z_{T'}}
(M'_{X/(X\times_{\mathbb S}X)^\sim,W'})
\cap [T']).$$
By the distinguished triangles
$$\to {\cal O}_{X_F}
\to 
L_{X/S}(\log D)
\to L_{X/S}(\log D/\log F)
\to 0,$$
$$\to M'_{X/(X\times_S X)^\sim,W'}\to 
L\varphi^*L_{X/S}(\log D)
\to
N_{T'/W'}\to,$$
$$\to M'_{X/(X\times_{\mathbb S} X)^\sim,W'}\to 
L\varphi^*L_{X/S}(\log D/\log F)
\to
N_{T'/W'}\to,$$
and by
$c_1({\cal O}_{X_F})\cap [T']=0$,
we obtain an equality
${c_d}^{T'}_{Z_{T'}}
(M'_{X/(X\times_SX)^\sim,W'})
\cap [T']
=
{c_d}^{T'}_{Z_{T'}}
(M'_{X/(X\times_{\mathbb S}X)^\sim,W'})
\cap [T']$.
Thus the assertion follows.
\qed

The following Proposition shows that
the localized intersection
product does not
depend on the choice of the base $S$.

\begin{pr}\label{pr4.3.a} 
Let $X$ be a a regular flat separated scheme of finite type over $S$ and $D\subset X$ be a divisor with normal crossings. 
Assume that $K$ is a finite extension 
of a complete discrete valuation field $K'$ and put
${\cal O}_{K'}=K'\cap {\cal O}_K$.  
Let $A$ be a closed subscheme of 
$(X\times_S X)^\sim$
satisfying the condition 
{\rm (A$'$)} after Definition
{\rm \ref{dflip}}
with respect to 
$S'={\rm Spec}\ {\cal O}_{K'}$.
Namely, we assume
that $\delta^{-1}(A) \cap \Sigma_{X/S'}$ is contained in $X_F$. 

Then, 
we have an inclusion
$\Sigma_{X/S}\subset \Sigma_{X/S'}$
and an equality
$$((\;,\Delta_X^{\log}))_{(X\times_{S'} X)^\sim}=
((\;,\Delta_X^{\log}))_{(X\times_S X)^\sim}$$
of maps 
${\rm Gr}^F_\bullet
G(A)\to 
{\rm Gr}^F_{\bullet-n}
G(\delta^{-1}(A)\cap \Sigma_{X/S'})$.
\end{pr}

{\it Proof}. The proof is similar to 
that of Proposition \ref{prbfS}. 
Let $W\subset A$ be an integral closed subscheme and let 
$\pi\colon W'\to W$ be as in the proof of Proposition \ref{prbfS}. 
We put $T'=
W'\times_{(X\times_S X)^\sim}\Delta_X^{\log}$, 
$Z'=\Sigma_{X/S'}$, and $d=\dim W_K+1$. 
Then, by the excess intersection formula \cite[Theorem 3.4.3]{KSI}, we have 
$$((W', \Delta^{\log}_X))_{(X\times_S X)^\sim}= \pi_*((-1)^dc^{T'}_{d,Z_{T'}}(M'_{X/(X\times_S X)^\sim,W'})\cap [T']),$$
$$((W', \Delta^{\log}_X))_{(X\times_{S'} X)^\sim}= \pi_*((-1)^dc^{T'}_{d,Z_{T'}}(M'_{X/(X\times_{S'} X)^\sim,W'}) \cap [T']).$$
Using the distinguished triangle
$\Omega^1_{S/S'}
\otimes^L{\cal O}_X\to 
L_{X/S'}(\log D)\to 
L_{X/S}(\log D)\to$,
we complete the proof
similarly as in the proof of
Proposition \ref{prbfS}. 
\qed

\newpage 
\section{Invariants of
wild ramification}\label{siw2}

We keep fixing
a complete discrete valuation
field $K$
with perfect residue
field $F$
of characteristic $>0$
and $S={\rm Spec}\ {\cal O}_K$.

In Section \ref{ssllpd},
we define invariants
of wild ramification
for
a finite \'etale morphism
$f\colon V\to U$
of regular flat separated
schemes of finite
type over $S$,
such that
the generic fiber
$V_K\to U_K$
is tamely ramified
with respect to
$K$ (Definition \ref{dftmT}).
The definition
uses the localized intersection product
with the log diagonal
constructed in 
Section \ref{sslipd}.
The definition is
extended to cover 
the case where $U$ and $V$
are not assumed regular
at the end of
Section \ref{ssexc}
as a consequence
of the excision formula,
Theorem \ref{thmexc}.
On the counterpart 
for a finite \'etale morphism
$f\colon V\to U$
of smooth separated
schemes of finite
type over $F$ defined
in \cite{KSA},
we also state
some complements.
In Section \ref{ssiw},
we establish
elementary properties
of the invariants
of wild ramification
defined in Section \ref{ssllpd}.
We define
the logarithmic different
and the Lefschetz classes
and derive
their basic properties
analogous to the classical
ones.

Before defining
the invariants in
the general case, 
we define and compute
the logarithmic different
and the Lefschetz class
using regular schemes
containing $U$ and $V$
as the complements
of divisors with
simple normal crossings
in Section \ref{ssdiff}.
We introduce
the target groups
where the invariants
of wild ramification
take values
as certain projective
limits with respect to
the system of compactifications
in Section \ref{sstg}.
We also introduce
in Theorem \ref{thmmap0}
a variant
that will be used
in the case where $K$ is
of characteristic $0$,
in Section \ref{sspc}.

\subsection{Logarithmic
different and
the Lefschetz class}\label{ssdiff}

Let $Y$ be a
regular flat 
separated scheme
of finite type
over $S={\rm Spec}\ {\cal O}_K$ 
and $V\subset Y$
be the complement
of a divisor $E$
with simple normal crossings.
Let $f\colon V\to U$
be a finite \'etale
morphism
of separated schemes
of finite type over $S$.
We consider the 
family ${\cal E}
=(E_i)_{i\in I}$
of irreducible
components of $E$
and we assume that
the closed subscheme
$\Sigma_{V/U}^{\cal E}
Y\subset Y$ 
(Definition \ref{dfSig}.1)
is supported
on the closed fiber.

By the assumption
that $V\to U$
is finite \'etale,
the diagonal
$\Delta_V$ is
an open and closed subscheme
of $V\times_UV$.
The closure $A$
of $(V\times_UV)
\setminus \Delta_V$
in the log product
$(Y\times_SY)^\sim
=(Y\times_SY)^\sim_{\cal E}$
satisfies the 
condition (A$'$)
after Definition \ref{dflip}
since 
$\Sigma_{V/U}^{\cal E}$
is the inverse image
$\delta^{-1}(A)$
by the log diagonal
$\delta\colon Y\to (Y\times_SY)^\sim$.

We also assume that
there exists a separated
scheme
$X$ of finite type over $S$
containing $U$
as the complement
of a Cartier divisor
and that $f\colon V\to U$
is extended to
a morphism $\bar f
\colon Y\to X$
over $S$ satisfying
$\bar f^{-1}(U)=V$.
Then, by Lemma \ref{lmEi}.2,
the same $A$
satisfies the condition (B)
in Proposition \ref{prprod}.
Thus, by applying 
the map
(\ref{eqprpr}),
we obtain
$$(((V\times_UV)
\setminus \Delta_V,
\Delta_Y^{\log}))_
{(Y\times_SY)^\sim}
\in F_0G(\Sigma_{V/U}^{\cal E}
Y).$$

\begin{df}\label{dfD}
Let $Y$ be a
regular flat
separated scheme
of finite type
over $S$
and $V\subset Y$
be the complement
of a divisor $E$ of $Y$
with simple normal crossings.
Let $f\colon V\to U$
be a finite \'etale
morphism
of separated schemes
of finite type over $S$.

We assume that
the closed subset
$\Sigma_{V/U}^{\cal E}
Y$ (Definition 
{\rm \ref{dfSig}.1})
defined for the 
family ${\cal E}
=(E_i)_{i\in I}$
of irreducible component of $E$
is supported
on the closed fiber.
We also assume that
there exists a separated
scheme
$X$ of finite type over $S$
containing $U$
as the complement
of a Cartier divisor
and that $f\colon V\to U$
is extended to
a morphism $\bar f
\colon Y\to X$
over $S$ satisfying
$\bar f^{-1}(U)=V$.

Then,
we define the 
{\rm the logarithmic different}
$D_{V/U,Y}^{\log}
\in F_0G(\Sigma_{V/U}^{\cal E}Y)$
by 
$$
D_{V/U,Y}^{\log}=
((V\times_UV\setminus
\Delta_V,\Delta_Y^{\log}))_
{(Y\times_SY)^\sim}.
$$
\end{df}

We compute the 
logarithmic different
explicitly using
regular models.
It will imply in 
particular 
(Corollary \ref{correxX}) 
that
if $U={\rm Spec}\ L$
and $V={\rm Spec}\ M$
for finite separable extensions
$L\subset M$ of $K$,
we have
\setcounter{equation}0
\begin{equation}
D_{V/U,Y}^{\log}
=
{\rm length}_{
{\cal O}_M}
\Omega^1_{{\cal O}_M/
{\cal O}_L}(\log/\log)
=
{\rm length}_{
{\cal O}_M}
\Omega^1_{{\cal O}_M/
{\cal O}_L}
-(e_{M/L}-1)
\label{eqDML}
\end{equation}
in 
${\mathbb Z}
=F_0G({\rm Spec}\
{\cal O}_M/{\mathfrak m}_M)$.
Recall that
$
{\rm length}_{
{\cal O}_M}
\Omega^1_{{\cal O}_M/
{\cal O}_L}$
is the classical different.

We consider a cartesian
diagram
\begin{equation}
\begin{CD}
V@>>> Y\\
@VfVV @VV{\bar f}V\\
U@>>> X
\end{CD}
\label{eqUVXY}
\end{equation}
of regular flat separated
schemes
of finite type over $S$.
Suppose that
$f\colon V\to U$
is finite \'etale and
that
$U=X\setminus D$ 
and $V=Y\setminus E$
are the complements
of divisors $D$ and $E$
with simple normal
crossings respectively.
Using the diagram
(\ref{eqUVXY}),
the logarithmic different
$D_{V/U,Y}^{\log}$
can be computed as follows.

We put $n=
\dim X_K+1=\dim Y_K+1$.
We consider 
the map
$\bar f^*
\Omega^1_{X/S}(\log D)
\to
\Omega^1_{Y/S}(\log E)$
of coherent 
${\cal O}_Y$-modules.
Let $\Sigma=\Sigma_{Y/X}
\subset Y$
be the closed subscheme
defined by
the annihilator
${\cal I}_\Sigma=
{\rm Ann}
({\rm Coker}
(\bar f^*
\Omega^1_{X/S}(\log D)
\to
\Omega^1_{Y/S}(\log E)))
\subset {\cal O}_Y.$
Since $Y$
is regular,
there exist a locally free
${\cal O}_Y$-module
${\cal V}$
and a surjection
${\cal V}\to 
\Omega^1_{Y/S}(\log E)$
by \cite[Corollaire 2.2.7.1]{SGA6}.
Hence, the
localized Chern class
$c_n(\Omega^1_{Y/S}(\log E)
-\bar f^*\Omega^1_{X/S}(\log D))_\Sigma
\cap [Y]
\in F_0G(\Sigma_{Y/X})$
is defined
in \cite[(3.24)]{KSA}
(cf.\ (\ref{eqck})).

The image of the
logarithmic
different
$D_{V/U,Y}^{\log}
\in F_0G(\Sigma_{V/U}^{\cal E}Y)$
in $F_0G(\Sigma_{Y/X})$
is computed using
the localized Chern class
$c_n(\Omega^1_{Y/S}(\log E)
-\bar f^*\Omega^1_{X/S}(\log D)
)_{\Sigma_{Y/X}}
\cap [Y]$ as follows.

\begin{pr}\label{prrexX}
Let $X$ and $Y$ be regular
flat separated schemes 
of finite type over $S$
and let $U=X\setminus D$ and
$V=Y\setminus E$ be
the complements of
divisors $D\subset X$
and $E\subset Y$
with simple normal
crossings.
Let $\bar f\colon
Y\to X$ be a morphism
over $S$
such that
$\bar f^{-1}(U)=V$
and the restriction
$f=\bar f|_V\colon
V\to U$
is finite \'etale.

We assume that
the support
$\Sigma=
\Sigma_{Y/X}\subset Y$
of the cokernel
$\Omega^1_Y(\log E)/
\bar f^*\Omega^1_X(\log D)$
is supported
on the closed fiber
$Y_F$.
We also assume that
there exists a
dense open subscheme
of $X$ smooth over $S$.

Then, we have
$\Sigma_{V/U}^{\cal E}Y
\subset \Sigma_{Y/X}$
and, for $n=\dim X_K+1$,
\setcounter{equation}0
\begin{equation}
D_{V/U,Y}^{\log}
=
(-1)^n
c_n(\Omega^1_{Y/S}(\log E)
-\bar f^*\Omega^1_{X/S}(\log D))_
\Sigma\cap ([Y])
\label{eqrexX}
\end{equation}
in $F_0G(\Sigma_{Y/X})$.
\end{pr}

{\it Proof.}
We consider the log products
$(X\times_SX)^\sim$
and 
$(Y\times_SY)^\sim$
with respect to $D$ and
$E$ respectively
and will apply
Proposition 
\ref{prrex}
to the commutative diagram
\begin{equation}
\begin{CD}
Y@>>> Y\\
@VVV @VVV\\
Y@>{\delta_Y}>> (Y\times_SY)^\sim\\
@V{\bar f}VV @VV{
(\bar f\times \bar f)^\sim}V\\
X@>{\delta_X}>> (X\times_SX)^\sim
\end{CD}
\label{eqXYs}
\end{equation}
where the upper square 
and the tall rectangle 
are cartesian.
We put 
$(Y\times_XY)^\sim
=
(Y\times_SY)^\sim
\times_{(X\times_SX)^\sim}
\Delta_X^{\log}$.
Since
the cokernel 
${\rm Coker}
(\bar f^*
\Omega^1_{X/S}(\log D))
\to
\Omega^1_{Y/S}(\log E))$
is the conormal
sheaf
$N_{Y/(Y\times_XY)^\sim}$,
the restriction
of the log diagonal map
$\delta_Y\colon
Y\to 
(Y\times_XY)^\sim$
to the complement
$\widetilde V
=Y\setminus \Sigma$
is an open immersion.
Hence the complement
$A=
(Y\times_XY)^\sim
\setminus
\Delta_{\widetilde V}
^{\log}$
is a closed subset
of $(Y\times_SY)^\sim$
such that $\delta^{-1}(A)=
\Sigma_{Y/X}$.
Since $A$ contains
$(V\times_UV)\setminus
\Delta_V$
as a subset,
we have an inclusion
$\Sigma_{V/U}^{\cal E}Y
\subset \Sigma_{Y/X}$.

We define
a bounded complex ${\cal C}$
of ${\cal O}_{(Y\times_SY)^\sim}$-modules
fitting in 
the distinguished triangle
$\to
{\cal C}
\to 
L(\bar f\times \bar f)^{\sim*}
{\cal O}_{\Delta_X^{\log}}
\to
{\cal O}_{\Delta_Y^{\log}}
\to $
as in (\ref{eqCV}).
We have $A\cap
(V\times_SV)=
(V\times_UV)\setminus \Delta_V$
and 
the restriction map
$F_nG(A)
\to 
F_nG((V\times_UV)\setminus \Delta_V)$
sends
$[{\cal C}]$
to
$[(V\times_UV)\setminus \Delta_V]$.
Hence, by Proposition
\ref{prprod},
the image of the
logarithmic different
$D_{V/U,Y}^{\log}$
by the map
$F_0G(\Sigma_{V/U}^{\cal E}Y)
\to
F_0G(\Sigma_{Y/X})$
is 
the localized
intersection product
$(([{\cal C}],
\Delta_Y^{\log}))_{
(Y\times_SY)^\sim}.$

In order to apply
Proposition
\ref{prrex} to the diagram
(\ref{eqXYs}),
we check that
its assumption is satisfied.
For a point $y$
of the closed fiber of $Y$,
we have an
open neighborhood $V'$ of $y$,
an open neighborhood $U'$ of 
$\bar f(y)$
and a cartesian diagram
$$\begin{CD}
V'@>>> Q\\
@VVV @VVV\\
U'@>>> P
\end{CD}$$
as in Lemma \ref{lmYQ}.
Then, in the diagram
$$\begin{CD}
V'@>>> (V'\times_SY)^\sim 
@>>> (Q\times_SY)^\sim \\
@VVV @VVV @VVV\\
U'@>>> (U'\times_SX)^\sim 
@>>> (P\times_SX)^\sim,
\end{CD}$$
the right square is cartesian
and 
the horizontal arrows
in the right square
are regular immersions
of codimension 1.
The compositions
of the horizontal arrows
are both sections
of smooth morphisms
of relative dimension $n$
and hence are
regular immersions
of codimension $n$.
Thus 
the condition
(\ref{prrex}.3)
is satisfied.

By the assumption
that there exists a
dense open subscheme
of $X$ smooth over $S$,
the excess conormal complexes
$M'_{Y/(Y\times_SY)^\sim,Y}=
M_{Y/(Y\times_SY)^\sim}=
L_{Y/S}(\log E)$
and 
$M'_{X/(X\times_SX)^\sim,Y}
=Lf^*M_{X/(X\times_SX)^\sim}=
Lf^*L_{X/S}(\log D)$
are quasi-isomorphic to
$\Omega^1_{Y/S}(\log E)$
and 
$f^*\Omega^1_{X/S}(\log D)$
respectively.
Applying Proposition \ref{prrex}.3
to the diagram
$$\begin{CD}
Y@>>>(Y\times_SY)^\sim\\
@VVV @VVV\\
X@>>>(X\times_SX)^\sim
\end{CD}$$
and to $T=W=Y$,
we obtain
$$(((f\times f)^*[\Delta_U]
-[\Delta_V],\Delta_Y^{\log}))_{
(Y\times_SY)^\sim}
=
[L\Lambda^nf^*
\Omega^1_{X/S}(\log D)
\to
L\Lambda^n
\Omega^1_{Y/S}(\log E)]
$$
in $F_0G({\Sigma}_{Y/X})$.
The right hand side
is equal to the image
of the localized
Chern class
$c_n(\Omega^1_{Y/S}(\log E)-
f^*\Omega^1_{X/S}(\log D))_{\Sigma}
\cap ([Y])$
in $F_0G(\Sigma_{Y/X})$
by Proposition \ref{prLL}.
\qed

\begin{cor}\label{correxX}
Let $L\subset M$ be
finite separable
extensions of $K$
and put $U={\rm Spec}\ L$
and $V={\rm Spec}\ M$.
Let $X$ and $Y$ be
the normalizations of $S$
in $U$ and $V$ respectively.
Then,
we have
$D_{V/U,Y}^{\log}
={\rm length}_{
{\cal O}_M}
\Omega^1_{{\cal O}_M/
{\cal O}_L}(\log/\log)
\in 
{\mathbb Z}$
{\rm (\ref{eqDML})}.
\end{cor}

For an automorphism
of a scheme
over $S$,
we define the Lefschetz
class as
the intersection
product of
the graph with the
log diagonal as follows.

\begin{df}\label{dfadm}
Let $X$ be a regular
flat separated scheme
of finite type over $S$
and $U=X\setminus D$
be the complement
of a divisor $D$
with simple normal crossings.
Let $\sigma$
be an automorphism
of $X$ over $S$
such that
$\sigma(U)=U$
and $U^\sigma=\emptyset$.

{\rm 1.}
Let $\widetilde
\Gamma_\sigma
\subset
(X\times_SX)^\sim$
be the schematic 
closure of
the graph
$\Gamma_\sigma
\subset
U\times_SU$
of the restriction
of $\sigma$.
We define
{\rm the logarithmic fixed part}
$X_{\log}^\sigma
\subset X$ by
$$X_{\log}^\sigma=
\Delta_X^{\log}
\times_{(X\times_SX)^\sim}
\widetilde \Gamma_\sigma.$$
We assume that
the intersection
$X_{\log}^\sigma
\cap \Sigma_{X/S}$
with the support $\Sigma_{X/S}$
of $\Omega^n_{X/S}(\log D)$
is supported in the closed fiber
set-theoretically.
We call 
the logarithmic intersection
product
$$((\widetilde
\Gamma_\sigma,
\Delta_X^{\log}))
_{(X\times_SX)^\sim}
\in F_0G(X_{\log}^\sigma
\cap \Sigma_{X/S})$$
{\rm the logarithmic Lefschetz class}.

{\rm 2.}
We say $\sigma$ is {\rm admissible}
if the following condition 
is satisfied:
For each irreducible component $D_i$
of $D$,
we have either $\sigma(D_i)=D_i$
or $\sigma(D_i)\cap D_i=\emptyset$.
\end{df}

We compute
the logarithmic Lefschetz class
using the Segre classes
\cite[4.2]{fulton},
under a slightly
weaker assumption
than in \cite[Lemma 5.4.8]{KSI}.

\begin{lm}\label{lmsig}
Let $X$ be a regular
flat separated scheme
of finite type over $S$
and $U=X\setminus D$
be the complement
of a divisor $D$
with simple normal crossings.
Let $\sigma$ 
be an automorphism
of $X$ over $S$
such that
$\sigma(U)=U$
and $U^\sigma=\emptyset$.
We assume $\sigma$
is admissible.
Let $D_1,\ldots,D_m$
be the irreducible components
of $D$ and
we put $\widetilde U=
X\setminus
\bigcup_{i:
\sigma(D_i)\cap D_i=\emptyset}D_i$.
Assume further
that there exists a dense open
scheme of $X$ smooth over $S$.
Then,

{\rm 1.}
Let 
$\gamma_{\widetilde U}:
\widetilde U\to X\times_SX$
be the restriction
of $\gamma
=(1,\sigma):
X\to X\times_SX$.
Then it induces a closed
immersion
$\tilde \gamma\colon
\widetilde U\to
(X\times_SX)^\sim$.
The image
$\widetilde 
\Gamma_\sigma=
\tilde \gamma(
\widetilde U)\subset
(X\times_SX)^\sim$
is the schematic closure
of $\Gamma_\sigma\subset
U\times_SU$.

{\rm 2.}
Assume that 
the generic fiber 
$X_{\log,K}^\sigma$
is empty
and let
$s(X_{\log}^\sigma,X)$
denote the Segre
class of
$X_{\log}^\sigma
\subset X$.
We put $n=\dim X_K+1$.
Then the log Lefschetz class
$((\widetilde \Gamma_\sigma,
\Delta_X^{\log}))
_{(X\times_SX)^\sim}
\in F_0G(X_{\log}^\sigma)$ 
is equal to the image of
$$\{c(\Omega^1_{X/S}(\log D))^*\cap 
s(X_{\log}^\sigma,X)\}_{\dim 0}=
\sum_{i=0}^{n}(-1)^i
c_i(\Omega^1_{X/S}(\log D))
s_{n-i}(X_{\log}^\sigma,X).$$
In particular,
if the logarithmic fixed part
$X_{\log}^\sigma$
is a Cartier divisor of $X$,
we have
$$((\Gamma_\sigma,
\Delta_X^{\log}))
_{(X\times_SX)^\sim}
=
\{c(\Omega^1_{X/S}(\log D))^*\cap 
(1+X_{\log}^\sigma)^{-1}\cap 
[X_{\log}^\sigma]\}_{\dim 0}.$$
\end{lm}

{\it Proof.}
1. We set $(X\times_SX)^0
=X\times_SX-\bigcup_{(i,j):D_i\cap D_j=\emptyset}
D_i\times D_j$.
By the definition of
$(X\times_SX)^\sim$,
we have ${\rm pr}_1^{-1}(D_i)
={\rm pr}_2^{-1}(D_i)$
in $(X\times_SX)^\sim$
for every irreducible
component $D_i$ of $D$.
Hence
$(X\times_SX)^\sim$
is a scheme over 
$(X\times_SX)^0$.
By the definition of 
$\widetilde U$,
it is the inverse image of
$(X\times_SX)^0
\subset X\times_SX$
by the map $\gamma:
X\to X\times_SX$.
Hence its restriction
$\gamma_{\widetilde U}\colon
\widetilde U\to 
(X\times_SX)^0$ 
is a closed immersion.

By the assumption
that $\sigma$ is admissible,
the map
$\gamma_{\widetilde U}\colon
\widetilde U\to 
(X\times_SX)^0$ 
induces a map 
$\tilde\gamma\colon
\widetilde U
\to (X\times_SX)^\sim$.
Since 
$\gamma_{\widetilde U}\colon
\widetilde U\to 
(X\times_SX)^0$ is a closed immersion,
the induced map
$\tilde\gamma\colon
\widetilde U
\to (X\times_SX)^\sim$
is also a closed immersion.

2. By the assumption that
$X_{\log,K}^\sigma$ 
is empty,
the underlying set 
of $X^\sigma_{\log}$
is a subset of the 
closed fiber $X_F$.
We apply \cite[Corollary 3.4.6]{KSI},
by taking $X\to (X\times_SX)^\sim\to X$
to be $V\to X\to S$ 
and $X^\sigma_{\log}\to 
\widetilde\Gamma_\sigma\to 
(X\times_SX)^\sim$ to be $T\to W\to X$
in \cite[Corollary 3.4.6]{KSI}.
Since $M_{X/(X\times_SX)^\sim}=
\Omega^1_{X/S}(\log D)$, we obtain
$((X,\widetilde
\Gamma_\sigma))_{(X\times_SX)^\sim}=
\{c(\Omega^1_{X/S}(\log D))^*\cap 
s(X_{\log}^\sigma,
\widetilde\Gamma_\sigma)\}_{\dim 0}$.
Since the open immersion
$\widetilde \Gamma_\sigma\to
X$ induces the identity on
$X_{\sigma}^{\log}$,
we have $s(X_{\log}^\sigma,\Gamma_\sigma)=
s(X_{\log}^\sigma,
\Delta_{\widetilde U})=
s(X_{\log}^\sigma,X)$.
Thus the assertion is proved.
\qed

\begin{cor}\label{corsi}
Assume further
that $\sigma$
is of finite order
and let $i$
be an integer prime
to the order of $\sigma$.
Then,
we have
$((\Gamma_\sigma,
\Delta_V^{\log}))=
((\Gamma_{\sigma^i},
\Delta_V^{\log}))$.
\end{cor}

{\it Proof.}
Since $X^{\sigma}_{\log}=
X^{\sigma^i}_{\log}$,
the assertion follows
from Lemma 
\ref{lmsig}.2.
\qed

For isolated singular points,
we have the following formula
similarly as
\cite[Lemma 3.4.14]{KSA}.

\begin{pr}\label{priso}
\setcounter{equation}0
Let $X$ be
a regular flat separated
scheme
of finite type
over $S$
and $\sigma$
be an automorphism
of $X$ over $S$.
Assume that
there exists a dense open
subscheme of $X$ smooth over $S$.
Let $x\in X$
be a closed point
in the closed fiber
and assume that
the fixed part
$X^{\sigma}$
is set-theoretically
equal to the set
$\{x\}$.

Let $f\colon
X'\to X$
be the blow-up at $x$
and $D$ be the exceptional
divisor.
Let $(X'\times_SX')^\sim$
denote the log product
with respect to $D$
and $\widetilde \Gamma_\sigma
\subset (X'\times_SX')^\sim$
be the proper transform
of the graph
$\Gamma_\sigma
\subset X\times_SX$
of $\sigma$.
Then, we have
\begin{equation}
f_*((\widetilde \Gamma_\sigma,
\Delta_{X'}))_{
(X'\times_SX')^\sim}
=[{\cal O}_{X^\sigma}]
-[x]
\label{eqiso}
\end{equation}
where
$[{\cal O}_{X^\sigma}]
={\rm length}\
{\cal O}_{X^\sigma}
\cdot[x].$
\end{pr}

{\it Proof.}
We apply Lemma \ref{lmsig}
to the automorphism
$\sigma$ of $X'$
admissible with respect 
to the exceptional divisor $D$.
By the exact sequence
$0\to 
\Omega^1_{D/F}
\to
\Omega^1_{X'/S}(\log D)
\otimes{\mathcal O}_D
\overset{\rm res}
\to {\mathcal O}_D\to 0$,
the total Chern class satisfies
$c(\Omega^1_{X'/S}(\log D))
=
c(\Omega^1_{D/F})
=
c({\mathcal O}(-1))^n
=(1-H)^n$ on $D$
where $H$ denote the class
of the hyperplane
of the projective space $D$.
Hence, by Lemma \ref{lmsig},
we obtain
$$((\widetilde \Gamma_\sigma,
\Delta_{X'}))_{
(X'\times_SX')^\sim}
=
\{(1+H)^ns(X_{\log}^{\prime \sigma},
X')\}_{\dim 0}.$$

Let ${\cal I}\subset {\mathcal O}_X$ 
and ${\cal J}\subset {\mathcal O}_{X'}$ 
denote
the ideal sheaf of $X^\sigma$
and the ideal sheaf of 
$X_{\log}^{\prime \sigma}$
respectively.
Then,
since ${\cal I}$
is generated by
$\sigma(t_i)-t_i$
for a minimal system $(t_i)$
of generators of
the maximal ideal ${\mathfrak m}_x$,
we compute
$f^* {\cal I}
={\cal J}\cdot {\cal I}_D$.
This means that
$X_{\log}^{\prime \sigma}$
is the residual scheme
\cite[Definition 9.2.1]{fulton}
to the Cartier divisor 
$D$ in the inverse image $f^*(X^\sigma)$.
Hence by \cite[Proposition 9.2]{fulton}, 
it implies that the Segre class satisfies
$$s(f^*(X^\sigma),X')_{\dim 0}
=
H^{n-1}+
\{(1+H)^ns(X_{\log}^{\prime \sigma},
X')\}_{\dim 0}$$
since the self intersection
$D\cdot D$ is $-H$.
Thus, we obtain
$$f_*((\widetilde \Gamma_\sigma,
\Delta_{X'}))_{
(X'\times_SX')^\sim}
=
f_*s(f^*(X^\sigma),X')_{\dim 0}
-
f_*H^{n-1}.$$
By
$f_*s(f^*(X^\sigma),X')_{\dim 0}
=
s(X^\sigma,X)_{\dim 0}
=[X^\sigma]$
(\cite[Proposition 4.2(2)]{fulton})
and
$f_*H^{n-1}=[x]$,
the assertion follows.
\qed

In the case
$X={\rm Spec}\
{\mathcal O}_L$
for a finite separable
extension $L$ of $K$,
we obtain the following.

\begin{cor}\label{corssw}
Let $L$ be
a finite separable
extension of $K$
and $\sigma$
be a non-trivial 
automorphism
of $L$ over $K$.
We put $X={\rm Spec}\
{\cal O}_L$
and let
$J_\sigma
\subset 
{\cal O}_L$
be the ideal
generated by
$\sigma(a)-a$
for $a\in {\cal O}_L$
and $\sigma(b)/b-1$
for $b\in {\cal O}_L$
and $b\neq 0$.
Then, we have
$$((\Gamma_\sigma,
\Delta_X^{\log}))_
{(X\times_SX)^\sim}
={\rm length}_{{\cal O}_L}
{\cal O}_L/
J_\sigma.$$
\end{cor}

\subsection{The target groups}\label{sstg}
\setcounter{equation}0
Let $f\colon V\to U$ 
be a finite \'etale morphism
of separated schemes
of finite type over $S={\rm Spec}\ {\cal O}_K$.
In this subsection,
we define
an abelian group
$F_0G(\partial_{V/U} W)$
and
a ${\mathbb Q}$-vector space
$F_0G(\partial_{V/U} W)_{\mathbb Q}$
for a separated scheme $W$
of finite type over $V$.
Assuming $U$ and $V$ are
regular,
for a finite \'etale 
scheme $V'$ over $V$,
invariants of
wild ramification
of $V'\to U$
will be defined
as elements
of the group
$F_0G(\partial_{V/U} V')_{\mathbb Q}$
in Section \ref{ssiw}.
Without assuming
the regularity
of $U$ and $V$,
the definition is
extended at the end of
Section \ref{ssexc}
as a consequence
of the excision formula,
Theorem \ref{thmexc}.

Let $f\colon V\to U$ 
be a morphism
of separated schemes
of finite type over $S$.
Recall that
an open immersion
$j\colon V\to Y$
is schematically dense
if the canonical map
${\mathcal O}_Y\to
j_*{\mathcal O}_V$
is an injection.
We define a category 
${\cal C}_{V\to U}$
of compactifications of
$f\colon V\to U$
as follows:
\begin{itemize}
\item
An object is
a morphism $\bar f\colon
Y\to X$
of proper schemes over $S$
such that $X$ and $Y$  
contain $U$ and $V$ respectively
as schematically dense
open subschemes
and that the diagram
\begin{equation}
\begin{CD}
V@>>> Y\\
@VfVV @VV{\bar f}V\\
U@>>> X
\end{CD}
\label{eqtg}
\end{equation}
is commutative.

\item
A morphism $(g,h)
\colon (\bar f'\colon
Y'\to X')
\to (\bar f
\colon Y\to X)$
is a pair of morphisms 
$g\colon X'\to X$
and 
$h\colon Y'\to Y$
of schemes
over $S$ extending the 
identities of $U$ and of $V$
such that
the diagram
\begin{equation}
\begin{CD}
Y'@>{f'}>>X'\\
@VhVV @VVgV\\
Y@>f>>X\\
\end{CD}
\label{eqtg2}
\end{equation}
is commutative.
\end{itemize}

\begin{lm}\label{lmCU}
Let $f\colon V\to U$ 
be a morphism
of separated schemes 
of finite type over $S$.

{\rm 1.}
The category 
${\cal C}_{V\to U}$
is cofiltered
and in particular non-empty.

{\rm 2.}
If $f\colon V\to U$
is finite flat,
then
the objects
$\bar f\colon
Y\to X$
such that
the diagram
{\rm (\ref{eqtg})}
is cartesian and that
$\bar f$ are finite flat
are cofinal
in the category 
${\cal C}_{V\to U}$.

{\rm 3.}
If $V$ is a $G$-torsor
over $U$ for a finite group $G$,
the objects
$\bar f\colon
Y\to X$
such that
the diagram
{\rm (\ref{eqtg})}
is cartesian, that
$\bar f$ are finite flat
and that the action of $G$
is extended to an action
on $Y$ over $X$
are cofinal
in the category 
${\cal C}_{V\to U}$.
\end{lm}

{\it Proof.}
1.
By Nagata's embedding
theorem \cite{nagata}, 
there exists a proper scheme
$X$ over $S$ containing $U$
as an open subscheme.
After replacing $X$
by the schematic closure
of $U$,
the open subscheme
$U$ is schematically
dense in $X$.
Further 
by Nagata's embedding
theorem \cite{nagata}, 
there exists a proper scheme
$Y$ over $X$ containing $V$
as an open subscheme.
After replacing $Y$
by the schematic closure of $V$
as above,
we obtain an object 
$Y\to X$
of ${\cal C}_{V\to U}$.

Let 
$Y\to X$ and $Y'\to X'$
be objects of ${\cal C}_{V\to U}$.
If there exists a map
$(Y\to X)\to (Y'\to X')$
of ${\cal C}_{V\to U}$,
it is unique since $V$ is
assumed schematically dense.
Let $X''$ 
be the schematic closure
of $U$ in
$X\times_SX'$
and $Y''$
be the schematic closure
of $V$ in
$X''\times_{X\times_SX'}
(Y\times_SY')$.
Then $Y''\to X''$
is an object of
${\cal C}_{V\to U}$
and there exist unique maps
$(Y''\to X'')
\to
(Y\to X)$
and
$(Y''\to X'')
\to
(Y'\to X')$

2.
Let $Y\to X$ be an object
of ${\cal C}_{V\to U}$.
Since $V$
is schematically dense
in $Y$,
the diagram (\ref{eqtg})
is cartesian.
Then 
it follows from \cite[5.7.10]{GR} that
there exists a blow-up
$X'\to X$ inducing an isomorphism
$U\times_XX'\to U$
such that
the proper transform $Y'$
of $Y$
is finite flat over $X'$.
After replacing $X'$ by
the schematic closure of $U$,
the immersion $U\to X$
is schematically dense.
Since $Y$ is flat over $X$,
the immersion $V\to Y$
is also schematically dense.

3.
Let $Y\to X$ be an object
of ${\cal C}_{V\to U}$.
By replacing $Y$ 
by the schematic closure
of the diagonal image of $V$
in the fibered product
$\prod_{\sigma \in G, X}Y$
over $X$,
we may assume that 
$Y$ carries
an action of $G$.
Then $Y'$ constructed
in the proof of 2.\
also carries
an action of $G$.
\qed

Recall that
for an object
$Y$ of the category
${\cal C}_{V/S}$
of compactifications of
$V$ over $S$,
the wild ramification locus
$\Sigma_{V/U}Y
\subset Y$
is defined
as a closed subscheme
in Definition 
{\rm \ref{dftmT}}.
Further, if
$V$ is schematically dense
in $Y$,
the closed subsets
$\Sigma_{V/U}Y$
form a projective system
by Lemma \ref{lmVYY'}.

\begin{df}\label{dfCU}
Let $f\colon V\to U$ 
be a finite \'etale morphism
of separated schemes 
of finite type over $S$.

{\rm 1.}
We define an abelian group
$F_0G(\partial_{V/U} U)$ and 
a ${\mathbb Q}$-vector
space
$F_0G(\partial_{V/U} U)_{\mathbb Q}$
as the inverse limits:
\begin{align*}
F_0G(\partial_{V/U} U)
&=\varprojlim_{(\bar f\colon
Y\to X)\in 
{\cal C}_{V\to U}}
F_0G(\bar f(\Sigma_{V/U}Y))
\\
F_0G(\partial_{V/U} U)_{\mathbb Q}
&=\varprojlim_{(\bar f\colon
Y\to X)\in 
{\cal C}_{V\to U}}
\bigl(F_0G(\bar f(\Sigma_{V/U}Y))
\otimes_{\mathbb Z}{\mathbb Q}\bigr).
\end{align*}
with respect
to the proper push-forward maps.

{\rm 2.}
Let $W$ be a separated
scheme of finite type
over $V$.
We define an abelian group
$F_0G(\partial_{V/U} W)$ and
a ${\mathbb Q}$-vector
space
$F_0G(\partial_{V/U} W)_{\mathbb Q}$
as the inverse limits:
\begin{align*}
F_0G(\partial_{V/U} W)
&=\varprojlim_{(Z\to Y)\in 
{\cal C}_{W\to V}}
F_0G(\Sigma_{V/U}Y\times_YZ)
\\
F_0G(\partial_{V/U} W)_{\mathbb Q}
&=
\varprojlim_{(Z\to Y)\in 
{\cal C}_{W\to V}}
\bigl(F_0G(\Sigma_{V/U}Y\times_YZ)
\otimes_{\mathbb Z}{\mathbb Q}\bigr)
\end{align*}
with respect
to the proper push-forward maps.
\end{df}

In the rest of this
subsection,
we will establish
properties for
$F_0G(\partial_{V/U} W)$.
The same proof also works for
$F_0G(\partial_{V/U} W)_{\mathbb Q}$.

For an object
$Z\to Y$ of 
${\cal C}_{W\to V}$,
we have a canonical map
$${\rm pr}_Z\colon
F_0G(\partial_{V/U} W)
\to 
F_0G(\Sigma_{V/U}Y\times_YZ).$$

Since we will assume
that the covering
of the generic fibers
$V_K\to U_K$
is tamely ramified
with respect to $K$
in the definition
of the invariants
in the next subsection,
the group
$F_0G(\partial_{V/U} W)$
is generated by
the classes
supported on the
closed fibers,
in practice.
The assumption
is always satisfied
if $K$ is of characteristic 0, 
by Corollary \ref{corptp}.

\begin{lm}\label{lmg!}
\setcounter{equation}0
Let 
\begin{equation}
\begin{CD}
U'@<<<V'@<<< W'\\
@VVV @VVV@VVgV\\
U@<<<V@<<< W\\
\end{CD}
\label{eqUVW}
\end{equation}
be a commutative diagram
of separated schemes
of finite type over $S$
such that the left
square is cartesian
and that
the map
$V\to U$ is finite \'etale.
Then, the push-forward maps
induce
\begin{equation}
g_!\colon
F_0G(\partial_{V'/U'} W')
\to 
F_0G(\partial_{V/U}W),\quad
g_!\colon
F_0G(\partial_{V'/U'} W')_{\mathbb Q}
\to 
F_0G(\partial_{V/U}W)_{\mathbb Q}.
\label{eqg!}
\end{equation}
\end{lm}

If $g$ is proper,
we write $g_*=g_!$
in (\ref{eqg!}).

{\it Proof.}
We define
a category ${\cal C}_{W'\to V'/W\to V}$
consisting of commutative diagrams
$$\begin{CD}
Y'@<<<Z'\\
@V{\bar h}VV@VV{\bar g}V\\
Y@<<< Z
\end{CD}$$
of schemes over $S$
compatible with
the right square of (\ref{eqUVW})
such that
$Z\to Y$ and
$Z'\to Y'$ are objects of
${\cal C}_{W\to V}$
and of
${\cal C}_{W'\to V'}$
respectively.
By Lemma \ref{lmVYY'},
for an object of
${\cal C}_{W'\to V'/W\to V}$
we have
$\bar h(\Sigma_{V'/U'}Y')
\subset 
\Sigma_{V/U}Y$
and the push-forward map
$\bar g_*\colon
F_0G(\Sigma_{V'/U'}Y'\times_{Y'}Z')
\to
F_0G(\Sigma_{V/U}Y\times_YZ)$
is defined.
Similarly as in Lemma \ref{lmCU},
the image of
${\cal C}_{W'\to V'/W\to V}$
in
${\cal C}_{W'\to V'}$
is cofinal.
Hence 
the map $g_!\colon
F_0G(\partial_{V'/U'} W')
\to 
F_0G(\partial_{V/U}W)$
is defined as the limit.
The map 
$g_!\colon
F_0G(\partial_{V'/U'} W')_{\mathbb Q}
\to 
F_0G(\partial_{V/U}W)_{\mathbb Q}$
is defined similarly.
\qed

\begin{lm}\label{lmff}
Let $f\colon V\to U$ 
be a finite \'etale morphism
of separated schemes 
of finite type over $S$
and let
$g\colon W'\to W$
be a finite flat
morphism of
separated schemes
of finite type over $V$.

{\rm 1.}
The pull-back maps
induce a map
$$g^*\colon
F_0G(\partial_{V/U} W)
\to
F_0G(\partial_{V/U} W').
$$

{\rm 2.}
Assume that
$g\colon W'\to W$
is of degree $d$.
Then, the composition
$g_*\circ g^*
\colon
F_0G(\partial_{V/U} W)
\to
F_0G(\partial_{V/U} W)$
is multiplication
by $d$.

{\rm 3.}
Assume that $W'$
is a $G$-torsor over $W$ for
a finite group $G$.
Then, the composition
$g^*\circ g_*\colon
F_0G(\partial_{V/U} W')
\to
F_0G(\partial_{V/U} W')$
is equal to 
$\sum_{\sigma\in G}\sigma^*$.
Consequently,
$g^*\colon
F_0G(\partial_{V/U} W
)_{\mathbb Q}
\to
F_0G(\partial_{V/U} W'
)_{\mathbb Q}$
is an isomorphism
to the $G$-fixed part.
\end{lm}

{\it Proof.}
1.
By Lemma \ref{lmCU}.2.,
it suffices to show the following:
Let 
$(h',h)\colon
((\bar g_1,{\rm id}_Y)\colon
(Z'_1\to Y)\to (Z_1\to Y))
\to 
((\bar g,{\rm id}_Y)\colon
(Z'\to Y)\to (Z\to Y))$
be a morphism 
in the category ${\cal C}_{W'\to V/ W\to V}$
defined in the proof of Lemma \ref{lmg!} such that
$\bar g_1\colon
Z'_1\to Z_1$ and
$\bar g\colon
Z'\to Z$
are finite flat
and that
the maps $W'\to W\times_{Z_1}Z'_1$
and
$W'\to W\times_ZZ'$ are isomorphisms.
Then,
the diagram
$$\begin{CD}
F_0G(\Sigma_{V/U}Y\times_YZ_1)
@>{\bar g_1^*}>>
F_0G(\Sigma_{V/U}Y\times_YZ'_1)\\
@V{h_*}VV @VV{h'_*}V\\
F_0G(\Sigma_{V/U}Y\times_YZ)
@>{\bar g^*}>>
F_0G(\Sigma_{V/U}Y\times_YZ')
\end{CD}$$
is commutative.

Since the diagram
is commutative if 
$Z'_1=Z'\times_ZZ_1$,
we may assume $h$ is the identity.
Let 
$z\in \Sigma_{V/U}Y\times_YZ$
be a closed point.
Then since the base changes of
the finite flat morphisms
$Z'\to Z$ and $Z'_1\to Z_1=Z$
to the henselization
are decomposed into
the disjoint unions of the spectra
of local rings,
the class
$\bar g^*[z]=
[\bar g^{-1}(z)]$
is equal to
$h'_*(\bar g_1^*[z])
=
h'_*([\bar g_1^{-1}(z)])$.
Hence the assertion follows.

2.
By Lemma \ref{lmCU}.2.,
it suffices to consider
objects
$\bar g\colon Z'\to Z$
of ${\cal C}_{W'\to W}$
such that
$\bar g$ is finite flat of degree $d$.
Then, for a closed point
$z\in \Sigma_{V/U}Y\times_YZ$,
we have
$\bar g_*\bar g^*[z]=d[z]$
and the assertion follows.

3.
Similarly as in the proof of 2,
by Lemma \ref{lmCU}.3.,
it suffices to consider
finite flat objects
$Z'\to Z$
of ${\cal C}_{W'\to W}$
such that the $G$-action
is extended.
Then, for a closed point
$z\in \Sigma_{V/U}Y\times_YZ'$,
we have
$\bar g^*\bar g_*[z]=
\sum_{\sigma\in G}[\sigma z]$
and the assertion follows.
\qed

Similarly,
the push-forward map
$$
f_*\colon
F_0G(\partial_{V/U} V)
\to
F_0G(\partial_{V/U} U)$$
and, 
if $V$ is a $G$-torsor over $U$,
the pull-back map
$$f^*\colon
F_0G(\partial_{V/U} U)
\to
F_0G(\partial_{V/U} V)
$$
are defined.
The following Lemma
is proved in the same way
as Lemma \ref{lmff}.

\begin{lm}\label{lmGal}
\setcounter{equation}0
Let $f\colon V\to U$
be a finite
\'etale morphism
of separated schemes
of finite type over $S$
and assume that $V$
is a $G$-torsor over $U$ for
a finite group $G$
of order $d$.

The composition
$f_*\circ f^*\colon
F_0G(\partial_{V/U} U)
\to
F_0G(\partial_{V/U} U)$
is multiplication by $d$
and the composition
$f^*\circ f_*\colon
F_0G(\partial_{V/U} V)
\to
F_0G(\partial_{V/U} V)$
is equal to 
$\sum_{\sigma\in G}\sigma^*$.
Consequently,
$f^*\colon
F_0G(\partial_{V/U} U
)_{\mathbb Q}
\to
F_0G(\partial_{V/U} V
)_{\mathbb Q}$
is an isomorphism
to the $G$-fixed part.
\end{lm}

Let 
$V'\to V\to U$ be
finite \'etale morphisms
of separated
schemes of finite type
over $S$.
Then, for an object
$g\colon Y'\to Y$
of ${\cal C}_{V'\to V}$,
we have an inclusion
$\Sigma_{V'/U}Y'
\subset
g^{-1}(\Sigma_{V/U}Y)$ 
by Lemma \ref{lmVYY'}.
Hence, for a separated
scheme $W$
of finite type over $V'$,
a canonical map
\begin{equation}
F_0G(\partial_{V/U} W)
\to F_0G(\partial_{V'/U} W)
\label{eqGal1}
\end{equation}
is defined.

Similarly,
let 
$V\to U\to U'$ be
finite \'etale morphisms
of separated
schemes of finite type
over $S$.
Then, 
for an object $Y$
of ${\mathcal C}_{V/S}$,
we have an inclusion
$\Sigma_{V/U}Y
\subset
\Sigma_{V/U'}Y$ 
by Lemma \ref{lmVYY'}.
Hence, for a separated
scheme $W$ of finite type over $V$,
a canonical map
\begin{equation}
F_0G(\partial_{V/U} W)
\to F_0G(\partial_{V/U'} W)
\label{eqGal2}
\end{equation}
is defined.

We introduce a variant.

\begin{df}\label{dfCV}
\setcounter{equation}0
Let $U$ 
be a separated scheme
of finite type over $S$
and ${\cal C}_{U/S}$
be the category of
compactifications
defined in the beginning of
Section {\rm \ref{sslim}}.

{\rm 1.}
We define an abelian group
$F_0G(\partial_F U)$ and a
${\mathbb Q}$-vector space
$F_0G(\partial_F U)_{\mathbb Q}$
as the inverse limits:
\begin{align*}
F_0G(\partial_F U)
&=\varprojlim_{X\in 
{\cal C}_{U/S}}
F_0G(X\times_SF)
\\
F_0G(\partial_F U)_{\mathbb Q}
&=\varprojlim_{X\in 
{\cal C}_{U/S}}
\bigl(F_0G(X\times_SF)
\otimes_{\mathbb Z}
{\mathbb Q}\bigr)
\end{align*}
with respect
to the proper push-forward maps.

{\rm 2.}
For a morphism
$f\colon V\to U$
of separated schemes
of finite type over $S$,
we define a map
$$f_!\colon
F_0G(\partial_F V)
\to 
F_0G(\partial_FU)$$
to be the limit
of the push-forward maps.

If $f$ is proper,
we write $f_*=f_!$.
\end{df}

Let $f\colon V\to U$
be a finite \'etale
morphism of
separated schemes 
of finite type over $S$
such that
the generic fiber
$f_K\colon V_K\to U_K$
is tamely ramified
with respect to $K$.
Then, 
the objects $Y$ 
of the category
${\mathcal C}_{V/S}$
of compactifications
of $V$ satisfying 
set-theortical inclusions
$\Sigma_{V/U}Y
\subset Y_F$
are cofinal
in ${\mathcal C}_{V/S}$.
Hence, 
we have a canonical map
\begin{equation}
F_0G(\partial_{V/U}V)
\to
F_0G(\partial_FV).
\label{eqUVF}
\end{equation}

\subsection{Definition
of invariants
of wild ramification}
\label{ssllpd}

In this subsection,
we define invariants
of wild ramification
without assuming
the regularity of compactification.

First, we recall the
existence
of an alteration.

\begin{lm}[{\cite[Theorem 6.5]{dJI}}]
\label{lmalt}
Let $X$ be a flat separated
scheme of finite type over 
$S={\rm Spec}\ {\cal O}_K$
and $U\subset X$ be a dense 
open subscheme.
Then, there exist
a scheme $Z$ over $S$
and a morphism
$\bar h\colon
Z\to X$ over $S$
satisfying the following
conditions:
\begin{itemize}
\item[{\rm (\ref{lmalt}.1)}]
The scheme $Z$
is regular flat separated
of finite type over $S$
and $W=\bar h^{-1}(U)$ 
is the complement
of a divisor $D$
with simple normal crossings.
\item[{\rm (\ref{lmalt}.2)}]
The morphism $\bar h:Z\to X$
is proper,
surjective
and generically finite.
\end{itemize}
\end{lm}

We give some sufficient
conditions
for simultaneous
good alterations
for a scheme $Y$
and for
a weakly semi-stable
scheme $Y'$ over $Y$.
This will be used in Proposition \ref{prcysd}.

\begin{cor}\label{lmcysd}
Let $Y$ be a flat separated
scheme of finite type
over $S$
and $V\subset Y$
be a dense open subscheme.
Let $Y'\to Y$
be a weakly semi-stable scheme
such that
the base change
$Y'_V=Y'\times_YV$
is smooth over $V$.
We assume that
either of the
following conditions
is satisfied:
\begin{itemize}
\item[{\rm (\ref{lmcysd}.1a)}]
\setcounter{equation}1
$Y'\to Y$
is a curve.
\item[{\rm (\ref{lmcysd}.1b)}]
There exist
a morphism
$\bar g_0\colon Z_0\to Y$
of schemes over $S$
satisfying the conditions
{\rm (\ref{lmalt}.1)}
and
{\rm (\ref{lmalt}.2)}
with $X$ and $\bar h$
replaced by
$Y$ and $\bar g_0$
and a log blow-up
$Z'_0\to Y'\times_YZ_0$
inducing an isomorphism
$Z'_0\times_{Z_0}W_0\to Y'\times_YW_0$
where $W_0=\bar g_0^{-1}(V)$.
Further, $Z'_0\to Z_0$ is
weakly strictly semi-stable
and satisfies the condition
{\rm (\ref{lmlbsn}.1)}
with $X\to S$
replaced by
$Z'_0\to Z_0$.
\end{itemize}

Then, there exist
a regular flat scheme $Z$
over $S$,
a proper surjective and
generically finite morphism
$\bar g\colon Z\to Y$
and a log blow-up
$Z'\to Y'\times_YZ$
satisfying the following
conditions:
The inverse image
$W=\bar g^{-1}(V)$
is the complement of
a divisor $D_Z$ with simple 
normal crossings,
the induced map
$Z'\times_ZW\to Y'\times_YW$
is an isomorphism,
the map $Z'\to Z$ is
weakly strictly semi-stable,
the scheme
$Z'$ is regular and
the divisor $Z'\times_ZD_Z$
has simple normal crossings.
\end{cor}

Only the case (\ref{lmcysd}.1a)
will be used in
the proof of
the conductor formula.

{\it Proof.}
First, we consider
the case where
$Y'$ is a curve over $Y$.
By replacing $Y$
by a finite covering
obtained by adjoining some square roots,
we may assume that $Y$
satisfies the condition
(\ref{lmD1}.1)
in Lemma \ref{lmD1}.
By Lemma \ref{lmalt},
there exist
a regular flat scheme $Z$
over $S$,
a proper surjective and
generically finite morphism
$\bar g\colon Z\to Y$
such that inverse image
$W=\bar g^{-1}(V)$
is the complement of
a divisor $D_Z$ with simple 
normal crossings.
Then, it suffices to apply
Lemmas \ref{lmD1}
and \ref{lmlbS}.

If (\ref{lmcysd}.1b)
is satisfied,
it suffices to
apply Lemma \ref{lmlbsn}.
\qed

\begin{pr}\label{prmap}
\setcounter{equation}0
Let $f\colon
V\to U$ be a finite
\'etale morphism of
regular flat
separated schemes
of finite type
over $S$.
Let $Y$ be a flat
separated scheme of finite type
over $S$ containing $V$
as an open subscheme
and let ${\cal D}
=(D_i)_{i\in I}$
be a finite family of
Cartier divisors
of $Y$ satisfying
$V=Y\setminus
\bigcup_{i\in I}D_i$
and 
$\Sigma_{V/U}^{\cal D}Y_K=
\emptyset$.
Let $A_{\cal D}
\subset (Y\times_SY)^\sim_{\cal D}$
denote the closure
$\overline{(V\times_UV\setminus \Delta_V)}$.

Let 
\setcounter{equation}0
\begin{equation}
\begin{CD}
Y@<{\bar g}<< Z
@>{\bar h}>> X\\
@A{\cup}AA @A{\cup}AA 
@AA{\cup}A\\
V@<g<< W
@>h>> U
\end{CD}
\label{eqmapd}
\end{equation}
be a cartesian diagram
of separated scheme of finite
type over $S$
satisfying the following conditions:
\begin{itemize}
\item[{\rm (\ref{eqmapd}a)}]
The scheme $Z$
is regular and flat
over $S$
and $W$ 
is the complement
of a divisor $D$
with simple normal crossings.
\item[{\rm (\ref{eqmapd}b)}]
The scheme $X$ contains
$U$ as the complement of
a Cartier divisor $B$
and we have
$h=f\circ g$.
\end{itemize}

Let 
$(Z\times_XZ)^\sim
\subset
(Z\times_SZ)^\sim$
be the log products
defined by
the family
$(D'_j)_j$
of irreducible
components of
the complement
$Z\setminus W$
and by $B\subset X$
{\rm (\ref{eqmapd}b)}.
Let
$A\subset 
(Z\times_SZ)^\sim$
denote the intersection
$(\bar g\times\bar g)^
{\sim -1}(A_{\cal D})
\cap (Z\times_XZ)^\sim$.

We put $n=\dim Z$.
Then, there exists a unique map
$$\begin{CD}
((\ ,\Delta_Z^{\log}))
\colon
{\rm Gr}^F_\bullet
G(W\times_UW\setminus
W\times_VW)
@>>>
{\rm Gr}^F_{\bullet-n}
G(\Sigma_{V/U}^{\cal D}
Y\times_YZ)
\end{CD}$$
that makes
the diagram
\begin{equation}
\xymatrix{
{\rm Gr}^F_\bullet
G(W\times_UW\setminus
W\times_VW)
\ar[rrr]^{\ \ \ 
((\ ,\Delta_Z^{\log}))}&&&
{\rm Gr}^F_{\bullet-n}
G(\Sigma_{V/U}^{\cal D}Y\times_YZ)\\
{\rm Gr}^F_\bullet
G(A)
\ar[u]^{\rm restriction}
\ar[rrru]_{((\ ,\Delta_Z^{\log}))
_
{(Z\times_SZ)^\sim}}
}
\label{eqprmap2}
\end{equation}
commutative.
\end{pr}

In the characterization of 
the map
$((\ ,\Delta_Z^{\log}))\colon
Gr^F_\bullet G(W\times_UW
\setminus
W\times_VW)\to
{\rm Gr}^F_{\bullet-n}
G(\Sigma_{V/U}^{\cal D}Y\times_YZ)$
in (\ref{eqprmap2}),
we may replace
$A$ by the closure of
$W\times_UW
\setminus
W\times_VW$.

Since the definition of 
the log product
$(Z\times_XZ)^\sim$
involves also the Cartier divisor
$B\subset X$,
it could be better denoted
by $(Z\times_{\mathbb X}Z)^\sim$.
However, since we 
always consider the log product
with $B$ as long as the base is
$X$,
we will use the notation
$(Z\times_XZ)^\sim$.

\noindent{\it Proof.}
We show that
the condition (A$'$)
after Definition \ref{dflip}
and (B) in
Proposition \ref{prprod}
are satisfied.
By the assumption that
the generic fiber
of $\Sigma_{V/U}^{\cal D}Y
=A_{\cal D} 
\cap \Delta_Y^{\log}$
is empty,
the intersection
$A\cap \Delta_Z^{\log}
\subset
(\bar g\times \bar g)^{\sim-1}
(A_{\cal D}
\cap \Delta_Y^{\log})$
is supported on the closed fiber.
Hence,
the condition (A$'$)
after Definition \ref{dflip}
is satisfied and the map
$((\ ,\Delta_Z^{\log}))_
{(Z\times_SZ)^\sim}
\colon
G(A)\to 
G(\Sigma_{V/U}^{\cal D}Y\times_YZ)$
is defined.

Let $(D'_j)_j$
be the irreducible
components of
$Z\setminus W$
and we put
$\bar h^*B=
\sum_jl_jD'_j$.
Since $W=\bar h^{-1}(U)$,
we have $l_j>0$
for each irreducible
component $D'_j$.
By Lemma \ref{lmEi}.2,
the intersection
$A\cap {\mathbf G}_{m,D'_j}
\subset
(Z\times_XZ)^\sim
\cap {\mathbf G}_{m,D'_j}$
is a closed subscheme
of
$\mu_{l_j,D'_j}$
for
each irreducible component
$D'_j$.
Hence the condition (B) in
Proposition \ref{prprod}
is also satisfied.
Since
the intersection
$A^\circ
=A\cap (W\times_UW)$
is equal to
$W\times_UW\setminus
W\times_VW$,
there exists a unique map
$((\ ,\Delta_Z^{\log}))
\colon
{\rm Gr}^F_\bullet
G(W\times_UW\setminus
W\times_VW)
\to 
{\rm Gr}^F_{\bullet-n}
G(\Sigma_{V/U}^{\cal D}Y\times_YZ)$
making the diagram
(\ref{eqprmap2}) commutative
by Proposition \ref{prprod}.
\qed

The map $((\ ,\Delta_Z^{\log}))$
is compatible with the pull-back as follows.

\begin{cor}\label{cormap00}
\setcounter{equation}0
We keep the notation in
Proposition {\rm \ref{prmap}}.
Further,
let $Z'$ be a regular separated
scheme of finite type over $S$
and $\pi\colon Z'\to Z$
be a morphism over $S$
such that $W'=\pi^{-1}(W)$
is the complement of
a divisor with simple normal crossings.

We assume $\dim Z'_K+1=\dim Z_K+1=n$.
Then, 
we have a commutative diagram
\begin{equation}
\begin{CD}
{\rm Gr}^F_\bullet
G(W\times_UW\setminus
W\times_VW)
@>{((\ ,\Delta_Z^{\log}))}>>
{\rm Gr}^F_{\bullet-n}
G(\Sigma_{V/U}^{\cal D}Y\times_YZ)\\
@V{(\pi\times\pi)^*}VV @VV{\pi^*}V\\
{\rm Gr}^F_{\bullet}
G(W'\times_UW'\setminus
W'\times_VW')
@>{((\ ,\Delta_{Z'}^{\log}))}>>
{\rm Gr}^F_{\bullet-n}
G(\Sigma_{V/U}^{\cal D}Y\times_YZ').
\end{CD}
\label{eqcormap2}
\end{equation}
\end{cor}

{\it Proof.}
It suffices to apply Lemma
\ref{corprod}.1.
\qed

We introduce a category of alterations.
Let $f\colon V\to U$
be a proper
morphism of 
separated reduced schemes
of finite type over $S$.
We define a category
${\mathcal A}_{V\to U}$
of alterations as follows:
\begin{itemize}
\item
An object is
a proper, surjective and generically finite
morphism $\bar g\colon
Z\to Y$
of proper schemes over $S$
such that $Y$ contains $V$ 
as a schematically dense
open subscheme
satisfying the following conditions:
\begin{itemize}
\item[{\rm (\ref{coralt}.1a)}]
The scheme $Z$
is regular and flat
over $S$
and $W=\bar g^{-1}(V)$ 
is the complement
of a divisor $D$
with simple normal crossings.
There exists a dense open subscheme
$V_0$ of $V$ such
that $\bar g^{-1}(V_0)\to V_0$
is finite flat of constant rank.
\item[{\rm (\ref{coralt}.1b)}]
There exists a proper scheme $X$
over $S$
containing $U$
as the complement of
a Cartier divisor $B$
and a cartesian diagram
\addtocounter{thm}{1}
\setcounter{equation}0
\begin{equation}
\begin{CD}
W@>\subset>>Z\\
@V{f\circ g}VV @VV{\bar h}V\\
U@>\subset>> X
\end{CD}
\label{eqZX}
\end{equation}
of schemes over $S$
where $g\colon W\to V$
is the restriction of 
$\bar g\colon Z\to Y$.
\end{itemize}

\item
A morphism $(\bar \pi,\varphi)
\colon (\bar g'\colon Z'\to Y')
\to (\bar g\colon Z\to Y)$
is a pair of 
a proper, surjective and generically finite
morphism $\bar\pi\colon Z'\to Z$
over $S$
and a morphism
$\varphi\colon Y'\to Y$
of schemes
over $S$ such that
the diagram
\begin{equation}
\begin{CD}
Z'@>{\bar g'}>>Y'@<\supset<<V\\
@V{\bar \pi}VV @VV{\varphi}V
@|\\
Z@>{\bar g}>>Y
@<\supset<<V\end{CD}
\label{eqalt2}
\end{equation}
is commutative
and that
there exists a dense open subscheme
$Z_0$ of $Z$ such
that $\bar \pi^{-1}(Z_0)\to Z_0$
is finite flat of constant rank.
\end{itemize}

\addtocounter{thm}{-1}
\begin{lm}\label{coralt}
\setcounter{equation}0
Let $f\colon V\to U$
be a proper,
surjective
and generically finite
morphism of separated schemes
of finite type over $S$.

{\rm 1.}
For an object $Y$ 
of the category 
${\cal C}_{V/S}$
of compactifications 
containing $V$ 
as a schematically dense
open subscheme,
there exists
an object
$\bar g\colon Z\to Y$
of the category
${\mathcal A}_{V\to U}$
of alterations.

{\rm 2.}
The category
${\mathcal A}_{V\to U}$
is cofiltered.
\end{lm}

{\it Proof.}
1. By Nagata's embedding
theorem \cite{nagata}, 
there exists
a proper scheme
$X$ over $S$ containing
$U$ as an open subscheme.
By replacing $X$
by a blow-up at
the complement 
$X\setminus U$ if necessary,
we may assume
$U\subset X$
is the complement of
a Cartier divisor $B$.
By replacing $Y$
by the closure
of the graph $\Gamma_f
\subset X\times_SY$,
we may assume
there exists a morphism
$\bar f\colon Y\to X$
such that
$\bar f^{-1}(U)=V$
and 
$\bar f|_V=f$.
Then, it suffices
to apply Lemma \ref{lmalt}
to the open immersion 
$V\to Y$
and to take some
disjoint union of connected component.

2.
It suffices to apply
Lemma \ref{lmalt}
to the open
immersion
$W\times_VW'
\to 
\overline{W\times_VW'}
\subset
Z\times_SZ'$
and to take some
disjoint union of connected component.
\qed

Note that 
the condition
(\ref{eqZX}b) is satisfied
if we have an object
$\bar f\colon Y\to X$
of ${\cal C}_{V\to U}$
such that
$X$ contains
$U$ as the complement
of a Cartier divisor $B$.
If $V\to U$
is a Galois covering,
such an object may be constructed as follows.

\begin{lm}\label{lmXYG}
Let $f\colon V\to U$
be a finite \'etale morphism
of regular separated
scheme of finite type
over $S$ and
$V\to Y$
be an open immersion 
of separated schemes of
finite type over $S$.
Let ${\mathcal D}=
(D_i)_{i\in I}$
be a finite family of Cartier
divisors of $Y$ such that
$V$ is the complement
of the union
$\bigcup_{i\in I}D_i$.
Assume that
$V$ is a $G$-torsor over $U$
and that the action of $G$
is extended to $Y$
and on ${\mathcal D}$.
Assume further that
the action of $G$
on $Y$ is admissible
in the sense that
the quotient $X=Y/G$
is defined as a scheme.

Then, 
the canonical map $\bar f\colon
Y\to X$
is finite,
the quotient
$X$ is separated of finite type
over $S$ and there
exists a Cartier divisor
$B$ of $X$ such that
the complement is $U$.
\end{lm}

{\it Proof.}
It suffices to show
the existence of $B$.
By the assumption
the sum $D=\sum_iD_i$
is stable by the action of $G$
and $V=Y\setminus D.$
The norm $B$ of $D$
is defined as a Cartier divisor of
$X$ since
${\mathcal O}_X
\to\bar f_*{\mathcal O}_Y$
is injective.
Since the inverse image of
the complement 
$X\setminus B$
is $V$,
we obtain
$U=X\setminus B$.
\qed

Let $f\colon V\to U$ be a finite
\'etale morphism of
regular separated
schemes of
finite type over $S$.
For a morphism
$g'\colon W\to V'$
of regular schemes of
finite type over $V$,
the pull-back map 
$$(g'\times g')^*
\colon
{\rm Gr}^F_nG(V'\times_UV'
\setminus
V'\times_VV')
\to 
Gr^F_nG(W\times_UW
\setminus
W\times_VW)$$
by
$g'\times g'\colon
W\times_SW\to
V'\times_SV'$
is defined by
Corollary \ref{cortdf}
and 
Lemma \ref{lmtfil}.

\begin{thm}\label{thmmap}
\setcounter{equation}0
Let $f\colon V\to U$ be a finite
\'etale morphism of
regular separated
schemes of
finite type over $S=
{\rm Spec}\ {\cal O}_K$
such that
$U_K\to V_K$
is tamely ramified
with respect to $K$.
Let $V'$ be a regular
flat scheme of finite
type over $S$
and $V'\to V$
be a proper morphism over $S$.
We assume that
$\dim V_K=\dim V'_K$
and put
$n=\dim V_K+1$.

Then, there exists a unique map
\setcounter{equation}0
\begin{equation}
\begin{CD}
((\ ,\Delta_{V'}))^{\log}
\colon {\rm Gr}^F_nG(V'\times_UV'
\setminus
V'\times_VV')
@>>>
F_0G(\partial_{V/U} V')_{\mathbb Q}.
\end{CD}
\label{eqmapG}
\end{equation}
satisfying the following property:

For an object $Y$ of 
${\cal C}_{V/S}$,
a finite family ${\cal D}$ of
Cartier divisors of $Y$
such that
$\Sigma_{V/U}Y
=\Sigma^{\cal D}_{V/U}Y$
(Definition 
{\rm \ref{dftmT}.1}),
an object $Y'\to Y$
of ${\cal C}_{V'\to V}$
and an object $\bar g'\colon Z\to Y'$
of ${\cal A}_{V'\to U}$
such that
$\bar g'$ is
generically
of constant degree $[W:V']$,
the diagram
\begin{equation}
\xymatrix{
{\rm Gr}^F_nG(V'\times_UV'
\setminus
V'\times_VV')
\ar[rr]^{\qquad ((\ ,\Delta_{V'}^{\log}))}
\ar[d]_{(g'\times g')^*}
&&
F_0G(\partial_{V/U} V')_{\mathbb Q}
\ar[d]^{{\rm pr}_{Y'}}
\\
Gr^F_nG(W\times_UW
\setminus
W\times_VW)
\ar[rrd]_{((\ ,\Delta_Z^{\log}))}
&&
F_0G(\Sigma_{V/U}Y\times_YY')_{\mathbb Q}
\\
&&
F_0G(\Sigma_{V/U}Y\times_YZ)
\ar[u]_{\frac1{[W:V']}
\bar g'_*}
}
\label{eqtmap}
\end{equation}
is commutative,
where $g'\colon W
=\bar g^{\prime -1}(V')\to V'$
is the restriction of
$\bar g'\colon Z\to Y'$
\end{thm}

{\it Proof.}
By the remark after
Definition \ref{dftmT},
there exist an object $Y$ of 
${\cal C}_{V/S}$ and
a finite family ${\cal D}$ of
Cartier divisors of $Y$
such that
$\Sigma_{V/U}Y
=\Sigma^{\cal D}_{V/U}Y$.
Further since ${\cal C}_{V'\to V}$
is non-empty,
there exist an object $Y'\to Y$
of ${\cal C}_{V'\to V}$
and an object $\bar g'\colon Z\to Y'$
of ${\cal A}_{V'\to U}$
by Lemma \ref{coralt}.1.
By the definition of
$F_0G(\partial_{V/U} V')_{\mathbb Q}$
as the projective limit,
it suffices to show that
the composition
of the lower maps in the diagram
(\ref{eqtmap})
is independent of
the choice of 
an object $\bar g'\colon Z\to Y'$
of ${\cal A}_{V'\to U}$
and that the compositions
form an inverse system
with respect to objects $Y'\to Y$
of ${\cal C}_{V'\to V}$.

The categories ${\cal C}_{V'\to V}$
and ${\cal A}_{V'\to U}$
are cofiltered by Lemmas
\ref{lmCU}.1 and \ref{coralt}.2.
Hence,
it suffices to show that the diagram
\begin{equation}
\begin{CD}
Gr^F_nG(W\times_UW
\setminus
W\times_VW)
@>{(\pi\times \pi)^*}>>
Gr^F_nG(W_1\times_UW_1
\setminus
W_1\times_VW_1)
\\
@V{((\ ,\Delta_Z^{\log}))}VV
@VV{((\ ,\Delta_{Z_1}^{\log}))}V
\\
F_0G(\Sigma_{V/U}Y\times_YZ)
@.
F_0G(\Sigma_{V/U}Y_1\times_{Y_1}Z_1)
\\
@V{\frac1{[W:V']}
\bar g_*}VV
@VV{\frac1{[W_1:V']}
\bar g_{1*}}V
\\
F_0G(\Sigma_{V/U}Y\times_YY')_{\mathbb Q}
@<{\varphi_*}<<
F_0G(\Sigma_{V/U}Y_1\times_{Y_1}Y'_1)_{\mathbb Q}
\end{CD}
\label{eqtmap2}
\end{equation}
is commutative for morphisms
$(\varphi,\psi)\colon
(Y'_1\to Y_1)
\to (Y'\to Y)$ of 
${\cal C}_{V'\to V}$
and
$(\bar \pi,\varphi)\colon
(\bar g_1\colon Z_1\to Y'_1)
\to (\bar g\colon Z\to Y')$ of 
${\cal A}_{V'\to U}$
where $\pi\colon W_1\to W$
denotes the restriction
of $\bar \pi\colon Z_1\to Z$
on the inverse image of $V'$.
By Corollary \ref{cormap00},
the pull-back
$\bar \pi^*\colon
F_0G(\Sigma_{V/U}Y\times_YZ)
\to
F_0G(\Sigma_{V/U}Y_1\times_{Y_1}Z_1)$
makes the upper half of
(\ref{eqtmap2})
into a commutative square.
On the other hand,
the push-forward
$\bar \pi_*\colon
F_0G(\Sigma_{V/U}Y_1\times_{Y_1}Z_1)
\to
F_0G(\Sigma_{V/U}Y\times_YZ)$
divided by the degree $[Z_1: Z]$
makes the lower half
into a commutative square.
The composition
$\bar\pi_*\circ 
\bar \pi^*$
is equal
to the multiplication
by $[R\bar \pi_*{\cal O}_{Z_1}]$
and induces
the multiplication
by $[Z_1:Z]={\rm rank}(R\bar \pi_*{\cal O}_{Z_1})$
on
$F_0G(\Sigma_{V/U}Y\times_YZ)$.
Hence, 
the assertion is proved.
\qed

\begin{df}\label{dfmapm}
Let the notation be as in
Theorem {\rm \ref{thmmap}}.
We call the map
\setcounter{equation}0
\begin{equation}
\begin{CD}
((\ ,\Delta_{V'}))^{\log}
\colon {\rm Gr}^F_nG(V'\times_UV'
\setminus
V'\times_VV')
@>>>
F_0G(\partial_{V/U} V')_{\mathbb Q}.
\end{CD}
\label{eqmapm}
\end{equation}
{\rm the logarithmic 
localized intersection
product with the diagonal}.
For an object
$Y'\to Y$ of ${\mathcal C}_{V'\to V}$,
we define
\begin{equation}
\begin{CD}
((\ ,\Delta_{Y'}))^{\log}
\colon {\rm Gr}^F_nG(V'\times_UV'
\setminus
V'\times_VV')
@>>>
F_0G(\Sigma_{V/U}Y
\times_YY')_{\mathbb Q}.
\end{CD}
\label{eqmapmY}
\end{equation}
as the composition of
{\rm (\ref{eqmapm})}
with the projection
$F_0G(\partial_{V/U} V')_{\mathbb Q}
\to 
F_0G(\Sigma_{V/U}Y
\times_YY')_{\mathbb Q}$
and call it also
the logarithmic 
localized intersection
product with the diagonal.
\end{df}

Since we assume
that
$V_K\to U_K$
is tamely ramified
with respect to
${\rm Spec}\ K$,
the target
group 
$F_0G(\partial_{V/U} V')_{\mathbb Q}$
is generated by
the classes
supported on the
closed fibers.
The assumption is
always satisfied if
$K$ is of characteristic 0,
by Corollary \ref{corptp}.

If $V$ is finite \'etale over $V$,
the graded piece
${\rm Gr}^F_nG
(V'\times_UV'
\setminus
V'\times_VV')$
is identified
with the free abelian group
$Z^0(V'\times_UV'
\setminus
V'\times_VV')$
generated by
the irreducible components
of $V'\times_UV'
\setminus
V'\times_VV'$.
Thus, in this case,
the maps (\ref{eqmapm})
and (\ref{eqmapmY})
define
\begin{equation}
\begin{CD}
((\ ,\Delta_{V'}))^{\log}
\colon Z^0(V'\times_UV'
\setminus
V'\times_VV')
@>>>
F_0G(\partial_{V/U} V')_{\mathbb Q},
\end{CD}
\label{eqmap}
\end{equation}
\begin{equation}
\begin{CD}
((\ ,\Delta_{Y'}))^{\log}
\colon Z^0(V'\times_UV'
\setminus
V'\times_VV')
@>>>
F_0G(\Sigma_{V/U}Y\times_YY')_{\mathbb Q}.
\end{CD}
\label{eqthmmap}
\end{equation}
Theorem \ref{thmmap}
implies that,
for an open and closed subscheme
$\Gamma$ of 
$V'\times_UV'
\setminus
V'\times_VV'$,
the logarithmic localized
intersection products
$((\Gamma ,\Delta_{Y'}))^{\log}$
for objects
$Y'\to Y$
of ${\cal C}_{V'\to V}$
such that $\Sigma_{V/U}^+Y_K=
\emptyset$
form a projective system
and defines an element
of
$F_0G(\partial_{V/U} V')_{\mathbb Q}
=
\varprojlim_{Y'\to Y}
F_0G(\Sigma_{V/U}Y\times_YY')_{\mathbb Q}$.

Keep assuming $V'\to V$
is finite \'etale and
let $Y'\to Y$
be an object of
${\cal C}_{V'\to V}$.
Assume that
$Y'$ is regular
and $V'\subset Y'$ is
the complement of
a divisor with simple normal crossings.
Assume further that
there exist a
proper scheme $X$
over $S$
containing $U$
as the complement
of a Cartier divisor
and a morphism
$Y'\to X$
extending $V'\to U$.
Then, the identity $Y'\to Y'$
is an object of ${\cal A}_{V'\to U}$
and, for a finite family
${\cal D}$ of Cartier divisors of $Y$
such that $\Sigma^{\cal D}_{V/U}Y_K
=\emptyset$,
the diagram
\begin{equation}
\begin{CD}
Z^0(V'\times_UV'
\setminus
V'\times_VV')
@>{((\ ,\Delta_{V'}))^{\log}}>>
F_0G(\partial_{V/U}V')
\\
@V{((\ ,\Delta_{Y'}^{\log}))}
VV
@VV{{\rm pr}_{Y'}}V
\\
F_0G(\Sigma_{V/U}^{\cal D}Y\times_Y Y')
_{\mathbb Q}
@<<<
F_0G(\Sigma_{V/U}Y\times_Y Y')
_{\mathbb Q}
\end{CD}
\end{equation}
is commutative.
Consequently, if
we assume resolution
of singularities
or if we assume
$\dim Y_K\le 1$,
we do not need
to introduce
denominator
and an integral version
\begin{equation}
((\ ,\Delta_{V'}))^{\log}_{\mathbb Z}
\colon
Z^0(V'\times_UV'
\setminus
V'\times_VV')
\to
F_0G(\partial_{V/U} V')
\label{eqmapZ}
\end{equation}
is defined
as the limit of
$((\ ,\Delta_{Y'}^{\log}))$.

Let $f\colon
V\to U$
be a finite \'etale
morphism
of smooth separated
schemes of finite
type over $F$
and $V'$ be
a finite \'etale scheme
over $V$.
Then, similarly as above,
slightly refining
\cite[Theorem 3.2.3]{KSA},
we define a map
\begin{equation}
\begin{CD}
(\ ,\Delta_{V'})^{\log}\colon
Z^0(V'\times_UV'
\setminus 
V\times_UV)
@>>>
CH_0(\partial_{V/U}V')_{\mathbb Q}.
\end{CD}
\label{eqmapF}
\end{equation}

We introduce a variant
of the map
(\ref{eqmap})
assuming $K$ is of
characteristic $0$.
This variant
is defined without
removing the diagonal
$\Delta_V
\subset V\times_UV$.

\begin{thm}\label{thmmap0}
Assume $K$ is of
characteristic $0$.
Let $V\to U$ be a finite
\'etale morphism of
smooth separated
schemes of
finite type over $S$
and let $Y$ be an object of 
${\cal C}_{V/S}$.
Then, there exists a unique map
\setcounter{equation}0
\begin{equation}
\begin{CD}
((\ ,\Delta_V))^{\log}
\colon Z^0(V\times_UV)
@>>>
F_0G(\partial_F V)_{\mathbb Q}.
\end{CD}
\label{eqmapK}
\end{equation}
satisfying the following property:

For an object $Y$ of 
${\cal C}_{V/S}$,
a finite family ${\cal D}$ of
Cartier divisors of $Y$
such that
$\Sigma^{\cal D}_{V/U}Y_K
=\emptyset$
and an object $\bar g\colon Z\to Y$
of ${\cal A}_{V\to U}$
such that
$\bar g$ is
generically
of constant degree $[W:V]$,
there exists a map 
$((\ ,\Delta_Z^{\log}))
\colon
Gr^F_nG(W\times_UW)\to
F_0G(Z\times_S{\rm Spec}\ F)$
that makes
the diagram
\begin{equation}
\xymatrix{
{\rm Gr}^F_nG(V\times_UV)
\ar[rr]^{\qquad ((\ ,\Delta_V^{\log}))}
\ar[d]_{(g\times g)^*}
&&
F_0G(\partial_F V)_{\mathbb Q}
\ar[d]^{{\rm pr}_{Y}}
\\
Gr^F_nG(W\times_UW)
\ar[rrd]^{((\ ,\Delta_Z^{\log}))}
&&
F_0G(Y\times_S{\rm Spec}\ F)_{\mathbb Q}
\\
Gr^F_nG((Z\times_XZ)^\sim)
\ar[rr]_{((\ ,\Delta_Z^{\log}))
_{(Z\times_SZ)^\sim}}
\ar[u]^{\rm restriction}
&&
F_0G(Z\times_S{\rm Spec}\ F)
\ar[u]_{\frac1{[W:V']}
\bar g_*}
}
\label{eqtmap0}
\end{equation}
commutative.
\end{thm}

By the assumption that
$K$ is of characteristic 0,
in the notation
of the proof of
Proposition \ref{prmap},
the generic fiber
$Z_K$ is smooth over $K$
and the log product
$(Z\times_XZ)^\sim
\subset
(Z\times_SZ)^\sim$
satisfies
the condition (A$'$)
after Definition \ref{dflip}
with $X$ replaced by $Z$.
Except this remark,
the proof is the same as
that of Theorem \ref{thmmap}.

We keep assuming that
$K$ is of characteristic $0$
and let
$f\colon V\to U$
be a finite \'etale
morphism of
regular flat
separated schemes 
of finite type over $S$.
Then, the maps
(\ref{eqmap})
and
(\ref{eqmapK})
are compatible
in the sense
that
the diagram
\begin{equation}
\begin{CD}
Z^0(V\times_UV
\setminus 
\Delta_V)
@>{((\ ,\Delta_V))^{\log}}>>
F_0G(\partial_{V/U} V)_{\mathbb Q}\\
@V{\cap}VV @VV
{\rm(\ref{eqUVF})}V\\
Z^0(V\times_UV)
@>{((\ ,\Delta_{V}))^{\log}}>>
F_0G(\partial_F V)_{\mathbb Q}
\end{CD}
\label{eqcmpt}
\end{equation}
is commutative.

\subsection{Elementary
properties
of the invariants
of wild ramification}
\label{ssiw}

The map (\ref{eqmap})
has 
the following compatibility.

\begin{pr}\label{prfun}
\setcounter{equation}0
Let $f\colon V\to U$
be a finite \'etale morphism
of regular schemes over $S$
such that
the generic fiber
$V_K\to U_K$
is tamely ramified with respect
to $K$.
Let $V'$
be a finite \'etale scheme over $V$.

{\rm 1.}
For a finite \'etale
morphism $g\colon V''\to V'$,
the diagram
\begin{equation}
\begin{CD}
Z^0(V'\times_UV'\setminus
V'\times_VV')
@>{((\ ,\Delta_{V'}))^{\log}}>>
F_0G(\partial_{V/U}V')
_{\mathbb Q}
\\
@V{(g\times g)^*}VV @VV{g^*}V\\
Z^0(V''\times_UV''\setminus
V''\times_VV'')
@>{((\ ,\Delta_{V''}))^{\log}}>>
F_0G(\partial_{V/U}V'')
_{\mathbb Q}
\end{CD}
\label{eqfun1}
\end{equation}

{\rm 2.}
Assume that
the generic fiber
$V'_K\to U_K$
is tamely ramified with respect
to $K$.
Then,
\begin{equation}
\begin{CD}
Z^0(V''\times_UV''\setminus
V''\times_VV'')
@>{((\ ,\Delta_{V''}))^{\log}}>>
F_0G(\partial_{V/U}V'')
_{\mathbb Q}
\\
@V{\rm can}VV 
@VV{\rm (\ref{eqGal1})}V\\
Z^0(V''\times_UV''\setminus
V''\times_{V'}V'')
@>{((\ ,\Delta_{V''}))^{\log}}>>
F_0G(\partial_{V'/U}V'')
_{\mathbb Q}
\end{CD}
\label{eqfun11}
\end{equation}
is commutative.

{\rm 3.}
For a finite \'etale
morphism $U\to U'$
such that
the generic fiber
$V_K\to U'_K$
is tamely ramified with respect
to $K$,
the diagram
\begin{equation}
\begin{CD}
Z^0(V'\times_UV'\setminus
V'\times_VV')
@>{((\ ,\Delta_{V'}))^{\log}}>>
F_0G(\partial_{V/U}V')
_{\mathbb Q}
\\
@V{\rm can}VV 
@VV{\rm (\ref{eqGal2})}V\\
Z^0(V'\times_{U'}V'\setminus
V'\times_VV')
@>{((\ ,\Delta_{V'}))^{\log}}>>
F_0G(\partial_{V/U'}V')
_{\mathbb Q}
\end{CD}
\label{eqfun2}
\end{equation}
is commutative.
\end{pr}

{\it Proof.}
1.
By Lemma \ref{lmff}.2.,
the map $g^*$
is injective.
Hence, by replacing
$V''$ if necessary,
we may assume $V''$
is a $G$-torsor over $V'$
for a finite group $G$.
Since an object
of ${\cal A}_{V''\to U}$
define an object
of ${\cal A}_{V'\to U}$,
the square with the arrow $g^*$
replaced by $|G|^{-1}g_*$
going the other way
is commutative
by the definition
of the map (\ref{eqmap}).
Since the images are
in the $G$-fixed part,
the assertion follows
from Lemma \ref{lmff}.3.

The rest is 
clear from the definition
and the remark after
Proposition \ref{prmap}.
\qed

We show a compatibility
with tame base change.

\begin{cor}\label{corfun}
\setcounter{equation}0
Let $f\colon V\to U$
be a finite \'etale morphism
of regular schemes over $S$
such that
the generic fiber
$V_K\to U_K$
is tamely ramified with respect
to $K$.
Let $g\colon U'\to U$
be a finite \'etale morphism
of regular schemes over $S$.
Let $V'\subset V\times_UU'$
be an open and closed subscheme
and $g'\colon V'\to V$
denote the projection.

Then,
if $g\colon U'\to U$ is
tamely ramified with respect
to $S$, the diagram
\begin{equation}
\begin{CD}
Z^0(V\times_UV\setminus
\Delta_V)
@>{((\ ,\Delta_V))^{\log}}>>
F_0G(\partial_{V/U}V)
_{\mathbb Q}
\\
@V{(g'\times g')^*}VV 
@VV{g^{\prime *}}V\\
Z^0(V'\times_{U'}V'\setminus
\Delta_{V'})
@>{{\rm can}\circ
((\ ,\Delta_{V'}))^{\log}}>>
F_0G(\partial_{V/U}V')
_{\mathbb Q}
\end{CD}
\label{eqtamap}
\end{equation}
is commutative.
\end{cor}

{\it Proof.}
By Lemma \ref{lmVYY'},
the canonical map
$F_0G(\partial_{V'/U'}V')
\to
F_0G(\partial_{V/U}V')$
is defined.
By Proposition \ref{prfun}.1
applied to
$V'\to V=V\to U$,
the diagram (\ref{eqtamap}) with
$Z^0(V'\times_{U'}V'\setminus
\Delta_{V'})$
replaced by
$Z^0(V'\times_UV'\setminus
V'\times_VV')$
is commutative.

By the assumption
that $g\colon U'\to U$
is tamely ramified with respect to $S$,
there exists an object
$X'$ of the category
${\cal C}_{U'/S}$ of compactifications
of $U'$ and a finite family ${\cal D}'$
of Cartier divisors $X'$
such that
$\Sigma_{U'/U}^{{\cal D}'}X'$
is empty.
For an object $Y'\to X'$ 
of ${\mathcal C}_{V'\to U'}$,
the closure of
the inverse image
$V'\times_UV'\setminus 
V'\times_{U'}V'$
of $U'\times_UU'\setminus \Delta_{U'}$
does not meet with
the log diagonal
$\Delta_{Y'}$
in the log product
$(Y'\times_SY')^\sim_{{\cal D}'}$
for the pull-back ${\cal D}'$
of ${\cal D}$.
Hence, we have
$((\Gamma ,\Delta_{V'}))^{\log}=0$
in $F_0G(\partial_{V'/U'}V')
_{\mathbb Q}$
if $\Gamma$
is in 
$V'\times_UV'\setminus 
V'\times_{U'}V'$.

By the assumption
$V'\subset V\times_UU'$,
we have
$(V'\times_{U'}V')\cap
(V'\times_VV')=\Delta_{V'}$.
Hence the assertion follows.
\qed

In the case where
$K$ is of characteristic 0,
the variant defined
in Theorem \ref{thmmap0}
satisfies properties
analogous to
Proposition \ref{prfun}
and Corollary \ref{corfun},
by the same proof.

The logarithmic different
and the logarithmic
Lefschetz class
defined in Section \ref{ssdiff}
are defined 
without assuming
the existence of
a regular model.

\begin{df}\label{dfdif}
\setcounter{equation}0
Let $f\colon V\to U$
be a finite \'etale
morphism
of regular flat
separated schemes
of finite type over $S$
such that
the generic fiber
$V_K\to U_K$
is tamely ramified with respect
to $K$.

{\rm 1.}
We call
\begin{equation}
D_{V/U}^{\log}=
((V\times_UV\setminus \Delta_V,
\Delta_V))^{\log}
\in F_0G(\partial_{V/U}V)
_{\mathbb Q}
\label{dfdifV}
\end{equation}
the logarithmic different
of $V$ over $U$.
We call
\begin{equation}
d_{V/U}^{\log}=
f_*D_{V/U}^{\log}
\in F_0G(\partial_{V/U}U)
_{\mathbb Q}
\label{dfdisV}
\end{equation}
{\rm the logarithmic discriminant}
of $V$ over $U$.

{\rm 2.}
Let $\sigma$
be an automorphism
of $V$ over $U$
such that
the fixed part $V^{\sigma}$
is empty
and let
$\Gamma_\sigma
\subset V\times_UV$
be the graph of $\sigma$.
We call
\begin{equation}
((\Gamma_\sigma,
\Delta_V))^{\log}
\in F_0G(\partial_{V/U}V)
_{\mathbb Q}
\label{dfLefV}
\end{equation}
{\rm the logarithmic Lefschetz
class}.
\end{df}

We show that
the log different
satisfies a chain rule
and that, 
for a Galois covering,
the logarithmic different
is the sum of
Lefschetz classes.

\begin{lm}\label{lmDs}
\setcounter{equation}0
Let $f\colon V\to U$
be a finite
\'etale morphism
of regular flat
separated schemes
of finite type over $S$
such that
the generic fiber
$V_K\to U_K$
is tamely ramified with respect
to $K$.

{\rm 1.}
Let $U'$ be a finite \'etale scheme
over $U$
such that
the generic fiber
$U'_K\to U_K$
is tamely ramified with respect
to $K$
and let
$g\colon V\to U'$ 
be a finite \'etale
morphism over $U$.
Then, 
we have
\begin{equation}
D^{\log}_{V/U}=
D^{\log}_{V/U'}+
g^*D^{\log}_{U'/U}
\label{eqDchn}
\end{equation}
in $F_0G(\partial_{V/U}V)
_{\mathbb Q}$.

{\rm 2.}
Assume that
$V$ is a $G$-torsor over $U$
for a finite group $G$.
Then we have
\begin{equation}
D_{V/U}^{\log}
=\sum_{\sigma\in G,\
\sigma \neq 1}
((\Gamma_\sigma,
\Delta_V))^{\log}
\label{eqDs}
\end{equation}
in $F_0G(\partial_{V/U}V)
_{\mathbb Q}$.

{\rm 3.}
Assume that
$V$ is a $G$-torsor
for a finite group $G$.
Let $N\subset G$
be a normal subgroup of
$G$
and let
$g\colon V\to V'$ 
be the corresponding
$N$-torsor.
Then, for an element
$\sigma'
\in G'=G/N,\ 
\neq 1$,
we have
$$g^*((\Gamma_{\sigma'},
\Delta_{V'}))^{\log}
=
\sum_{\sigma\in G,
\sigma\mapsto 
\sigma'}
((\Gamma_{\sigma},
\Delta_V))^{\log}$$
in $F_0G(\partial_{V/U}V)
_{\mathbb Q}$.
\end{lm}

{\it Proof.}
1.
It follows from
$V\times_UV
\setminus \Delta_V=
(V\times_{U'}V
\setminus \Delta_V)
\amalg
(g\times g)^{-1}
(U'\times_UU'
\setminus \Delta_{U'})$
and
Proposition \ref{prfun}.1
applied to $U'\times_UU'
\setminus \Delta_{U'}$.

2.
Clear from
$V\times_UV\setminus \Delta_V=
\coprod_{\sigma\in G,\
\sigma \neq 1}
\Gamma_\sigma$.

3.
It follows from
$(g\times g)^{-1}
(\Gamma_{\sigma'})
=
\coprod_{\sigma
\mapsto \sigma'}
\Gamma_{\sigma}$
and
Proposition \ref{prfun}.1.
\qed

\begin{cor}\label{corchn}
\setcounter{equation}0
Let the notation
be as in 
Lemma {\rm \ref{lmDs}.1}.
Let $f'\colon U'
\to U$ denote the morphism
and assume that the map
$g\colon V\to U'$
is of constant degree
$[V:U']$.
Then, 
for the discriminants defined
in Definition {\rm \ref{dfdif}.1},
we have
\begin{equation}
d^{\log}_{V/U}=
f'_*d^{\log}_{V/U'}+
[V:U']\cdot d^{\log}_{U'/U}
\label{eqdchn}
\end{equation}in $F_0G(\partial_{V/U}U)
_{\mathbb Q}$.
\end{cor}

{\it Proof.}
We take
the push-forward
of
(\ref{eqDchn}).
Then, similarly as Lemma \ref{lmff}.2,
we obtain (\ref{eqdchn}).
\qed

\begin{cn}\label{cnsi}
Let $f\colon V\to U$
be a finite \'etale
morphism
of regular flat
separated schemes
of finite type over $S$
such that
the generic fiber
$V_K\to U_K$
is tamely ramified with respect
to $K$.
Let $\sigma$
be an automorphism
of $V$ over $U$
such that
the fixed part $V^{\sigma}$
is empty.

Then, for an integer
$i$ prime to the
order of $\sigma$,
we have
$$((\Gamma_\sigma,
\Delta_V))^{\log}=
((\Gamma_{\sigma^i},
\Delta_V))^{\log}.$$
\end{cn}

\begin{pr}\label{prsi}
Conjecture {\rm \ref{cnsi}}
is true if $\dim V_K\le 1$.
\end{pr}

{\it Proof.} 
By the resolution of singularity for two dimensional schemes, 
regular objects $Y$ of
the category ${\mathcal C}_{V/S}$
are cofinal in ${\mathcal C}_{V/S}$.
Further,
the regular objects $Y$
such that the action of $\sigma$
is extended to an admissible
action on $Y$
are cofinal in ${\mathcal C}_{V/S}$.
Hence, the assertion of Conjecture
\ref{cnsi} follows from Corollary \ref{corsi}.
\qed

\newpage 
\section{Formulas for invariants of
wild ramification}
\label{sfml}

In this section,
we establish formulas
for the invariants
of wild ramification
defined in the
previous section.
We state and prove
the results for
the map (\ref{eqmap}).
However,
they also hold for
the map (\ref{eqmapF})
by the same argument.

We prove an excision formula
Theorem \ref{thmexc}
and a blow-up formula
Proposition \ref{prblup} 
in Section \ref{ssexc}.
We establish some
preliminary formulas
in Section \ref{ssdiv}.
In Section \ref{ssss},
we prove a formula 
Proposition \ref{prcysd} for
some semi-stable families
applying the log Lefschetz trace formula
Theorem \ref{thmlLTF},
which will play a crucial role
in the proof of the conductor formula.

\subsection{Divisors and
projective space bundles}
\label{ssdiv}
The results
in this subsection
will be used
in the proof of
the blow-up formula
and of
the excision formula
in the next subsection.
In Propositions \ref{prU1}
and \ref{prU2},
we compute the log localized
intersection product
of some classes
supported on the inverse image
of a divisor.
In Proposition
\ref{prpbdl} and 
Lemma \ref{lmpbdl},
we give formulas
for a projective space
bundle.

We keep the notation
that
$f\colon
V\to U$ denotes
a finite \'etale morphism
of regular flat separated
schemes of finite
type over $S$
and $n=\dim V_K+1$
and the assumption that
the generic fiber
$V_K\to U_K$
is tamely ramified
with respect to
$K$ (Definition \ref{dftmT}).
Although we also state
corresponding
formulas for finite
\'etale morphism 
$f\colon
V\to U$ 
of smooth separated
schemes of finite
type over $F$,
the proof is similar and easier
and will be omitted.

We prepare to state
Proposition \ref{prU1}.
We consider
a cartesian diagram
\begin{equation}
\begin{CD}
V@>f>> U\\
@VgVV @VVV\\
V_0@>{f_0}>> U_0
\end{CD}
\label{eqU0}
\end{equation}
of regular flat
separated schemes
of finite type over $S$
where $f_0$ is finite \'etale
and the vertical arrows are proper.
We assume
that the generic fiber
$f_{0,K}\colon
V_{0,K}\to U_{0,K}$
is tamely ramified with
respect to $K$.
If the map $g\colon 
V\to V_0$
is birational,
the diagram
\begin{equation}
\begin{CD}
{\rm Gr}^F_nG
(V\times_{U_0}V
\setminus
V\times_{V_0}V)
@>{((\ ,\Delta_V))^{\log}}>>
F_0G(\partial_{V_0/U_0} V)_{\mathbb Q}\\
@A{(g\times g)^*}AA
@VV{g_*}V\\
{\rm Gr}^F_nG
(V_0\times_{U_0}V_0
\setminus
\Delta_{V_0})
@>{((\ ,\Delta_{V_0}))^{\log}}>>
F_0G(\partial_{V_0/U_0} V_0)_{\mathbb Q}
\end{CD}
\label{eqmap00}
\end{equation}
is commutative.

\begin{pr}\label{prU1}
\setcounter{equation}0
Let $f\colon V\to U$
be a finite \'etale morphism
of regular flat separated
schemes of finite
type over $S$
and $n=\dim V_K+1$.
Suppose that
we have a cartesian
diagram
{\rm (\ref{eqU0})}
such that
$V_{0,K}\to U_{0,K}$
is tamely ramified
with respect to $K$.

Let $U_1\subset U$ 
be a regular divisor
and 
$i$ denote the immersion
$U_1\to U$ and its base changes.
Assume that
either $U_1$ is a scheme over
$K$
or a scheme over $F$.
We put
$V_1=V\times_UU_1$.
Then, for
$\Gamma_1
\in {\rm Gr}^F_nG
(V_1\times_{U_0}V_1
\setminus 
V_1\times_{V_0}V_1)$, 
we have
\begin{equation}
((\Gamma_1,
\Delta_V))^{\log}
=
\begin{cases}
-
i_*((\Gamma_1
\cdot c_1({\rm pr}_2^*N_{U_1/U}),
\Delta_{V_1}))^{\log}
&\text{ if }
V_1=V_{1,K}\\
-
i_*(\Gamma_1
\cdot c_1({\rm pr}_2^*N_{U_1/U}),
\Delta_{V_1})^{\log}
&\text{ if }
V_1=V_{1,F}
\end{cases}
\end{equation}
in 
$F_0G
(\partial_{V_0/U_0}
V)_{\mathbb Q}$.
\end{pr}

{\it Proof.}
Let $Y\to Y_0$ be an object
of ${\cal C}_{V\to V_0}$
and let
${\cal D}_0$ be
a finite family
of Cartier divisors
such that
$\Sigma_{V_0/U_0}
^{{\cal D}_0}Y_0=
\Sigma_{V_0/U_0}Y_0.$
Let $Y_1\subset Y$ be
the closure of $V_1$.
By replacing $Y$ by
the blow-up at $Y_1$
if necessary,
we may assume
$Y_1$ is a divisor of $Y$.
We take
an object $Z\to Y$ of
${\cal A}_{V\to U}$.
Let
$g^*V_1
=\sum_je_jW_j$
be the decomposition
by irreducible components
and let
$Z_j$ denote the closure 
of $W_j$.
Replacing $Z$
if necessary,
we may assume that
$\sum_jZ_j$ 
has simple normal crossings
and meets
$D=Z\setminus W$
transversely.

Let $(Y\times_SY)^\sim$
be the log product
with respect to
the pull-back
${\cal D}$
of ${\cal D}_0$
and 
$(Z\times_SZ)^\sim$
be the log product
with respect to
the divisor $D$
with simple normal crossings.
We consider the map
$(\bar g\times\bar g)^\sim
\colon
(Z\times_SZ)^\sim
\to
(Y\times_SY)^\sim$
and its restriction
$g\times g
\colon
W\times_SW
\to
V\times_SV$.
For an irreducible component
$Z_j$, let $i_j\colon 
Z_j\to Z$
be the closed immersion
and $\bar g_j\colon Z_j\to Y_1$
be the restrictions of
$\bar g\colon Z\to Y$.
Let
$g_j\colon W_j
=W\cap Z_j\to V_1$
be the restriction 
of $\bar g_j$.
The intersection 
$D_j=Z_j\cap D$
is a divisor with simple normal crossings.

Let
$(Z_i\times_S Z_j)^\sim$
denote the fiber product
$(Z_i\times_S Z_j)
\times_{Z\times_S Z}
(Z\times_S Z)^\sim$.
Define 
$(g_i\times g_j)^*
\colon
{\rm Gr}^F_nG(V_1
\times_{U_0}V_1
\setminus
V_1
\times_{V_0}V_1)
\to
{\rm Gr}^F_nG(W_i
\times_{U_0}W_j\setminus
W_i\times_{V_0}W_j)$
as the pull-back by
$g_i\times g_j
\colon
W_i \times_SW_j
\to
V_1\times_SV_1$.
Then, by Corollary \ref{cortorE}.1,
we obtain
\begin{equation}
(g\times g)^*
(\Gamma_1)
=
\sum_{i,j}
e_ie_j\cdot
(g_i\times g_j)^*
(\Gamma_1)
\label{eqij}
\end{equation}
in ${\rm Gr}^F_nG(
W_1\times_{U_0}W_1
\setminus
W_1\times_{V_0}W_1)$.

Let $A_0\subset
(Y_0\times_SY_0)
^\sim_{{\cal D}_0}$
be the closure of
$V_0\times_{U_0}V_0
\setminus \Delta_{V_0}$
and let
$A\subset (Z\times_SZ)^\sim$
be the intersection
of the pull-back of $A_0$
with $(Z\times_XZ)^\sim$.
For each $i,j$,
we put
$A_{ij}=A
\cap (Z_i\times_SZ_j)^\sim
\subset 
(Z\times_XZ)^\sim$.
We have
$A\cap (W\times_SW)
=(W\times_{U_0}W)
\setminus 
(W\times_{V_0}W)$
and
$A_{ij}\cap (W_i\times_SW_j)
=(W_i\times_{V_0}W_j)
\setminus 
(W_i\times_{U_0}W_j)$.
We take
$\Gamma_{ij}
\in {\rm Gr}^F_nG
(A_{ij})$
lifting
$(g_i\times g_j)^*
(\Gamma_1)$.
Then, by (\ref{eqij}),
we have
\begin{equation}
((\Gamma_1,
\Delta^{\log}_Y))
=
\frac1{[W:V]}
\sum_{i,j}e_ie_j\cdot
\bar g_*
((\Gamma_{ij},
\Delta_Z^{\log}))_
{(Z\times_SZ)^\sim}.
\label{eqGij}
\end{equation}

We continue the proof
assuming
$U_1=U_{1,K}$.
The proof of 
the other case
$U_1=U_{1,F}$
is similar and omitted.
Let ${\cal I}_{Z_j}
\subset {\cal O}_Z$
be the invertible ideal
defining $Z_j\subset Z$.
For each $i,j$,
we show
\begin{equation}
\bar g_*((\Gamma_{ij},
\Delta^{\log}_Z))
_{(Z\times_SZ)^\sim}
=
-
\bar g_*
i_{i*}\left(((\Gamma_{ii},
\Delta^{\log}_{Z_i}))
_{(Z_i\times_SZ_i)^\sim}
\cdot c_1({\cal I}_{Z_j})
\right)
\label{eqZij}
\end{equation}
in $F_0G
(\Sigma_{V_0/U_0}Y_0\times_{Y_0}Y)_{\mathbb Q}$.
If $i=j$,
the equality
(\ref{eqZij})
follows from
Lemma \ref{lmX1pp}.1.
(In the case
$U_1=U_{1,F}$,
we apply
Lemma \ref{lmX1pp}.2.)

We assume $Z_i\neq Z_j$.
We put $Z_{ij}=
Z_i\times_Z Z_j,
W_{ij}=
W_i\times_WW_j,
A_{i,ij}=
A\cap (Z_i\times_SZ_{ij})^\sim$
and let 
$g_{ij}\colon 
W_{ij}\to V_1$
be the restriction of $g$.
Then, the immersions
$(Z_i\times_SZ_j)^\sim
\to
(Z\times_SZ)^\sim$ 
and
$i_{ij}\colon
Z_{ij}\to Z$ are
regular immersions of
codimension 2.
Hence by Lemma \ref{lmfm1},
we obtain
$$((\Gamma_{ij},
\Delta^{\log}_Z))
_{(Z\times_SZ)^\sim}
=
i_{ij*}
((\Gamma_{ij},
\Delta^{\log}_{Z_{ij}}))
_{(Z_i\times_SZ_j)^\sim}.$$
Further, the immersion
$(Z_i\times_SZ_{ij})^\sim
\to
(Z_i\times_SZ_j)^\sim$
is a regular immersion of
codimension 1.
Hence by Lemma \ref{lmfm1},
we obtain
$$((\Gamma_{ij},
\Delta^{\log}_{Z_{ij}}))
_{(Z_i\times_SZ_j)^\sim}
=
\bigl(\bigl(
(\Gamma_{ij},
(Z_i\times_SZ_{ij})^\sim)
_{(Z_i\times_SZ_j)^\sim}
,
\Delta^{\log}_{Z_{ij}}\bigr)\bigr)
_{(Z_i\times_SZ_{ij})^\sim}.$$
Since both
$(\Gamma_{ij},
(Z_i\times_SZ_{ij})^\sim)
_{(Z_i\times_SZ_j)^\sim}$
and
$(\Gamma_{ii},
(Z_i\times_SZ_{ij})^\sim)
_{(Z_i\times_SZ_i)^\sim}
\in {\rm Gr}^F_{n-1}G(A_{i,ij})$
are liftings of
$(g_i\times g_{ij})^*
(\Gamma_1)$,
we have
$$\bigl(\bigl(
(\Gamma_{ij},
(Z_i\times_SZ_{ij})^\sim)
_{(Z_i\times_SZ_j)^\sim}
,
\Delta^{\log}_{Z_{ij}}\bigr)\bigr)
_{(Z_i\times_SZ_{ij})^\sim}
=
\bigl(\bigl(
(\Gamma_{ii},
(Z_i\times_SZ_{ij})^\sim)
_{(Z_i\times_SZ_i)^\sim}
,
\Delta^{\log}_{Z_{ij}}\bigr)\bigr)
_{(Z_i\times_SZ_{ij})^\sim}$$
similarly as Proposition \ref{prprod}.
By applying the associativity
Lemmas \ref{lmfm2} and \ref{lmfm4}
to the diagram
$$\begin{CD}
(Z_i\times Z_{ij})^\sim @<<< Z_{ij}\\
@VVV @VVV\\
(Z_i\times Z_i)^\sim @<<< Z_i,
\end{CD}$$
we obtain
$$
\bigl(\bigl(
(\Gamma_{ii},
(Z_i\times_SZ_{ij})^\sim)
_{(Z_i\times_SZ_i)^\sim}
,
\Delta^{\log}_{Z_{ij}}\bigr)\bigr)
_{(Z_i\times_SZ_{ij})^\sim}
=
\bigl(
((\Gamma_{ii},
\Delta^{\log}_{Z_i}))
_{(Z_i\times_SZ_i)^\sim},
Z_{ij}
\bigr)_{Z_i}.$$
Thus, we obtain
$$((\Gamma_{ij},
\Delta^{\log}_Z))
_{(Z\times_SZ)^\sim}
=
i_{ij*}
\bigl(
((\Gamma_{ii},
\Delta^{\log}_{Z_i}))
_{(Z_i\times_SZ_i)^\sim},
Z_{ij}
\bigr)_{Z_i}
=
-i_{i*}\left(
((\Gamma_{ii},
\Delta^{\log}_{Z_i}))
_{(Z_i\times_SZ_i)^\sim}
\cdot c_1({\cal I}_{Z_j})\right)$$
and 
the equality
(\ref{eqZij}) is proved.

Therefore, the sum in
the right hand side
of (\ref{eqGij})
is equal to
\begin{equation}
-
\sum_ie_i\cdot
i_*
\bar g_{i*}\left(
((\Gamma_{ii},
\Delta_{Z_i}^{\log}))_
{(Z_i\times_SZ_i)^\sim}
\cdot
\textstyle{\sum_je_j}
c_1({\cal I}_{Z_j})
\right).
\label{eqpar}
\end{equation}
Since $g^*V_1=
\sum_je_jW_j$
as a divisor of $W$,
the restriction
$\sum_je_j
c_1({\cal I}_{Z_j})|_{W_i}$
is equal to
$c_1(\bar g_i^*N_{Y_1/Y})|_{W_i}.$
Hence by Proposition 
\ref{prprod},
we have
\begin{eqnarray*}
&
((\Gamma_{ii},
\Delta_{Z_i}^{\log}))_
{(Z_i\times_SZ_i)^\sim}
\cdot
\textstyle{\sum_je_j}
c_1({\cal I}_{Z_j})
=
((\Gamma_{ii}
\cdot
\textstyle{\sum_je_j}
c_1({\rm pr}_2^*
{\cal I}_{Z_j}),
\Delta_{Z_i}^{\log}))_
{(Z_i\times_SZ_i)^\sim}\\
&
=
((\Gamma_{ii}
\cdot{\rm pr}_2^*
c_1(\bar g_i^*N_{Y_1/Y}),
\Delta_{Z_i}^{\log}))_
{(Z_i\times_SZ_i)^\sim}.
\end{eqnarray*}
If $Z_i\to Y_1$ is surjective,
it defines an object of
${\mathcal A}_{V_1\to U_1}$
and  we have
$$
\bar g_{i*}
((\Gamma_{ii}
\cdot{\rm pr}_2^*
c_1(\bar g_i^*N_{Y_1/Y}),
\Delta_{Z_i}^{\log}))_
{(Z_i\times_SZ_i)^\sim}=
[W_i\colon V_1]
\cdot
((\Gamma_1
\cdot{\rm pr}_2^*
c_1(N_{Y_1/Y}),
\Delta_{Y_1}^{\log})).
$$
We show that,
if $Z_i\to Y_1$
is not surjective,
the left hand side is $0$.
By replacing
$Z_i$ by an alteration,
we may assume that
there is an object
$Z_0\to Y_1$ of
${\mathcal A}_{V_1\to U_1}$
and that
$Z_i\to Y_1$
factors through $\pi\colon Z_i
\to Z_0$.
Since $\pi_*
\pi^*$
is multiplication by
${\rm rank}\
R\pi_*{\mathcal O}_{Z_i}=0$
on $F_0G(\Sigma_{V_0/U_0}Y
\times_YZ_0)$,
the assertion follows
from Corollary \ref{cormap00}.
Thus, by
$\sum_ie_i[W_i:V_1]
=
[W:V]$,
(\ref{eqpar})
is equal to
the right hand side of (\ref{prU1}.1).
\qed

We consider a regular
divisor $U_1\subset U$
as in Proposition
\ref{prU1}.
Let $(U\times_UU)^\sim$
denote the log product
with respect to
the Cartier divisor $U_1$.
It is the union of
$U$ with
$E={\mathbf G}_{m,U_1}$
meeting at $U_1$.
It is canonically 
identified
with the fiber product
$(U\times_SU)^\sim
\times_{U\times_SU}
\Delta_U$.
The closed subscheme
$E={\mathbf G}_{m,U_1}
\subset (U\times_UU)^\sim$ 
is the inverse
image of
$U_1\subset U=\Delta_U$
by the canonical map
$(U\times_UU)^\sim\to U$,
as in Lemma \ref{lmEi}.1.
Let $(V\times_UV)^\sim
=(V\times_UV)
\times_U
(U\times_UU)^\sim$
denote the log product
with respect to
$V_1=V\times_UU_1$.
Then, we define the localized
log intersection product
\setcounter{equation}0
\addtocounter{thm}1
\begin{equation}
\begin{CD}
((\ ,\Delta_V))^{\log}
\colon {\rm Gr}^F_nG
((V\times_UV)^\sim
\setminus
(\Delta_V\times_U
(U\times_UU)^\sim))
@>>>
F_0G(\partial_{V/U} V)_{\mathbb Q}
\end{CD}
\label{eqmapE}
\end{equation}
similarly as follows,
in order to state
Proposition \ref{prU2}.

We consider an
object $Y$ of ${\cal C}_{V/S}$
and a finite family ${\cal D}$
of Cartier divisors of
$Y$ such that
$\Sigma_{U/V}Y=
\Sigma^{\cal D}_{U/V}Y$.
Let $Y_1\subset Y$ denote
the closure of
$V_1$ and
${\cal D}_1$ denote
the restriction of ${\cal D}$.
By replacing $Y$ by
the blow-up at $Y_1$
if necessary,
we assume
$Y_1$ is a divisor of $Y$.
We also consider an object
$\bar g\colon Z\to Y$ of ${\cal A}_{V\to U}$.
We assume that
it also define
an object
of ${\cal A}_{(V\setminus V_1)\to 
(U\setminus U_1)}$.
In particular,
$W_0=\bar g^{-1}(V\setminus V_1)$
is the complement
of a divisor 
$D'\supset D$ of $Z$
with simple normal crossings.
Let $Z\to X$ be a proper morphism
of schemes over $S$
such that $X$ contains
$U$
as the complement of
a Cartier divisor $B$.

We define the log products
$(Y\times_SY)^\sim$
and 
$(Y\times_SY)^\approx$
with respect to
${\cal D}$ 
and to 
the union ${\cal D}'$ of
${\cal D}$ with $Y_1$
respectively.
Similarly, we define the log products
$(Z\times_SZ)^\sim$
and 
$(Z\times_SZ)^\approx$
with respect to
$D$ and to $D'$.
They form commutative
diagrams 
$$
\begin{CD}
(Y\times_SY)^\approx
@<{(\bar g\times
\bar g)^\approx}<<
(Z\times_SZ)^\approx\\
@VVV@VVV\\
(Y\times_SY)^\sim
@<<<
(Z\times_SZ)^\sim
\end{CD}\qquad
\begin{CD}
(V\times_SV)^\sim
@<{(g\times g)^\sim}<<
(W\times_SW)^\sim\\
@VqVV@VVV\\
V\times_SV
@<<<
W\times_SW
\end{CD}
$$
where the right square is
obtained by taking
the base change to $V\times_SV$.

We define a
closed subset
$A\subset 
(Y\times_SY)^\approx$
to be the closure of
$(V\times_UV)^\sim
\setminus
(\Delta_V\times_U
(U\times_UU)^\sim)
\subset
(V\times_SV)^\sim$
and put
$A_Z=
(\bar g\times
\bar g)^{\approx-1}(A)
\cap 
(Z\times_XZ)^\approx$.
Here
$(Z\times_XZ)^\approx$
denote the log product
defined with respect to
$D'$ and $B$.
For $\Gamma
\in {\rm Gr}^F_nG
((V\times_UV)^\sim
\setminus
(\Delta_V\times_U
(U\times_UU)^\sim))$,
we take an element
$\widetilde \Gamma
\in {\rm Gr}^F_nG(A_Z)$
lifting
the pull-back
$(g\times g)^{\sim*}
(\Gamma)$
by
$(g\times g)^{\sim}\colon
(W\times_SW)^\sim
\to
(V\times_SV)^\sim$.
Then,
$((\Gamma,
\Delta_Y^{\log}))$
is defined as
$\bar g_*
((\overline \Gamma,
\Delta_Z^{\log}))
_{(Z\times_SZ)^\approx}$
divided by
$[Z:Y]$.
The map
(\ref{eqmapE})
is defined as
the projective limit.

By the associativity
Lemma \ref{lmfm4}
applied to
$(Z\times_SZ)^\sim
\gets
(Z\times_SZ)^\approx
\gets Z$,
we obtain a commutative
diagram
\begin{equation}
\begin{CD}
{\rm Gr}^F_nG
(V\times_UV
\setminus
\Delta_V)
@>{((\ ,\Delta_V))^{\log}}>>
F_0G(\partial_{V/U} V)_{\mathbb Q}\\
@V{q^*}VV
@|\\
{\rm Gr}^F_nG
((V\times_UV)^\sim
\setminus
(\Delta_V\times_U
(U\times_UU)^\sim))
@>{((\ ,\Delta_V))^{\log}}>>
F_0G(\partial_{V/U} V)_{\mathbb Q}
\end{CD}
\label{eqmapE0}
\end{equation}
for the pull-back
$q^*$ by the projection
$q\colon 
(V\times_SV)^\sim
\to V\times_SV$.

We put
$U_0=U\setminus U_1,
V_0=V\setminus V_1$ and
we consider the diagram
\begin{equation}
\begin{CD}
{\rm Gr}^F_nG
(V\times_UV
\setminus
\Delta_V)
@.
\\
@VVV @.\\
{\rm Gr}^F_nG
((V\times_UV)^\sim
\setminus
(\Delta_V\times_U
(U\times_UU)^\sim))
@>{((\ ,\Delta_V))^{\log}}>>
F_0G(\partial_{V/U} V)_{\mathbb Q}
\\
@VVV @VVV\\
{\rm Gr}^F_nG
(V_0\times_{U_0}V_0
\setminus
\Delta_{V_0})
@>{((\ ,\Delta_{V_0}))^{\log}}>>
F_0G(\partial_{V/U} V_0)_{\mathbb Q}
\end{CD}
\label{eqUUsim}
\end{equation}
where the 
upper left vertical arrow
is induced by
the immersion
$V\times_UV\to
(V\times_UV)^\sim$
and the lower 
left vertical arrow is
the restriction.
Since the log product
$(Z\times_XZ)^\sim$
is defined with respect to $B$
whose complement is $U$,
the closed subset
$A_Z$ need {\em not} satisfy
the condition (B) in
Proposition \ref{prprod}
with respect to
the components of
$D'\setminus D$.
Consequently,
the square is {\em not} necessarily commutative.
However,
the compositions with the upper vertical
arrow form a commutative
diagram
since the image
of the composition
$V\times_UV\to
(V\times_UV)^\sim
\to 
(X'\times_{X'}X')^\sim$
lies in the image
of the log diagonal
$X'\to (X'\times_{X'}X')^\sim$
for a compactification $X'$ 
of $U$ containing $U_0$
as the complement of
a Cartier divisor $B'$
extending $U_1$
and the log product $(X'\times_{X'}X')^\sim$
defined with respect to $B'$.

\addtocounter{thm}{-1}
\begin{pr}\label{prU2}
Let $f\colon V\to U$
be a finite \'etale morphism
of regular flat separated
schemes of finite
type over $S$
and $n=\dim V_K+1$.
We assume that
$V_K\to U_K$
is tamely ramified
with respect to $K$.

Let $U_1\subset U$ 
be a regular divisor
and 
$i\colon
U_1\to U$ denote the immersion.
Assume either
$U_1$ is a scheme over
$K$ or a scheme over $F$.
Let
$(V\times_SV)^\sim$
be the log product
with respect to
the Cartier divisor 
$V_1=V\times_UU_1$
of $V$.
Let $q_1\colon E
={\mathbf G}_{m,U_1}\to U_1$
be the projection
and $q_1\colon E\times_{U_1}
(V_1\times_{U_1}V_1)
\to V_1\times_{U_1}V_1$
also denote
the base change.
We regard
$E\times_{U_1}
(V_1\times_{U_1}V_1)$
as a closed subscheme
of $(V\times_SV)^\sim$
as above.

Then, for
$\Gamma_1
\in {\rm Gr}^F_{n-1}G
(V_1\times_{U_1}V_1
\setminus \Delta_{V_1})$, 
the product
{\rm (\ref{eqmapE})}
satisfies
\begin{equation}
((q_1^*\Gamma_1,
\Delta_V))^{\log}
=
\begin{cases}
i_*((\Gamma_1,
\Delta_{V_1}))^{\log}
&\text{ if }
V_1=V_{1,K}\\
i_*(\Gamma_1,
\Delta_{V_1})^{\log}
&\text{ if }
V_1=V_{1,F}
\end{cases}
\end{equation}
in 
$F_0G
(\partial_{V/U}
V)_{\mathbb Q}$.
\end{pr}

{\it Proof.}
We keep the notation
in the definition of
(\ref{eqmapE}) above.
We put
$g^*V_1
=\sum_{j\in J}
e_jW_j$.
For each irreducible component
$W_j$, let $Z_j$
be the closure,
$i_j\colon 
Z_j\to Z$
be the closed immersion
and $\bar g_j\colon Z_j\to Y_1$
be the restrictions of
$\bar g\colon Z\to Y$.
Let $g_j\colon W_j\to V_1$
be the restriction 
of $\bar g_j$.
Let ${\cal D}_Z=
(D_k)_{k\in I}$
and ${\cal D}'_Z=
(D_k)_{k\in I'}$
be the families
of irreducible components
of $D=Z\setminus W
\subset 
D'=Z\setminus W_0$
indexed by $I\subset 
I'=I\amalg J$
respectively.

For $j\in J$, let 
$(Z_j\times_S Z_j)^\sim$
and 
$(Z_j\times_S Z_j)^\approx$
be the log product
with respect to
the families
${\cal D}_j=
(D_k\cap Z_j)_{k\in I}$
and 
${\cal D}'_j=
(D_k\cap Z_j)_{k\in I',D_k\neq Z_j}$
respectively.
Let $E_j\subset
(Z\times_S Z)^\approx$
denote the inverse
image of
$Z_j\subset Z$
by either of
the projections
$(Z\times_S Z)^\sim
\to Z$.
Then,
the canonical map
$(Z_j\times_S Z_j)^\approx
\to
(Z_j\times_S Z_j)^\sim$
is of finite tor-dimension by
Corollary \ref{cortdf}
and 
$E_j$ is flat over
$(Z_j\times_S Z_j)^\approx$
by Lemma \ref{lmEi}.
Hence,
the canonical map
$\bar q_j\colon
E_j
\to 
(Z_j\times_S Z_j)^\sim$
is of finite tor-dimension.
Let
$q_j\colon
E_j^\circ
\to 
W_j\times_S W_j$
be the base change
of $\bar q_j$.
We consider the commutative diagram
$$\begin{CD}
E_1=(V\times_SV)^\sim
\times_{V\times_SV}
(V_1\times_SV_1)
@<<< E_j^\circ @>{\subset}>>E_j\\
@V{q_1}VV@V{q_j}VV@VV{\bar q_j}V\\
V_1\times_SV_1
@<{g_j\times g_j}<<
W_j\times_SW_j
@>{\subset}>>
(Z_j\times_SZ_j)^\sim
\end{CD}$$
where the right horizontal
arrows are open immersions.

Let 
$(g_j\times g_j)^*\colon
{\rm Gr}^F_{n-1}G
(V_1\times_{U_1}V_1
\setminus \Delta_{V_1})
\to
{\rm Gr}^F_{n-1}G
(W_j\times_{U_1}W_j
\setminus 
W_j\times_{V_1}W_j)$
denote
the pull-back by
$g_j\times g_j
\colon
W_j \times_SW_j\to
V_1\times_SV_1$.
Then, by Corollary \ref{cortorE}.2,
we obtain
\begin{equation}
(g\times g)^{\sim *}
(q_1^*\Gamma_1)
=
\sum_j
e_j\cdot
q_j^*(g_j\times g_j)^*
(\Gamma_1)
\label{eqqi}
\end{equation}
in ${\rm Gr}^F_nG(
(W_1
\times_{U_0}W_1)^\sim
\setminus
(W_1
\times_{V_0}W_1)^\sim)$.

Let 
$A_{\cal D}\subset
(Y\times_SY)^\sim_{\cal D}$
be the closure of
$V\times_UV\setminus \Delta_V$
and let
$A\subset
(Z\times_SZ)^\approx$
be the intersection
of the pull-back of $A_{\cal D}$ with
$(Z\times_XZ)^\approx
=
(Z\times_SZ)^\approx
\times_
{(X\times_SX)^\sim}X$.
For each $Z_j$,
let $A_j\subset
(Z_j\times_SZ_j)^\sim$
be the intersection
of the pull-back of $A$ with
$(Z_j\times_XZ_j)^\sim
=
(Z_j\times_SZ_j)^\sim
\times_
{(X\times_SX)^\sim}X$.
We have
$A\cap (W\times_SW)^\sim
=(W\times_VW)^\sim
\setminus 
(W\times_UW)^\sim$
and
$A_j\cap (W_j\times_SW_j)
=(W_j\times_{U_1}W_j)
\setminus 
(W_i\times_{V_1}W_j)$.
We take
$\Gamma_j
\in {\rm Gr}^F_nG(A_j)$
lifting
$(g_j\times g_j)^*
(\Gamma_1)$.
Then, by (\ref{eqqi}),
we have
\begin{equation}
((q_1^*\Gamma_1,
\Delta_Y))^{\log}
=
\frac1{[W:V]}
\sum_je_j\cdot
\bar g_*
((\bar q_j^*\Gamma_j,
\Delta_Z^{\log}))_
{(Z\times_SZ)^\approx}.
\label{eqGi}
\end{equation}

We continue the proof
assuming
$U_1=U_{1,K}$.
The proof of 
the other case
$U_1=U_{1,F}$
is similar using
Lemma \ref{lmX1p}.2.
By Lemma \ref{lmX1p}.1,
we have
$((\bar q_j^*\Gamma_j,
\Delta_Z^{\log}))_
{(Z\times_SZ)^\approx}
=
i_{j*}
(((\bar q_j^*\Gamma_j,
\Delta^{\log}_{Z_j}))
_{(Z_j\times_SZ_j)^\sim},
1_{Z_j})
_{{\mathbf G}_{m, Z_j}}$.
Since $\bar q_j\colon
E_j\to 
(Z_j\times_SZ_j)^\sim$
is of finite tor-dimension,
we have
$((\bar q_j^*\Gamma_j,
\Delta^{\log}_{Z_j}))
_{(Z_j\times_SZ_j)^\sim}
=
\bar q_j^*((\Gamma_j,
\Delta^{\log}_{Z_j}))
_{(Z_j\times_SZ_j)^\sim}$
by the associativity
Lemmas  \ref{lmfm2}
and  \ref{lmfm4}
where
$\bar q_j\colon
{\mathbf G}_{m,Z_j}
\to Z_j$ in the right hand side
denotes the restriction of
$\bar q_j\colon
E_j
\to (Z_j\times_SZ_j)^\sim$.
Thus we obtain
$((\bar q_j^*\Gamma_j,
\Delta_Z^{\log}))_
{(Z\times_SZ)^\sim}
=
i_{j*}((\Gamma_j,
\Delta^{\log}_{Z_j}))
_{(Z_j\times_SZ_j)^\sim}$
and
$$\bar g_*
((\bar q_j^*\Gamma_j,
\Delta_Z^{\log}))_
{(Z\times_SZ)^\approx}
=
\bar g_{j*}
((\Gamma_j,
\Delta_{Z_j}^{\log}))_
{(Z_j\times_SZ_j)^\sim}$$
for the right hand side
of (\ref{eqGi}).
The right hand side is equal to
$[W_j:V_1]\cdot
{\rm pr}_{Y_1}
((\Gamma_1,
\Delta_{V_1}))^{\log}$
if $Z_j\to Y_1$
is generically finite
and is $0$ if otherwise
similarly as at the end of
the proof of Proposition \ref{prU1}.
Therefore,
the sum of the right hand side
of (\ref{eqGi})
is equal to
$$
\sum_je_j
[W_j:V_1]\cdot
{\rm pr}_{Y_1}
((\Gamma_1,
\Delta_{V_1}))^{\log}.$$
Since
$\sum_je_j
[W_j:V_1]=
[W:V]$,
the assertion follows.
\qed

\begin{pr}\label{prpbdl}
\setcounter{equation}0
{\rm 1.}
Let $f\colon V\to U$
be a finite \'etale morphism
of regular flat separated
schemes of finite
type over $S$
and $n=\dim V_K+1$.
We assume that
$V_K\to U_K$
is tamely ramified
with respect to $K$.
Let ${\cal E}$
be a locally free 
${\cal O}_U$-module of rank $c$,
$p\colon
P={\mathbf P}({\cal E})\to U$
be the 
associated
${\mathbf P}^{c-1}$-bundle
and $P_V=P\times_UV$
be the base change.
Let $\Gamma\subset
V\times_UV\setminus \Delta_V$
be  an open and closed
subscheme and
we regard
$\Gamma_P=\Gamma\times_UP$
as an open and closed
subscheme of
$P_V\times_PP_V
\setminus \Delta_{P_V}
=
(V\times_UV\setminus \Delta_V)
\times_UP$.

Then, we have
\setcounter{equation}0
\begin{equation}
p_*((\Gamma_P,\Delta_{P_V}))^{\log}
=
c\cdot
((\Gamma,\Delta_V))^{\log}
\label{eqpbdl}
\end{equation}
in $F_0G(\partial_{V/U} V)
_{\mathbb Q}$.

{\rm 2.}
Let the assumption
be the same as in {\rm 1.}\
except that we assume
$f\colon V\to U$
is a finite \'etale morphism
of smooth separated
schemes of finite
type over $F$
and $n=\dim V$.

Then, we have
\begin{equation}
p_*(\Gamma_P,
\Delta_{P_V})^{\log}
=c\cdot
(\Gamma,\Delta_V)^{\log}
\label{eqpbdl2}
\end{equation}
in $CH_0(\partial_{V/U} V)
_{\mathbb Q}$.
\end{pr}

\noindent{
\it Proof.}
1. By the flattening
theorem \cite{GR},
there exist
a proper scheme
$X$ over $S$
containing
$U$ as a dense
open subscheme
and 
a locally free
${\cal O}_X$-module
${\cal E}_X$ of rank $c$
extending ${\cal E}$.
Replacing $X$ by a blow-up,
we may assume 
$U$ is the complement
of a Cartier divisor $B$.

Let $Y$ be an object
of ${\cal C}_V$
and ${\cal D}$
be a finite family
of Cartier divisors
of $Y$ such that
$\Sigma_{U/V}Y=
\Sigma_{U/V}^{\cal D}Y$.
Let $Z\to Y$
be an object
of ${\cal A}_{V\to U}$.
Let $P_Z
=Z\times_XP_X$
denote the base change
of the projective
space bundle $P_X=
{\mathbf P}({\cal E}_X)$.
Let 
$A_{\cal D}\subset
(Y\times_SY)^\sim_{\cal D}$
be the closure of
$V\times_UV\setminus \Delta_V$
and let
$A\subset
(Z\times_SZ)^\sim$
denote the intersection
$(\bar g\times\bar g)^{\sim -1}
(A_{\cal D})
\cap (Z\times_XZ)^\sim$
as in Proposition \ref{prmap}.
We regard
$A_P=A\times_XP_X$
as a closed subscheme
of
$(P_Z\times_SP_Z)^\sim
=(Z\times_SZ)^\sim
\times_{X\times_SX}(P_X\times_SP_X)$
by the diagonal maps
$X\to X\times_SX$ and
$P_X\to P_X\times_SP_X$.
Then, 
$((\Gamma,\Delta_V))^{\log}$
is defined using the image
of
$((\widetilde \Gamma,
\Delta_Z^{\log}))_
{(Z\times_SZ)^\sim}$
by taking
a lifting 
$\widetilde \Gamma
\in F_nG(A)$
of $(g\times g)^*\Gamma$.
The product
$((\Gamma_P,
\Delta_{P_V}))^{\log}$
is defined using the image
of
$((p^*\widetilde \Gamma,
\Delta_{P_Z}^{\log}))_
{(P_Z\times_SP_Z)^\sim}$
where $p^*\colon
F_nG(A)\to
F_{n+c-1}G(A_P)$
denotes the pull-back.

We apply the associativity
formula,
Lemma \ref{lmfm2},
to $A_P\to
(P_Z\times_SP_Z)^\sim
\gets
P_Z\times_ZP_Z
\gets
P_Z.$
Since a projective space bundle
$P_Z$ is smooth over $Z$,
the diagonal
$P_Z\to 
P_Z\times_ZP_Z$
is a regular immersion
and hence is of
finite tor-dimension.
By applying
Lemma \ref{lmfm2},
we obtain
\begin{equation}
((\Gamma_P,\Delta_{P_Z}^{\log}))_
{(P_Z\times_SP_Z)^\sim}
=
(((\Gamma_P,\Delta_Z^{\log}))_
{(Z\times_SZ)^\sim},
\Delta_{P_Z})_
{P_Z\times_ZP_Z}.
\label{eqPZ}
\end{equation}
Since the projection
$p_Z\colon
P_Z\to Z$ is smooth,
we have
$((\Gamma_P,\Delta_Z^{\log}))_
{(Z\times_SZ)^\sim}
=p_Z^*((\Gamma,\Delta_Z^{\log}))_
{(Z\times_SZ)^\sim}$.
Since
$(\Delta_{P_Z},
\Delta_{P_Z})_
{P_Z\times_ZP_Z}
=
(-1)^{c-1}c_{c-1}
(\Omega^1_{P_Z/Z})$,
the right hand side
of (\ref{eqPZ}) is
equal to
$$p_Z^*((\Gamma,\Delta_Z^{\log}))_
{(Z\times_SZ)^\sim}
\cdot
(-1)^{c-1}c_{c-1}
(\Omega^1_{P_Z/Z})$$
by the excess intersection formula
for the usual intersection product.
Since
$\deg (-1)^{c-1}c_{c-1}
(\Omega^1_{P_Z/Z})=c$,
by the projection formula,
we obtain
$$p_{Z*}
((\Gamma_P,\Delta_{P_Z}^{\log}))_
{(P_Z\times_SP_Z)^\sim}
=
c\cdot ((\Gamma,\Delta_Z^{\log}))_
{(Z\times_SZ)^\sim}$$
and the assertion follows.

We also omit
the similar and easier proof
of 2.
\qed

\begin{lm}\label{lmpbdl}
{\rm 1.}
Let $f\colon V\to U$,
${\cal E}$,
$p\colon
P={\mathbf P}({\cal E})\to U$,
$\Gamma
\subset V\times_UV\setminus 
\Delta_V$ etc.\
be the same as in
Proposition 
{\rm\ref{prpbdl}}.
We consider
the pull-back
$(p\times p)^*\Gamma$
as an open and
closed subscheme
of $P_V\times_UP_V\setminus 
P_V\times_VP_V$.
For an integer $m$,
we put $c_{c,m}=
\deg c_{c-1}(\Omega^1_{
{\mathbf P}^{c-1}}(m)).$

Then, we have
\setcounter{equation}0
\begin{equation}
p_*(((p\times p)^*
\Gamma\cdot c_{c-1}(
{\rm pr}_1^*
\Omega^1_{P/U}(m)
\otimes
{\rm pr}_2^*{\cal O}(m'))
,\Delta_{P_V}))^{\log}
=
c_{c,m+m'}\cdot
((\Gamma,\Delta_V))^{\log}
\label{eqbdl}
\end{equation}
in $F_0G(\partial_{V/U} V)
_{\mathbb Q}$.

{\rm 2.}
Let the assumption
be the same as in {\rm 1.}\
except that
$f\colon V\to U$
is a finite \'etale
morphism of
separated smooth
schemes of finite type over $F$.

Then, we have
\begin{equation}
p_*((p\times p)^*
\Gamma\cdot c_{c-1}(
{\rm pr}_1^*
\Omega^1_{P/U}(m)
\otimes
{\rm pr}_2^*{\cal O}(m'))
,\Delta_{P_V})^{\log}
=c_{c,m+m'}\cdot
(\Gamma,\Delta_V)^{\log}
\label{eqbdl2}
\end{equation}
in $CH_0(\partial_{V/U} V)
_{\mathbb Q}$.
\end{lm}

An elementary calculation
shows 
\begin{equation}
c_{c,m}=
\begin{cases}
\dfrac 1m((m-1)^c-(-1)^c)
&
\text{ for }m\neq 0\\
(-1)^{c-1}c
&
\text{ for }m= 0.
\end{cases}
\label{eqcnm}
\end{equation}

\noindent{
\it Proof.}
1. We keep the notation
$X,{\cal E}_X,P_Z$,
etc.\ in the proof
of Proposition
\ref{prpbdl} above.
Since $P_Z$
is smooth over $Z$,
similarly as above,
we obtain
\begin{eqnarray*}
&&(((p\times p)^*
\Gamma\cdot 
c_{c-1}(
{\rm pr}_1^*
\Omega^1_{P/U}(m)
\otimes
{\rm pr}_2^*{\cal O}(m'))
,\Delta^{\log}_{P_Z}))_
{(P_Z\times_SP_Z)^\sim}\\
&=&
(((p\times p)^*
\Gamma
,\Delta^{\log}_{P_Z}))
_{(P_Z\times_SP_Z)^\sim}
\cdot 
c_{c-1}(\Omega^1_{P_Z/Z}(m+m'))\\
&=&
p_Z^*
((\Gamma
,\Delta_Z))
_{(Z\times_SZ)^\sim}
\cdot 
c_{c-1}(\Omega^1_{P_Z/Z}(m+m')).
\end{eqnarray*}
By applying the 
projection formula,
we obtain
(\ref{eqbdl}).

We also omit
the similar and easier proof
of 2.
\qed

\subsection{Excision formula}
\label{ssexc}

We keep the notation
$f\colon V\to U$
etc.\ as in the
previous subsection.
We prove the excision formula,
Theorem \ref{thmexc}.
We begin with
the following 
blow-up formula.

\begin{pr}\label{prblup}
\setcounter{equation}0
Let $f\colon V\to U$
be a finite \'etale morphism
of regular flat separated
schemes of finite
type over $S$
such that
$V_K\to U_K$
is tamely ramified with respect to
${\rm Spec}\ K$.
We put $n=\dim U_K+1$.
Let $U_1\subset U$
be a regular closed subscheme
of codimension $c$
and $p\colon U'\to U$
be the blow-up at $U_1$.
Assume either
$U_1$ is a scheme over
$K$ or a scheme over $F$.
We consider the
cartesian diagram
$$\begin{CD}
V'@>p>>V@<i<<V_1\\
@VVV @VVV@VVV\\
U'@>p>>U@<i<<U_1
\end{CD}$$
where we use the same letters
to denote the base change.

Let $\Gamma\subset
V\times_UV\setminus \Delta_V$
be an open and closed
subscheme and
we regard
$\Gamma'=\Gamma\times_UU'$
and
$\Gamma_1=\Gamma\times_UU_1$
as open and closed
subschemes of
$V'\times_{U'}V'
\setminus \Delta_{V'}$
and of
$V_1\times_{U_1}V_1
\setminus \Delta_{V_1}$
respectively.

Then, we have
\begin{equation}
((\Gamma,\Delta_V))^{\log}
=p_*((\Gamma',\Delta_{V'}))^{\log}
+
(c-1)
\begin{cases}
i_*((\Gamma_1,\Delta_{V_1}))^{\log}
&\text{ if } U_1=U_{1,K}\\
i_*(\Gamma_1,\Delta_{V_1})^{\log}
&\text{ if } U_1=U_{1,F}
\end{cases}
\label{eqblup}
\end{equation}
in $F_0G(\partial_{V/U} V)_{\mathbb Q}$.
\end{pr}

\noindent{
\it Proof.}
We consider the pull-back
$(p\times p)^*\Gamma$
by $p\times p\colon
V'\times_S V'\to
V\times_S V$.
We have
\begin{equation}
((\Gamma,
\Delta_V))^{\log}
=p_*
(((p\times p)^*\Gamma,
\Delta_{V'}))^{\log},
\label{eqUU'}
\end{equation}
by the commutative diagram
(\ref{eqmap00}).
We compute
$(p\times p)^*\Gamma$.
Note that
$p_1\colon V'_1
\to V_1$
is a ${\mathbf P}^{c-1}$-bundle.
The morphism $U'\to U$
of regular scheme
is locally of complete intersection.
Since
$(p\times p)^*[\Delta_U]
=[{\cal O}_{U'}
\otimes^L_{{\cal O}_U}
{\cal O}_U
\otimes^L_{{\cal O}_U}
{\cal O}_{U'}]
=\sum_i(-1)^i
[{\cal T}or^{{\cal O}_U}_i
({\cal O}_{U'},
{\cal O}_{U'})]$,
by applying
Corollary \ref{corbp}
to the blow-up $U'\to U$,
we obtain
\begin{equation}
(p\times p)^*[\Delta_U]
-[\Delta_{U'}]
=\sum_{i=1}^{c-1}
(-1)^{i-1}
\sum_{j=1}^i
[{\rm pr}_1^*
\Omega_{U'_1/U_1}^i(i)
\otimes 
{\rm pr}_2^*
N_{U'_1/U'}^{\otimes -j}]
\label{eqSigma}
\end{equation}
in ${\rm Gr}^F_nG(U'\times_UU')$,
where
${\rm pr}_k\colon
U'_1\times_{U_1}U'_1
\to U'_1$
denote the projections.
Let $\Sigma$
denote the right hand
side of (\ref{eqSigma}).
We will use the computation
\begin{eqnarray}
&&
\Sigma\cdot 
c_1({\rm pr}_2^*
N_{U'_1/U'})
=
\sum_{i=0}^{c-1}
(-1)^i
[{\rm pr}_1^*
\Omega_{U'_1/U_1}^i(i)]
\otimes 
\bigl([{\rm pr}_2^*
N_{U'_1/U'}^{\otimes -i}]
-[{\cal O}_{U'_1\times_{U_1}U'_1}]
\bigr)
\label{eqSig}
\\
&=&
(-1)^{c-1}\cdot
\Bigl(c_{c-1}(
{\rm pr}_1^*
\Omega_{U'_1/U_1}^1(1)
\otimes 
{\rm pr}_2^*
{\cal O}(-1))
-
c_{c-1}(
{\rm pr}_1^*
\Omega_{U'_1/U_1}^1(1))
\Bigr)
\nonumber
\end{eqnarray}
in ${\rm Gr}^F_{n-1}G(U'_1\times_{U_1}U'_1)$
that follows from
$N_{V'_1/V'}=
{\cal O}(1)$.

By
$(p\times p)^*[\Delta_U]
-[\Delta_{U'}]=\Sigma$,
we have
$(p\times p)^*\Gamma
=\Gamma'+
(p_1\times p_1)^*\Gamma_1
\cdot \Sigma$
since $\Gamma$ is flat
over $U$.
Thus, by
(\ref{eqUU'}),
we obtain
$$((\Gamma,
\Delta_V))^{\log}
=p_*
((\Gamma',
\Delta_{V'}))^{\log}
+
p_*(((p_1\times p_1)^*\Gamma_1
\cdot \Sigma,
\Delta_{V'}))^{\log}.$$
Let $i'\colon
V'_1\to V'$ 
denote the immersion.
If $U_1=U_{1,K}$,
then Proposition \ref{prU1}
gives us
\begin{equation}
(((p_1\times p_1)^*\Gamma_1\cdot
\Sigma,\Delta_{V'}))^{\log}
=
-
i'_*
(((p_1\times p_1)^*\Gamma_1\cdot
\Sigma\cdot
c_1({\rm pr}_2^*
N_{U'_1/U'}),
\Delta_{V_1}))^{\log}.
\label{eqpbdlx}
\end{equation}
If $U_1=U_{1,F}$,
we replace the double parentheses
in the right hand side
by a single parantheses.
By substituting
(\ref{eqSig}) in
(\ref{eqpbdlx})
and applying
Lemma \ref{lmpbdl}.1
to the projective space bundle
$p_1\colon V'_1
\to V_1$,
we obtain
\begin{equation}
p_*(((p_1\times p_1)^*\Gamma_1\cdot
\Sigma,\Delta_{V'}))^{\log}
=
i_*(-1)^{c-1}(c_{c,0}-c_{c,1})
((\Gamma_1,
\Delta_{V_1}))^{\log}
\end{equation}
in $F_0G
(\Sigma_{V/U}V)_{\mathbb Q}$.
Thus, the assertion follows
from $c_{c,0}-c_{c,1}=
(-1)^{c-1}(c-1)$
(\ref{eqcnm}).
\qed

We state and prove
an {\it excision} formula.

\begin{thm}\label{thmexc}
Let $f\colon V\to U$
be a finite \'etale morphism
of regular flat separated
schemes of finite
type over $S$
such that
$V_K\to U_K$
is tamely ramified with respect to
${\rm Spec}\ K$.
Let $U_1\subset U$ be 
a regular closed subscheme
and we consider the cartesian diagram
$$\begin{CD}
V_1@>j>> V @<i<< V_0&=V\setminus V_1\\
@VVV @VVV @VVV&\\
U_1@>j>> U @<i<< U_0&=U\setminus U_1
\end{CD}$$
For an open and closed subscheme $\Gamma$
of $V\times_UV$,
we put
$\Gamma_0=
\Gamma\times_UU_0$ and
$\Gamma_1=
\Gamma\times_UU_1$.

\setcounter{equation}0
Then, we have
\begin{equation}
((\Gamma,\Delta_V))^{\log}=
j_!((\Gamma_0,\Delta_{V_0}))^{\log}+
\begin{cases}
i_*((\Gamma_1,\Delta_{V_1}))^{\log}
&
\text{ if $U_1$ is flat over $S$ }\\
i_*(\Gamma_1,\Delta_{V_1})^{\log}
&
\text{ if $U_1=U_{1,F}$ }
\end{cases}
\label{eqexc}
\end{equation}
in $F_0G(\partial_{V/U} V)_{\mathbb Q}$.
\end{thm}

In the right hand side,
$j_!((\Gamma_0,\Delta_{V_0}))^{\log}$
and
$i_*((\Gamma_1,\Delta_{V_1}))^{\log}$
are elements
of $F_0G(\partial_{V/U} V)_{\mathbb Q}$
defined as the image of the map (\ref{eqg!}).

{\it Proof.}
By a standard devissage,
we may assume
either
$U_1=U_{1,K}$ or
$U_1=U_{1,F}$.
By Propositions 
\ref{prpbdl} and \ref{prblup},
it suffices to prove
the case where $U_1$
is a divisor of $U$.
We put $V_1=V\times_UU_1$.
Let $(V\times_SV)^\sim$
denote the log product
$(V\times_SV)^\sim_{V_1}$.
We consider the pull-back
$q^*\Gamma$
by the projection
$q\colon
(V\times_SV)^\sim
\to
V\times_SV$.
In the notation
of (\ref{eqmapE}),
we have
\begin{equation}
((\Gamma,
\Delta_V))^{\log}
=
((q^*\Gamma,
\Delta_V))^{\log},
\label{eqUsim}
\end{equation}
by the commutative diagram
(\ref{eqmapE0}).
Let $\widetilde \Gamma
\subset
(V\times_SV)^\sim$
denote the proper transform
of $\Gamma$
and we put
$\Gamma_1=
\Gamma\times_VV_1$.
We also have
\begin{equation}
((\Gamma_0,
\Delta_{V_0}))^{\log}
=
((\widetilde \Gamma,
\Delta_V))^{\log}
\end{equation}
by the commutative part
of the diagram
(\ref{eqUUsim}).

Let $q_1\colon
E\to 
V_1\times_SV_1$
be the base change of $q$
and $q_1^*\Gamma_1$ be the
pull-back.
Since $q\colon(U\times_SU)^\sim
\to U\times_SU$ is locally of complete intersection by Corollary \ref{cortdf}
and since $\Gamma$ is flat over $U$, 
we have an equality
\begin{equation}
q^*\Gamma=
\widetilde \Gamma
+
\tilde i_*(q_1^*\Gamma_1)
\label{eqdq}
\end{equation}
in $G((V\times_UV)^\sim
\setminus
(U\times_UU)^\sim)$
by Corollary \ref{cortlp}.
By Proposition \ref{prU2},
we have
$((q_1^*\Gamma_1,\Delta_V))^{\log}
=((\Gamma_1,\Delta_{V_1}))^{\log}$
if $U_1=U_{1,K}$
and
$((q_1^*\Gamma_1,\Delta_V))^{\log}
=(\Gamma_1,\Delta_{V_1})^{\log}$
if $U_1=U_{1,F}$.
Thus the assertion is proved.
\qed

Similarly and more easily,
we have the following analogue
of Theorem \ref{thmexc}.

\begin{pr}\label{prexc}
\setcounter{equation}0
Let $f\colon V\to U$
be a finite \'etale morphism
of smooth separated
schemes of finite
type over $F$.
Let $U_1\subset U$ be a smooth closed subscheme
and $U_0=U\setminus U_1$
be the complement.
For an open and closed subscheme $\Gamma$
of $V\times_UV$,
we put
$\Gamma_0=
\Gamma\times_UU_0$ and
$\Gamma_1=
\Gamma\times_UU_1$.

Then, we have
\begin{equation}
(\Gamma,\Delta_V)^{\log}=
(\Gamma_0,\Delta_{V_0})^{\log}+
(\Gamma_1,\Delta_{V_1})^{\log}
\label{eqexcF}
\end{equation}
in $CH_0(\partial_{V/U} 
V)_{\mathbb Q}$.
\end{pr}

We generalize the
definition
of the map
(\ref{eqmap})
for not necessarily regular $U$.
For a noetherian scheme
$X$,
let $\Gamma(X,{\mathbb Z})$
be the ${\mathbb Z}$-module
of ${\mathbb Z}$-valued
locally constant
functions on $X$.

\begin{cor}\label{corexc}
\setcounter{equation}0
For every finite \'etale 
morphism $f\colon V\to U$
of separated
schemes of finite
type over $S$
such that 
$V_K\to U_K$
is tamely ramified
with respect to $K$,
there exists a 
unique way to attach
a morphism
$$
((\ ,\Delta_V))^{\log}\colon
\Gamma(V\times_UV
\setminus \Delta_V,{\mathbb Z})
\to 
F_0G(\partial_{V/U} V)_{\mathbb Q}$$
satisfying the following properties:

{\rm (1)} If $U$ is 
regular and flat
of dimension $n$
over $S$,
it is the composition
$$\begin{CD}
\Gamma(V\times_UV
\setminus \Delta_V,
{\mathbb Z})
@>>>
Gr^F_nG(V\times_UV
\setminus \Delta_V)
@>{((\ ,\Delta_V))^{\log}}>>
F_0G(\partial_{V/U} V)_{\mathbb Q}
\end{CD}$$
where the first arrow
is the natural isomorphism.

If $U$ is smooth
of dimension $n$
over $F$,
it is the composition
$$\begin{CD}
\Gamma(V\times_UV
\setminus \Delta_V,
{\mathbb Z})
@>>>
CH_n(V\times_UV
\setminus \Delta_V)
@>{(\ ,\Delta_V)^{\log}}>>
F_0G(\partial_{V/U} V)_{\mathbb Q}
\end{CD}$$
where the first arrow
is the natural isomorphism.

{\rm (2)} Assume $U=\coprod_iU_i$
is a finite decomposition 
by regular subschemes.
Let $j_i\colon U_i\to U$
denote the immersion
and put $V_i=V\times_UU_i$
for each $i$.
Then, the diagram
\begin{equation}
\begin{CD}
\Gamma(V\times_UV
\setminus \Delta_V,
{\mathbb Z})
@>{((\ ,\Delta_V))^{\log}}>>
F_0G(\partial_{V/U} V)_{\mathbb Q}\\
@V{(j_i^*)_i}VV@AA{\sum j_{i!}}A\\
\bigoplus_i 
\Gamma(V_i\times_{U_i}V_i
\setminus \Delta_{V_i},
{\mathbb Z})
@>{\bigoplus_i((\ ,\Delta_{V_i}))^{\log}}>>
\bigoplus_i F_0G(\partial_{V_i/U_i} V_i)_{\mathbb Q}
\end{CD}
\label{eqswclgen}
\end{equation}
is commutative.
\end{cor}

{\it Proof.}
The uniqueness is a consequence
of the existence of a finite partition
by regular subschemes.
To show the existence,
it suffices to compare the maps defined by 
taking
partitions by regular subschemes.
By taking a common refinement,
it is reduced to verify the following.
Let $U$ be a regular scheme
and $U=\coprod_iU_i$
be a finite stratification
by regular subschemes.
Then, the maps defined in (1)
make the diagram 
(\ref{eqswclgen}) commutative.
It follows from 
Theorem \ref{thmexc}
and Proposition \ref{prexc}
by the induction on the maximum of 
the codimensions of $U_i$ in $U$.
\qed

By the same argument,
we have the following variant.

\begin{cor}\label{corexc0}
\setcounter{equation}0
Assume $K$ is of
characteristic $0$.
For every finite \'etale 
morphism $f\colon V\to U$
of separated
schemes of finite
type over $K$,
there exists a 
unique way to attach
a morphism
$$
((\ ,\Delta_V))^{\log}\colon
\Gamma(V\times_UV,{\mathbb Z})
\to 
F_0G(\partial_F V)_{\mathbb Q}$$
satisfying the following properties:

{\rm (1)} If $U$ is 
regular and flat
of dimension $n$
over $S$,
it is the composition
$$\begin{CD}
\Gamma(V\times_UV,
{\mathbb Z})
@>>>
Gr^F_nG(V\times_UV)
@>{((\ ,\Delta_V))^{\log}}>>
F_0G(\partial_F V)_{\mathbb Q}
\end{CD}$$
where the first arrow
is the natural isomorphism.

{\rm (2)} Assume $U=\coprod_iU_i$
is a finite decomposition 
by smooth subschemes.
Let $j_i\colon U_i\to U$
denote the immersion
and put $V_i=V\times_UU_i$
for each $i$.
Then, the diagram
\begin{equation}
\begin{CD}
\Gamma(V\times_UV,
{\mathbb Z})
@>{((\ ,\Delta_V))^{\log}}>>
F_0G(\partial_F V)_{\mathbb Q}\\
@V{(j_i^*)_i}VV@AA{\sum j_{i!}}A\\
\bigoplus_i 
\Gamma(V_i\times_{U_i}V_i,
{\mathbb Z})
@>{\bigoplus_i((\ ,\Delta_{V_i}))^{\log}}>>
\bigoplus_i F_0G(\partial_F V_i)_{\mathbb Q}
\end{CD}
\label{eqswclgen0}
\end{equation}
is commutative.
\end{cor}

\subsection{A semi-stable case}\label{ssss}

In this subsection,
we establish a crucial step
in the proof of 
the conductor formula.
Namely, in Proposition \ref{prcysd},
we compare the log localized
intersection products
$((T,\Delta_V))^{\log}$
and
$((\Gamma,\Delta_{V'}))^{\log}$
for a morphism $f\colon
V'\to V$
using the Lefschetz
trace formula Theorem \ref{thmlLTF},
assuming among other things that
$f$ is extended
to a weakly semi-stable
morphism 
of compactifications.

We consider 
a commutative diagram
\setcounter{thm}1
\setcounter{equation}0
\begin{equation}
\begin{CD}
U'@<<< V'\\
@VVV@VVfV\\
U@<<<V
\end{CD}
\label{eqcysd}
\end{equation}
of separated schemes 
of finite type over $S$
where
the horizontal
arrows are finite \'etale
and the vertical arrows
are smooth.
We further consider
a commutative diagram
\begin{equation}
\begin{CD}
\Gamma @>>> V'\times_{U'}V'\\
@VVV @VVV\\
T@>>> V\times_UV.
\end{CD}
\label{eqcysdT}
\end{equation}
where $T\subset
V\times_UV$
and 
$\Gamma \subset
V'\times_{U'}V'$
are open and
closed subschemes.

Let $V_T^{\prime (1)}$
and $V_T^{\prime (2)}$ 
denote the base change
$V'\times_VT$
with respect to
the first and the second projections
respectively.
We identify
the fiber product
$V_T^{\prime (1)}\times_T
V_T^{\prime (2)}$
with an open
and closed subscheme
$(V'\times_UV')
\times_
{V\times_UV}T$
of
$V'\times_UV'$.
Then, (\ref{eqcysdT})
implies that
$\Gamma$ is a closed
subscheme of
$V_T^{\prime (1)}\times_T
V_T^{\prime (2)}$.

We compare
the elements
$((T,\Delta_V))^{\log}$
and
$f_!((\Gamma,\Delta_{V'}))^{\log}$
of
$F_0G(\partial_FV)_{\mathbb Q}$,
assuming
that we have 
a commutative diagram
\begin{equation}
\begin{CD}
U'@>{\subset}>> X'
&\ \supset B'\\
@AAA@AAA\\
V'@>{\subset}>> Y'
&\ \supset D\\
@VfVV@VV{\bar f}V\\
V@>{\subset}>>Y
&\ \supset E_j
\end{CD}
\label{eqcysd2}
\end{equation}
of separated schemes
of finite type
over $S$ satisfying 
the following
conditions:

\begin{itemize}
\item[(\ref{ssss}.1.4)]
The schemes
$Y,Y'$ and $X'$
are proper over $S$
and 
$Y$ is the disjoint
union of irreducible
components.
The scheme $V$
is the complement
in $Y$ 
of a finite family
${\cal E}
=(E_j)_{j\in J}$
of Cartier divisors
and $U'$ is
the complement of
a Cartier divisor $B'$ of $X'$.
The morphism $\bar f\colon
Y'\to Y$ is proper 
weakly semi-stable
of relative dimension $d$
such that
$Y'_V=Y'\times_YV\to V$
is smooth.
The subscheme $D$ 
is a divisor
of $Y'$ over $Y$
with simple normal crossings relatively to $Y$
and
$V'=Y'_V\setminus D_V$.
\end{itemize}

Let $D_1,\ldots,D_m$
be the irreducible components
of $D$.
Let $Y^{\prime (1)}_T,
D_1^{(1)},\ldots,
D_m^{(1)}$
denote the base change
of $Y',D_1,\ldots,D_m$ over $Y$
by the composition
$T\to V\times_UV\to V\to Y$
of the first projection.
Similarly, we define
$Y^{\prime (2)}_T,
D_1^{(2)},\ldots,
D_m^{(2)}$ 
as the base change
using the second projections.
Let
$(Y^{\prime (1)}_T\times_T
Y^{\prime (2)}_T)^\sim$
denote the log product
with respect to the families
of Cartier divisors
$(D_1^{(1)},\ldots,
D_m^{(1)})$ and
$(D_1^{(2)},\ldots,
D_m^{(2)})$.

Let $K'$ be
a finite extension of $K$
and 
$\gamma'\colon
{\rm Spec}\ K'\to T$
be a morphism over $K$.
By the valuative criterion,
the compositions 
${\rm Spec}\ K'\to T
\to Y$
with the two projections
are extended to
$S'={\rm Spec}\
{\cal O}_{K'}\to Y$.
Let $Y^{\prime (1)}_{S'}$
and 
$Y^{\prime (2)}_{S'}$
denote the base change
with respect to
the two morphisms
respectively.

We regard $Y$
as a log scheme
with the log structure
defined by the finite family
of Cartier divisors
${\cal E}
=(E_j)_{j\in J}$.
We assume that
the restrictions of
$S'\to Y$
to the closed point
$s'\in S'$
define the same log point
$s'\to Y$.
Then, 
we have a canonical
isomorphism
$\iota_{s'}\colon
Y^{\prime (1)}_{s'}
\to
Y^{\prime (2)}_{s'}$.
Hence,
if the second projection
$\Gamma\to V^{\prime(2)}_T$
is proper,
the alternating sum
${\rm Tr}((\gamma^{\prime*}\Gamma)^*
\colon
H^*_c(V'_{\bar K'},
{\mathbb Q}_\ell))$
(\ref{eqTr*})
is defined 
for a prime number
$\ell$ invertible on $S$.

\begin{pr}\label{prcysd}
\setcounter{equation}0
Let the notations be
as in {\rm (\ref{ssss}.1.1)--(\ref{ssss}.1.4)}
and assume either of
{\rm (\ref{lmcysd}.1a)}
or
{\rm (\ref{lmcysd}.1b)}
is satisfied.
Let $T\subset
V\times_UV\setminus \Delta_V$
be an open and
closed subscheme and
let $\widetilde \Gamma$
be a closed subscheme of
$(Y^{\prime (1)}_T\times_T
Y^{\prime (2)}_T)^\sim$
flat over $T$
such that
$\Gamma=
\widetilde \Gamma\cap
(V_T^{\prime (1)}\times_T
V_T^{\prime (2)})$
is an open subscheme
of $V'\times_{U'}V'$.
We regard $Y$
as a log scheme
with the log structure
defined by 
${\cal E}
=(E_j)_{j\in J}$.
Let ${\cal E}'$
be the finite family of Cartier
divisors of $Y'$ consisting of
the pull-back of
${\cal E}$
and the irreducible components of
$D$.
We assume that the generic fibers
$\Sigma_{V/U}^{\cal E}Y_K$
and 
$\Sigma_{V'/U'}^{{\cal E}'}Y'_K$
are empty.

Then, there exist
a finite family
$(K_i)_{i\in I}$
of finite extensions
of $K$,
a family
$(\gamma_i\colon
{\rm Spec}\ K_i\to T
)_{i\in I}$
of morphisms over $S$
and rational numbers
$(r_i)_{i\in I}$
satisfying 
the following properties:

Let $s_i\in S_i
={\rm Spec}\ 
{\cal O}_{K_i}$ denote
the closed point
for $i\in I$.
Then,
for each $i\in I$,
the log points 
$\bar \gamma_i\colon
s_i\to Y$
defined by the unique maps
$S_i\to Y$
extending the composition
${\rm Spec}\ K_i\to 
T\to Y$ with
the first 
and the second projections
are equal to 
each other.
Further,
for a prime number 
$\ell$
invertible on $S$,
we have
${\rm Tr}((\gamma_i^*\Gamma)^*
\colon
H^*_c(V'_{\bar K_i},
{\mathbb Q}_\ell))
\in {\mathbb Q}$
and,
for the logarithmic product
{\rm (\ref{eqthmmap})},
\begin{align}
((T,\Delta_Y^{\log}))
&=\sum_ir_i[\bar \gamma_i(s_i)],
\label{eqccT}\\
\bar f_*
((\Gamma,\Delta_{Y'}^{\log}))
&=
\sum_ir_i
{\rm Tr}((\gamma_i^*\Gamma)^*
\colon
H^*_c(V'_{\bar K_i},
{\mathbb Q}_\ell))
\cdot [\bar \gamma_i(s_i)]
\label{eqcc}
\end{align}
in
$F_0G(\Sigma_{V/U}Y)_{\mathbb Q}$.
\end{pr}

{\it Proof.}
We take an object
$\bar g\colon Z\to Y$
of the category ${\cal A}_{V\to U}$
of alterations.
Since the conditions
{\rm (\ref{lmcysd}.1a)}
and 
{\rm (\ref{lmcysd}.1b)}
are stable by the base change,
we may assume that
we have a log blow-up
$Z'\to Z\times_YY'$
as in the conclusion
of Lemma {\rm \ref{lmcysd}}.
Hence,
the inverse image $W'=V'\times_{Y'}Z'$
is the complement $D'$ of
a divisor with simple normal crossings. 
By the assumption on
the upper square in
(\ref{eqcysd2}),
$\bar g'\colon Z'\to Y'$ defines
an object
of ${\cal A}_{V'\to U'}$.

Let $g\colon W\to V$
be the restriction of $\bar g\colon
Z\to Y$.
We put 
$(g\times g)^*(T)=
\sum_jm_j T_j
\in F_nG((W\times_SW)
\times_
{V\times_SV}T)$
and, for each $j$,
let $\bar T_j
\subset
(Z\times_SZ)^\sim$
be the schematic closure.
Then, we have
\begin{equation}
((T,\Delta_Y^{\log}))
=
\frac1{[W:V]}
\sum_jm_j\cdot
\bar g_*((\bar T_j,
\Delta^{\log}_Z))
_{(Z\times_SZ)^\sim}.
\label{eqcsT}
\end{equation}

We define the
log products
$(Y\times_SY)^\sim$
and
$(Y'\times_SY')^\sim$
with respect to ${\cal E}$
and ${\cal E}'$.
The log product
$(Y_T^{\prime (1)}\times_T
Y_T^{\prime (2)})
^{\sim}$ defined with respect
to the pull-backs
of the irreducible components
of $D$
is canonically identified with
$(Y'\times_SY')^\sim
\times_{(Y\times_SY)^\sim}T$.
We define
$(Z'\times_SZ')^\sim$
to be the log product
with respect to
the irreducible components
of a divisor $Z'\setminus (Z'\times_ZW)$
with simple normal crossings
and
the pull-backs of
the irreducible components of $D\subset Y'$.
The inverse image
of $W\times_SW$
by 
$(Z'\times_SZ')^\sim
\to(Z\times_SZ)^\sim$
is canonically identified with
$(Y'\times_SY')^\sim
\times_{(Y\times_SY)^\sim}
(W\times_SW)$.

For $T_j$ as above,
both
$(Y_T^{\prime (1)}\times_T
Y_T^{\prime (2)})
^{\sim}
\times_TT_j$
and 
$(Z'\times_SZ')^\sim
\times_{(Z\times_SZ)^\sim}
T_j$
are identified with
$(Y'\times_SY')^\sim
\times_{(Y\times_SY)^\sim}T_j$.
We regard
$\widetilde\Gamma_j=
\widetilde\Gamma\times_TT_j$
as a closed subscheme
of 
$(Y_T^{\prime (1)}\times_T
Y_T^{\prime (2)})
^{\sim}
\times_TT_j
=
(Z'\times_SZ')^\sim
\times_{(Z\times_SZ)^\sim}
T_j$.
By the assumption
that
$\widetilde \Gamma
\subset 
(Y_T^{\prime (1)}\times_T
Y_T^{\prime (2)})
^{\sim}$ is 
flat over $T$
and by the flattening
theorem
\cite{GR},
there exists a proper
modification
$\bar q_j\colon
\bar T'_j\to
\bar T_j$ for
each $\bar T_j$
satisfying the following
conditions:
The map
$\bar q_j\colon
\bar T'_j\to
\bar T_j$ induces
the identity on 
the dense open subscheme $T_j$
and 
the schematic closure 
$\bar \Gamma'_j$ of
$\widetilde\Gamma_j
\subset
(Y_T^{\prime (1)}\times_T
Y_T^{\prime (2)})
^{\sim}
\times_TT_j
=
(Z'\times_SZ')^\sim
\times_{(Z\times_SZ)^\sim}
T_j$
in the base change
$(Z'\times_SZ')^\sim
\times_{(Z\times_SZ)^\sim}
\bar T'_j$
is flat over $\bar T'_j$.
Then, similarly as (\ref{eqcsT}),
we have
\begin{equation}
((\Gamma,\Delta_{Y'}^{\log}))
=
\frac1
{[W:V]}\sum_jm_j\cdot
\bar g'_*\bar q_{j*}((\bar\Gamma'_j,
\Delta^{\log}_{Z'}))
_{(Z'\times_SZ')^\sim}.
\label{eqcsG}
\end{equation}

For each $j$, we put
\begin{equation}
((\bar T'_j,
\Delta^{\log}_Z))
_{(Z\times_SZ)^\sim}
=\sum n_i[s'_i]
\label{eqs'i}
\end{equation}
in
$F_0G(\bar T'_j\times_
{(Z\times_SZ)^\sim} 
\Delta^{\log}_Z)$.
Since $s'_i\in \bar T'_j$
is a closed point
and $T_j$ is dense in 
$\bar T'_j$,
there exist 
a discrete valuation field
$K_i$ and 
a map
$\gamma_i\colon
S_i={\rm Spec}\ 
{\cal O}_{K_i}
\to \bar T_j$ extending
${\rm Spec}\ K_i
\to T_j$
such that $s'_i$
is the image of the closed point
$S_i$ of $S_i$.
Since the image of $s_i$ in
$(Y\times_S Y)^\sim$ 
is in the log diagonal,
the log points 
$s_i\to Y$
defined by the two projections
are equal
to each other.
Thus, by (\ref{eqcsT}),
we obtain (\ref{eqccT}).

We prove the equality
(\ref{eqcc}).
We fix $j$ and
let
$p_j\colon
\bar\Gamma'_j
\to \bar T'_j$
denote the projection.
First we show
\begin{equation}
p_{j*}((\bar\Gamma'_j,
\Delta^{\log}_{Z'}))
_{(Z'\times_SZ')^\sim}
=
\sum_i n_i
\deg
(\bar \Gamma'_{j,s_i},
\Delta^{\log}_{Z'_{s_i}})
_{(Z'_{s_i}\times_{s_i}
Z'_{s_i})^\sim}
\cdot 
[s'_i]
\label{eqGT}
\end{equation}
in 
$F_0G(\bar T'_j\times_
{(Z\times_SZ)^\sim}
\Delta^{\log}_Z)$.
Since $Z'\to Z$
is log smooth,
the morphism
$(Z'\times_SZ')^\sim
\to (Z\times_SZ)^\sim$
is smooth
and
the log diagonal map
$Z'\to (Z'\times_ZZ')^\sim
=
(Z'\times_SZ')^\sim
\times_{ (Z\times_SZ)^\sim}Z$
is a regular immersion.
By applying the associativity
Lemma \ref{lmfm2}
to
$(Z'\times_SZ')^\sim
\gets
(Z'\times_ZZ')^\sim
\gets Z'$, 
we obtain
$$
((\bar\Gamma'_j,
\Delta^{\log}_{Z'}))
_{(Z'\times_SZ')^\sim}
=
\bigl(((\bar\Gamma'_j,
(Z'\times_ZZ')^\sim))
_{(Z'\times_SZ')^\sim},
\Delta^{\log}_{Z'}\bigr)
_{(Z'\times_ZZ')^\sim}.$$
Since 
$(Z'\times_SZ')^\sim$
is smooth over
$(Z\times_SZ)^\sim$,
it is tor-independent with
$Z$.
Hence by applying
the associativity
Lemma \ref{lmfm4}
to $(Z\times_SZ)^\sim
\gets(Z'\times_SZ')^\sim
\gets\bar\Gamma'_j$
as $X\gets X'\gets W$,
we obtain
$$((\bar\Gamma'_j,
(Z'\times_ZZ')^\sim))
_{(Z'\times_SZ')^\sim}
=
((\Delta_Z^{\log},\bar\Gamma'_j))
_{(Z\times_SZ)^\sim}.$$
Since $p_j\colon
\bar\Gamma'_j
\to \bar T'_j$
is flat,
by further applying the associativity
Lemma \ref{lmfm2}
to
$(Z\times_SZ)^\sim
\gets \bar T'_j\gets
\bar\Gamma'_j$,
we obtain
$$
((\Delta_Z^{\log},\bar\Gamma'_j))
_{(Z\times_SZ)^\sim}
=p_j^*((\bar T'_j,
\Delta^{\log}_Z))
_{(Z\times_SZ)^\sim}.
$$
Thus, 
we obtain
$$
((\bar\Gamma'_j,
\Delta^{\log}_{Z'}))
_{(Z'\times_SZ')^\sim}
=
(p_j^*((\bar T'_j,
\Delta^{\log}_Z))
_{(Z\times_SZ)^\sim},
\Delta^{\log}_{Z'})
_{(Z'\times_ZZ')^\sim}.$$
Substituting
(\ref{eqs'i}),
we see that
the right hand side is equal to
$$\sum_in_i
(p_j^*[s'_i],
\Delta^{\log}_{Z'})
_{(Z'\times_ZZ')^\sim}
=
\sum_in_i
(\bar \Gamma'_{j,s_i},
\Delta^{\log}_{Z'_{s_i}})
_{(Z'_{s_i}\times_{s_i}
Z'_{s_i})^\sim}
$$
in 
$F_0G(\bar \Gamma'_j\times_
{(Z\times_SZ)^\sim}
\Delta^{\log}_Z)$.
Thus the equality
(\ref{eqGT}) is proved.

For the right hand side
of (\ref{eqGT}),
we show
\begin{equation}
\deg
(\bar \Gamma'_{j,s_i},
\Delta^{\log}_{Z'_{s_i}})
_{(Z'_{s_i}\times_{s_i}
Z'_{s_i})^\sim}
=
{\rm Tr}((\gamma_i^*\Gamma)^*
\colon
H^*_c(V'_{\bar K_i},
{\mathbb Q}_\ell))
\label{eqTFcs}
\end{equation}
by applying Theorem \ref{thmlLTF}
to the base changes
$Z'\times_ZS_i\to S_i
={\rm Spec}\ {\cal O}_{K_i}$.
To apply it,
we verify that the assumptions are satisfied.
The log blow-up of the product
$(Y^{\prime (1)}_T\times_T
Y^{\prime (2)}_T)'$
with respect to the families
of Cartier divisors
$(D_1^{(1)},\ldots,
D_m^{(1)})$ and
$(D_1^{(2)},\ldots,
D_m^{(2)})$
contains the log product
$(Y^{\prime (1)}_T\times_T
Y^{\prime (2)}_T)^\sim$
as the complement of
the proper transforms
$(D^{(1)}_T\times_TY^{(2)}_T)'$
and
$(Y^{(1)}_T\times_TD^{(1)}_T)'$.
Let 
$\Gamma'$
be the closure of
$\widetilde \Gamma$
in the product
$(Y^{\prime (1)}_T\times_T
Y^{\prime (2)}_T)'$.
We show
\begin{equation}
\Gamma'\cap
(D^{(1)}_T\times_TY^{(2)}_T)'
=
\Gamma'\cap
(Y^{(1)}_T\times_TD^{(1)}_T)'.
\label{eqGY'}
\end{equation}

Let $(Y'_U\times_UY'_U)'$
denote the log product
with respect to the families
of Cartier divisors
$(D_{1,U},\ldots,
D_{m,U})$.
We have an open immersion
$(Y^{\prime (1)}_T\times_T
Y^{\prime (2)}_T)'
\to(Y'_U\times_UY'_U)'$
as the base change of
$T\to V\times_UV$.
Since $\Gamma
\subset V'\times_{U'}V'$,
the image of
$\Gamma'$
by
$(Y^{\prime (1)}_T\times_T
Y^{\prime (2)}_T)'
\subset
(Y'_U\times_UY'_U)'
\to
X'_U\times_UX'_U$
is in the diagonal
$X'_U$.
Since $U'\subset X'$
is assumed to be the complement
of a Cartier divisor $B'$,
we have 
$\Gamma'\cap
(D_U\times_UY'_U)'
=
\Gamma'\cap
(Y'_U\times_UD_U)'$
by \cite[Proposition 1.1.6.2]{KSA}.
Thus, we obtain (\ref{eqGY'}).
Consequently, the base change to $K_i$
satisfies the inclusion (\ref{eqG'}).

We construct
a commutative diagram
(\ref{eqPT}) of monoids
for the two base changes
$Z'\times_ZS_i\to S_i
={\rm Spec}\ {\cal O}_{K_i}$
satisfying the
condition (P) loc.\ cit.
Since the log structures
of $Z$ and $Z'$ are
defined by divisors
with simple normal crossings,
we have a commutative diagram
\begin{equation}
\begin{CD}
{\mathbb N}^m
@>>>
{\mathbb N}^{m'}\\
@VVV@VVV\\
\Gamma(Z,\bar {\mathcal M}_Z)
@>>>
\Gamma(Z',\bar {\mathcal M}_{Z'})
\end{CD}
\label{eqPZZ'}
\end{equation}
of morphisms of monoids,
locally lifted to charts.
Since the image of
the closed point by the 
composition
$S_i={\rm Spec}\ {\cal O}_{K_i}
\to \bar T_j\to (Z\times_SZ)^\sim$
lies in the log diagonal,
the compositions with
the two propositions 
define the same log points
$s_i\to Z$. Hence
there exists one morphism
${\mathbb N}^m
\to 
{\mathbb N}$
of monoids that makes
the diagram
\begin{equation}
\begin{CD}
{\mathbb N}
@<<<
{\mathbb N}^m\\
@VVV@VVV\\
\Gamma(S_i,\bar {\mathcal M}_{S_i})
@<<<
\Gamma(Z,\bar {\mathcal M}_Z).
\end{CD}
\label{eqPZSi}
\end{equation}
commutative for
the compositions with
the two projections.
We define a monoid $P$
to be the saturated sum
${\mathbb N}+^{\rm sat}_
{{\mathbb N}^m}
{\mathbb N}^{m'}$.
Then, the commutative diagrams
(\ref{eqPZZ'})
and (\ref{eqPZSi})
induces a morphism
$P\to 
\Gamma(Z'\times_ZS_i,\bar {\mathcal M}_{Z'\times_ZS_i})$,
locally lifted to charts.
It defines
a commutative diagram
(\ref{eqPT}) of monoids
satisfying the
condition (P).
Thus we may apply Theorem \ref{thmlLTF} and
we obtain (\ref{eqTFcs}).

Therefore the equality
(\ref{eqcc}) follows
from (\ref{eqcsG}),
(\ref{eqGT}) and (\ref{eqTFcs}).
\qed

By the same argument,
in the case where $K$ is
of characteristic $0$,
we have the following.

\begin{pr}\label{prcysd0}
\setcounter{equation}0
Assume that $K$
is of characteristic $0$
and let the notations be
as in {\rm (\ref{ssss}.1.1)--(\ref{ssss}.1.4)}.
We assume that
either {\rm (\ref{lmcysd}.1a)}
or
{\rm (\ref{lmcysd}.1b)}
is satisfied.
Let $T\subset
V\times_UV$
be an open and
closed subscheme and
let $\widetilde \Gamma$
be a closed subscheme
of the log product
$(Y^{\prime (1)}_T\times_T
Y^{\prime (2)}_T)^\sim$
flat over $T$
such that
$\Gamma=
\widetilde \Gamma\cap
(V_T^{\prime (1)}\times_T
V_T^{\prime (2)})$
is an open subscheme
of $V'\times_{U'}V'$.
We regard $Y$
as a log scheme
with the log structure
defined by a finite family
of Cartier divisors
${\cal E}
=(E_j)_{j\in J}$
satisfying $V=
Y\setminus
\bigcup_{j\in J}E_j$.

Then, there exist
a finite family
$(K_i)_{i\in I}$
of finite extensions
of $K$,
a family
$(\gamma_i\colon
{\rm Spec}\ K_i\to T
)_{i\in I}$
of morphisms over $S$
and rational numbers
$(r_i)_{i\in I}$
satisfying the following properties:

Let $s_i\in S_i
={\rm Spec}\ 
{\cal O}_{K_i}$ denote
the closed point
for $i\in I$.
Then,
for each $i\in I$,
the log points 
$\bar \gamma_i\colon
s_i\to Y$
defined by the unique maps
$S_i\to Y$
extending the composition
${\rm Spec}\ K_i\to 
T\to Y$ with
the first 
and the second projections
are equal to 
each other.
Further
for a prime number 
$\ell$
invertible on $S$,
we have
\begin{align}
((T,\Delta_Y^{\log}))
&=\sum_ir_i[\bar \gamma_i(s_i)],
\label{eqccT0}\\
\bar f_*
((\Gamma,\Delta_{Y'}^{\log}))
&=
\sum_ir_i
{\rm Tr}((\gamma_i^*\Gamma)^*
\colon
H^*_c(V'_{\bar K_i},
{\mathbb Q}_\ell))
\cdot [\bar \gamma_i(s_i)]
\label{eqcc0}
\end{align}
in
$F_0G(Y\times_SF)_{\mathbb Q}$.
\end{pr}

\newpage 
\section{The Swan class
and a conductor formula}\label{sSw}

We keep the notation
that $K$ 
is a complete discrete 
valuation ring
and
$S={\rm Spec}\ {\cal O}_K$
as in the previous sections.
We fix a prime number
$\ell$
different from
the characteristic $p$
of the residue field $F$
of $K$.
First, we define
the Swan character
classes for 
Galois coverings in
Section \ref{ssSwch}.
We define the Swan class
of a locally constant
$\bar {\mathbb F}_\ell$-sheaf
in Section \ref{ssdSw}.
We extend the definition
of the Swan class
to a constructible
sheaf 
in Section \ref{ssSct}
using the excision
formula Proposition \ref{prExc}.2.,
assuming $K$
is of characteristic 0.
We prove
a conductor formula 
for some relative curves
in Section \ref{sscfs}
and derive the general case
in Section \ref{sspc}.
In an equal characteristic case,
more elementary proof is found in
\cite[Corollaries 5.12, 5.13]{Tsu}.

In this paper,
we state and prove
results for
$\overline{\mathbb F}_\ell$-sheaves.
The corresponding
results for
$\overline{\mathbb Q}_\ell$-sheaves
are obtained simply by
taking reduction modulo $\ell$.

\subsection{Swan character
classes}\label{ssSwch}

We define the Swan character
class for a Galois covering.

\begin{df}\label{dfSwch}
Let $U$ be a separated
regular flat
scheme of finite
type over $S$
and 
$f\colon V\to U$
be a finite \'etale $G$-torsor
for a finite group $G$
such that
the generic fiber
$V_K\to U_K$
is tamely ramified
with respect to $K$.
Then, for
an element 
$\sigma\in G$,
we define
{\rm the Swan character
class} $s_{V/U}(\sigma)
\in F_0G(\partial_{V/U}V)
_{\mathbb Q}$
by 
\begin{equation}
s_{V/U}(\sigma)
=\begin{cases}
D_{V/U}^{\log}
\quad
&\text{ for }\sigma=1\\
-((\Gamma_\sigma,\Delta_V))^{\log}
&\text{ for }\sigma\neq1.
\label{eqSwch}
\end{cases}
\end{equation}
\end{df}

By Corollaries
\ref{correxX} and \ref{corssw},
for a finite Galois
extension $L$
of $K$ of Galois
group $G$
and $U={\rm Spec}\ K,
V={\rm Spec}\ L$,
we have
\begin{equation}
s_{V/U}(\sigma)
=\begin{cases}
{\rm length}_{{\cal O}_L}
\Omega^1_{{\cal O}_L/
{\cal O}_K}(\log/\log)
\quad
&\text{ for }\sigma=1\\
-
{\rm length}_{{\cal O}_L}
{\cal O}_L/J_\sigma
&\text{ for }\sigma\neq1.
\label{eqswch}
\end{cases}
\end{equation}
in $F_0G(\partial_{V/U}V)
={\mathbb Z}$,
where $J_\sigma$
is the ideal of
${\cal O}_L$
generated by
$\sigma(a)-a$
for $a\in {\cal O}_L$
and
$\sigma(b)/b-1$
for $b\in L^\times$.

\begin{lm}\label{lmSwch}
Let the notation be
as in Definition
{\rm \ref{dfSwch}}.
Then, the following hold:

{\rm 1.}
We have
$$\sum_{\sigma\in G}
s_{V/U}(\sigma)=0.$$

{\rm 2.}
Let $H$ be a subgroup
of $G$ and
$g\colon
V\to U'$
be the corresponding
$H$-torsor.
Then, for $\sigma\in H$,
we have
$$s_{V/U}(\sigma)=
\begin{cases}
s_{V/U'}(1)+
g^*D_{U'/U}^{\log}
\quad&\text{ if }
\sigma=1\\
s_{V/U'}(\sigma)
\quad&\text{ if }
\sigma\neq1.
\end{cases}$$

{\rm 3.}
Let $N$ be a normal subgroup
of $G$ and
$G'=G/N$ be the quotient.
Let $g\colon V\to V'$
be the corresponding
$N$-torsor.
Then, for $\sigma'\in G'$,
we have
$$g^*s_{V'/U}(\sigma')=
\sum_{\sigma\in G,
\bar \sigma=\sigma'}
s_{V/U}(\sigma).$$
\end{lm}

{\it Proof.}
1. Clear from
Lemma \ref{lmDs}.2.

2. For $\sigma=1$,
it follows from
Lemma \ref{lmDs}.1.
For $\sigma\neq 1$,
it is clear from
the definition.

3. Clear from
Lemma \ref{lmDs}.3.
\qed

\begin{cor}\label{corswch}
If the order
of $\sigma\in G$
is not a power
of $p$,
we have
$$s_{V/U}(\sigma)=0.$$
\end{cor}

{\it Proof.}
By Lemma \ref{lmSwch}.2,
we may assume
that $G$ is
the cyclic group generated by
$\sigma$.
Assume the order
of $\sigma$
is not a power of $p$.
Let
$N\subset G$
be the $p$-Sylow
subgroup
and $U'\to U$
be the corresponding
$G'=G/N$-torsor.
Then, since
the order of $G'$
is prime to $p$,
the finite \'etale morphism
$U'\to U$
is tamely ramified
with respect to $S$
by Corollary \ref{corptp}.
Hence, it follows from
Corollary \ref{corfun}
applied to $V'=V$.
\qed

For an element
$\sigma\in G$
of order a power
of $p$
and an integer $i$
prime to $p$,
Conjecture \ref{cnsi}
predicts
$$s_{V/U}(\sigma)=
s_{V/U}(\sigma^i).$$

\begin{cor}\label{corIy}
Let the notation
be as in Definition
{\rm \ref{dfSwch}}.
Let $X$ be a
normal proper scheme over $S$
containing $U$
as a dense open subscheme
and let $Y$
be the normalization
of $X$ in $V$.
Let $\sigma\in G$
be an element
not contained
in any conjugate
of a $p$-Sylow group
of the inertia group
$I_{\bar y}\subset G$
for any geometric 
point $\bar y$ of $Y$.
Then, we have
$$s_{V/U}(\sigma)=0.$$
\end{cor}

{\it Proof.}
By  Corollary
\ref{corswch},
it suffices to
consider the case
where the order of $\sigma$
is a power of $p$.
As in the proof of Corollary
\ref{corswch},
we may assume that
$G$ is the cyclic group
generated by $\sigma$.
Let $N\subset G$ be the
unique maximal proper subgroup
generated by $\sigma^p$.
Let $V'\to U$
be the corresponding
$G'=G/N$-torsor
and let
$Y'$ be normalization
of $X$ in $V'$.
Then, by the assumption,
the inertia group
at every geometric point
is a subgroup of $N$ and
hence $Y'\to X$ is \'etale.
Hence
the assertion follows from
Corollary \ref{corfun}.
\qed

\subsection{Swan class
of a locally constant sheaf}
\label{ssdSw}

We briefly recall the Brauer
trace of an $\ell$-regular
element \cite{Il}.
Let $G$ be a finite group
and $\ell$ be a prime number.
An element $\sigma\in G$
is called an $\ell$-{\it regular
element}
if the order of $\sigma$
is prime to $\ell$.
Let $G^{(\ell)}$
denote the subset
of $G$ consisting
of $\ell$-regular
elements.
For an element $\sigma$
of a pro-finite group,
we say that $\sigma$
is $\ell$-regular
if it is a projective
limit of $\ell$-regular
elements.

Let $M$ be an
$\overline
{\mathbb F}_\ell$-vector space of
finite dimension $n$
and $\sigma$
be an automorphism of
$M$ of order prime to $\ell$.
Then 
the Brauer trace ${\rm Tr}^{\rm Br}
(\sigma:M)
\in {\mathbb Z}_\ell^{\rm ur}$
is defined as follows.
Let $\alpha_1,\ldots,
\alpha_n
\in
\bar {\mathbb F}_\ell^{\times}$
be the eigenvalues
of $\sigma$ on $M$
counted with multiplicities
and 
$\tilde \alpha_1,\ldots,
\tilde \alpha_n\in
{\mathbb Z}_\ell^{{\rm ur}\times}$
be the liftings
of finite orders
prime to $\ell$.
Then 
the Brauer trace
is defined by
${\rm Tr}^{\rm Br}
(\sigma:M)=\sum_{i=1}^n
\tilde \alpha_i$.

Let $f\colon V\to U$ 
be a finite \'etale 
$G$-torsor for
a finite group $G$
such that
the generic fiber
$V_K\to U_K$
is tamely ramified
with respect to $K$.
By Corollary \ref{corswch},
we have
$s_{V/U}(\sigma)=0$,
if the order 
of $\sigma$ is not
a power of $p$.
In the following,
let $G_{(p)}$
denote the subset
of $G$ consisting
of elements 
of order a power of $p$.
For $\ell\neq p$,
we have
$G_{(p)}
\subset G^{(\ell)}$.
We put
$F_0G(\partial_{V/U}V)
_{{\mathbb Q}(\zeta_{p^\infty})}
=
F_0G(\partial_{V/U}V)
_{\mathbb Q}
\otimes_{\mathbb Q}
{\mathbb Q}(\zeta_{p^\infty}),
F_0G(\partial_{V/U}V)
_{{\mathbb Z}[\zeta_{p^\infty}]}
=
F_0G(\partial_{V/U}V)
\otimes_{\mathbb Z}
{\mathbb Z}[\zeta_{p^\infty}]$
etc.

\begin{df}\label{dfSwn}
\setcounter{equation}0
Let $U$ be a regular flat
separated scheme 
of finite type over $S=
{\rm Spec}\ {\cal O}_K$
and let ${\cal F}$ be
a locally constant constructible
$\overline{\mathbb F}_\ell$-sheaf
on $U$.
Let $f\colon V\to U$ 
be a finite \'etale 
$G$-torsor for
a finite group $G$
such that
$f^*{\cal F}$
is a constant sheaf on $V$.
We assume that
the generic fiber
$V_K\to U_K$
is tamely ramified
with respect to $K$.
Let $M$ be the 
$\overline{\mathbb F}_\ell$-representation
of $G$ corresponding to ${\cal F}$.

Then, we define
{\rm the Swan class}
${\rm Sw}_{V/U}{\cal F}\in 
F_0G(\partial_{V/U}V)_{\mathbb Q(\zeta_{p^\infty})}$
by
\begin{equation}
{\rm Sw}_{V/U}{\cal F}
=
\sum_{\sigma \in G_{(p)}}
{\rm Tr}^{\rm Br}(\sigma:M)
\cdot
s_{V/U}(\sigma).
\label{eqswan}
\end{equation}
\end{df}

By Lemma \ref{lmSwch}.1
and
${\rm Tr}^{\rm Br}(1:M)
=\dim M$,
the defining
equality (\ref{eqswan})
is equivalent to
the following:
\begin{equation}
{\rm Sw}_{V/U}{\cal F}
=
\sum_{\sigma \in G_{(p)},
\neq 1}
(\dim M-{\rm Tr}^{\rm Br}(\sigma:M))
\cdot
((\Gamma_{\sigma},\Delta_V))^{\log}.
\label{eqSwn}
\end{equation}
Recall that in \cite{KSA},
the Swan class
is defined similarly
for a locally constant sheaf
on a smooth scheme
over a perfect field
and is called the naive Swan class
and is denoted by ${\rm Sw}'$.
Modifying the notation,
we remove ``$'$''.
If we assume
Conjecture \ref{cnsi}
asserting that
$s_{V/U}(\sigma)
=s_{V/U}(\sigma^i)$
for an integer $i$ prime to $p$,
the Swan class
${\rm Sw}_{V/U}{\cal F}$ 
is in fact defined in the subspace
$F_0G(\partial_{V/U}V)
_{\mathbb Q}
\subset
F_0G(\partial_{V/U}V)
_{\mathbb Q(\zeta_{p^\infty})}$.

If we assume a strong form
of resolution of
singularity,
the Swan character
class is defined integrally
(\ref{eqmapZ}) and hence
the Swan class
${\rm Sw}_{V/U}{\cal F}$
is defined as an element
of
$F_0G(\partial_{V/U}V)
_{\mathbb Z[\zeta_{p^\infty}]}$.
Further, if we assume
Conjecture \ref{cnsi},
the Swan class
${\rm Sw}_{V/U}{\cal F}$ 
is in fact defined 
integrally in the subgroup
$F_0G(\partial_{V/U}V)
\subset
F_0G(\partial_{V/U}V)
_{\mathbb Z[\zeta_{p^\infty}]}$.
When, we emphasize
that it is defined integrally,
we write
${\rm Sw}_{V/U}^{\mathbb Z}{\cal F}$ 
and call it the integral
Swan class.
Note that
Conjecture \ref{cnsi}
itself is a consequence
of a strong form
of equivariant resolution
of singularities.

Similarly as \cite[Lemma 4.3.10]{KSA},
we have the following
analogue of \cite[Th\'eor\`eme 2.1]{Il}.

\begin{pr}
{\rm (cf.\ \cite[Corollaire 3.4]{Vidal})}
Let $U$ be a regular flat
separated scheme 
of finite type over $S=
{\rm Spec}\ {\cal O}_K$
and ${\cal F}_1$ 
and ${\cal F}_2$ 
be locally constant
constructible sheaves of
$\overline{\mathbb F}_\ell$-modules
on $U$.
Let $f\colon V\to U$
be a finite \'etale
$G$-torsor for
a finite group $G$
such that
$f^*{\cal F}_1$ 
and $f^*{\cal F}_2$ 
are constant on $V$
and that
the generic fiber
$V_K\to U_K$
is tamely ramified
with respect to $K$.

Let $X$ be a proper
normal scheme
over $S$
containing
$U$ as a dense open subscheme.
Assume that,
for every geometric
point $\bar x$ of
$X$,
the restriction
to a $p$-Sylow
subgroup of the
inertia group $I_{\bar x}$
of the representations
$M_1$ and $M_2$
of $G$
corresponding to
${\cal F}_1$ and
${\cal F}_2$ are
isomorphic.
Then, we have
${\rm Sw}_{V/U}{\cal F}_1=
{\rm Sw}_{V/U}{\cal F}_2$.
\end{pr}

{\it Proof.}
It follows from 
(\ref{eqswan}) and
Corollary
\ref{corIy}.
\qed

\begin{lm}\label{lmSwU}
Let $U,V, G$ and ${\cal F}$
be as in Definition {\rm \ref{dfSwn}}.
Let $f'\colon
V'\to U$ be a
finite \'etale $G'$-torsor
for a finite group $G'$
such that $f^{\prime *}{\cal F}$
is a constant sheaf on $V'$.
We assume that
$V'_K\to U_K$
is tamely ramified with respect to
${\rm Spec}\ K$.
Let $g\colon V'\to V$
be a morphism over $U$
compatible with a group
homomorphism
$G'\to G$.
Then, we have
$${\rm Sw}_{V'/U}{\cal F}
=g^*{\rm Sw}_{V/U}{\cal F}
$$
in $F_0G(\partial_{V/U}V')
_{\mathbb Q(\zeta_{p^\infty})}$.
\end{lm}

{\it Proof.}
It follows from
the definition
and Lemma
\ref{lmSwch}.3.
\qed

By Lemma \ref{lmSwU},
the Swan class
${\rm Sw}_{V/U}{\cal F}$
is $G$-invariant
and hence
$\frac 1{|G|}f_*{\rm Sw}_{V/U}{\cal F}
\in
F_0G(\partial_{V/U}U)
_{\mathbb Q(\zeta_{p^\infty})}$
is independent
of the choice
of a Galois covering $V$
trivializing ${\cal F}$.
Thus the following
definition makes sense.

\begin{df}\label{dfSwU}
Let $U$ be a regular flat
separated scheme 
of finite type over $S=
{\rm Spec}\ {\cal O}_K$
and ${\cal F}$ be a locally constant
constructible sheaf of
$\overline{\mathbb F}_\ell$-modules
on $U$.

{\rm 1.}
We say that
${\cal F}$ is {\rm tamely ramified
on the generic fiber}
if there exists a finite \'etale
surjective morphism
$f\colon V\to U$
such that
$f^*{\cal F}$
is constant on $V$
and that
the generic fiber
$V_K\to U_K$
is tamely ramified
with respect to $K$.

{\rm 2.}
Assume that
${\cal F}$ is tamely ramified
on the generic fiber
and that $U$ is connected.
Let $f\colon V\to U$
be a finite \'etale
$G$-torsor for
a finite group $G$
such that
$f^*{\cal F}$
is constant on $V$
and that
${\cal F}$ corresponds
to a faithful $\bar {\mathbb F}_\ell$-representation of $G$.

Then, we put
$$F_0G(\partial_{\cal F}U)_{{\mathbb Q}(\zeta_{p^\infty})}
=
F_0G(\partial_{V/U}U)_
{{\mathbb Q}(\zeta_{p^\infty})}$$
and define {\rm the Swan class}
$${\rm Sw}_U{\cal F}
=\frac 1{|G|}f_*
{\rm Sw}_{V/U}{\cal F}$$
to be the image
by the isomorphism
from the $G$-fixed part
$\frac 1{|G|}f_*
\colon
F_0G(\partial_{V/U}V)^G_{\mathbb Q(\zeta_{p^\infty})}
\to 
F_0G(\partial_{V/U}U)_{\mathbb Q(\zeta_{p^\infty})}=
F_0G(\partial_{\cal F}U)_{{\mathbb Q}(\zeta_{p^\infty})}$.
If $U$ is not connected,
we define the Swan class
componentwise.
\end{df}

For the $G$-torsor $V$ in Definition
\ref{dfSwU}.2,
we have 
${\rm Sw}_{V/U}{\cal F}
=f^*{\rm Sw}_U{\cal F}$
by Lemma \ref{lmGal}.2.
If $K$ is of characteristic $0$,
every locally constant sheaf on $U$
is tamely ramified on the generic fiber.
In the case where $U={\rm Spec}\ K$
and ${\cal F}$ is wildly ramified,
the Swan class
${\rm Sw}_U{\cal F}
\in 
F_0G(\partial_{\cal F}U)_{\mathbb Q(\zeta_{p^\infty})}
=\mathbb Q(\zeta_{p^\infty})$
is nothing but
the Swan conductor
$${\rm Sw}_K{\cal F}
=\frac1{|G|}
\sum_{\sigma\in G_{(p)}}
{\rm Tr}^{\rm Br}(\sigma
\colon M)\cdot
f_*s_{L/K}(\sigma)$$
by (\ref{eqswch}),
known to be an integer $\ge 1$.

The Swan classes 
satisfy the following 
additivity and 
the excision formula.

\begin{pr}\label{prExc}
Let $U$ be a regular flat
separated scheme 
of finite type over $S=
{\rm Spec}\ {\cal O}_K$
and ${\cal F}$ be a locally constant
constructible sheaf of
$\overline{\mathbb F}_\ell$-modules
on $U$.
We assume that
${\cal F}$ 
is tamely ramified 
on the generic fiber.

{\rm 1.}
For an exact sequence
$0\to {\cal F}'
\to {\cal F}
\to {\cal F}''\to 0$
of locally constant
constructible sheaves of
$\overline{\mathbb F}_\ell$-modules
on $U$,
we have
\setcounter{equation}0
\begin{equation}
{\rm Sw}_U{\cal F}=
{\rm Sw}_U{\cal F}'+
{\rm Sw}_U{\cal F}''
\label{eqaddS}
\end{equation}
in 
$F_0G(\partial_{\cal F}U)_{\mathbb Q(\zeta_{p^\infty})}$.

{\rm 2.}
Let
$U_1\subset U$
be a regular closed subscheme
and $U_0=U\setminus U_1$
be the complement.
For the immersions
$i\colon U_1\to U$
and
$j\colon U_0\to U$,
we have
\begin{equation}
{\rm Sw}_U{\cal F}
=
j_!{\rm Sw}_{U_0}{\cal F}|_{U_0}
+
i_*{\rm Sw}_{U_1}{\cal F}|_{U_1}.
\label{eqexcS}
\end{equation}
in 
$F_0G(\partial_{\cal F}U)_{\mathbb Q(\zeta_{p^\infty})}$.
\end{pr}

{\it Proof.}
1.
Clear from (\ref{eqswan}).

2.
It follows from
(\ref{eqswan})
and the excision formula
Theorem \ref{thmexc}.
\qed

For a smooth $\bar {\mathbb Q}_\ell$-sheaf
${\cal F}$ on $U$,
its Swan class
${\rm Sw}_U{\cal F}$ 
is defined as the Swan class
${\rm Sw}_U\bar {\cal F}$ 
of the reduction 
$\bar {\cal F}
={\cal F}_0/\lambda {\cal F}_0$
modulo $\ell$.
Though the $\bar {\mathbb F}_\ell$-sheaf
$\bar {\cal F}$
itself depend on the choice of a lattice
${\cal F}_0$,
is defined as the Swan class
${\rm Sw}_U\bar {\cal F}$ 
is well-defined by
Proposition \ref{prExc}.1.


We prove an induction formula
for the Swan classes.
The following
results are regarded as the relative conductor formula in the case of relative dimension $0$.

\begin{pr}\label{prind}
Let $U$ be a regular flat
separated scheme 
of finite type over $S$
and $f\colon V\to U$
be a finite \'etale
$G$-torsor for
a finite group $G$
such that
the generic fiber
$V_K\to U_K$
is tamely ramified
with respect to $K$.
Let $H\subset G$
be a subgroup
and $g\colon V\to U'$
be the corresponding
$H$-torsor.
Let ${\cal F}$ be a locally constant
constructible sheaf of
$\overline{\mathbb F}_\ell$-modules
on $U'$
such that
$g^*{\cal F}$
is a constant sheaf
on $V$ and let
$h\colon U'\to U$
denote the canonical map.

Let $T\subset G$
be a complete set of
representatives
of $G/H$.
Then, 
we have
\setcounter{equation}0
\begin{equation}
{\rm Sw}_{V/U}
h_*{\cal F}=
\sum_{\tau\in T}
\tau_*({\rm Sw}_{V/U'}{\cal F}
+{\rm rank}\ {\cal F}
\cdot g^*D^{\log}_{U'/U})
\label{eqind}
\end{equation}
in
$F_0G(\partial_{V/U}V)
_{\mathbb Q(\zeta_{p^\infty})}$.
\end{pr}

{\it Proof.}
Let $M$ be the representation
of $H$ corresponding to ${\cal F}$.
Then, we have
${\rm Tr}^{\rm Br}(\sigma:{\rm Ind}_H^GM)=
\sum_{\tau\in T,
\tau^{-1}\sigma
\tau\in H}
{\rm Tr}^{\rm Br}(\tau^{-1}\sigma
\tau:M)$.
Hence,
by the definition of
the Swan class,
the left hand side
of (\ref{eqind}) is equal to
\begin{align}
\sum_{\sigma\in G_{(p)}}
{\rm Tr}^{\rm Br}(\sigma:{\rm Ind}_H^GM)
\cdot
s_{V/U}(\sigma)
=&
\sum_{\sigma\in G_{(p)}}
\sum_{\tau\in T,
\tau^{-1}\sigma
\tau\in H}
{\rm Tr}^{\rm Br}(\tau^{-1}\sigma
\tau:M)
\cdot
s_{V/U}(\sigma)
\label{eqind3}
\\
=&
\sum_{\tau\in T}
\tau_*
\left(\sum_{\rho\in H_{(p)}}
{\rm Tr}^{\rm Br}(\rho:M)
\cdot
s_{V/U}(\rho)\right).
\nonumber
\end{align}
Thus, it follows from
Lemma \ref{lmSwch}.2.
\qed

\begin{cor}\label{corind}
Let $f\colon U\to V$ be a 
finite \'etale morphism
of regular flat
schemes of finite type over $S$
such that
the generic fiber
$U_K\to V_K$
is tamely ramified
with respect to $K$.
Let ${\cal F}$
be a locally constant
constructible
sheaf of
$\overline {\mathbb F}_\ell$-modules
on $U$
such that
there exists
a finite \'etale
morphism 
$g\colon U'\to U$ over $S$
such that
$g^*{\cal F}$
is constant on $U'$ and
that $U'_K\to V_K$
is tamely ramified
with respect to $K$.

Then, 
we have
\setcounter{equation}0
\begin{equation}
{\rm Sw}_Vf_*{\cal F}=
f_*{\rm Sw}_U{\cal F}
+{\rm rank}\ {\cal F}\cdot d_{U/V}^{\log}.
\label{eqindd}
\end{equation}
In particular,
for ${\cal F}=
\overline{\mathbb F}_\ell$,
we obtain
\begin{equation}
{\rm Sw}_Uf_*\overline{\mathbb F}_\ell=
d_{U/V}^{\log}.
\label{eqSw1}
\end{equation}
\end{cor}

{\it Proof.}
It follows immediately
from (\ref{eqswan}),
(\ref{dfdisV}),
Proposition \ref{prind}
and the remark on the Galois closure
after Definition \ref{dftmT}.
\qed

We expect the following
generalization of the Hasse-Arf theorem
to hold.

\begin{cn}\label{cnHA}
The Swan class
${\rm Sw}_U{\cal F}$
is in the image of the map
$F_0G(\partial_{V/U}U)
\to F_0G(\partial_{V/U}U)_{\mathbb Q(\zeta_{p^\infty})}$.
\end{cn}

Note that Conjecture \ref{cnHA}
is much stronger than
the statement that
the Swan class
${\rm Sw}_{V/U}{\cal F}$
is in the image of the map
$F_0G(\partial_{V/U}V)
\to F_0G(\partial_{V/U}V)_{\mathbb Q(\zeta_{p^\infty})}$,
which is, as we have seen above,
a consequence of 
a strong form
of equivariant
resolution of
singularities.
We will prove Conjecture
\ref{cnHA} in the case 
$\dim U_K\le 1$
later at Corollary \ref{corint}.

Conjecture \ref{cnHA}
is related to
the following
conjecture of Serre.

\begin{cn}[{\cite[Section 6]{Sear}}]\label{cnSe}
Let $A$ be a regular
local noetherian ring
and $G$ be a
finite group of
automorphisms of $A$.
Assume that
the fixed part
$A^G$ is noetherian
and that
for every $\sigma\in G,
\sigma\neq 1$,
the quotient
$A/I_\sigma$
by the ideal 
$I_\sigma=
\langle \sigma(a)-a;
a\in A\rangle$
is of finite length.
Then the ${\mathbb Z}$-valued
function
$a_G$ of $G$
defined by
$$a_G(\sigma)=
\begin{cases}
{\rm length}
\ A/I_\sigma
\quad&\text{ if }
\sigma\neq 1\\
-\sum_{\tau
\in G,\tau\neq 1}
a_G(\tau)
\quad&\text{ if }
\sigma= 1
\end{cases}$$
is a character
of $G$.
\end{cn}

We prove Conjecture \ref{cnSe}
in the case where $\dim A=2$ at
the end of Section \ref{ssrk1}.

\begin{lm}\label{lmSe}
Assume that the fraction
field of $A$ is of characteristic $0$
and that the residue
field $F$ of $A$ is 
of characteristic $p>0$.
Then,
Conjecture {\rm \ref{cnHA}}
for $U$ such that
$n=\dim U_K+1$
implies 
Conjecture {\rm \ref{cnSe}}
for $A$ of dimension $n$.
\end{lm}

{\it Proof.}
First, we consider
the following special case.
Let $Y$ be a regular
flat separated
scheme of finite type
over $S={\rm Spec}\
{\cal O}_K$
and $y\in Y$ be
a closed point
in the closed fiber
as in Proposition
\ref{priso}.
Let $G$ be a finite
group of automorphisms
of $Y$ over $S$
such that,
for every $\sigma\in G,
\sigma\neq 1$,
the fixed part
$Y^\sigma$
is equal to
$\{y\}$
set-theoretically.
We assume that
the quotient
$X=Y/G$ exists as a scheme
of finite type
over $S$.
We put $V=Y\setminus \{y\}$
and $f\colon V
\to U=V/G\subset X$.

We show that
Conjecture \ref{cnHA}
for $f\colon V\to U$
and $G$
implies
Conjecture \ref{cnHA}
for $A={\cal O}_{Y,y}$
and $G$.
We consider the image
$s_G(\sigma)
\in {\mathbb Q}$
of
$s_{V/U}(\sigma)
\in F_0G(\partial_
{V/U}V)_{\mathbb Q}$
by
$F_0G(\partial_
{V/U}V)_{\mathbb Q}\to
F_0G(\Sigma_
{V/U}Y)_{\mathbb Q}
\to
F_0G(\{y\})_{\mathbb Q}
={\mathbb Q}$.
Conjecture \ref{cnHA}
implies
that the function
$s_G(\sigma)$
is a character of
$G$.
By Proposition
\ref{priso},
we have
$a_G=r_G-u_G+s_G$
where
$r_G$ and $u_G$
denote the
characters
of the regular and the unit
representations of $G$
respectively.
Hence, the assertion
is proved in this case.

We reduce
the general case
to the special
case above
similarly 
as in the proof of
\cite[Lemma (5.3)]{Artin}.
By replacing $A$ by 
the completion,
we may assume $A$
and hence $A^G$
are complete.
Let $C$
be a complete valuation
ring such that
$p$ is a prime of $C$
and the residue field
$k$ is the same
as that of
$A^G$.
Then, 
by \cite[Chapitre 0, Th\'eor\`eme 19.8.8 (ii)]{EGA41}
there exists
a finite injection
$C[[t_1,\ldots,t_{n-1}]]
\to A^G.$
Let $W(\bar k)$
be the ring of Witt
vectors
and we take a 
local homomorphism
$C\to W$.
By replacing
$A$ by a factor of
the completion
$A\hat \otimes_CW$,
we may assume
$k$ is algebraically
closed
and $C=W$.
Then, the rest
of the argument
is the same as
that in the proof of
\cite[Lemma (5.3)]{Artin}
by replacing
$k[[t_1,\ldots,t_{n-1}]]$
by
$W[[t_1,\ldots,t_{n-1}]]$
and
$k\{t_1,\ldots,t_{n-1}\}$
by the strict localization
$W\{t_1,\ldots,t_{n-1}\}$.
\qed

Conjecture \ref{cnsi} implies
the weaker statement 
that
the Swan class
${\rm Sw}_U{\cal F}$
is in the image of the map
$F_0G(\partial_{V/U}U)_{\mathbb Q}
\to F_0G(\partial_{V/U}U)_{\mathbb Q(\zeta_{p^\infty})}$.
Since we do not know the proof of
Conjecture
\ref{cnsi} in general,
we make the following definition.
Let 
${\rm pr}_{{\mathbb Q}(\zeta_{p^\infty})/
{\mathbb Q}}:
F_0G(\partial_{V/U}U)_{\mathbb Q(\zeta_{p^\infty})}
\to
F_0G(\partial_{V/U}U)_{\mathbb Q}$
be the projection
induced by
$\varinjlim_n
\displaystyle{
\frac1{[\mathbb Q(\zeta_{p^n}):
\mathbb Q]}{\rm Tr}_{
\mathbb Q(\zeta_{p^n})/
\mathbb Q}}$.

\begin{df}\label{dfSwQ}
The rationalized Swan class
${\rm Sw}^{\mathbb Q}_U{\cal F}\in
F_0G(\partial_{V/U}U)_{\mathbb Q}$
is the image of
the Swan class
${\rm Sw}_U{\cal F}$
by the projection
${\rm pr}_{{\mathbb Q}(\zeta_{p^\infty})/
{\mathbb Q}}:
F_0G(\partial_{V/U}U)_{\mathbb Q(\zeta_{p^\infty})}
\to
F_0G(\partial_{V/U}U)_{\mathbb Q}$.
\end{df}

In \cite{KSA},
the rationalized Swan class
is called the Swan class
and is denoted by ${\rm Sw}$.
By modifying the notation there
we write the Swan class
by ${\rm Sw}$
and
the rationalized Swan class
by ${\rm Sw}^{\mathbb Q}$.

\subsection{Conductor formula
for a relative curve}\label{sscfs}

Let $f\colon U'\to U$
be a smooth morphism
of separated regular flat
schemes
of finite type over $S$
and ${\cal F}$ be
a locally constant constructible
sheaf of $\overline{\mathbb F}_\ell$-modules
on $U'$.
Let
$\pi'\colon V'\to U'$
be a finite \'etale morphism
such that
the pull-back
$\pi^{\prime*}{\cal F}$
is a constant sheaf
on $V'$.
We assume that
the following conditions
are satisfied:
\begin{itemize}
\item[(\ref{sscfs}.0.1)]
There exists a proper
smooth scheme $X'$
over $U$
containing $U'$
as the complement
$U'=X'\setminus D$
of a divisor $D$
with simple normal crossing relatively to 
$U$.
The finite \'etale morphism
$\pi'\colon V'\to U'$
is tamely ramified
with respect to $X'$.
\end{itemize}
Then, 
by \cite{App},
the higher direct images
$R^qf_!
{\mathbb F}_\ell$
and $R^qf_!{\cal F}$
are
locally constant
sheaves on $U$.
Further,
for the alternating
sum of ranks,
we have 
\setcounter{equation}1
\begin{equation}
{\rm rank}\
Rf_!{\cal F}=
{\rm rank}\ {\cal F}\cdot
{\rm rank}\ Rf_!
\overline {\mathbb F}_\ell.
\label{eqrank}
\end{equation}

We assume that
the locally constant sheaves
$R^qf_!{\cal F}$
on $U$ are tamely ramified
on the generic fiber
and that ${\cal F}$
on $U'$ is tamely ramified
on the generic fiber.
Then the Swan class
$${\rm Sw}_URf_!{\cal F}
=\sum_q(-1)^q
{\rm Sw}_UR^qf_!{\cal F}$$
$\in F_0G(\partial_
{\cal F}U)_{\mathbb Q(\zeta_{p^\infty})}$
is also defined
as the alternating sum
and the Swan class
${\rm Sw}_{U'}{\cal F}$
is defined.
A conductor formula
should express
the Swan class
${\rm Sw}_URf_!{\cal F}$
in terms of
the class
${\rm Sw}_{U'}{\cal F}$.
We prove a conductor
formula for a relative curve
under a certain
assumption
in Corollary \ref{corcfs}
and 
in general
in Section \ref{sspc},
assuming $K$ is of characteristic 0.

First, we give
a general formalism
Proposition \ref{prcfs2}
to prove a conductor formula.
We consider 
a commutative diagram
\begin{equation}
\begin{CD}
U'@<{\pi'}<< V'\\
@VfVV@VV{g}V\\
U@<{\pi}<<V
\end{CD}
\label{eqtrcf}
\end{equation}
of regular flat separated schemes 
of finite type over $S$
satisfying the 
following condition:
\begin{itemize}
\item[(\ref{sscfs}.0.4)]
The horizontal
arrows are finite \'etale,
$V$ is a 
$G$-torsor over $U$
and $V'$ is a 
$G'$-torsor over $U'$
for finite groups
$G$ and $G'$.
The arrow $g$
is compatible with
a morphism
$\varphi\colon
G'\to G$
of finite groups
in the sense that,
for $\sigma\in G'$
and $\tau=\varphi(\sigma)\in G$,
we have
$g\circ \sigma=\tau\circ g$.
\end{itemize}

\setcounter{equation}0
\begin{lm}\label{lmtrcf}
We consider a
commutative diagram
{\rm (\ref{eqtrcf})}
of separated regular flat
schemes
of finite type over $S$
satisfying the conditions
{\rm (\ref{sscfs}.0.1)}
and {\rm (\ref{sscfs}.0.4)}.
Assume that the higher direct image
$R^qg_!
{\mathbb F}_\ell$
is a constant sheaf on $V$
for every $q\ge0$.

Let $\sigma\in G'$ be
an element and we put
$\tau=\varphi(\sigma)
\in G$.
Let $\bar \eta$
be a geometric point of $V$
and $\bar \tau$
be an automorphism
of $\bar \eta$
compatible with $\tau$. 
Let $\sigma^*
\circ \bar \tau^*$
denote the automorphisms
of
$H^q_c(V'_{\bar \eta},{\mathbb Q}_\ell)$
and
$H^q_c(V'_{\bar \eta},{\mathbb F}_\ell)$
defined by the pull-back by
$\sigma\times \bar \tau$
on $V'_{\bar \eta}=
V'\times_V{\bar \eta}$.

{\rm 1.}
The automorphism
$\sigma^*
\circ \bar \tau^*$
of 
$H^q_c(V'_{\bar \eta},
{\mathbb F}_\ell)$
is independent
of the choice of
$\bar \tau$.

{\rm 2.}
Assume 
that $\sigma$
and $\bar \tau$
are $\ell$-regular.
Let 
$\sigma^*$
denote the automorphism
$\sigma^*
\circ \bar \tau^*$
of $H^*_c(V'_{\bar \eta},
{\mathbb F}_\ell)$
independent of
$\bar \tau$.
Then,
the alternating sum 
${\rm Tr}
(\sigma^*\circ \bar \tau^*
\colon
H^*_c(V'_{\bar \eta},
{\mathbb Q}_\ell))$
is independent
of 
$\bar \tau$
and is equal 
to the alternating sum
${\rm Tr}^{\rm Br}(\sigma^*,
H^*_c(V'_{\bar \eta},
{\mathbb F}_\ell))$
of Brauer traces.
\end{lm}

{\it Proof.}
1.
By the assumption
that 
$R^qg_!{\mathbb F}_\ell$
is constant,
the action of
$\sigma\circ \bar \tau$
on 
$H^q_c(V'_{\bar \eta},{\mathbb F}_\ell)$
is independent
of the lifting $\bar \tau$
of $\tau$.

2.
Since 
$\sigma$
and $\bar \tau$
are assumed $\ell$-regular,
the actions of
$\sigma^*
\circ \bar \tau$
on 
$H^*_c(V'_{\bar \eta},{\mathbb Q}_\ell)$
and on
$H^*_c(V'_{\bar \eta},
{\mathbb F}_\ell)$
are of finite order prime-to-$\ell$.
Hence the assertion follows.
\qed

\begin{lm}\label{lmBr}
We consider a
commutative diagram
{\rm (\ref{eqtrcf})}
of separated regular flat
schemes
of finite type over $S$
satisfying the conditions
{\rm (\ref{sscfs}.0.1)}
and {\rm (\ref{sscfs}.0.4)}.
Let  ${\cal F}$ be
a locally constant constructible
sheaf of $\overline{\mathbb F}_\ell$-modules
on $U'$ such that
the pull-back $\pi^{\prime*}{\cal F}$
is a constant sheaf on $V'$.
Assume that the higher direct image
$R^qg_!
{\mathbb F}_\ell$
and
the pull-back
$\pi^*R^qf_!{\cal F}$
are constant sheaves on $V$
for every $q\ge 0$.

Let $\tau\in G$
be an $\ell$-regular
element.
Let $\bar \xi$ be
a geometric point of $U$
and lift it to
a geometric point
$\bar \eta$ of $V$.
Let $M$ be the
$\bar {\mathbb F}_\ell
$-representation of $G'$
corresponding to ${\cal F}$.
Then, we have
\setcounter{equation}0
\begin{equation}
{\rm Tr}^{\rm Br}(\tau,
H^*_c(U'_{\bar \xi},{\mathcal F}))
=
\frac 1{|G'|}
\sum_{
\begin{array}c
\scriptstyle{
\sigma\in G^{\prime(\ell)},\
\rho\in G;}\\
\scriptstyle{
\varphi(\sigma)
=\rho\tau \rho^{-1}}
\end{array}}
{\rm Tr}^{\rm Br}
(\sigma^*:
H^*_c(V'_{\bar \eta},
{\mathbb F}_\ell))
\cdot
{\rm Tr}^{\rm Br}(\sigma^*:M).
\label{eqBr}
\end{equation}
\end{lm}

{\it Proof.}
By the assumption that
$\pi^*R^qf_!{\cal F}$
are constant on $V$,
we may regard
$H^q_c(U'_{\bar \xi},{\cal F})$
as an
$\bar {\mathbb F}_\ell$-representation
of $G$
and the Brauer trace
${\rm Tr}^{\rm Br}(\tau,
H^*_c(U'_{\bar \xi},{\mathcal F}))$
is defined.

We consider
$H^q_c(V'_{\bar \xi},\bar {\mathbb F}_\ell)$
as an
$\bar {\mathbb F}_\ell$-representation
of $G'$ by
the pull-back action of $G'$
on $V'_{\bar \xi}
=V'\times_U{\bar \xi}$.
By \cite[Lemma 2.2]{Il},
the alternating sum
$H^*_c(V'_{\bar \xi},\bar {\mathbb F}_\ell)$
defines an element
$[H^*_c(V'_{\bar \xi},\bar {\mathbb F}_\ell)]$
of the Grothendieck group
$P_\ell(G')$
of the exact category of
finitely generated
projective $\bar {\mathbb F}_\ell[G']$-modules
(\cite[Partie III 1.3]{SeRF}).

Since $\tau$
is assumed $\ell$-regular,
we have a decomposition
$H^*_c(V'_{\bar \xi},\bar {\mathbb F}_\ell)
=\bigoplus_{
\chi\in \widehat{\langle\tau\rangle}}
H^*_c(V'_{\bar \xi},\bar {\mathbb F}_\ell)_{\chi}$
by characters 
of the cyclic subgroup
$\langle\tau\rangle
\subset G$
and an equality
$[H^*_c(V'_{\bar \xi},\bar {\mathbb F}_\ell)]
=\sum_{
\chi\in \widehat{\langle\tau\rangle}}
[H^*_c(V'_{\bar \xi},
\bar {\mathbb F}_\ell)_{\chi}]$
in $P_\ell(G')$.
Similarly as 
\cite[Lemma 2.2]{Il},
we obtain
\begin{equation}
\dim
H^*_c(U'_{\bar \xi},{\mathcal F})_\chi
=\dim 
([H^*_c(V'_{\bar \xi},\bar {\mathbb F}_\ell)_{\chi}]
\cdot [M])^{G'}.
\label{eqchi}
\end{equation}
Thus, we have
\begin{equation}
{\rm Tr}^{\rm Br}(\tau,
H^*_c(U'_{\bar \xi},{\mathcal F}))
=\sum_{
\chi\in \widehat{\langle\tau\rangle}}
\chi(\tau)\cdot
\dim 
([H^*_c(V'_{\bar \xi},\bar {\mathbb F}_\ell)_{\chi}]
\cdot [M])^{G'}.
\label{eqtau}
\end{equation}
Similarly as \cite[Lemma 2.3]{Il},
the right hand side
of (\ref{eqtau})
is equal to
\begin{eqnarray*}
&&\sum_{
\chi\in \widehat{\langle\tau\rangle}}
\chi(\tau)\cdot
\frac 1{|G'|}
\left(\sum_{\sigma \in G^{\prime(\ell)}}
{\rm Tr}^{\rm Br}(\sigma^*:
H^*_c(V'_{\bar \xi},\bar {\mathbb F}_\ell)_{\chi})
\cdot
{\rm Tr}^{\rm Br}(\sigma^*:M)\right)
\\
&=&
\frac 1{|G'|}
\sum_{\sigma \in G^{\prime(\ell)}}
\left(\sum_{
\chi\in \widehat{\langle\tau\rangle}}
\chi(\tau)\cdot
{\rm Tr}^{\rm Br}(\sigma^*:
H^*_c(V'_{\bar \xi},\bar {\mathbb F}_\ell)_{\chi})
\right)\cdot
{\rm Tr}^{\rm Br}(\sigma^*:M)
\\
&=&
\frac 1{|G'|}
\sum_{\sigma \in G^{\prime(\ell)}}
{\rm Tr}^{\rm Br}
(\sigma^*\times \bar\tau^*:
H^*_c(V'_{\bar \xi},
{\mathbb F}_\ell))
\cdot
{\rm Tr}^{\rm Br}(\sigma^*:M).
\end{eqnarray*}

For $\rho\in G$,
let $V'_{\rho(\bar \eta)}$
denote the geometric fiber
of $V'\to V$
by the composition
of $\bar \eta\to V$
with $\rho\colon V\to V$.
Then, the geometric fiber
$V'_{\bar \xi}
=V'\times_U{\bar \xi}$
is the disjoint union
$\coprod_{\rho\in G}
V'_{\rho(\bar \eta)}$
and we have
$$
{\rm Tr}^{\rm Br}(\sigma^*\times \bar\tau^*:
H^*_c(V'_{\bar \xi},{\mathbb F}_\ell))
=
\sum_{\rho\in G;
\
\varphi(\sigma)
=\rho\tau \rho^{-1}}
{\rm Tr}^{\rm Br}
(\sigma^*\circ \bar\tau^*:
H^*_c(
V'_{\rho(\bar \eta)},
{\mathbb F}_\ell)).$$
We have
${\rm Tr}^{\rm Br}
(\sigma^*\circ \bar\tau^*:
H^*_c(
V'_{\rho(\bar \eta)},
{\mathbb F}_\ell))=
{\rm Tr}^{\rm Br}
(\sigma^*:
H^*_c(V'_{\bar \eta},
{\mathbb F}_\ell))$
since
$R^qg_!{\mathbb F}_\ell$
are assumed constant on $V$.
\qed

\begin{pr}\label{prcfs2}
\setcounter{equation}0
We consider a
commutative diagram
{\rm (\ref{eqtrcf})}
of separated regular flat
schemes
of finite type over $S$
satisfying the conditions
{\rm (\ref{sscfs}.0.1)}
and {\rm (\ref{sscfs}.0.4)}.
Assume that $U,U',V$ and $V'$
are irreducible
and that $V_K\to U_K$
and $V'_K\to U'_K$
are tamely ramified
with respect to $K$.

Let $\ell$ be
a prime number invertible
on $S$ and
${\cal F}$ be
a locally constant constructible
sheaf of $\overline{\mathbb F}_\ell$-modules
on $U'$.
Assume further
that the pull-backs
$\pi^*R^qf_!\pi'_*
{\mathbb F}_\ell$
and $\pi^*R^qf_!{\cal F}$
are constant on $V$.
Let $\bar \xi$ be
a geometric point of $U$
and $\bar \eta$ be
a geometric point of $V$
above $\bar \xi$.
For an $\ell$-regular element
$\sigma\in G'$,
define
${\rm Tr}^{\rm Br}
(\sigma^*:
H^*_c(V'_{\bar \eta},
{\mathbb F}_\ell))$
as in Lemma {\rm \ref{lmtrcf}}.
Then, we have
the following.

{\rm 1.}
Assume that,
for each {\rm non-trivial}
$\ell$-regular element
$\sigma\in G'$
and
$\tau=\varphi(\sigma)
\in G$,
we have
\begin{equation}
{\rm Tr}^{\rm Br}
(\sigma^*:
H^*_c(V'_{\bar \eta},
{\mathbb F}_\ell))
\cdot
((\Gamma_\tau,
\Delta_V))^{\log}
=
g_!((\Gamma_\sigma,
\Delta_{V'}))^{\log}
\label{eqtrcd}
\end{equation}
in $F_0G(\partial_{U/V}V)
_{\mathbb Q}$.
Then, we have
\begin{equation}
{\rm Sw}_URf_!{\cal F}-
{\rm rank}\ {\cal F}\cdot
{\rm Sw}_URf_!\overline
{\mathbb F}_\ell
=
g_!{\rm Sw}_{U'}{\cal F}
\label{eqpcs}
\end{equation}
in $F_0G(\partial_{V/U}U)
_{\mathbb Q(\zeta_{p^\infty})}$.

{\rm 2.}
Assume $K$ is of characteristic $0$
and suppose
$U$ and hence $V,U',V'$ 
are schemes over $K$.
Assume that,
for each 
$\ell$-regular element
$\sigma\in G'$
and
$\tau=\varphi(\sigma)\in G$,
we have
\begin{equation}
{\rm Tr}^{\rm Br}
(\sigma^*:
H^*_c(V'_{\bar \eta},
{\mathbb F}_\ell))
\cdot
((\Gamma_\tau,
\Delta_V))^{\log}
=
g_!((\Gamma_\sigma,
\Delta_{V'}))^{\log}
\label{eqtrcd0}
\end{equation}
in $F_0G(\partial_FV)
_{\mathbb Q}$.
Then, we have
\begin{equation}
{\rm Sw}_URf_!\overline
{\mathbb F}_\ell
=
-f_!
((\Delta_{U'},\Delta_{U'}))^{\log}
+
\chi_c(U'_{\bar \xi})
\cdot
((\Delta_U,\Delta_U))^{\log}
\label{eqpcs1}
\end{equation}
in $F_0G(\partial_FU)
_{\mathbb Q(\zeta_{p^\infty})}$.
\end{pr}

{\it Proof.}
1.
Let $M$ denote
the representation of $G'$
corresponding to
the locally constant 
$\overline {\mathbb F}_\ell$-sheaf
${\cal F}$
on $U'$.
By (\ref{eqrank})
and by the definition
of the Swan class
(\ref{eqSwn}),
the equality (\ref{eqpcs})
is equivalent to the following:
\begin{align}
&
\frac1{|G|}
\sum_{\tau\in G^{(\ell)},
\tau\neq 1}
\Bigl(
{\rm Tr}^{\rm Br}(\tau,
H^*_c(U'_{\bar \xi},{\cal F}))
-{\rm rank}\ {\cal F}\cdot
{\rm Tr}^{\rm Br}(\tau,
H^*_c(U'_{\bar \xi},{\mathbb F}_\ell))
\Bigr)
\cdot
\pi_*
((\Gamma_\tau,
\Delta_{V}))^{\log}
\label{eqSwBr}
\\
=\ &
\frac1{|G'|}
\sum_{\sigma\in G^{\prime(\ell)},
\sigma\neq 1}
({\rm Tr}^{\rm Br}(\sigma,M)-
\dim M)\cdot
f_!\pi'_*
((\Gamma_{\sigma},
\Delta_{V'}))^{\log}.
\nonumber
\end{align}
Substituting 
(\ref{eqBr}),
we see that the left hand side
of (\ref{eqSwBr}) is
\begin{align}
&
\frac1{|G|}
\sum_{\tau\in G^{(\ell)},
\tau\neq 1}
\frac 1{|G'|}
\sum_{
\begin{array}c
\scriptstyle{
\sigma\in G^{\prime(\ell)},
\sigma\neq 1\
\rho\in G;}\\
\scriptstyle{
\varphi(\sigma)
=\rho\tau \rho^{-1}}
\end{array}}
{\rm Tr}^{\rm Br}
(\sigma^*:
H^*_c(V'_{\bar \eta},
{\mathbb F}_\ell))
\cdot
({\rm Tr}^{\rm Br}(\sigma,M)-
\dim M)
\cdot
\pi_*
((\Gamma_\tau,
\Delta_{V}))^{\log}
\label{eqSwBr3}
\\
=&
\frac1{|G'|}
\sum_{\sigma\in G^{\prime(\ell)},
\sigma\neq 1}
({\rm Tr}^{\rm Br}(\sigma,M)-
\dim M)
\cdot
\frac 1{|G|}
\sum_{\rho\in G}
{\rm Tr}^{\rm Br}
(\sigma^*:
H^*_c(V'_{\bar \eta},
{\mathbb F}_\ell))
\cdot
\pi_*((\Gamma_{\rho^{-1}
\varphi(\sigma)\rho},
\Delta_V))^{\log}.
\nonumber
\end{align}
By the assumption
(\ref{eqtrcd}),
each term in
the second summation
in the second line of
(\ref{eqSwBr3}) is 
$$
\pi_*\rho_*g_!
((\Gamma_\sigma,
\Delta_{V'}))^{\log}
=
f_!\pi'_*
((\Gamma_\sigma,
\Delta_{V'}))^{\log}.$$
Thus the equality
(\ref{eqSwBr}) 
follows.

2.
Similarly as above,
the left hand
side of the
equality (\ref{eqpcs1})
is equal to
\begin{equation}
-
\frac1{|G|}
\sum_{\tau\in G^{(\ell)}}
{\rm Tr}^{\rm Br}(\tau,H^*_c(U'_{\bar \xi},
{\mathbb F}_\ell))\cdot
\pi_*
((\Gamma_\tau,
\Delta_V))^{\log}
+
\frac1{|G|}
\sum_{\tau\in G^{(\ell)}}
\dim H^*_c(U'_{\bar \xi},
{\mathbb F}_\ell)
\cdot
\pi_*
((\Gamma_\tau,
\Delta_{V}))^{\log}
\label{eqSwBr2}
\end{equation}
Substituting 
(\ref{eqBr}),
we see that
the first term in
(\ref{eqSwBr2}) is equal to
$$
\frac 1{|G'|}
\sum_{\sigma \in G^{\prime(\ell)}}
{\rm Tr}^{\rm Br}(\sigma:1)
\cdot
\frac1{|G|}
\sum_{\rho\in G}
{\rm Tr}^{\rm Br}
(\sigma^*:
H^*_c(V'_{\bar \eta},
{\mathbb F}_\ell))
\cdot
\pi_*((\Gamma_{\rho^{-1}
\varphi(\sigma)\rho},
\Delta_V))^{\log}.
$$
Hence, similarly as above,
by the assumption
(\ref{eqtrcd0}),
it is further equal to
$$
\frac 1{|G'|}
\sum_{\sigma \in G^{\prime(\ell)}}
f_!\pi'_*
((\Gamma_\sigma,
\Delta_{V'}))^{\log}
=
f_!
\frac 1{|G'|}
\pi'_*\pi^{\prime*}
((\Delta_{U'},
\Delta_{U'}))^{\log}
=
f_!
((\Delta_{U'},
\Delta_{U'}))^{\log}.$$
Similarly,
the second term
in (\ref{eqSwBr2})
is equal to
$\chi_c(U'_{\bar \xi})$times
$\frac 1{|G|}
\sum_{\tau \in G^{(\ell)}}
\pi_*((\Gamma_\tau,
\Delta_V))^{\log}
=
((\Delta_U,
\Delta_U))^{\log}$
and
the equality
(\ref{eqpcs1})
follows.
\qed

We derive
from the crucial step 
Proposition \ref{prcysd}
a sufficient condition
on a relative curve
for the assumptions 
(\ref{eqtrcd}) 
and (\ref{eqtrcd0}) of
Proposition \ref{prcfs2}
to be satisfied.

\setcounter{equation}0
\begin{pr}\label{prcfs}
Let $f\colon U'\to U$
be a smooth morphism
of relative dimension $1$
of separated regular flat schemes
of finite type over $S$.
Let $\ell$ be
a prime number invertible
on $S$ and
${\cal F}$ be
a locally constant constructible
sheaf of $\overline{\mathbb F}_\ell$-modules
on $U'$.
We consider a commutative
diagram {\rm (\ref{eqtrcf})}
satisfying the conditions
{\rm (\ref{sscfs}.0.1)}
and
{\rm (\ref{sscfs}.0.4)}.
Assume that the pull-back
$\pi^{\prime*}
{\cal F}$
is constant on $V'$
and that
$\pi^*
R^qf_!\pi'_*
{\mathbb F}_\ell$
is constant on $V$
for every $q\ge 0$.

We assume that
$V$ and $V'$
are integral.
Let $\eta$ denote
the generic point
of $V$ 
and let
$\bar \eta$ be a
geometric point above $\eta$.
Let $\sigma\in G'$
and $\tau=\varphi(\sigma)
\in G$ be
{\rm non-trivial}
$\ell$-regular elements.

We also consider a commutative
diagram
\begin{equation}
\begin{CD}
V'@>{\subset}>> Y'
&\ \supset D\\
@VVgV
@VV{\bar g}V\\
V@>{\subset}>>Y
\end{CD}
\label{eqtrcf2}
\end{equation}
of separated schemes 
of finite type over $S$
and a finite family ${\cal E}$ of
Cartier divisors $Y$
satisfying
the following condition:
\begin{itemize}
\item[{\rm (\ref{prcfs}.2)}]
The schemes $Y$ and $Y'$ are proper over $S$.
There exist
$d$ sections 
$(t_k\colon Y\to Y')_k$
such that the pair
$(Y',(t_k))$
is a $d$ pointed stable
curve of genus $g$
over $Y$
and that
$D\subset Y'$
is the disjoint union of
the sections $t_k(Y)$. 
The open subscheme $V
\subset Y$
is the complement
of the union of ${\cal E}$.
The restriction 
$Y'_V=Y'\times_YV\to V$
is smooth 
and 
$V'\subset Y'_V$
is the complement
of $D_V$.
The action of 
$G$ on $V$ is extended
to an admissible action on $Y$
and on ${\cal E}$,
in the sense that the quotient
$X=Y/G$ exists as a scheme.

Let ${\cal E}'$
denote the family of Cartier divisors
of $Y'$ defined as the union of
the pull-back of ${\cal E}$
and the sections
$(t_k(Y))_k$.
Then, $\Sigma_{V/U}^{\cal E}Y_K$
and $\Sigma_{V'/U'}^{{\cal E}'}Y'_K$
are empty.
\addtocounter{equation}1
\end{itemize}
Then, we have
\begin{equation}
\bar g_*
((\Gamma_{\sigma},
\Delta_{Y'}^{\log}))
=
{\rm Tr}^{\rm Br}(\sigma^*,
H^*_c(V'_{\bar \eta},
{\mathbb F}_\ell))
\cdot
((\Gamma_\tau,
\Delta_Y^{\log}))
\label{eqpc2}
\end{equation}
in
$F_0G(\Sigma_{V/U}Y)_{\mathbb Q}$,
\end{pr}

{\it Proof.}
Since the smooth compactification
$Y'_V\supset V$
is unique,
the action of
$G'$ on $V'$ is
extended uniquely to that on $Y'_V$
compatible with
the action of $G$ on $V$.
Further, since
an extension $Y'\supset Y'_V$
is unique,
the action of
$G'$ on $Y'_V$ is
extended uniquely to that on $Y'$
compatible with
the admissible action of $G$ on $Y$.
Since $Y'\to Y$
is projective
and $G'$ acts on 
the relatively ample sheaf
$\Omega^1_{Y'/Y}(\log D)$,
the action of $G'$ on $Y'$
is admissible in the sense
that the quotient $X'=Y'/G'$
exists as a scheme.
Hence, by Lemma \ref{lmXYG},
the quotients
$X=Y/G$ and $X'=Y'/G'$
contains $U$ and $U'$
as the complements
of Cartier divisors respectively.
In particular,
we have a commutative
diagram (\ref{eqcysd2})
satisfying (\ref{ssss}.1.4).

Let $T\subset
V'\times_VV'$
be the graph 
$\Gamma_\tau$ of $\tau$.
The base changes
$Y_T^{\prime(1)}$ and
$Y_T^{\prime(2)}$ 
defined as in Proposition \ref{prcysd}
are isomorphic to $Y'_V$.
We define the log product
and the log blow-up
$(Y_T^{\prime(1)}\times_T
Y_T^{\prime(2)})^\sim
\subset
(Y_T^{\prime(1)}\times_T
Y_T^{\prime(2)})'$
with respect to the log structure
defined by the sections
$(t_k(T))_k$.
Since the automorphism
$\sigma$ of $Y'_V=Y_T^{\prime(1)}$ 
permutes
the sections $(t_k(T))_k$
of $Y_T^{\prime(1)}$,
the intersection
$t_k(T)\cap
\sigma(t_k(T))$
is either empty or
equal to a divisor $t_k(T)$
of $Y_T^{\prime(1)}$.
Hence by the universality of blow-up,
the closed immersion
$\gamma=(1,\sigma)\colon
Y_T^{\prime(1)}\to 
Y_T^{\prime(1)}\times_T
Y_T^{\prime(2)}$
is uniquely lifted
to a closed immersion
$\gamma'\colon
Y_T^{\prime(1)}
\to 
(Y_T^{\prime(1)}\times_T
Y_T^{\prime(2)})'$.
Let $\Gamma'\subset
(Y_T^{\prime(1)}\times_T
Y_T^{\prime(2)})'$
denote the image of
$\gamma'$
and
$\widetilde \Gamma
=\Gamma'\cap
(Y_T^{\prime(1)}\times_T
Y_T^{\prime(2)})^\sim$
be the intersection.
Then, 
$\Gamma'$ and hence
$\widetilde \Gamma$
are flat over $T$.

We regard $Y$
as a log scheme 
with the log structure
defined by ${\cal E}$
and define the log product
$(Y\times_SY)^\sim$.
By applying
Proposition
\ref{prcysd} to 
$Y'\to Y$,
we obtain
a finite family
$(K_i)_{i\in I}$
of finite extensions
of $K$,
a family 
$(\gamma_i\colon
{\rm Spec}\
K_i\to T)_{i\in I}$
of morphisms over $S$
extended to
$(\bar\gamma_i\colon
{\rm Spec}\
{\cal O}_{K_i}\to 
(Y\times_SY)^\sim)_{i\in I}$
such that the image of the closed
points
$\bar \gamma_i(s_i)$
are in the log diagonal
$\Delta_Y^{\log}
\subset (Y\times_SY)^\sim$
and 
a family $(r_i)_{i\in I}$
of rational numbers
satisfying
$((T,\Delta_Y^{\log}))
=\sum_ir_i[\bar \gamma_i(s_i)]$
and 
\begin{equation}
\bar f'_*
((\Gamma,\Delta_{Y'}^{\log}))
=
\sum_ir_i
{\rm Tr}((\gamma_i^*\Gamma)^*
\colon
H^*_c(V'_{\bar K_i},
{\mathbb Q}_\ell))
\cdot [\bar \gamma_i(s_i)]
\end{equation}
Thus, it suffices to show
\begin{equation}
{\rm Tr}^{\rm Br}(\sigma^*,
H^*_c(V'_{\bar \eta},
{\mathbb F}_\ell))
=
{\rm Tr}((\gamma_i^*\Gamma)^*
\colon
H^*_c(V'_{\bar K_i},
{\mathbb Q}_\ell))
\label{eqccf}
\end{equation}
for each $i$.

We regard the composition
$s_i\to (Y\times_SY)^\sim
\to Y$ as a morphism
of log scheme.
Since $\bar \gamma_i(s_i)
\in (Y\times_SY)^\sim$
are in the log diagonal,
the composition
$s_i\to Y\overset \tau\to Y$
of the morphisms
of log schemes is
the same as the original morphism
$s_i\to Y$ of log schemes.
Let $\bar s_i$ be a log geometric point
above $s_i$ and let
$\widetilde Y_{\bar s_i}$
be the log strict localization.
Let $\tilde \tau$
be the automorphism
of $\widetilde Y_{\bar s_i}$
induced by $\tau$.
We lift the generic geometric point
$\bar \eta\to Y$
to a geometric point
$\bar \eta\to 
\widetilde Y_{\bar s_i}$
dominating the generic point
$\tilde \eta\in
\widetilde Y_{\bar s_i}$.
Let $k_0$
be the fixed subfield of
$\kappa(\tilde \eta)$
by an automorphism
$\tilde \tau$ of order prime to $\ell$.
We take an 
$\ell$-regular lifting 
$\bar \tau
\in {\rm Gal}(
\bar\eta/k_0)$
of 
$\tilde \tau
\in 
{\rm Gal}(
\tilde \eta/k_0)
=\langle\tilde \tau\rangle$.
By applying Lemma \ref{lmtrcf}.2.,
we obtain
$${\rm Tr}^{\rm Br}(\sigma^*,
H^*_c(V'_{\bar \eta},
{\mathbb F}_\ell))
=
{\rm Tr}(\sigma^*\circ \bar \tau^*,
H^*_c(V'_{\bar \eta},
{\mathbb Q}_\ell)).$$

We deduce
(\ref{eqccf})
from Proposition \ref{prcc}.
We show that the assumption
of Proposition \ref{prcc}
is satisfied.
The intersection 
$\widetilde \Gamma
\cap
(V_T^{\prime (1)}\times_T
V_T^{\prime (2)})$
is the graph
$\Gamma_\sigma$
of $\sigma$.
Hence, the second projection
$\widetilde \Gamma
\cap
(V_T^{\prime (1)}\times_T
V_T^{\prime (2)})=\Gamma_\sigma
\to V_T^{\prime (2)}$
is proper
and
$\Gamma_\sigma$ is flat over $T$.
Thus $\Gamma_\sigma$ satisfies
the conditions 
in Proposition \ref{prcc}.
We define a map
$Y_{\bar s_i}\to S_0$
to a regular noetherian scheme
satisfying 
the condition (\ref{prcc}.2)
for $Y_{\bar t}\to S$
in the notation there.
We consider the map
$Y\to 
\bar {\cal M}_{g,d}$
to the moduli
of $d$ pointed stable curves
of genus $g$
defined by 
the pointed stable
curve $(X,(t_k))$.
Let
$S_0$ be the strict localization
of $\bar {\cal M}_{g,d}$
at the geometric point $\bar s_i$.
Then, the map
$Y_{\bar s_i}\to S_0$
satisfies the
condition
(\ref{prcc}.2).

Since $\bar \tau$
is compatible with
$\tilde \tau$,
we may apply
Proposition \ref{prcc}
and we obtain
$${\rm Tr}(\sigma^*\circ \bar \tau^*,
H^*_c(V'_{\bar \eta},
{\mathbb Q}_\ell))
=
{\rm Tr}((\gamma_i^*\Gamma)^*
\colon
H^*_c(U'_{\bar K_i},
{\mathbb Q}_\ell)).
$$
Thus, the equality (\ref{eqccf})
is proved.
\qed

Similarly,
if $K$ is of
characteristic $0$,
the same argument
gives us the following variant.

\begin{pr}\label{prcfs0}
We assume
$K$ is of
characteristic $0$.
Let the assumption
be the same as
in Proposition {\rm\ref{prcfs}}
except
that $U$ is a scheme over $K$
and that we do not
assume $\sigma$
or $\tau$ be non-trivial.
Then, the equality 
{\rm (\ref{eqpc2})}
holds in
$F_0G(Y_F)_{\mathbb Q}$.
\end{pr}

We derive a conductor
formula for relative
curves
from 
Proposition \ref{prcfs2}
assuming ${\rm char}\ K=0$.

\begin{cor}\label{corcfs}
Assume that
$K$ is of characteristic $0$.
Let $f\colon U'\to U$
be a smooth morphism
or relative dimension $1$
of separated regular flat schemes
of finite type over $S$
and let $\ell$
be a prime number
invertible on $S$.
We assume that
$U$ and $U'$ are connected.
Let
$\bar f\colon X'\to U$
be a proper smooth curve
with geometrically connected
fibers of genus $g$
and let
$D$ be a divisor
of $X'$
finite \'etale of degree $d$
over $U$
such that
$U'=X'\setminus D$
and $2g-2+d>0$.

{\rm 1.}
\setcounter{equation}0
Let ${\cal F}$
be a locally constant constructible
sheaf 
of $\overline{\mathbb F}_\ell$-modules
on $U'$.
Then, there exists
a finite \'etale covering
$\pi\colon V\to U$
such that
we have an equality
\begin{equation}
{\rm Sw}_URf_!{\cal F}-
{\rm rank}\ {\cal F}\cdot
{\rm Sw}_URf_!\overline
{\mathbb F}_\ell
=
f_!{\rm Sw}_{U'}{\cal F}
\label{eqpcsc}
\end{equation}
in $F_0G(\partial_{V/U}U)
_{\mathbb Q(\zeta_{p^\infty})}$.

{\rm 2.}
Assume that
the schemes $U$ and $U'$
are schemes over $K$.
Then, we have
\begin{equation}
{\rm Sw}_URf_!\overline
{\mathbb F}_\ell
=-
f_!
((\Delta_{U'},\Delta_{U'}))^{\log}
+
\chi_c(U'_{\bar \xi})
\cdot
((\Delta_U,\Delta_U))^{\log}
\label{eqpcsc1}
\end{equation}
in $F_0G(\partial_FU)
_{\mathbb Q(\zeta_{p^\infty})}$
for a geometric point $\bar\xi$ of $U$.
\end{cor}

{\it Proof.}
By the assumption
${\rm char}\ K=0$,
the locally constant sheaf
${\cal F}$
on $U'$
and
the locally constant sheaves
$R^qf_!{\cal F}$
and
$R^qf_!\overline
{\mathbb F}_\ell$
on $U$
are 
tamely ramified on the generic fiber.

We define a commutative
diagram {\rm (\ref{eqtrcf})}
satisfying the condition
{\rm (\ref{sscfs}.0.4)}.
Since we assume $K$
is of characteristic $0$,
the sheaf ${\cal F}$ is
tamely ramified along $D$
by Abhyankar's lemma
\cite[Proposition 5.5]{SGA1}.
Hence, we may take a $G'$-torsor
$\pi\colon V'\to U'$
for a finite group $G'$
that is
tamely ramified along $D$
such that
the pull-back
$\pi^*{\cal F}$
is a constant sheaf
on $U'$.
Let $Y'$ be the normalization
of $X'$ in $U'$.
Then, since $V'\to U'$
is tamely ramified along $D$,
the proper curve
$Y'$ is smooth over $U$
and $V'\subset Y'$
is the complement of
a divisor $D'$
finite \'etale over $U$
by Lemma \ref{lmSn} below.

Since $Y'\to U$ is proper smooth,
its Stein factorization $D''\to U$ is
finite \'etale
\cite[Remarque (7.8.10)]{EGA3}.
Let $\pi'\colon V\to U$ be
a finite \'etale $G$-torsor
trivializing
the Stein factorization
of the finite \'etale coverings 
$D'$ and $D''$.
Replacing $V'$
by a connected component
of $U'\times_UV$
and $G'$ by the stabilizer
in $G'\times H$,
we obtain
a commutative diagram
(\ref{eqtrcf})
satisfying the condition
(\ref{sscfs}.0.4).

The proper smooth curve $Y'\to V$
has geometrically connected fibers
of genus $g'$
and $D'_V$ is the union of
$d'$ disjoint sections.
By the assumption
that $2g-2+d>0$,
we have
$2g'-2+d'>0$.
Hence,
with an ordering
of sections
$V\to D'$,
the pair $(Y',D')$
defines a $d'$ pointed smooth 
stable curve of genus $g'$.
Further replacing $V$ if necessary,
we may assume that
$\pi^{\prime*}R^qf_!{\mathbb F}_\ell$,
$\pi^{\prime*}R^qf_!{\cal F}$
and
$\pi^{\prime*}R^qf_!
{\mathbb Z}/n{\mathbb Z}$
are constant on $V$
for some integer $n\ge 3$
invertible on $S$.

Let $Y$ be a proper scheme
over $S$ containing $V$
as the complement
of a family ${\cal E}$
of Cartier divisors
such that
$\Sigma_{V/U}Y=
\Sigma^{\cal E}_{V/U}Y$.
Let $V\to 
{\cal M}_{g',d',n}$
be the morphism
to the moduli space
defined by the $d'$ pointed
smooth curve $Y'$ of
genus $g'$ over $V$.
By replacing $Y$
by the schematic closure
of the graph of the map
$V\to {\cal M}_{g',d',n}$ in 
$Y\times_{\mathbb Z}
\bar {\cal M}_{g',d',n}$,
we may assume that
$V\to 
{\cal M}_{g',d',n}$ is
extended to a morphism
$Y\to \bar
{\cal M}_{g',d',n}$.
Further by replacing $Y$
if necessary,
we may and do assume that
the action of $G$
on $Y$ is admissible
in the sense that the
quotient $Y/G$ exists as a scheme
and that ${\cal E}$
carries an action of $G$.
The pull-back of the universal family
by the map
$Y\to \bar
{\cal M}_{g',d',n}$
is a pointed stable
curve over $Y$
and satisfies
the condition
(\ref{prcfs}.2).

Thus, the assumptions
in Propositions
\ref{prcfs} and
\ref{prcfs0} are satisfied.
By Propositions
\ref{prcfs} and
\ref{prcfs0},
the assumptions
(\ref{eqtrcd}) and
(\ref{eqtrcd0})
in Proposition
\ref{prcfs2}
are satisfied respectively.
Thus the assertion follows.
\qed

\begin{lm}\label{lmSn}
Let $S$ be a normal scheme
and $X$ be a smooth curve
over $S$.
Let $D$ be a divisor of 
$X$ \'etale over $S$
and $U=X\setminus D$ 
be the complement.
Let $V\to U$ 
be a finite \'etale morphism
tamely ramified along $D$
and $Y$ be the normalization
of $X$ in $V$.
Then, $Y$ is smooth over
$S$ and
$V$ is the complement
of a divisor $E$ of $Y$
\'etale over $S$.
\end{lm}

{\it Proof.}
Let $\bar x$
be a geometric point
of $X$ and 
$t$ be a function on 
a neighborhood
defining $D$.
Let $\bar y$
be a geometric point
of $Y$ above $\bar x$.
Then, by Abhyankar's lemma
\cite[Proposition 5.5]{SGA1},
$Y$ is \'etale locally isomorphic
to 
$X[T]/(T^n-t)$
for an integer $n\ge 1$
invertible at $\bar x$
on a neighborhood of $\bar y$.
Hence the assertion follows.
\qed

\subsection{Swan class
of a constructible sheaf}
\label{ssSct}

In the rest of
this section,
we assume that the characteristic of $K$
is 0.
We define the Swan class for
a constructible sheaf
on a scheme over $K$.

For a separated scheme $U$ 
of finite type over $K$,
let 
$K(U,\overline{\mathbb F}_\ell)$
be the Grothendieck group of
constructible $\overline{\mathbb F}_\ell$-sheaves
on the \'etale site of $U$.
More precisely,
it is the quotient
of the free abelian group
generated by the isomorphism classes
$[{\cal F}]$ of
constructible 
$\overline{\mathbb F}_\ell$-sheaves
${\cal F}$
on the \'etale site of $U$
divided by the relations
$[{\cal F}]=
[{\cal F}']+[{\cal F}'']$
for exact sequences
$0\to {\cal F}'\to 
{\cal F}\to {\cal F}''\to 0$.

\begin{lm}\label{lmKU}
The abelian group
$K(U,\overline{\mathbb F}_\ell)$
is generated by the classes
$[i_!{\cal F}]$
where $i:Z\to U$
is a locally closed immersion of
smooth subscheme
and ${\cal F}$
is a locally constant constructible sheaf of
$\overline{\mathbb F}_\ell$-modules
on $Z$.
The relations are given by
\setcounter{equation}0
\begin{equation}
[i_!{\cal F}]=
[i_!{\cal F}']+
[i_!{\cal F}'']
\label{eqKUx}
\end{equation}
for exact sequences
$0\to {\cal F}'\to 
{\cal F}\to {\cal F}''\to 0$
of locally constant constructible 
$\overline{\mathbb F}_\ell$-modules
on $Z$ and 
\begin{equation}
[i_!{\cal F}]=
[i_{0!}{\cal F}|_{Z_0}]+
[i_{1!}{\cal F}|_{Z_1}]
\label{eqKUc}
\end{equation}
for smooth locally closed subschemes
$Z_1\subset Z$
where $i_0:Z_0=Z\setminus Z_1\to U$
and $i_1:Z_1\to U$
are the immersions.
\end{lm}

{\it Proof.}
We consider
the free abelian group
generated by the classes
$[i_!{\cal F}]$
where $i:Z\to U$
are locally closed immersions of
smooth subschemes
and ${\cal F}$
are locally constant constructible
$\overline{\mathbb F}_\ell$-modules
on $Z$.
Let $K'$ denote
its quotient
by the relations
(\ref{eqKUx}) and
(\ref{eqKUc}).
Clearly, we have
a canonical map $K'\to 
K(U,\overline
{\mathbb F}_\ell)$.
The inverse
is defined as follows.

For a constructible sheaf
${\cal F}$ on $U$,
there exists a finite partition
$U=\coprod_iU_i$ by 
smooth schemes
such that ${\cal F}|_{U_i}$
are locally constant.
It follows from
(\ref{eqKUc}) that
the sum
$\sum_i[{\cal F}|_{U_i}]$
is independent of
the partition.
Thus,
the class $[{\cal F}]
=\sum_i[{\cal F}|_{U_i}]
\in K'$
is well-defined.
Further,
the equalities
(\ref{eqKUx}) and
(\ref{eqKUc})
implies 
$[{\cal F}]=
[{\cal F}']+[{\cal F}'']$
for exact sequences
$0\to {\cal F}'\to 
{\cal F}\to {\cal F}''\to 0$.
Thus, the map
$K(U,\overline
{\mathbb F}_\ell)
\to K'$ is well-defined
and is the inverse of
the map
$K'\to K(U,\overline
{\mathbb F}_\ell)$ above.
\qed

\begin{pr}\label{prSwan}
For separated schemes $U$ 
of finite type over $K$,
there exists a unique way to attach
morphisms
\setcounter{equation}0
\begin{equation}
{\rm Sw}_U\colon
K(U,\overline{\mathbb F}_\ell)
\to
F_0G(\partial_F U)
_{\mathbb Q(\zeta_{p^\infty})}
\label{eqSwan}
\end{equation}
satisfying the following properties:

{\rm (1)} If $U$ is smooth over $K$
and if
${\cal F}$
is a locally constant constructible
sheaf of
$\overline{\mathbb F}_\ell$-modules
on $U$,
we have
${\rm Sw}_U([{\cal F}])=
{\rm Sw}_U{\cal F}$.

{\rm (2)} For an immersion
$i\colon U'\to U$,
the diagram
\begin{equation}
\begin{CD}
K(U,\overline{\mathbb F}_\ell)
@>{{\rm Sw}_U}>>
F_0G(\partial_F U)_{\mathbb Q(\zeta_{p^\infty})}\\
@A{i_!}AA @AA{i_!}A\\
K(U',\overline{\mathbb F}_\ell)
@>{{\rm Sw}_{U'}}>>
F_0G(\partial_F U')_{\mathbb Q(\zeta_{p^\infty})}
\end{CD}
\label{eqSwan2}
\end{equation}
is commutative.
\end{pr}

{\it Proof.}
By (1) and (2),
the map ${\rm Sw}_U$
is characterized
by ${\rm Sw}_U([i_!{\cal F}])
=i_!{\rm Sw}_{U'}{\cal F}$
for a locally constant
constructible sheaf
${\cal F}$ on
a regular subscheme $U'$ and
the immersion $i\colon
U'\to U$.
Hence the uniqueness follows
from Lemma \ref{lmKU}.
Further by Lemma \ref{lmKU}, 
the existence follows 
from Proposition \ref{prExc}.
\qed

We define a modification of
the map
${\rm Sw}_U$
with stronger compatibility
for push-forward.
For a separated scheme $U$
of finite type over $K$,
let ${\rm Const}(U)$ be the ${\mathbb Z}$-module
of constructible ${\mathbb Z}$-valued
functions on $U$.
For a 
constructible $\overline{\mathbb F}_\ell$-sheaf
${\cal F}$ on $U$,
let ${\rm rank}\ {\cal F}
\in {\rm Const}(U)$
be the constructible
function defined by
${\rm rank}\ {\cal F}(x)=
\dim {\cal F}_{\bar x}$.
Let ${\rm rank}\colon
K(U,\overline{\mathbb F}_\ell)
\to {\rm Const}(U)$ be the homomorphism
sending the class $[{\cal F}]$
to the function ${\rm rank}\ {\cal F}$.

Similarly as Proposition \ref{prSwan}, 
we have a map
${\rm Ch}_U\colon
{\rm Const}(U)
\to
F_0G(\partial_F U)_{\mathbb Q}$
characterized as follows.

\begin{pr}\label{lmC}
For separated schemes $U$ 
of finite type over $K$,
there exists a unique way to attach
morphisms
\setcounter{equation}0
\begin{equation}
{\rm Ch}_U\colon
{\rm Const}(U)
\to
F_0G(\partial_FU)_{\mathbb Q}
\label{eqC}
\end{equation}
satisfying the following properties:

{\rm (1)} If $U$ is smooth over $K$,
for
the constant function
$1_U$,
we have
${\rm Ch}_U(1_U)=
((\Delta_U,\Delta_U))^{\log}$.

{\rm (2)} For an immersion
$i\colon U'\to U$,
the diagram
\begin{equation}
\begin{CD}
{\rm Const}(U)
@>{{\rm Ch}_U}>>
F_0G(\partial_FU)_{\mathbb Q}\\
@A{i_!}AA @AA{i_!}A\\
{\rm Const}(U')
@>{{\rm Ch}_{U'}}>>
F_0G(\partial_FU')_{\mathbb Q}
\end{CD}\label{eqCh}
\end{equation}
is commutative.
\end{pr}

{\it Proof.}
It follows from
the excision formula
Theorem \ref{thmexc}.
\qed

\begin{df}\label{dfSwb}
Let $U$ be a separated scheme 
of finite type over $K$.
For a constructible 
$\overline{\mathbb F}_\ell$-sheaf 
${\cal F}$ on $U$,
we define the total Swan class
$\overline{\rm Sw}_U{\cal F}
\in F_0G(\partial_F U)_{
{\mathbb Q}(\zeta_{p^{\infty}})}$
by
\setcounter{equation}0
\begin{equation}
\overline{\rm Sw}_U{\cal F}=
{\rm Sw}_U{\cal F}-
{\rm Ch}_U({\rm rank}\ {\cal F}).
\label{eqSwanb}
\end{equation}
\end{df}

\begin{cor}\label{corSwb}
{\rm 1.}
Assume $U$
is smooth and ${\cal F}$ is
a locally constant constructible
$\overline{\mathbb F}_\ell$-sheaf
of constant rank
on $U$.
Let 
$f\colon V\to U$ be
a finite \'etale $G$-torsor
for a finite group $G$
such that
$\pi^*{\cal F}$ is a constant
sheaf.
Then, we have
\setcounter{equation}0
\begin{equation}
\overline {\rm Sw}_U{\cal F}=
-\frac 1{|G|}
\sum_{\sigma \in G}
{\rm Tr}^{\rm Br}(\sigma:M)
\cdot
f_*((\Gamma_{\sigma},\Delta_V))^{\log}.
\label{eqSwb}
\end{equation}

{\rm 2.}
For separated schemes $U$ 
of finite type over $K$,
the collection of the maps
$\overline{\rm Sw}_U\colon
K(U,\overline{\mathbb F}_\ell)
\to
F_0G(\partial_F U)
_{\mathbb Q(\zeta_{p^\infty})}$ 
is characterized by the
following properties:

{\rm (1)} 
Under the assumption in {\rm 1.},
we have
\begin{equation}
\overline {\rm Sw}_U([{\cal F}])=
-\frac 1{|G|}
\sum_{\sigma \in G}
{\rm Tr}^{\rm Br}(\sigma:M)
\cdot
f_*((\Gamma_{\sigma},\Delta_V))^{\log}.
\label{eqSwbb}
\end{equation}

{\rm (2)} For an immersion
$i\colon U'\to U$,
the diagram
\begin{equation}
\begin{CD}
K(U,\overline{\mathbb F}_\ell)
@>{\overline {\rm Sw}_U}>>
F_0G(\partial_F U)_{\mathbb Q(\zeta_{p^\infty})}\\
@A{i_!}AA @AA{i_!}A\\
K(U',\overline{\mathbb F}_\ell)
@>{\overline {\rm Sw}_{U'}}>>
F_0G(\partial_F U')_{\mathbb Q(\zeta_{p^\infty})}
\end{CD}\label{eqSwim}
\end{equation}
is commutative.
\end{cor}

{\it Proof.}
1. We have
\begin{equation}
\overline{\rm Sw}_U{\cal F}=
{\rm Sw}_U{\cal F}-
{\rm rank}\ {\cal F}\cdot
((\Delta_U,\Delta_U))^{\log}.
\label{eqSwbr}
\end{equation}
Thus, the equality
(\ref{eqSwb}) follows from
definition (\ref{eqSwn})
of the Swan class
and the equality
$|G|\cdot
((\Delta_U,\Delta_U))^{\log}=
\sum_{\sigma\in G}
f_*((\Gamma_\sigma,\Delta_V))^{\log}$.

2. We define the map 
$\overline{\rm Sw}_U$
by
$\overline{\rm Sw}_U
={\rm Sw}_U
-{\rm Ch}_U\circ {\rm rank}$.
Then, the commutative diagram
(\ref{eqSwim}) follows
from (\ref{eqSwan2})
and (\ref{eqCh}).
The uniqueness 
is clear from Lemma \ref{lmKU}.
\qed

\subsection{Conductor formula}
\label{sspc}

We keep the assumption
that $K$ is of characteristic $0$.
We show that
the diagram (\ref{eqSwim}) is commutative
for arbitrary morphisms
over $K$.
Changing the notation,
$f\colon U\to V$
denotes an arbitrary morphism
over $K$ of
separated schemes of
finite type over $K$.

\begin{thm}\label{thmcf}
Let $f\colon U\to V$ be a morphism
of separated schemes
of finite type over $K$.
Then, 
the diagram
\setcounter{equation}0
\begin{equation}
\begin{CD}
K(U,\overline{\mathbb F}_\ell)
@>{\overline {\rm Sw}_U}>>
F_0G(\partial_FU)_{\mathbb Q(\zeta_{p^\infty})}\\
@V{f_!}VV @VV{f_!}V\\
K(V,\overline{\mathbb F}_\ell)
@>{\overline {\rm Sw}_V}>>
F_0G(\partial_FV)_{\mathbb Q(\zeta_{p^\infty})}
\end{CD}\label{eqSwf}
\end{equation}
is commutative.
\end{thm}

{\it Proof.}
It suffices to show
the equality 
\begin{equation}
\overline{\rm Sw}_Vf_![{\cal F}]=
f_!\overline{\rm Sw}_U[{\cal F}]
\label{eqpc}
\end{equation}
for a constructible sheaf
${\cal F}$ on $U$.
We prove this by induction
on the dimensions of $U$ and $V$.
By a standard devissage
using the excision formula,
it suffices to show that
there exist dense open subschemes
$U'\subset U,V'\subset V$
such that $f(U')\subset V'$
and that we have
$$\overline{\rm Sw}_{V'}f|_{U'!}[{\cal F}|_{U'}]=
f_!\overline{\rm Sw}_{U'}[{\cal F}|_{U'}].$$
Hence, we may assume 
the following condition is satisfied.
\begin{itemize}
\item
The sheaf ${\cal F}$
is locally constant
and the scheme $V$ is smooth.
\end{itemize}
The formula
(\ref{eqpc})
is compatible
with the composition of morphisms.
Hence, by the induction on relative dimension,
we may further assume the following.
\begin{itemize}
\item[(\ref{thmcf}.3)]
The morphism
$f\colon U\to V$
is smooth of relative dimension $\le 1$.
\end{itemize}

Since we are allowed to shrink $V$,
we may assume that
$V$ is connected and that
there exists a proper smooth
curve $X$ over $V$
containing $U$
as the complement
of a divisor $D$ finite
\'etale over $V$.
Since the formula
(\ref{eqpc})
is proved for a finite \'etale
morphism
in Corollary \ref{corind},
by replacing $V$ 
by the Stein factorization
of $X\to V$,
we may assume
that the relative dimension is 1
and that
the geometric fibers of
$X\to V$ is connected.
Further shrinking $U$ and $V$,
we may replace the
condition (\ref{thmcf}.3)
by the following.
\begin{itemize}
\item[(\ref{thmcf}.4)]
There exist 
a proper smooth and geometrically connected curve
$\bar f\colon X\to V$
of genus $g$
and an open immersion $U\to X$
such that $U$ is the complement
of a divisor $D\subset X$
finite and \'etale over $V$
of degree $d$ such that
$2g-2+d>0$.
\end{itemize}
Then, applying Corollary \ref{corcfs},
we obtain the equalities
(\ref{eqpcsc}) and
(\ref{eqpcsc1}).
The equality (\ref{eqpc}) follows
from them
together with (\ref{eqSwbr}).
\qed

We derive some consequences
of Theorem \ref{thmcf}.

\begin{cor}\label{corUVsm}
Let $f\colon U\to V$
be a smooth morphism
of smooth separated schemes
of finite type over $K$.
Assume that
$R^qf_!{\mathbb F}_\ell$
is locally constant for
every $q\ge 0$.

{\rm 1.}
Let ${\cal F}$ be
a constructible
sheaf of
$\overline{\mathbb F}_\ell$-modules
of constant rank
on $U$.
Assume that
$R^qf_!{\cal F}$ is 
locally constant for
every $q\ge 0$.
Then, we have
\setcounter{equation}0
\begin{equation}
{\rm Sw}_VRf_!{\cal F}=
f_!
{\rm Sw}_U{\cal F}
+
{\rm rank}\ {\cal F}\cdot
{\rm Sw}_VRf_!
\overline {\mathbb F}_\ell
\label{eqmain1}
\end{equation}
in $F_0(\partial_FV)_
{{\mathbb Q}(\zeta_{p^\infty})}$.

{\rm 2.}
We have
\begin{equation}
{\rm Sw}_VRf_!
\overline {\mathbb F}_\ell
=
{\rm rank}\ Rf_!
\overline {\mathbb F}_\ell
\cdot
((\Delta_V,\Delta_V))^{\log}
-
f_!
((\Delta_U,\Delta_U))^{\log}
\label{eqmain2}
\end{equation}
in $F_0(\partial_FV)_
{{\mathbb Q}(\zeta_{p^\infty})}$.
\end{cor}

{\it Proof.}
1.
Since ${\rm rank}\ Rf_!{\cal F}
=
{\rm rank}\ {\cal F}\cdot
{\rm rank}\ Rf_!
\overline {\mathbb F}_\ell$,
it suffices to apply
Theorem \ref{thmcf}
to $[{\cal F}]-{\rm rank}\ {\cal F}
\cdot [{\mathbb F}_\ell]$.

2.
It suffices to apply
Theorem \ref{thmcf}
to $[{\mathbb F}_\ell]$.
\qed

If there exist
a proper smooth 
scheme
$\bar f\colon X\to V$
and a divisor $D$ of $X$
with normal crossings relatively to $V$
such that
$U$ is the complement
$X\setminus D$,
the assumption of Corollary
\ref{corUVsm} is satisfied.
Further, if $d$ denotes
the relative dimension 
of $X$ over $V$
we have
$${\rm rank}\ Rf_!
\overline {\mathbb F}_\ell
=
(-1)^d{\rm deg}\ c_d(\Omega^1_{X/V}(\log D)).$$

In particular,
for $V={\rm Spec}\ K$,
we have the following.

\begin{cor}\label{corUsm}
Let $U$ be a smooth
separated scheme of finite type
over $K$
and ${\cal F}$
be a smooth $\overline {\mathbb F}_\ell$-sheaf
of constant rank on $U$.
Then, we have
\setcounter{equation}0
\begin{eqnarray}
{\rm Sw}_KR\Gamma_c(U_{\overline K},{\cal F})
&=&
{\rm deg}\ {\rm Sw}_U{\cal F}+
{\rm rank}\ {\cal F}\cdot
{\rm Sw}_KR\Gamma_c(U_{\overline K},
\overline {\mathbb F}_\ell),
\label{eqcfK}\\
{\rm Sw}_KR\Gamma_c(U_{\overline K},
\overline {\mathbb F}_\ell)
&=&
-
{\rm deg}((\Delta_U,\Delta_U))^{\log}.
\label{eqcfKc}
\end{eqnarray}
\end{cor}

The equality
(\ref{eqcfKc}) implies
the conductor formula of Bloch
in the case proved in
\cite{KSI} as follows.
We assume that
$U$ is proper smooth over $K$
and 
$X$ is a proper regular flat scheme
over $S={\rm Spec}\ {\cal O}_K$
such that $U=X_K$
and the reduced closed fiber
$X_{F, {\rm red}}$
is a divisor with simple normal crossings.
Then, by Proposition \ref{prbfS}
and 
\cite[(5.4.2.6)]{KSI},
we have
$((\Delta_U,\Delta_U))^{\log}
_{
(U\times_SU)^\sim}=
((\Delta_U,\Delta_U))^{\log}
_{
(U\times_{\mathbb S}U)^\sim}
=
(-1)^dc_d(
\Omega^1_{X/S}(\log/\log))_{X_F}.$
Thus, in this case,
the equality
(\ref{eqcfKc})
is equivalent to
\cite[Theorem 6.2.5]{KSI}
and hence
to the conductor formula
of Bloch \cite{bloch}.
This proof of
the conductor formula
of Bloch uses the same tools
including the
localized intersection product.
However,
the excision formula allows us
to reduce the proof to
relative curves.

\newpage 
\section{A computation
in the rank 1 case}\label{srk1}

We state Conjecture \ref{cncF}
comparing
the Swan class
of a sheaf 
of rank 1
with the cycle class
defined in 
\cite[Section 5.1]{Kato}
and prove it in Theorem \ref{thmrk1}
assuming $\dim U_K\le 1$.
Using it,
we prove the integrality conjecture
Conjecture \ref{cnHA}
under the
assumption $\dim U_K\le 1$.
In Section \ref{sskum} and 
in the second half of
Section \ref{ssrk1},
we will assume that $K$
is of characteristic $0$.

\subsection{Ramification
of characters}\label{ssch}
\setcounter{equation}0

We briefly recall the
theory of ramification
of characters of Galois groups
in \cite{deg1}.
For a field $K$, 
let $X_K$ denote the 
dual group
$H^1(K,{\mathbb Q}/{\mathbb Z})
=H^2(K,{\mathbb Z})$
of the abelian quotient
$G_K^{\rm ab}$ of
the absolute Galois group
$G_K={\rm Gal}(\bar K/K)$.
The cup-product
defines a canonical pairing
\addtocounter{thm}1
\begin{equation}
(\ ,\ )_K\colon
X_K\times K^\times
=
H^2(K,{\mathbb Z})
\times
H^0(K,{\mathbf G}_m)
\to Br(K)=
H^2(K,{\mathbf G}_m).
\label{eqcp}
\end{equation}

Assume $K$ is a henselian discrete
valuation field
and let $F$ be the residue field
of characteristic $p>0$.
In this subsection,
we drop the assumption
that $F$ is perfect.
We briefly recall the definition
of the exact sequence
$$0\to 
\Omega_F\to
\Omega_F(\log)
\overset{\rm res}\to F\to 0$$
of $F$-vector spaces.
A canonical map $d\log \colon
{\cal O}_K^\times\to 
\Omega_F=\Omega_{F/F^p}$ 
is defined by $a\mapsto 
\bar a^{-1}d\bar a$.
The $F$-vector space $\Omega_F(\log)$ is defined as
the amalgamate sum of
$\Omega_F$
with $F\otimes_{\mathbb Z}K^\times$
over $F\otimes_{\mathbb Z}{\cal O}_K^\times$
with respect to the map
$d\log \colon
{\cal O}_K^\times\to 
\Omega_F$ 
and the inclusion
${\cal O}_K^\times\to K^\times$.
The valuation
$K^\times\to {\mathbf Z}$
induces the residue map
${\rm res}\colon
\Omega_F(\log)\to F$.
The map $d\log \colon
{\cal O}_K^\times\to 
\Omega_F$ is canonically extended to
$d\log \colon
K^\times\to 
\Omega_F(\log)$.

We identify
a character $\chi
\in X_F$ 
with the corresponding
unramified character
$\chi\in X_K$
and regard
$X_F$ as a subgroup
of $X_K$.
For $a\in F$,
let $\chi_a
\in X_F$ 
be the character
defined by
the Artin-Schreier
equation $T^p-T=a$.
We define a map
$\chi\colon F\to X_K$
by sending $a\in F$
to $\chi_a\in X_F\subset
X_K$.
In \cite[(1.4)]{deg1},
it is shown that
there exists a unique map
$\lambda_K\colon
\Omega_F(\log)
\to 
Br(K)$ that makes the diagram
\begin{equation}
\begin{CD}
F\times K^\times
@>{(a,b)\mapsto a\cdot d\log b}>>
\Omega_F(\log)\\
@V{\chi\times 1}VV
@VV{\lambda_K}V\\
X_K\times K^\times
@>{(\ ,\ )_K}>>
Br(K)
\end{CD}
\label{eqOB}
\end{equation}
commutative.

The main construction in 
\cite[Definition (2.1)]{deg1}
is the increasing filtration
$F_\bullet$
of $X_K$ indexed by
$r\in {\mathbb N}$.
We have
$X_K=\bigcup_{r\ge0}
F_rX_K$.
The subgroup $F_0X_K$
consists of
the characters
at most tamely ramified.
For $r\ge 1$,
we put $U^r_K=
1+{\mathfrak m}^r_K
\subset K^\times$.
Then, the pairing 
$(\ ,\ )_K\colon
X_K\times K^\times$ maps
$F_rX_K\times U^r_K$
for $r\ge 1$
and
$F_0X_K\times K^\times$
to ${\rm Im}\ \lambda_K
\subset Br(K)$.
For an extension $L$ of
henselian discrete valuation field 
such that ${\cal O}_K=K\cap
{\cal O}_L$ and
${\mathfrak m}_K
{\cal O}_L=
{\mathfrak m}^e_L$,
the canonical map
$X_K\to X_L$
sends 
$F_rX_K$ to $F_{re}X_L$.

For $r\ge 1$,
we put
${\rm Gr}^F_rX_K=
F_rX_K/F_{r-1}X_K$.
A canonical injection
\begin{equation}
{\rm rsw}_{r,K}\colon
Gr^F_rX_K
\to {\rm Hom}_F(
{\mathfrak m}_K^r/
{\mathfrak m}_K^{r+1},
\Omega_F(\log))
\label{eqrswi}
\end{equation}
is defined
in \cite[Corollary (5.2)]{deg1}.
It is characterized
by the following
properties:
\begin{itemize}
\item[(1)]
For
$\chi\in F_rX_K$ and 
$c\in {\mathfrak m}^r_K$,
we have
\begin{equation}
(\chi,1-c)_K=\lambda_K(
{\rm rsw}_{r,K}(\chi)(\bar c)).
\label{eqrswc}
\end{equation}
\item[(2)]
Let $L$ be an arbitrary
extension of
henselian discrete valuation field
such that ${\cal O}_K=K\cap
{\cal O}_L$ and
${\mathfrak m}_K
{\cal O}_L=
{\mathfrak m}^e_L$.
Let $F_L$ denote
the residue field of $L$.
Then, the diagram
\begin{equation}
\begin{CD}
Gr^F_rX_K
@>{{\rm rsw}_{r,K}}>>
Hom_F(
{\mathfrak m}^r_K/
{\mathfrak m}^{r+1}_K,
\Omega_F(\log))
\\
@VVV @VVV\\
Gr^F_{er}X_L
@>{{\rm rsw}_{er,L}}>>
Hom_{F_L}
({\mathfrak m}^{er}_L/
{\mathfrak m}^{er+1}_L,
\Omega_{F_L}(\log))
\end{CD}
\label{eqrswL}
\end{equation}
is commutative.
\end{itemize}
For an element
$\chi\in F_rX_K
\setminus F_{r-1}X_K$,
the injection
\begin{equation}
{\rm rsw}_{r,K}(\chi)\colon
{\mathfrak m}^r_K/
{\mathfrak m}^{r+1}_K
\longrightarrow
\Omega_F(\log)
\label{eqrsw}
\end{equation}
is called the 
refined Swan conductor
of $\chi$ and will be
denoted by
${\rm rsw}\ \chi$.

We compute the refined Swan conductor
of a Kummer character of degree $p$
explicitly.
Assume that
$K$ is of characteristic $0$
and the residue field $F$
is of characteristic $p$.
Assume further that $K$
contains a primitive $p$-th root
$\zeta_p$ of 1.
We identify ${\mathbb Z}/p{\mathbb Z}
=\mu_p$ by $\zeta_p$
and the $p$-torsion
part $X_K[p]=H^1(K,{\mathbb Z}/p{\mathbb Z})$
with $K^\times/K^{\times p}=
H^1(K,\mu_p)$
by the isomorphism
$\theta\colon
K^\times/K^{\times p}
\to X_K[p]$
of Kummer theory.
For $a\in K^\times$,
let $\theta_a\in X_K[p]$
denote the corresponding
character.

We put $z=\zeta_p-1$.
Then, we have
$z^p+pz\equiv (z+1)^p-1=0
\bmod pz^2$
and ${\rm ord}\ z^p
={\rm ord}\ pz
>{\rm ord}\ pz^2$.
Hence,
for an element $a\in {\cal O}_K$,
the reduction of the Kummer equation
$(1-zt)^p=1-az^p$
gives the Artin-Schreier equation
$t^p-t=\bar a$
and the unramified character
$\chi_{\bar a}\in X_K[p]$
is identified with
$1-az^p\in 
K^\times/K^{\times p}$.
In particular,
we have
$1+z^p{\mathfrak m}_K
\subset K^{\times p}$.
Consequently,
we have a commutative diagram
$$\begin{CD}
K^\times/K^{\times p}
@>{\theta}>> X_K[p]\\
@A{\bar a\mapsto 1-az^p}AA
@AAA\\
F@>{\chi}>>X_F[p].
\end{CD}$$

\begin{pr}\label{prkum}
Let $K$ be a henselian
discrete valuation field
of mixed characteristic $(0,p)$
containing a primitive 
$p$-th root $\zeta_p$
of $1$.
We put $e'=
p\cdot {\rm ord}_K(\zeta_p-1)
={\rm ord}_Kz^p$.
We define a decreasing filtration
$F^\bullet$ on $K^\times/
K^{\times p}$
by $F^m(K^\times/
K^{\times p})={\rm Image}\ U_K^m$
for $m\ge 1$ and
by $F^0(K^\times/
K^{\times p})=K^\times/
K^{\times p}$.

{\rm 1. (\cite[Proposition 4.1]{deg1})}
The isomorphism
$\theta \colon
K^\times/
K^{\times p}\to
X_K[p]$
induces an isomorphism
\setcounter{equation}0
\begin{equation}
F^m(K^\times/
K^{\times p})\to 
F_rX_K[p]
\label{eqkum}
\end{equation}
for $0\le m=e'-r\le e'$.
In particular,
we have
$F_{e'}X_K[p]=X_K[p]$.

{\rm 2.}
For $a\in K^\times$
such that $d\log a\neq 0$
in $\Omega_F(\log)$,
the map
\begin{equation}
\begin{CD}
{\rm rsw}_{e',K}(\theta_a)
\colon
{\mathfrak m}^{e'}_K
/{\mathfrak m}^{e'+1}_K
\to
\Omega_F(\log)
\end{CD}
\label{eqrswk1}
\end{equation}
sends $c\cdot z^p$ to
$-c\cdot d\log a$
for $c\in {\cal O}_K$.

{\rm 3.}
For $1\le m=e'-r<e'$
and 
$a=1-b\in U^m_K,
\notin U^{m+1}_K$,
the map
\begin{equation}
\begin{CD}
{\rm rsw}_{r,K}(\theta_a)
\colon
{\mathfrak m}^r_K
/{\mathfrak m}^{r+1}_K
\to
\Omega_F(\log)
\end{CD}
\label{eqrswk}
\end{equation}
sends $c\cdot z^p/b$ to
$c\cdot d\log b$
for $c\in {\cal O}_K$.
\end{pr}

{\it Proof.}
We identify the $p$-torsion
part $Br(K)[p]=
H^2(K,\mu_p)$
with $H^2(K,\mu_p^{\otimes 2})$
by $\zeta_p$.
Then, for $a,b\in K^\times$,
the cup-product
$(\theta_a,b)_K
\in Br(K)$
is identified
with the Galois symbol
$\{a,b\}$
defined as
$\theta_a\cup \theta_b
\in H^2(K,\mu_p^{\otimes 2})$.
Let $a\in K^\times$
and let $L$ be an arbitrary
extension of henselian discrete
valuation field.
Then, for $c\in {\cal O}_L$,
we have
\begin{equation}
(\theta_a,1-z^pc)_L=
\{a,1-z^pc\}=
-\{1-z^pc,a\}
=-(\chi_{\bar c},a)_L
=-
\lambda_L(
\bar c\cdot d\log a).
\label{eqe'}
\end{equation}
By the characterization
of ${\rm rsw}_{r,K}$,
the equality
(\ref{eqe'}) implies that
the map
${\rm rsw}_{r,K}$
is the zero-map
for $r>e'$.
Hence, by
the injectivity of
${\rm rsw}_{r,K}$,
we have
${\rm Gr}^F_rX_K[p]=0$
for $r>e'$.
Thus by
$X_K=\bigcup_{r\ge0}
F_rX_K$,
we obtain 
$F_{e'}X_K[p]=X_K[p]$.
Now, the equality 
(\ref{eqe'})
implies the assertion 2.

To show the remaining assertions,
we use the following
elementary lemma
on the symbol map.

\begin{lm}[{\cite[Lemma 6]{residue}}]\label{lm1x}
For $b,c\in K^\times\setminus \{1\}$,
we have
\setcounter{equation}0
\begin{equation}
\{1-b,1-c\}
=\{1-bc,-b\}+
\{1-bc,1-c\}-
\{1-bc,1-b\}.
\label{eq1x}
\end{equation}
\end{lm}
{\it Proof of Lemma.}
Since $\{x,y\}=0$
for $x,y\in K^\times$
satisfying $x+y=1$,
we have
$$\left\{
1-b,1-\frac{1-bc}{1-b}\right\}=
\left\{
1-bc,1-\frac{1-bc}{1-b}\right\}.$$
Since
$\displaystyle{
1-\frac{1-bc}{1-b}
=\frac{-b(1-c)}{1-b}}$,
the right hand side is equal to that
of (\ref{eq1x}).
Further, since
$\displaystyle{
\frac1{1-b}
+\frac{-b}{1-b}=1}$,
the left hand side is equal to that
of (\ref{eq1x}).
\qed

We go back to the proof of Proposition
\ref{prkum}.
Let $1\le m=e'-r< e'$ be
an integer,
$b\in {\mathfrak m}_K^m$
and put $a=1-b$.
Let $L$ be an 
arbitrary extension
of discrete valuation
field
and $c\in {\mathfrak m}_K^{-m}
{\cal O}_L$ 
be an arbitrary
element.
Since
$U_L^{e'_L+1}
\subset L^{\times p}$
for $e'_L={\rm ord}_Lz^p$,
we have
$\{U_L^{e'_L},U_L^1\}=0$
by Lemma \ref{lm1x}.
Hence,
by $b,z^pc\in {\mathfrak m}_L$,
we have 
\begin{align}
(\theta_a,1-z^pc)_L&=
\{1-b,1-z^pc\}=
\{1-bcz^p,-b\}
\label{eqm}
\\
&=(\chi_{\overline{bc}},
-b)_L
=
\lambda_L(\overline{bc}\cdot
d\log (-b))
=
\lambda_L(
\overline{bc}\cdot
d\log b)
\nonumber
\end{align}
further 
by Lemma \ref{lm1x}.

Similarly as above,
the equality (\ref{eqm})
together with
the characterization
and the injectivity
of ${\rm rsw}_{n,K}$
shows that
$\theta$ maps
$U^m_K$
to
$F_rX_K[p]$,
by induction on $1\le 
m=e'-r\le e'$.
Further
the equality (\ref{eqm})
implies the assertion 3.
Hence, the composition
of the map 
\begin{align*}
{\rm Gr}_F^m(
K^\times/
K^{\times p})
&=
\begin{cases}
{\mathfrak m}_K^m/
{\mathfrak m}_K^{m+1}
&
\text{ if }p\nmid m\\
{\mathfrak m}_K^m/
({\mathfrak m}_K^n)^p
{\mathfrak m}_K^{m+1}
&
\text{ if }m=np\\
\end{cases}
\\
&
\begin{CD}
@>{{\rm Gr}\theta}>>
{\rm Gr}^F_r
X_K
@>{{\rm rsw}_{r,K}}>>
Hom_F({\mathfrak m}_K^r/
{\mathfrak m}_K^{r+1},
\Omega_F(\log))\end{CD}
\end{align*}
sends
$1-b$ to the map
$c\mapsto bc/z^p\cdot d\log b$
and is
an injection.
Since
$\theta\colon
K^\times/
K^{\times p}\to
X_K[p]$ is an isomorphism, 
we conclude
$\theta(U^m_K)
=F_rX_K[p]$.
\qed

\subsection{Kummer covering
of degree $p$}\label{sskum}
\setcounter{equation}0

We apply the
theory recalled in
the previous section
to the
following geometric situation.
Let $K$ be
a complete discrete valuation
field of characteristic
$0$
with perfect residue field
$F$ of characteristic $p>0$.
Let $X$ be a regular
flat separated scheme
of finite type
over $S={\rm Spec}\
{\cal O}_K$
and $D$ be a divisor with
simple normal crossings.
Let $D_1,\ldots,D_n$
be the irreducible components
of $D$. For an irreducible
component $D_i$,
let $K_i$ be the fraction
field of the completion
$\widehat {\cal O}_{X,\xi_i}$
of the local ring at
the generic point $\xi_i$ of $D_i$.
The residue field $F_i=\kappa(\xi_i)$
of the complete discrete
valuation field $K_i$
is the function field of $D_i$.
The fiber
$\Omega^1_{X/S}(\log D)_{\xi_i}
\otimes_{{\cal O}_{X,\xi_i}}
F_i$
is identified
with the $F_i$-vector
space
$\Omega_{F_i}(\log)$
in the notation 
of the last subsection.

Let $\chi\in
H^1(U,{\mathbb Q}/{\mathbb Z})$
be a character.
Then, for each local field $K_i$,
the restriction
defines a character
$\chi_i\in X_{K_i}$. 
By the ramification
theory recalled in Section \ref{ssch},
the Swan conductor
$r_i=
{\rm sw}_{K_i}(\chi_i)\ge 0$
is defined for each $K_i$.
We define 
the Swan divisor of $\chi$ by 
$D_\chi=
\sum_{i=1}^mr_iD_i$.
Let 
$E=
\sum_{i,r_i>0}D_i$
be the support of $D_\chi$.
For each irreducible component $D_i$
such that $r_i>0$,
the refined Swan conductor
${\rm rsw}_{K_i}(\chi_i)$
defines a non-zero map
$${\cal O}_X(-D_{\chi})_{\xi_i}
\otimes_{{\cal O}_{X,\xi_i}}
F_i
={\mathfrak m}_{K_i}^{r_i}/
{\mathfrak m}_{K_i}^{r_i+1}
\to 
\Omega^1_{X/S}(\log D)_{\xi_i}
\otimes_{{\cal O}_{X,\xi_i}}
F_i
=\Omega_{F_i}(\log).$$
In \cite[Theorem (7.1),
Proposition (7.3)]{deg1},
it is shown that
there exists an ${\cal O}_E$-linear injection
\begin{equation}
{\rm rsw}\ \chi\colon
{\cal O}_X(-D_{\chi})
\otimes_{{\cal O}_X}
{\cal O}_E
\to 
\Omega^1_{X/S}(\log D)
\otimes_{{\cal O}_X}
{\cal O}_E
\label{eqrswE}
\end{equation}
inducing
${\rm rsw}_{K_i}(\chi_i)$
at each generic point.

\begin{df}\label{dfcln}
We say that $\chi$
is {\rm clean} with respect to $X$
if the map
${\rm rsw}\ \chi\colon
{\cal O}_X(-D_{\chi})
\otimes_{{\cal O}_X}
{\cal O}_E
\to 
\Omega^1_{X/S}(\log D)
\otimes_{{\cal O}_X}
{\cal O}_E$
is a locally splitting injection.

Assume $\chi$
is clean with respect to $X$.
Then, we say that $\chi$
is s-{\rm clean} with respect to $X$
if,
for each irreducible
component $D_i$ of $E$,
the composition
$$\begin{CD}
{\cal O}_X(-D_{\chi})
\otimes_{{\cal O}_X}
{\cal O}_{D_i}
@>{{\rm rsw}_i(\chi)}>>
\Omega^1_{X/S}(\log D)
\otimes_{{\cal O}_X}
{\cal O}_{D_i}
@>{{\rm res}}>>
{\cal O}_{D_i}
\end{CD}$$
is either an isomorphism
or the zero-map,
depending
on $D_i$.
\end{df}

It is conjectured in \cite{Kato}
that
for any $\chi$,
there exists a proper
modification
$X'$ of $X$
such that
$\chi$ is clean
with respect to $X'$,
see Lemma \ref{lmscl2}.3.

We compute the Swan divisor
$D_\chi$
and the map
${\rm rsw}\ \chi$
(\ref{eqrswE})
for a Kummer character $\chi$
of order $p$
explicitly.

\begin{lm}\label{lmFma}
Let $A$ be a regular
local ring
such that
the fraction field
is of characteristic $0$
and
the residue field
is of characteristic 
$p>0$.
Let $t_1,\ldots,t_n
\in A$
be a part of regular
system of parameters.
Assume that
$A$ contains a primitive
$p$-th root $\zeta_p$ of
$1$ and that
$t_1,\ldots,t_n$
divide $p$.
Let $K_i$
be the fraction
field of the completion
${\cal O}_{K_i}$
of the discrete
valuation ring
$A_{(t_i)}$
for each $i$.

Let $m=(m_1,\ldots,m_n)$
be a family of integers
satisfying
$0\le m_i\le e'_i=
p\cdot {\rm ord}_{K_i}(\zeta_p-1)$
and let $F^mA^\times$
denote the subgroup
$1+t_1^{m_1}
\cdots t_n^{m_n}A$
for $m\neq 0$
and $F^0A^\times
=A^\times$.
Then, the 
inverse image 
of 
$\bigoplus_i
F^{m_i}(K_i^\times
/K_i^{\times p})
\subset
\bigoplus_i
K_i^\times
/K_i^{\times p}$
by the canonical map
$A^\times/
A^{\times p}
\to
\bigoplus_i
K_i^\times
/K_i^{\times p}$
is the image of
$F^mA^\times$.
\end{lm}

{\it Proof.}
First, we show the
case where $n=1$.
We prove it by
induction on $m=m_1$.
It is obvious for $m=0$.
Assume $m=1$
and that the image
of $a\in A^\times$
in $K_1^\times/
K_1^{\times p}$
is in $F^1(K_1^\times/
K_1^{\times p}).$
Let $F_1$ denote
the residue field of $K_1$.
Then, we have
$\bar a \in F_1^{\times p}$.
We put
$\bar a=b^p$
for $b\in F_1^\times$.
Since $b$ is integral
over the normal ring
$A/t_1A$,
we have $b\in A/t_1A$
and $b\in (A/t_1A)^\times$.
Take a unit
$c\in A^\times$ lifting $b$.
Then, 
$a/c^p$
is in $F^1A^\times
=1+t_1A$.

Assume $m\ge 1$
and that
the image of
$a\in F^mA^\times
=1+t_1^mA$
is in 
$F^{m+1}(K_1^\times/
K_1^{\times p}).$
If $p\nmid m$,
we have $a\in 1+
t_1^{m+1}
{\cal O}_{K_1}$.
Since
$A\cap
t_1^{m+1}
{\cal O}_{K_1}
=
t_1^{m+1}A$,
we have
$a\in F^{m+1}A^\times$.
Assume $p|m$
and we put
$a=1+t_1^mb$.
Then, we have
$\bar b\in F_1^p$.
Similarly as above,
there is an element
$c\in A$
such that
$b\equiv c^p\bmod t_1$.
Then,
$a/(1+t^{m/p}c)^p$
is in $F^{m+1}A^\times
=1+t_1^{m+1}A$.

We prove the general case.
Assume that the image of
$a\in A^\times$
is in 
$\bigoplus_i
F^{m_i}(K_i^\times/
K_i^{\times p}).$
Then, for each $i$,
there exists $a_i\in A^\times$
such that
$a/a_i^p\in F^{m_i}A^{\times}$.
Let $m'_i\ge 0$ be
the smallest integer
satisfying $p\cdot
m'_i\ge m_i$.
Then, since
the $p$-th power map
$(A/t_i^{m'_i})^\times
\to
(A/t_i^{m_i})^\times$
is injective,
the class
$\bar a_i
\in (A/t_i^{m'_i})^\times$
is uniquely determined
by the condition
$a/a_i^p\in F^{m_i}A^{\times}$.
Further, for $i\neq j$,
the $p$-th power map
$(A/(t_i^{m'_i},t_j^{m'_j}))^\times
\to
(A/(t_i^{m_i},t_j^{m_j}))^\times$
is also injective.
Hence, there exists
a unique element
$b\in 
(A/t_1^{m'_1}\cdots
t_n^{m'_n}A)^\times$
satisfying
$b\equiv a_i
\bmod t_i^{m'_i}$.
Let $c\in A^\times$
be a unit lifting $b$.
Then,
we have
$a/c^p$
is in $F^mA^\times$.
\qed

\begin{cor}\label{corkuma}
Let $\chi\in
H^1(U,{\mathbb Z}/p{\mathbb Z})$
be a character of
order $p$
and $x\in D$
be a point.
Let $D_1,\ldots,D_n$
be the irreducible components
of $D$ containing $x$
and put
$D_\chi=
\sum_ir_iD_i$
on a neighborhood of $x$.
We put $A={\cal O}_{X,x}$
and, for each irreducible
component $D_i$,
we put
$e'_i=
p\cdot {\rm ord}_{D_i}(\zeta_p-1)$
and $m_i=e'_i-r_i$.

{\rm 1.}
On a neighborhood of $x$,
there exists
an element $a\in 
\Gamma(U,{\cal O}_U^\times)$
such that
$\chi$ is defined by
$t^p=a$ and satisfying
one of
the following conditions:
\begin{itemize}
\item[{\rm (\ref{corkuma}.1)}]
${\rm ord}_{D_i}a$
is prime to $p$
for at least one
$D_i$.
\item[{\rm (\ref{corkuma}.2)}]
$a$ is a unit at $x$
and its image in $A^\times$
is in $F^mA^\times$
for $m=(m_1,\ldots,m_n)$
in the notation of Lemma
{\rm \ref{lmFma}}.
\end{itemize}

{\rm 2.}
Assume
$D_\chi
=
\sum_ie'_iD_i$
and let $a$ be as in {\rm 1.}
Then, the map 
$${\rm rsw}\ \chi
\colon
{\cal O}_E(-D_\chi)
\to
\Omega^1_{X/S}(\log D)
\otimes
{\cal O}_E$$
sends $z^p$ to $-d\log a$
where $z=\zeta_p-1$.

{\rm 3.}
Assume
$D_\chi
<
\sum_ie'_iD_i$.
Then, the 
condition
{\rm (\ref{corkuma}.2)}
holds.
Let $b$
be a basis
of the invertible sheaf
${\cal O}_X(-\sum_im_iD_i)$
on a neighborhood of $x$.
We put 
$a=1-bc\in {\cal O}_X^\times$
where $c\in {\cal O}_X$ 
on a neighborhood of $x$
as in {\rm (\ref{corkuma}.2)}.
Then,
the map 
$${\rm rsw}\ \chi
\colon
{\cal O}_E(-D_\chi)
\to
\Omega^1_{X/S}(\log D)
\otimes
{\cal O}_E$$
sends $a\cdot z^p/b$ 
to $c\cdot d\log b+dc$.
\end{cor}

{\it Proof.}
1.
Let $a$ be a rational function
on $X$
such that
$\chi$ is defined by
$t^p=a$ on the generic point.
The regular local ring
$A={\cal O}_{X,x}$ 
is a UFD.
Hence, 
we may assume 
$0\le {\rm ord}_ya<p$
for every discrete valuation
defined by a point $y\in X$
of codimension 1,
after dividing $a$ by 
the $p$-th power of a rational function
and shrinking $X$
if necessary.
For a point $y\in X$
of codimension 1,
if the valuation of $a$
at $y$ is not divisible by $p$,
then $\chi$ is ramified at $y$.
Hence, 
$a$ is a unit on $U$.
Further,
if the condition
(\ref{corkuma}.1)
is not satisfied,
then $a$
is a unit at $x$.
Then, by
Lemma \ref{lmFma},
after dividing $a$ by 
the $p$-th power of a rational function,
the condition
(\ref{corkuma}.2)
is satisfied.

2. 3.
Clear from Proposition \ref{prkum}
and the equality
$-da/b=
c\cdot d\log b+dc$.
\qed

We give a condition
for character $\chi$
of order $p$ to be
clean.

\begin{pr}\label{prcln}
Let $\chi\in
H^1(U,{\mathbb Z}/p{\mathbb Z})$
be a character of
order $p$
and $x\in D$
be a point.
Assume that
$\chi$ is not tamely ramified
at $x$
and that $K$ is of characteristic $0$ 
and contains
a primitive $p$-th root $\zeta_p$ of $1$.
Let $D_1,\ldots,D_n$
be the irreducible components
of $D$ containing $x$
and let $C$
denote
the intersection 
$D_1\cap \cdots\cap D_n$.
We put $D_\chi=
\sum_ir_iD_i>0$
as in Corollary {\rm \ref{corkuma}}.

For each irreducible 
component $D_i$,
we put
$e'_i=p\cdot
{\rm ord}_{D_i}(\zeta_p-1)$
and $m_i=e'_i-r_i$.
Let $t_i\in \Gamma(X,{\cal O}_X)$
be an element defining $D_i$.

{\rm 1.}
Assume $D_\chi=
\sum_ie'_iD_i$.
Then,
$\chi$ is clean
at $x$
if and only if
on a neighborhood of $x$,
there exists
an element $a\in 
\Gamma(U,{\cal O}_U^\times)$
such that
$\chi$ is defined by
$t^p=a$ and satisfying
one of
the following conditions:
\begin{itemize}
\item[{\rm (\ref{prcln}.1)}]
${\rm ord}_{D_i}a$
is prime to $p$
for at least one
$D_i$.
\item[{\rm (\ref{prcln}.2)}]
$a$ is a unit at $x$
and $da|_C$ has
no zero at $x$.
\end{itemize}

{\rm 2.}
Assume $D_\chi<
\sum_ie'_iD_i$.
Then,
$\chi$ is clean
at $x$
if and only if,
on a neighborhood of $x$,
there exists
an element $a\in 
\Gamma(U,{\cal O}_U^\times)$
such that
$\chi$ is defined by
$t^p=a$ and satisfying
one of
the following conditions:
\begin{itemize}
\item[{\rm (\ref{prcln}.3)}]
$a=1-
u\cdot
t_1^{m_1}\cdots t_n^{m_n}$
for a unit $u$ at $x$
and,
for at least one $D_i$,
the integer $m_i$
is prime to $p$.
\item[{\rm (\ref{prcln}.4)}]
$a=1-
c\cdot
t_1^{m_1}\cdots t_n^{m_n}$
for a regular function
$c$ at $x$
such that $dc|_C$ has
no zero at $x$.
\end{itemize}
\end{pr}

{\it Proof.}
We have an exact sequence
$0\to \Omega^1_C
\to 
\Omega^1_{X/S}(\log D)
\otimes_{{\cal O}_X}
{\cal O}_C
\to
\bigoplus_i
{\cal O}_C
\to 0$
and $(d\log t_i)_i$
defines a splitting.
Hence, if
the condition
(\ref{corkuma}.1)
holds,
then $\chi$
is clean
and we have
$D_\chi=
\sum_ie'_iD_i$.
Thus,
it suffices
to consider the case
where
(\ref{corkuma}.2) holds.

1. Assume
$D_\chi=
\sum_ie'_iD_i$
and hence $m_i=0$
for every $i$.
Then,
$\chi$ is clean
if and only $da|_C$
has no zero at $x$.

2.
Assume
$D_\chi<
\sum_ie'_iD_i$
and hence $m_i>0$
for some $i$.
Then, by 
Corollary \ref{corkuma}.1.,
there exists 
a regular function $c$ at $x$
such that
$\chi$ is defined
by $t^p=a$ for
$a=1-c\cdot
t_1^{m_1}\cdots t_n^{m_n}$.
By Corollary \ref{corkuma}.3.\
and by the local splitting above,
$\chi$ is clean
if and only if
either
$c\cdot (m_i)_i\in 
\bigoplus_i{\cal O}_C$
or
$dc\in 
\Omega^1_C$
has no zero at $x$.
The condition that
$c\cdot (m_i)_i\in 
\bigoplus_i{\cal O}_C$
has no zero at $x$
means that
$c$ is a unit at $x$
and one of $m_i$
is prime to $p$.
The second condition
is equivalent to
(\ref{prcln}.4).
\qed

\begin{cor}\label{corscl}
Let the assumption be as in 
Proposition {\rm \ref{prcln}}.

{\rm 1.}
Assume $D_\chi=
\sum_ie'_iD_i$.
Then,
$\chi$ is clean
at $x$
if and only if
$\chi$ is s-clean at $x$.

{\rm 2.}
Assume $D_\chi<
\sum_ie'_iD_i$.
Then,
$\chi$ is s-clean
at $x$
if and only if,
on a neighborhood of $x$,
there exists
an element $a\in 
\Gamma(U,{\cal O}_U^\times)$
such that
$\chi$ is defined by
$a$ and satisfying
either 
the condition
{\rm (\ref{prcln}.3)}
or 
the following condition:
\begin{itemize}
\item[{\rm (\ref{prcln}.4$'$)}]
$a=1-
c\cdot
t_1^{m_1}\cdots t_n^{m_n}$
for a regular function
$c$ at $x$
such that $dc|_C$ has no zero at $x$
and,
for every $D_i$,
the integer $m_i$
is divisible by $p$.
\end{itemize}
\end{cor}

{\it Proof.}
1.
If the condition
(\ref{prcln}.1)
or
(\ref{prcln}.2)
is satisfied,
then $\chi$
is s-clean at $x$.

2.
If the condition
(\ref{prcln}.3)
is satisfied,
then $\chi$
is s-clean at $x$.
Assume
the condition
(\ref{prcln}.4)
is satisfied.
Then, in the notation
of the proof of
Proposition \ref{prcln}.2,
$\chi$
is s-clean at $x$
if and only if
either
$c\cdot (m_i)_i\in 
\bigoplus_i{\cal O}_C$
has no zero at $x$
or
$c|_C\cdot (m_i)_i=0$.
The first condition
is equivalent to
(\ref{prcln}.3).
By the condition
(\ref{prcln}.4),
we have
$c|_C\neq 0$.
Hence, the second
condition
$c|_C\cdot (m_i)_i=0$
is equivalent to
that
$m_i$ is divisible
by $p$ for every $i$.
\qed

We recall the main
result from 
\cite{Kato} and
prove a complement.

\begin{lm}\label{lmscl2}
Let the assumption be as in 
Proposition {\rm \ref{prcln}}.
Assume $\dim X_K+1=2$
and let $\Sigma\subset \Sigma_s
\subset D$
be the sets of points
where ${\cal F}$
is {\rm not} clean 
and {\rm not} $s$-clean 
with respect to $X$
respectively.

{\rm 1.}
The subsets $\Sigma$
and $\Sigma_s$
consist of finitely
many closed points of $D$.

{\rm 2. \cite[Remark 4.13]{Kato}}
Let $x$
be a closed point of $D$
and $f\colon X'\to X$
be the blow-up at $x$.
If $\chi$ is clean at $x$,
then
$\chi$ is clean
on a neighborhood of
$f^{-1}(x)$.

{\rm 3. \cite[Theorem 4.1]{Kato}}
There exists a
successive blow-up
$f\colon X'\to X$
at $\Sigma$
such that 
$\chi$ is clean with
respect to $X'$.

{\rm 4.}
There exists a
successive blow-up
$f\colon X'\to X$
at $\Sigma_s$
such that 
$\chi$ is s-clean with
respect to $X'$.
\end{lm}

{\it Proof.}
1.
Clear from the definition.

2.
3.
See 
\cite[Remark 4.13]{Kato}
and 
\cite[Theorem 4.1]{Kato}
respectively.

4.
By 3.,
we may assume 
$\chi$ is clean with respect
to $X$.
Let $x\in \Sigma_s$.
If $x$ is a singular
point of $D$,
then
$\chi$ is $s$-clean at $x$.
Hence,
$x$ is a smooth point of $D$
and we may assume $D$
is irreducible.
We put $D_\chi=rD$
and prove the assertion
by induction on $r>0$.
By \cite[Corollary 4.9]{Kato}
and \cite[Theorem (8.1)]{deg1},
we have $e=1$
in the notation
of \cite[Corollary 4.9]{Kato}
for the blow-up $X'\to X$
at $x$.
Hence,
we have $r'=r-e=r-1<r$
and the assertion follows
by induction.
\qed

We expect that
Lemma \ref{lmscl2}.3
holds in arbitrary dimension.

For an s-clean character
of order $p$,
we give a local description
of the normalization
in the corresponding
cyclic covering
of degree $p$.

\begin{pr}[cf.\ {\cite[Lemma 1]{TS}}]
\label{prkumY}
Let $X$ be a regular flat separated
scheme of finite of finite
type over ${\cal O}_K$ and
$U=X\setminus D$
be the complement of
a divisor $D$ with simple normal crossing.
Assume that
$K$ is of characteristic $0$ and contains
a primitive $p$-th root $\zeta_p$ of $1$.

Let $\chi\in
H^1(U,{\mathbb Z}/p{\mathbb Z})$
be a character of
order $p$.
Let $V\to U$ be the cyclic
Galois covering of order $p$
trivializing $\chi$
and $Y$ be
the integral closure
of $X$ in $V$.

Assume that
$\chi$ is s-clean
with respect to $X$.
Then,
there exists
an fs log structure 
${\cal M}_Y$ on $Y$
such that
$(Y,{\cal M}_Y)$ is log regular {\rm \cite{KT}}
and that
$V$ is the maximum
open subscheme
where ${\cal M}_Y$ is
trivial.
We have
${\cal M}_Y=j_*{\cal O}_V^\times
\cap {\cal O}_Y$
where $j\colon V\to Y$
denotes the open immersion.
\end{pr}

{\it Proof.}
The proof is similar
to that of
\cite[Lemma 1]{TS}.
Since the question
is local, 
we may assume that
$X={\rm Spec}\ A$
is affine
and that
the log structure
on $X$ is defined by the
chart $P=
\prod_i{\mathbb N}e_i
\to A$ sending
the basis
$e_i$ to $t_i$
defining irreducible
components $D_i$ of $D$.
By Corollary \ref{corscl},
it suffices to
consider each case 
(\ref{prcln}.1)--(\ref{prcln}.3) 
and (\ref{prcln}.4$'$).
We take the notation
in Proposition
\ref{prcln}.

(\ref{prcln}.1) 
Assume $a=u\prod_i
t_i^{m_i}$
where $u\in A^\times$
and $p\nmid m_i$
for at least one $i$.
We put $Q_0=
P\times {\mathbb Z}e_u,\
e_t=\frac 1p(\sum_im_ie_i+e_u)$
and $Q_1=Q_0+
\langle e_t\rangle$.
Let $Q\subset Q_1^{\rm gp}$ 
be the saturation
of $Q_1$.
We put
$B_1=A[t]/(t^p-u
\prod_it_i^{m_i})$
and define
a monoid homomorphism
$Q_1\to B_1$
by $e_u\mapsto u$
and $e_t\mapsto t$.
We define $B=
B_1\otimes_{A[Q_1]}A[Q],\
Y={\rm Spec}\ B$
and define a log
structure ${\cal M}_Y$ on $Y$ 
by $Q\to B$.

We show that $Y$ is log regular
and is equal to
the normalization of
$X$ in $V$.
By the assumption
that there exists
$m_i$ prime to $p$,
the quotient
$Q_1^{\rm gp}/Q_1^\times=
(P^{\rm gp}\times {\mathbb Z}e_t)/
{\mathbb Z}e_u$
is torsion free.
The quotient
$B_1/I_1$ by the ideal
$I_1\subset B_1$
generated by 
$Q_1\setminus Q_1^\times$
is equal to
$\overline A=
A/(t_i; i=1,\ldots,n)$
and is regular.
Since $B_1$ is flat over $A$,
we have
$\dim B_1=\dim A
=\dim \overline A
+n
=\dim B_1/I_1
+{\rm rank}\
Q_1^{\rm gp}/Q_1^\times$.
Hence, similarly as
the proof of \cite[Claim]{TS},
the log scheme $Y$
is regular
by \cite[Proposition (12.2)]{KT}.
Since the normal scheme 
$Y$ is finite
over $X$
and $Y\times_XU=V$,
it is the normalization
of $X$ in $V$.

(\ref{prcln}.2) 
Assume $a\in A^\times$
and $da|_C$ has no zero.
We put
$B=A[t]/(t^p-a)$,
$\bar A=A/(t_1,\ldots,t_n)$
and
$\bar B=\bar A[t]/(t^p-\bar a)$.
By the assumption
that $d\bar a$ is non-vanishing,
the ring
$\bar B$ is regular
and hence
$B=A[t]/(t^p-a)$
is the normalization
on a neighborhood of $x$.
The open subscheme
$V$ of $Y={\rm Spec}\ B$
is the complement
of a divisor with
simple normal crossing
defined by
$(t_1\cdots t_n)$.

(\ref{prcln}.3) 
Let $u\in A^\times$ be a unit,
$0\le m_i\le e'_i$ be
integers and
we put $b=u\prod_it_i^{m_i}$
and $a=1-b$.
We assume
that $p\nmid m_i<e'_i$
for at least one $i$.
We put $z=\zeta_p-1$ as above
and $t=1-z/s$.
Then, by an elementary
computation,
the equation
$t^p=a$
is equivalent to
$$s^p=\frac 1b(s^p-(s-z)^p).$$
We define a polynomial
$f\in {\cal O}_K[S]$
of degree $p-1$
by $f=(S^p-(S-z)^p)/z^p$.
We have $f\equiv 1-S^{p-1}
\bmod z$.
We put
$B_1=A[s]/(s^p-\frac{z^p}b\cdot f(s))$.
We put $r_i=
e'_i-m_i\ge 0$
and
define a unit
$w\in A^\times$ by 
$z^p/b=w\cdot 
\prod_it_i^{r_i}$.
The assumption
$m_i<e'_i$
for at least one $i$
means that
$z^p/b=w\cdot 
\prod_it_i^{r_i}$
is in the maximal ideal at $x$.
Hence, shrinking $X$
if necessary,
we may assume $f(s)$ 
is a unit of $B$.

We put $Q_0=
P\times {\mathbb Z}e_w,\
e_s=\frac 1p
(\sum_ir_ie_i+e_w)$
and $Q_1=Q_0+
\langle e_s\rangle$.
Let $Q\subset Q_1^{\rm gp}$ 
be the saturation
of $Q_1$.
We define
a monoid homomorphism
$Q_1\to B_1$
by $e_w\mapsto w\cdot f(s)$
and $e_s\mapsto s$.
We define $B=
B_1\otimes_{A[Q_1]}A[Q],\
Y={\rm Spec}\ B$
and define a log
structure ${\cal M}_Y$ on $Y$ 
by $Q\to B$.
Then similarly
as in the case
(\ref{prcln}.1),
we see that
$Y$ is log regular
and is equal to
the normalization of
$X$ in $V$.

(\ref{prcln}.4$'$)
Let $c\in A$,
$0\le m_i\le e'_i$ be
integers and
we put $b=c\prod_it_i^{m_i}$
and $a=1-b$.
We assume
that $m_i<e'_i$
for at least one $i$
and $p$ divides
$m_i=p\cdot m'_i$ for
every $i$.
We also assume
that $dc|_C$
is not 0 at $x$.
We put $t=1-
\prod_it_i^{m'_i}s$.
Then, the equation
$t^p=a$
is equivalent to
$(1-
\prod_it_i^{m'_i}s)^p
=1-\prod_it_i^{m_i}c$.
We define a 
polynomial
$g\in A[S]$
of degree $p$ by
$g=(S-\prod_it_i^{-m'_i})^p
+
\prod_it_i^{-m_i}$.
Then,
as in the case (\ref{prcln}.2),
$B=A[s]/(g(s)+c)$
is the normalization
of $A$.
The open subscheme
$V$ of $Y={\rm Spec}\ B$
is the complement
of a divisor with
simple normal crossing.

Since $(Y,{\cal M}_Y)$
is log regular,
we have
${\cal M}_Y=j_*{\cal O}_V^\times
\cap {\cal O}_Y$.
\qed

We consider
the sheaf of differential
forms
$\Omega^1_{Y,{\cal M}_Y/S}$
defined with respect to
the log structure
${\cal M}_Y$ and the 
trivial log structure on $S$.

\begin{cor}[cf.\ {\cite[Lemma 1]{TS}}]
\label{corTS1}
\setcounter{equation}0
Let the notation
and the assumption
be as in Proposition
{\rm \ref{prkumY}}
and $\sigma$ be
a generator of
${\rm Gal}(V/U)$.
Let ${\cal I}_\sigma$
denote the ideal sheaf
defining the
fixed part $Y^{\sigma}
\subset Y$.
Then, we have
the following.

{\rm 1.}
For each geometric point
$\bar y$ of the fixed
part $Y^\sigma$,
the action of
$\sigma$ on the stalk
${\cal M}_{Y,\bar y}/
{\cal O}^\times_{Y,\bar y}$
is trivial.

{\rm 2.}
We define an ideal sheaf
${\cal J}_\sigma$
of ${\cal O}_Y$
to be that
generated by
${\cal I}_\sigma$
and 
$\sigma(b)/b-1$
for $b\in {\cal M}_Y$.
Then, 
${\cal J}_\sigma$
is an invertible ideal
defined by an effective
Cartier divisor
$D_\sigma$.
Further, we have
\begin{equation}
\pi^*D_\chi=
pD_\sigma.
\label{eqDcs}
\end{equation}

{\rm 3.}
We put $D'=(p-1)D_\sigma$.
Then, the coherent
${\cal O}_Y$-module
$\Omega^1_{Y/X}
(\log/\log)$ defined by
$\Omega^1_{Y/X}
(\log/\log)=
{\rm Coker}(
\pi^*\Omega^1_{X/S}
(\log D)\to
\Omega^1_{Y,{\cal M}_Y/S})$
is an invertible
${\cal O}_{D'}$-module.

{\rm 4.}
Define an ${\cal O}_Y$-linear
surjection
$\varphi_\sigma\colon
\Omega^1_{Y,{\cal M}_Y/S}
\to {\cal J}_\sigma/
{\cal J}_\sigma^2$
by
$\varphi_\sigma(da)
=\sigma(a)-a$
and
$\varphi_\sigma(d\log b)
=\sigma(b)/b-1$.
Then, it induces
an isomorphism
$\Omega^1_{Y/X}
(\log/\log)
\otimes_{{\cal O}_Y}
{\cal O}_{D_\sigma}
\to {\cal J}_\sigma/
{\cal J}_\sigma^2=
{\cal O}_{D_\sigma}
(-D_\sigma)$.

{\rm 5.}
We put $E_Y=
(\pi^*E)_{\rm red}$
where $E=
\sum_{i,r_i>0}D_i$
denotes the support of $D_\chi$.
Then, the sequence
\begin{equation}
\begin{CD}
0\to
{\cal O}_{E_Y}
(-\pi^*D_\chi)
@>{{\rm rsw}\
\chi}>>
\pi^*\Omega^1_X(\log D)
\otimes_{{\cal O}_Y}
{\cal O}_{E_Y}\\
@>>>
\Omega^1_{Y,{\cal M}_Y}
\otimes_{{\cal O}_Y}
{\cal O}_{E_Y}
@>{\varphi_\sigma}>>
{\cal O}_{E_Y}
(-D_\sigma)
\to 0
\end{CD}
\label{eqTS1}
\end{equation}
is exact.
\end{cor}

{\it Proof.}
Since the assertions are local,
it suffices
to consider the
cases
(\ref{prcln}.1)--(\ref{prcln}.3) 
or (\ref{prcln}.4$'$)
respectively.

1.
Let $\pi\colon
Y\to X$ be the canonical
map.
Since the canonical map
$\bar {\cal M}_{X,\pi(\bar y)}
\to 
\bar {\cal M}_{Y,\bar y}$
induces an isomorphism
$\bar {\cal M}_{X,\pi(\bar y)}
^{\rm gp}\otimes{\mathbb Q}
\to 
\bar {\cal M}_{Y,\bar y}
^{\rm gp}\otimes{\mathbb Q}$,
the assertion follows.

2.\ and 3.
We take the notation
in the proof of Proposition
\ref{prkumY}.

(\ref{prcln}.1)
The ideal
${\cal J}_\sigma$
is generated by
$\sigma(x)/x-1$
for $x\in Q$.
Hence, it is generated by $\zeta_p-1$.
We have
$D_\sigma=
{\rm div}(\zeta_p-1)$
and
$D_\chi=p\cdot
{\rm div}(\zeta_p-1)$
by Proposition \ref{prkumY}.

The ${\cal O}_Y$-module
$\Omega^1_{Y/X}
(\log/\log)$
is generated by
$d\log t$
and the annihilator is
$(p)$.
We have
${\rm div}\ p=
(p-1)\cdot
{\rm div}(\zeta_p-1)$.

(\ref{prcln}.2)
The ideal
${\cal J}_\sigma$
is generated by
$\sigma(t)-t
=(\zeta_p-1)t$.
Since $t$ is a unit,
it is generated by $\zeta_p-1$.
We also have
$D_\sigma=
{\rm div}(\zeta_p-1)$
and
$D_\chi=p\cdot
{\rm div}(\zeta_p-1)$
by Proposition \ref{prkumY}.

The ${\cal O}_Y$-module
$\Omega^1_{Y/X}
(\log/\log)$
is generated by
$dt$
and the annihilator is
$(p)$.
We have
${\rm div}\ p=
(p-1)\cdot
{\rm div}(\zeta_p-1)$.

(\ref{prcln}.3)
The ideal 
${\cal J}_\sigma$
is generated by
$\sigma(s)/s-1$.
By $s=z/(1-t)$,
we have
$\sigma(s)/s-1
=(1-t)/(1-\zeta_pt)-1
=(\zeta_p-1)t/(1-\zeta_pt)
=zt/(1-t-zt)$.
Since $1-t=z/s$,
it is further equal to
$zt/(z/s-zt)=st/(1-st)$.
Since $1-st$ and $t$
are unit,
the ideal ${\cal J}_\sigma$
is generated by $(s)$.
Since $(s^p)=(z^p/b)$,
we have
$p\cdot D_\sigma=
\pi^*D_\chi$.

We put
$g=
S^p-\frac {z^p}bf(S)
\in A[S]$.
Then, the 
${\cal O}_Y$-module
$\Omega^1_{Y/X}
(\log/\log)$
is generated by
$d\log s$ and 
the relation is
given by
$p\cdot d\log s=
d\log f(s)$
and 
$g'(s)\cdot s\cdot
d\log s=0$.
Since $g'(s)=
ps^{p-1}-\frac {z^p}bf'(s)
=
s^{p-1}(p-s\cdot f'(s)/f(s))$,
the annihilator is
$(p-s\cdot f'(s)/f(s))$.
Since $g'(s)=
\prod_{i=1}^{p-1}(s-\sigma^i(s))$
and
${\rm div}(1-\sigma^i(s)/s)
=D_\sigma$
for each $i$,
we have
${\rm div} 
(p-s\cdot f'(s)/f(s))=
(p-1)\cdot D_\sigma$.

(\ref{prcln}.4$'$)
The ideal 
${\cal J}_\sigma$
is generated by
$\sigma(s)-s$.
By $t=1-
\prod_it_i^{m'_i}s$,
we have
$\sigma(s)-s
=(\zeta_p-1)t/\prod_it_i^{m'_i}$.
Since $t$ is a unit,
the ideal ${\cal J}_\sigma$
is generated by $(z/
\prod_it_i^{m'_i})$.
Thus, we have
$p\cdot D_\sigma=
\pi^*D_\chi$.

The ${\cal O}_Y$-module
$\Omega^1_{Y/X}
(\log/\log)$
is generated by
$ds$
and the annihilator is
$(g'(s))$.
Since $g'(s)
=\prod_{i=1}^{p-1}
(s-\sigma^i(s))$
and
${\rm div}(s-\sigma^i(s))
=D_\sigma$
for each $i$,
we have
${\rm div}\ g'(s)=
(p-1)\cdot D_\sigma$.

4.
By the assertion 3.,
the map
$\varphi_\sigma\colon
\Omega^1_{Y/X}
(\log/\log)
\otimes_{{\cal O}_Y}
{\cal O}_{D_\sigma}
\to {\cal O}_{D_\sigma}
(-D_\sigma)$
is a surjection
of invertible
${\cal O}_{D_\sigma}$-modules
and hence is
an isomorphism.

5.
It follows from 4.\ that
the sequence
(\ref{eqTS1}) is exact at
$\Omega^1_{Y,{\cal M}_Y}
\otimes_{{\cal O}_Y}
{\cal O}_{E_Y}$
and at
${\cal O}_{E_Y}
(-D_\sigma)$.
Hence, the kernel
of the map in the middle
is an invertible 
${\cal O}_{E_Y}$-module.
By the assumption
that $\theta$ is clean,
the map ${\rm rsw}\ \theta$
is a locally split
injection.
The composition
${\cal O}_{E_Y}
(-\pi^*D_\chi)
\to
\pi^*\Omega^1_X(\log D)
\otimes_{{\cal O}_Y}
{\cal O}_{E_Y}
\to
\Omega^1_{Y,{\cal M}_Y}
\otimes_{{\cal O}_Y}
{\cal O}_{E_Y}$
is the 0-map
by (\ref{eqrswL}).
Hence,
the assertion follows.
\qed

\subsection{A computation
in the rank 1 case}\label{ssrk1}
\setcounter{equation}0

Let $X$ and $U=
X\setminus D$
over $S={\rm Spec}\
{\cal O}_K$
be as in the previous
subsection.
We briefly recall the definition of
the 0-cycle class 
$c_{\cal F}$ in \cite{Kato}
for a smooth sheaf ${\cal F}$
of $\bar {\mathbb F}_\ell$-vector
spaces
of rank 1 on $U$.
Let $D_1,\ldots,D_n$
be the irreducible components of $D$
and let $E=\sum_{r_i>0}D_i
\subset X$
be the support
of the Swan divisor
$D_{\chi}=
\sum_ir_iD_i$.
We put $n=\dim X_K+1$.
The divisor
$E$ is supported on
the closed fiber
$X_F$.
Hence,
the coherent
${\cal O}_E$-module
$\Omega^1_{X/S}(\log D)
\otimes_{{\cal O}_X}
{\cal O}_E$
is locally free
of rank $n$
and the bivariant Chern class
$c(\Omega^1_{X/S}(\log D)
\otimes_{{\cal O}_X}
{\cal O}_E)$
is defined
as an operator
$CH^*(E\to E)$.

Assume $\chi$ is clean
with respect to $X$.
Then, we define
the 0-cycle class $c_\chi\in CH_0(E)$
by
\begin{eqnarray}
c_\chi&=&
\{c(\Omega^1_{X/S}(\log D)
\otimes_{{\cal O}_X}
{\cal O}_E)^*\cap
(1+D_\chi)^{-1}\cap D_\chi\}_{\dim 0}
\label{eqswchi}\\
&=&
(-1)^{n-1}\sum_{i=1}^m
r_i\cdot
c_{n-1}({\rm Coker}({\rm rsw}_i(\chi)))\cap [D_i].\nonumber
\end{eqnarray}
By \cite[Theorem 5.2]{Kato},
the cycle classes $c_\chi$
define an element
of $F_0G(\partial_{V/U}U)$
for the finite \'etale
Galois covering $V\to U$
trivializing $\chi$,
if $\dim U_K\le 1$.

We fix an isomorphism
$\overline {\mathbb F}_\ell^\times
\to 
{\mathbb Q}/{\mathbb Z}[\frac1\ell]
\subset
{\mathbb Q}/{\mathbb Z}$.
For a character
$\chi\in H^1(U,
{\mathbb Q}/{\mathbb Z}[\frac1\ell])$
of order prime to $\ell$,
let ${\cal F}_\chi$
denote
the corresponding
locally constant
constructible sheaf of
$\overline {\mathbb F}_\ell$-vector spaces
of rank 1 on $U$.

\begin{cn}\label{cncF}
Let $X$ be a regular
flat separated scheme
of finite type over $S$
and $U=X\setminus D$
be the complement
of a divisor with simple
normal crossings.
Let $f\colon V\to U$
be the \'etale cyclic
covering trivializing
$\chi$. 
Assume that
$\chi$ is clean with respect to $X$
and that $\chi$ is
tamely ramified on the generic fiber.

Then, we have
\setcounter{equation}0
\begin{equation}
{\rm Sw}_{V/U} {\cal F}_{\chi}
=f^*c_{\chi}
\label{eqcF}
\end{equation}
in $F_0G(\partial_{V/U}V)_
{{\mathbb Q}(\zeta_{p^\infty})}.$
\end{cn}
We prove a refinement of
Conjecture \ref{cncF}
assuming $\dim U_K\le 1$
in Theorem \ref{thmrk1}
at the end of this section.
Similarly as \cite[Lemma 5.1.2]{KSA},
Conjecture \ref{cncF}
implies 
Conjecture \ref{cnHA},
by Brauer induction.

We show that
the class $c_\chi$
also satisfies
an excision formula.

\begin{lm}\label{lmexcc}
Let $X$ be a regular
flat separated scheme
of finite type over $S$
and $U=X\setminus D$
be the complement
of a divisor with simple
normal crossings.
Let $X_1$ be 
a regular divisor
meeting $D$
transversely.
We put $U_0=
X\setminus (D\cap X_1)$
and $U_1=U\cap X_1$.

Let $\chi$ be a character
on $U$
and let $\chi_0=
\chi|_{U_0}$
and $\chi_1=
\chi|_{U_1}$
be the restrictions.
Assume that
both $\chi$ and
$\chi_0$
are clean with respect to
$X$.
Then 
$\chi_1$
is also clean with respect to
$X_1$
and we have
$$c_\chi=
c_{\chi_0}+c_{\chi_1}.$$
\end{lm}

{\it Proof.}
The union
$D'=D\cup X_1$
is a divisor
of $X$ 
with simple normal crossings
and
the intersection
$D_1=D\cap X_1$
is a divisor
of $X_1$ 
with simple normal crossings.
Let $D_i$
be an irreducible component
of $E$
and we put
$C_i=D_i\cap X_1$.
The image of the map
$\Omega^1_{X/S}
(\log D)
\otimes_{{\cal O}_X}
{\cal O}_{C_i}
\to
\Omega^1_{X/S}
(\log D')
\otimes_{{\cal O}_X}
{\cal O}_{C_i}$
is canonically
identified
with
$\Omega^1_{X_1/S}
(\log D_1)
\otimes_{{\cal O}_{X_1}}
{\cal O}_{C_i}$.
Hence, if
both $\chi$ and $\chi_0$
are clean,
then $\chi$
is strongly
clean on a neighborhood of
$C_i$
in the terminology
of \cite[Definition (7.4)]{deg1}.
Thus, by
\cite[Theorem (9.1)]{deg1},
$\chi_1$ is strongly
clean with respect to
$X_1$
and 
$D_{\chi_1}$
is the pull-back of $D_{\chi}$.

We put $E_1=E\cap X_1$.
Then, by the exact sequence
$0\to
\Omega^1_{X/S}(\log D)
\otimes_{{\cal O}_X}
{\cal O}_E
\to
\Omega^1_{X/S}(\log D')
\otimes_{{\cal O}_X}
{\cal O}_E
\to 
{\cal O}_{E_1}\to 0$,
the difference
$c_\chi-
c_{\chi_0}$
is equal to
\begin{align*}
&\{c(\Omega^1_{X/S}(\log D)
\otimes_{{\cal O}_X}
{\cal O}_E)^*\cap
(1+D_\chi)^{-1}\cap D_\chi\}_{\dim 0}
\\&\qquad
-
\{c(\Omega^1_{X/S}(\log D')
\otimes_{{\cal O}_X}
{\cal O}_E)^*\cap
(1+D_\chi)^{-1}\cap D_\chi\}_{\dim 0}
\nonumber
\\
=&\
\{c(\Omega^1_{X/S}(\log D)
\otimes_{{\cal O}_X}
{\cal O}_E)^*\cap
(1+D_\chi)^{-1}\cap(X_1\cap D_\chi)
\}_{\dim 0}.
\end{align*}
By $D_{\chi_1}=
X_1\cap D_\chi$
and by
$c(\Omega^1_{X/S}(\log D))^*\cap
X_1=
c(\Omega^1_{X_1/S}(\log D_1))^*$,
the right hand side
is equal to
$c_{\chi_1}$.
\qed

Let $U'\to U$
be a finite \'etale morphism
tamely ramified with
respect to $X$
and let $X'$
be the normalization
of $X$ in $U'$.
Then, $X'$ has a natural
log structure
such that
$U'$
is the maximum open
subscheme
where the log structure
is trivial
and the map
$X'\to X$
is log \'etale with respect
to this log structure.
By taking a regular
proper
subdivision of the
associated fan
\cite[Section 10]{KT},
we obtain
a log blow-up
$X''\to X'$
such that
$X''$ contains
$U'$
as the complement
of a divisor with
simple normal crossings.

\begin{lm}\label{lmtmc}
Let $X$ be a regular
flat separated scheme
of finite type over $S$
and $U=X\setminus D$
be the complement
of a divisor with simple
normal crossings.
Assume $\dim X_K=1$.
Let $g\colon
U'\to U$
be a finite \'etale
morphism
tamely ramified
with respect to $X$.
Let $X'$ be a log blow-up
of the normalization of
$X$
and $\bar g\colon
X'\to X$
be the canonical map.

Let $\chi$ be a character
on $U$
and let $\chi'$
be the pull-back to
$U'$.
Assume that
$\chi$ is clean with respect to
$X$.
Then 
$\chi'$
is also clean with respect to
$X'$
and we have
$$c_{\chi'}=
\bar g^*c_{\chi}.$$
\end{lm}

{\it Proof.}
We may assume that
the subdivision 
defining $X'$
induces a subdivision
of the fan associated
to $X$ and defines
a log blow-up $X_1\to X$.
The induced map
$X'\to X_1$ is
finite.
Since $\bar g\colon
X'\to X$
is log \'etale,
the map
$\bar g^*\Omega^1_{X/S}
(\log D)
\to
\Omega^1_{X'/S}
(\log D')$
is an isomorphism.
At each singular
point of $D$,
$\chi$ is strongly clean
with respect to $X$.
Hence, by \cite[Theorem (8.1)]{deg1},
the divisor
$D_{\chi_1}$
is the pull-back of
$D_\chi$
and the divisor
$D_{\chi'}$
is also the pull-back of
$D_\chi$.
Further by \cite[Theorem (8.1)]{deg1},
$\chi'$ is clean with respect
to $X'$
and we have
$c_{\chi'}=
\bar g^*c_\chi$.
\qed

In the rest of this section,
we assume that $K$ is of characteristic $0$.
The case where $K$ is of characteristic $p>0$
is studied similarly as in \cite{KSA} and
in \cite{TS}.

We first show that
the computation
in the previous subsection
implies
Conjecture \ref{cncF}
for a character of
order $p$
under a slightly stronger
assumption.
Since $(Y,{\cal M}_Y)$
in Proposition \ref{prkumY}
is log regular,
by \cite{KT},
there exists a
log blow-up
$\widetilde Y\to Y$
satisfying the following
property:
The map
$\widetilde Y\to Y$
induces an
isomorphism over
$V$ and 
the scheme
$\widetilde Y$
is regular and
contains
$V$ as the complement
of a divisor $D_{\widetilde Y}$
with simple normal
crossings.
We regard
$\widetilde Y$
as a log scheme
with the log structure
defined by $D_{\widetilde Y}$.
Then, the map
$\widetilde Y\to Y$
is log \'etale.

\begin{pr}\label{prth}
Let $X$ be a regular
flat separated scheme
of finite type over 
$S={\rm Spec}\ {\cal O}_K$
and $U=X\setminus D_X$
be the complement
of a divisor with simple
normal crossings.
Assume that $K$
contains a primitive $p$-th
root $\zeta_p$ of $1$.
Let $\theta$
be a character of order $p$
s-clean with respect to $X$.

Let $\widetilde Y\to Y$
be a log blow-up as above
and assume that
the action of
$G={\rm Gal}(V/U)
\simeq
{\mathbb Z}/p{\mathbb Z}$
is extended
to an action on 
$\widetilde Y$.
Let $\sigma$
be a generator
of $G$
and assume
that $\sigma$
is an admissible
automorphism 
(Definition {\rm \ref{dfadm}.2})
of $\widetilde Y$.

{\rm 1.}
We put
$c_\sigma=
((\Gamma_\sigma,
\Delta_{\widetilde Y}))_
{(\widetilde Y
\times_S\widetilde Y)^\sim}
\in F_0G(\partial_
{V/U}\widetilde Y)$.
Then, we have
$$c_\sigma=
\{c^*(\Omega^1_{\widetilde Y/S}
(\log D_{\widetilde Y}))
\cdot (1+D_\sigma)^{-1}
\cdot D_\sigma\}_{\dim 0}.$$

{\rm 2.}
Let $\pi\colon
\widetilde Y\to X$
be the canonical map.
Then, we have
\setcounter{equation}0
\begin{equation}
{\rm Sw}^{\mathbb Z}
_{V/U,\widetilde Y} 
{\cal F}_{\theta}
=\pi^*
c_{\theta}
\label{eqth}
\end{equation}
in $F_0G(\partial_
{V/U}\widetilde Y)$.
\end{pr}

{\it Proof.}
1.
This is 
Lemma \ref{lmsig}.2.

2.
Since $\widetilde Y
\to Y$ is log \'etale,
we have an exact sequence
$0\to 
{\cal O}_{E_{\widetilde Y}}
(-\tilde \pi^*D_\theta)
\to
\Omega_{X/S}
(\log D_X)\otimes {\cal O}_
{E_{\widetilde Y}}
\to
\Omega_{\widetilde Y/S}
(\log D_{\widetilde Y})
\otimes {\cal O}_
{E_{\widetilde Y}}
\to
{\cal O}_{E_{\widetilde Y}}
(-D_\sigma)\to 0$
by Corollary \ref{corTS1}.2.
Hence the ratio of the
total Chern classes
$c(\Omega_{\widetilde Y/S}
(\log D_{\widetilde Y})
\otimes {\cal O}_
{E_{\widetilde Y}})\cdot
c(\Omega_{X/S}
(\log D_X)
\otimes {\cal O}_E))^{-1}$
is equal to
$(1-D_\sigma)
(1-\pi^*D_\theta)^{-1}$.
Thus, by
the equality
$\pi^*D_\theta=pD_\sigma$
(\ref{eqDcs}),
we have
\begin{eqnarray*}
\pi^* c_\theta
&=&
\pi^*
\{c^*(\Omega^1_{X/S}(\log D_X))
\cdot (1+\pi^*D_\theta)^{-1}
\cdot
D_\theta\}_{\dim 0}
\\
&=&
\{c^*(\Omega^1_{\widetilde Y/S}
(\log D_{\widetilde Y}))
\cdot (1+\pi^*D_\sigma)^{-1}
\cdot
p\cdot D_\sigma\}_{\dim 0}.
\end{eqnarray*}
By 1.,
the right hand
side is equal to
$p\cdot c_\sigma$.

Since $c_{\sigma^i}=
c_\sigma$
for every $i\in 
({\mathbb Z}
/p{\mathbb Z})^\times$,
the integral Swan class
${\rm Sw}^{\mathbb Z}
_{V/U,\widetilde Y} 
{\cal F}_{\theta}$
is equal to
$p\cdot
s_{V/U}(\sigma)=
-p\cdot c_\sigma$.
Thus the assertion follows.
\qed

We recall an induction
step from \cite{TS},
which will be used
in the proof of
Theorem \ref{thmrk1} below.

\begin{lm}[cf.\ {\cite[Lemma 2]{TS}}]
\label{lmTS2}
Let $X$ be a regular scheme
and $U$ be the complement
of a divisor with simple normal crossings.
Let $\chi,\theta\in
H^1(U,{\mathbb Q}/{\mathbb Z})$
be characters
clean with respect to $X$.
Assume $\theta$
is of order $p$
and s-clean.
Let $V\to U$
be the cyclic covering
of degree $p$
trivializing $\theta$
and $Y$ be the 
normalization of $X$ in
$V$ with the
log structure defined by $V$.
Let $\widetilde Y\to Y$
be the log blow-up
defined by
a regular proper subdivision
of the fan of $Y$
and $\pi\colon 
\widetilde Y\to X$
be the canonical map.

Assume $\chi$ is clean
with respect to $X$ 
and the pull-back
$\chi'=\pi^*\chi$ is clean
with respect to $\widetilde Y$.
Assume further that
the following condition
is satisfied:
\begin{itemize}
\item[{\rm (\ref{lmTS2}.1)}]
For $D_\chi=\sum_ir_iD_i$
and
$D_\theta=\sum_is_iD_i$,
the condition
$r_i=0$ implies $s_i=0$,
and the condition
$r_i>0$ implies $r_i>s_i$.
\end{itemize}
Then, we have
\setcounter{equation}1
\begin{equation}
\pi^*c_\chi=
c_{\chi'}
+
D^{\log}_{U_1/V,Y_1}
\label{eqTS2}
\end{equation}
in $CH_0(E_{\widetilde Y})$.
\end{lm}

The proof is the same
as \cite[Lemma 2]{TS}
by using the exact sequence
(\ref{eqTS1})
and we omit it.

\begin{cor}\label{corTS2}
We keep the notation
and the assumptions
in Lemma {\rm \ref{lmTS2}}
except that we do not assume
{\rm (\ref{lmTS2}.1)}.
Assume further
that $\chi$ is of order
$n=mp$
and $\theta=m\cdot \chi$.
Assume that
the Swan class 
${\rm Sw}^{\mathbb Z}
_{V/U_1,Y} {\cal F}_{\chi'}$
is defined integrally and
that 
we have 
\setcounter{equation}0
\begin{equation}
{\rm Sw}^{\mathbb Z}
_{V/U_1,Y} {\cal F}_{\chi'}
=\bar g^*c_{\chi'}
\label{eqTS2ih}
\end{equation}
Then, 
the Swan class 
${\rm Sw}^{\mathbb Z}
_{V/U,Y}
{\cal F}_{\chi}$
is also defined integrally and
we have 
\setcounter{equation}1
\begin{equation}
{\rm Sw}^{\mathbb Z}
_{V/U,Y}
{\cal F}_{\chi}
=\bar f^*c_{\chi}.
\label{eqTS2i}
\end{equation}
\end{cor}

{\it Proof.}
By the definition
of the Swan class
and by $s_{V/U}(\sigma^i)=
s_{V/U}(\sigma)$
for an integer prime to $p$,
we obtain
$${\rm Sw}^{\mathbb Z}
_{V/U,Y}
{\cal F}_{\chi}-
D^{\log}_{V/U,Y}
=
{\rm Sw}^{\mathbb Z}
_{V/U_1,Y}
{\cal F}_{\chi_1}-
D^{\log}_{V/U_1,Y}.$$
By $\theta=m\cdot \chi$,
the assumption 
{\rm (\ref{lmTS2}.1)}
is satisfied.
Hence, by applying
Lemma \ref{lmTS2},
we obtain
$\bar f^*c_\chi=
\bar g^*c_{\chi_1}
+
\bar g^*D^{\log}_{U_1/V,Y_1}$.
Thus, it follows
from the assumption
(\ref{eqTS2ih})
and the chain rule
$D^{\log}_{U/V,Y}=
D^{\log}_{U/U_1,Y}
+
\bar g^*D^{\log}_{U_1/V,Y_1}$.
\qed

In the rest of the paper,
we consider the case
$\dim U_K=1$.
In this case,
the strong form 
of resolution
of singularity
is known and
consequently
the Swan class
${\rm Sw}^{\mathbb Z}
_{V/U}{\cal F}$
is defined integrally
as an element
of $F_0G(\partial_{V/U}V)$.

\begin{thm}\label{thmrk1}
Assume $K$ is of
characteristic $0$.
Let $U$ be a regular
flat separated scheme
of finite type over $S$
such that
$\dim U_K=1$.
Let ${\cal F}=
{\cal F}_\chi$
be a locally constant
constructible
sheaf of 
$\overline
{\mathbb F}_\ell$-vector
spaces of rank $1$
and $\chi
\in H^1(U,{\mathbb Q}/
{\mathbb Z})$
be the corresponding
character.
Let $f\colon
V\to U$ be the cyclic
covering trivializing $\chi$.
Then, we have 
\setcounter{equation}0
\begin{equation}
[K(\zeta_p):K]
\cdot
{\rm Sw}^{\mathbb Z}_{V/U}
{\cal F}_{\chi}
=
[K(\zeta_p):K]
\cdot
f^*c_{\chi}
\label{eqrk1}
\end{equation}
in $F_0G(\partial_{V/U}V)$.
\end{thm}

{\it Proof.}
Let $\chi'$ be the
$p$-primary part of $\chi$
and $V'\to U$
be the cyclic covering
trivializing $\chi'$ and
let
$\pi\colon
V\to V'$ be
the canonical map.
Then,
since ${\rm Sw}^{\mathbb Z}
_{V/U}{\cal F}_{\chi}
=\pi^*
{\rm Sw}^{\mathbb Z}
_{V'/U}{\cal F}_{\chi'}$
and $c_{\chi}
=c_{\chi'}$,
we may assume 
that the order of
$\chi$ is a power of $p$.

We show that
we may assume $U=U_K$.
Let $X$ be a proper
regular flat scheme
over $S$ containing
$U$ as the complement
of a divisor $D$ with
simple normal crossings.
By blowing up 
some closed points
in the closed fiber
of $X$,
we may assume
that $U_K\subset X$
is the complement
of a divisor with
simple normal crossings.
We show the claim
by the induction
on the number
of irreducible
components of $U_F$.

If the number is 0,
then $U_F$ is empty and
there is nothing to prove.
Let $C$ be an
irreducible component
of $U_F$.
Let $\chi_0$
be the restriction of
$\chi$ to $U_0=U
\setminus C$.
By blowing up $X$
at the boundary of $C$,
we may assume
that both
$\chi$ and 
$\chi_0$
are clean with respect to $X$.
Then, by the excision formulas 
Proposition \ref{prExc}.2
and
Lemma \ref{lmexcc},
the equality
(\ref{eqcF})
for $\chi_0$
is equivalent to
that for $\chi$.
Thus, by the induction,
the claim is proved.

We assume $U=U_K$.
By taking
the base change to
$K(\zeta_p)$
and by applying
Corollary \ref{corfun}
and Lemma \ref{lmtmc},
we may assume that
$K$ contains 
a primitive $p$-th
root of 1.

Assume $\chi$ is
of order $p$.
Then, by Lemma \ref{lmscl2}.4,
we may assume $\chi$
is $s$-clean with respect to $X$
by replacing $X$ by
a blow-up.
Then, it follows
from Proposition \ref{prth}.

Assume $\chi$ is
of order $p^n$
and
we prove the assertion
by induction on $n\ge 1$.
Similarly as above,
we may assume
that
$\theta=p^{n-1}\chi$
is s-clean.
Let $U_1\to U$
be the cyclic covering
of degree $p$
and $g\colon
V\to U_1$
be the canonical map.
Let $Y_1$ be the normalization
of $X$ in $U_1$ 
and $\widetilde Y_1\to Y_1$
be a blow-up as in Corollary
\ref{corTS2}.
Let $\pi\colon 
\widetilde Y_1\to X$
be the canonical map.
Then, further similarly as above,
we may assume that
$\chi'=\pi^*\chi$
is clean with respect to 
$\widetilde Y$.
Then, by the induction 
hypothesis,
we have
${\rm Sw}^{\mathbb Z}
_{V/U_1,\widetilde Y}
{\cal F}_{\chi'}
=g^*c_{\chi'}.$
Thus it follows from
Corollary \ref{corTS2}.
\qed


We deduce the integrality of 
the Swan classes
and the conjecture of Serre
under the assumption 
$\dim U_K\le 1$.

\begin{cor}\label{corint}
{\rm 1.}
Let $U$ be a regular flat separated
scheme of finite type over ${\cal O}_K$.
If $\dim U_K\le 1$,
Conjecture {\rm \ref{cnHA}}
is true.

{\rm 2.}
For a regular local ring $A$
of $\dim A\le 2$,
Conjecture {\rm \ref{cnSe}}
is true.
\end{cor}

{\it Proof.}
1.
By the Brauer induction and
by the induction formula
Proposition \ref{prind},
we may assume 
${\rm rank}\ {\cal F}=1$.
Then, it follows from
Theorem \ref{thmrk1}.

2. 
Since the positive characteristic
case is proved in \cite{Artin},
it follows from
1 and Lemma \ref{lmSe}.
\qed

As in the classical ramification
theory,
our proof of the integrality
Conjecture {\rm \ref{cnHA}}
is by the reduction to
the rank 1 case
using Brauer induction.

\noindent 
K.\ Kato: Department of Mathematics, 
University of Chicago, 
Chicago 60637 Illnois, USA\
kkato@math.uchicago.edu\\
T.\ Saito: Department of Mathematical Sciences, 
University of Tokyo, 
Tokyo 153-8914, Japan\
t-saito@ms.u-tokyo.ac.jp

\end{document}